\documentclass[10pt]{amsbook}

\let\svref\ref
\usepackage{amsmath}
\usepackage{amsthm}
\usepackage{amssymb}
\usepackage{pstricks}

\newfont{\rmm}{cmr10 scaled 1000}
\newfont{\itt}{cmsl10 scaled 1000}

\theoremstyle{plain}

\newtheorem{theo}{Theorem}[section]
\newtheorem{lemm}[theo]{Lemma}
\newtheorem{prop}[theo]{Proposition}
\newtheorem{coro}[theo]{Corollary}

\theoremstyle{definition}

\newtheorem{defi}[theo]{Definition}
\newtheorem{rema}[theo]{Remark}

\newcommand{\bq}{\begin{equation}}

\usepackage{ifthen}

\newcounter{lemma}[section]

\newcounter{tempcounter}

\newcommand{\mylabel}[1]{\label{#1}\setcounter{tempcounter}{\thechapter}%
\addtocounter{tempcounter}{-1}\refstepcounter{tempcounter}\label{ch#1}}

\newcommand{\myref}[1]
{
\ifthenelse{ \equal{\svref{ch#1}} {\thechapter}  }
{\ref{#1}}
{\ifthenelse{ \equal{\svref{ch#1}} {0} }
{\ref{#1} of  Introduction}
{\ref{#1} of Chapter ~\ref{ch#1}}}}
\newcommand{\lb}{\label}
\newcommand{\mlb}{\mylabel}

\newcommand{\mrf}[1]{\myref{#1}}

\newcommand{\RC}{Rokhlin Conjecture}
\newcommand{\RRR}{{\mathbf R}}
\newcommand{\CCC}{{\mathbf{C}}}
\newcommand{\ZZZ}{{\mathbf{Z}}}
\newcommand{\KKK}{{\mathbf{K}}}
\newcommand{\NNN}{{\mathbf{N}}}

\newcommand{\PP}{{\mathcal P}}

\newcommand{\GG}{{\mathcal G}}

\newcommand{\HH}{{\mathcal H}}
\renewcommand{\AA}{{\mathcal A}}
\newcommand{\gA}{{\mathfrak{A}}}

\newcommand{\pr}{\partial}

\renewcommand{\a}{\alpha}
\renewcommand{\b}{\beta}

\newcommand{\ti}{\times}
\newcommand{\bk}{\backslash}

\newcommand{\bede}{\begin{defi}}
\renewcommand{\beth}{\begin{theo}}
\newcommand{\bele}{\begin{lemm}}
\newcommand{\bepr}{\begin{prop}}
\newcommand{\bere}{\begin{rema}}

\newcommand{\beeq}{\begin{equation}}
\newcommand{\bega}{\begin{gather}}
\newcommand{\been}{\begin{enumerate}}
\newcommand{\beal}{\begin{aligned}}

\newcommand{\enth}{\end{theo}}
\newcommand{\enle}{\end{lemm}}
\newcommand{\enpr}{\end{prop}}
\newcommand{\enre}{\end{rema}}

\newcommand{\enga}{\end{gather}}
\newcommand{\eneq}{\end{equation}}
\newcommand{\enen}{\end{enumerate}}
\newcommand{\enal}{\end{aligned}}

\newcommand{\e}{\epsilon}
\newcommand{\krest}{\begin{picture}(14,14)
\put(00,04){\line(1,0){14}}
\put(00,02){\line(1,0){14}}
\put(06,-4){\line(0,1){14}}
\put(08,-4){\line(0,1){14}}
\end{picture}     }

\newcommand{\sbs}{\subset}
\newcommand{\pom}{polynomial}
\newcommand{\pol}{polyhedron}
\newcommand{\poms}{polynomials}
\newcommand{\pols}{polyhedra}

\newcommand{\sut}{such that}
\newcommand{\G}{\Gamma}
\newcommand{\D}{\Delta}

\begin{document}
\title{Proof of the Rokhlin's Conjecture on Arnold's
surfaces}

\author{F.Nicou}
\address{Universit\'e de Rennes I, Campus de Beaulieu, F-35042 Rennes Cedex,
France}
\email{nicou@maths.univ-rennes1.fr}
\begin{abstract}
In this paper we prove that Arnold Surfaces of
all real algebraic curves of even degree with non-empty real part
are standard (Rokhlin's Conjecture).
There is an obvious connection with classification of Arnold Surfaces
up to isotopy of $S^4$ and Hilbert's Sixteen Problem on the
arrangements of connected real components of curves.
First, we  consider some $M$-curves, i.e curves  of a prescribed degree
having the greatest possible number of connected real components,
and prove that Arnold surfaces of these
curves are standard.
Afterwards,
we define a procedure of modification "perestroika"
of these $M$-curves which allows to
prove the Rokhlin's Conjecture.
\end{abstract}
\maketitle

\tableofcontents

\part{\bf Introduction}
\mlb{p:intro}
A real algebraic curve $\AA$ of even degree is by definition a
homogeneous
real polynomial $F(x_0,x_1,x_2)$ of degree $m=2k$ in three variables.
The set
$$\CCC \AA=\{ (x_0:x_1: x_2 ) \in \CCC P^2 \mid F(x_
0,x_1,x_2 )=0\}$$
is called {\it complex point set }
of the curve.
The complex conjugation $conj$ of $\CCC P^2$ induces an antiholomorphic
involution on $\CCC \AA$ with fixed point  set
 $\RRR \AA = \CCC \AA \cap \RRR P^2$, which is the set of real
points of the curve.
Given a real algebraic curve $\AA$ with non-singular
complex set point $\CCC \AA$,
using the fact that $\CCC  P^2/conj $ is diffeomorphic to $S^4$ \cite{ArS},
one can  associate to
the triple
$(\CCC  P^2,conj,\CCC\AA)$
 two pairs $(S^4,\gA)$ where $\gA$ is a
closed connected
surface embedded
in $S^4$ called {\it  Arnold Surface}.
 These surfaces were introduced in
\cite{Ar}. The present  paper is devoted to
studying the topology of these pairs.
 The Arnold Surfaces are constructed as follows.
The real part $\RRR \AA$ is a union of disjoint circles embedded
in $\RRR  P^2$
which divide $\RRR  P^2$ into two  parts:\\
 $\RRR  P^2_+= \{ x \in \RRR  P^2
\mid F(x_0,x_1,x_2)  \geq  0\}$ and
 $\RRR  P^2_-= \{ x \in \RRR  P^2
\mid F(x_0,x_1,x_2)  \leq  0\}$, which have
 common boundary $\RRR \AA$.
Gluing $\RRR  P^2_+$,  respectively $\RRR  P^2_-$, to
the  quotient  $\CCC \AA/conj$
along their common boundary $\RRR \AA$,
 we get two closed surfaces in $S^4$:
$$\gA_\pm=(\CCC \AA/conj)\cup(\RRR  P^2)_\pm$$
One of these surfaces is obviously non orientable. Changing the
 sign of $F$
if necessary, we can always assume that the non-orientable
 part of $\RRR  P^2$ lies
entirely in $\RRR  P^2_-$. Then  $\gA_-$ is non-orientable.
  Note that these surfaces are not smooth along
$\RRR \AA$ but can be easily smoothed
since $\CCC \AA$ is nowhere tangent to $\RRR  P^2$.
In this paper we consider only the smoothed Arnold surfaces and we
keep for them
the same name "Arnold surfaces", since no confusion is possible.
We shall say that an embedded surface
$S\sbs S^4$ is standard if it is a connected sum of
standard tori and standard $\RRR P^2$.
Moreover, when it does not lead to  confusion,
we shall call curve a  real algebraic
curve with non-singular complex point set.
\vskip0.1in
{\bf Rokhlin's conjecture:} {\it Let $\AA$ be a real algebraic curve with
 non-empty real part $\RRR \AA$.
\footnote
{In the case of empty real part, $\gA_-$ is trivially non-standard.}
Then both Arnold surfaces are standard.}
\vskip0.1in The proof of the Rokhlin's Conjecture for the surface $\gA_+$
associated to the maximal nest curves can be extracted from the
paper  \cite{Ak}
of S.Akbulut.
(Recall that the maximal nest curve is a curve whose real part is
, up to isotopy in $\RRR P^2$,  $k$ circles
linearly ordered by inclusion.)
In the paper \cite{Fi},
S.Finashin
proved Rokhlin's Conjecture for
$L$ curves of even degree.
(Recall that an $L$ curve
is a perturbation of a curve
which splits into a union of $2k$ real lines in general position.)
\vskip0.1in

Here is a brief description of the methods of our proof.

Let $k$ be a positive integer. Consider an algebraic curve  of degree $2k$
whose real part is described up to isotopy in $\RRR P^2$ as
$ (k-1)(2k-1) +1$ disjoint circles embedded in $\RRR P^2$
such that:
$ \frac {(k-1)(k-2)}{2}$ circles lie inside one circle and
$ \frac {k(k-3)}{2} $ circles lie outside this circle where the inside
(resp, outside) of a circle design respectively the part of $\RRR P^2$
homeomorphic to a disc (resp,  a M\"{o}bius band).

  Such a curve was initially obtained in 1876 from Harnack's construction
based on perturbations of singularities \cite{Har}.
  It can be as well obtained by the Patchworking procedure
due to O.Viro. We recall the Patchworking procedure in the
 Preliminary Section.
 The construction of Harnack curve using the
 Patchworking method was given by  I.Itenberg
in \cite{It}.
We recall these construction
of the so-called Harnack curve
in the Chapter \ref{ch:PatchHar}.

Our proof starts in the Chapter \ref{ch:Mope}.
In the Chapter \ref{ch:ArnHar},
we prove that Arnold surfaces for
Harnack curves are standard. The proof makes
use of a desingularization procedure.

Afterwards (sections \mrf{su:typ1}, \mrf{su:typ2})
we define a procedure of modification
of the Harnack curve which allows to obtain
(up  to conj-equivariant isotopy) all real algebraic curves
with non-empty real part.
In the chapter  \ref{ch:arn}
we prove that a curve $\AA$ obtained as the result of
such modification has
standard Arnold surfaces if $\RRR \AA$ is non-empty.

\vskip0.1in
{\bf Acknowledgement}
\vskip0.1in

I am very grateful to Ilia Itenberg for posing the problem
and for useful conversations during the course of this work.

\newpage

\chapter*{Preliminaries-Patchworking Process}
\lb{ch:Patchproc}

\subsection{Introduction}
\mlb{su:Preliminaries}
\vskip0.1in
The T-curves are constructed by a combinatorial procedure
due to O.Viro \cite{Vi}.
They are naturally introduced via  the procedure of
patchworking polynomials.

Here we introduce some notations and definitions.
\vskip0.1in
In what follows, $\KKK$ denotes either
the real field $\RRR$ or the complex field $\CCC$;
$\ZZZ$ denotes the  ring of integers, $\NNN$ the set of positive integers,
$\RRR_+$ the set of positive reals.

\vskip0.1in

A linear combination of products
$\{ x^iy^j \mid
(i,j) \in \ZZZ ^2 \} $ with coefficients from $\KKK$ is called
an {\it $L$-polynomial} over $\KKK$ in two variables.
L-polynomials over $\KKK$ in two variables
form a ring isomorphic to the ring of regular function of the variety
$(\KKK^*)^2 $. The variety $(\KKK^*)^2$ is known as an algebraic torus
over $\KKK $.\\
Let us recall some properties of $(\KKK ^*)^2$.\\
Let $l$ be the map $(\KKK^*)^2 \to \RRR^2$ defined by:
$l(x,y)=(ln(\vert x \vert),ln(\vert y \vert)).$
Let $U_{\KKK} =\{ x \in \KKK \mid \vert x \vert =1 \}$ and  $a$ be the map
$(\KKK^*)^2 \to U_{\KKK}^2$ defined by :
$a(x,y)= (\frac {x}{\vert x \vert}, \frac {y}{\vert y \vert})$.
The map
$la=(l,a): (\KKK^*)^2 \to \RRR^2 \times U_{\KKK}^2$ is a diffeomorphism.\\
Being an abelian group $(\KKK^*)^2$ acts on itself by translations.\\
Let us recall some translations involved into this action.
For $(i,j) \in \RRR^2$ and $t>0$,
denote by $qh_{(i,j),t}:(\KKK ^*)^2 \to (\KKK ^*)^2 $
 the translation defined by formula :
$qh_{(i,j),t} (x,y)=(t^{i} x,t^{j} y)$.
For $(x',y') \in U_{\KKK}^2$ denote by $S_{(x',y')}$ the translation
$(\KKK ^*)^2 \to (\KKK ^*)^2$ defined by formula :
$S_{(x',y')}(x,y)=(x'.x,y'.y)$.

\vskip0.1in
Any usual real polynomial
$f=f(x,y)=\sum_{(i,j)\in\NNN^2}a_{i,j}x^iy^j $
can be considered as a
Laurent polynomial with coefficients from $\KKK $.
\vskip0.1in
Let
$f=f(x,y)=\sum_{(i,j)\in\NNN^2}a_{i,j}x^iy^j $
be a real polynomial in two variables.
Denote by $\D(f)$ and call {\it Newton polyhedron} of $f$ the convex hull of
$\{w\in\NNN^2\mid a_w\not= 0\}$ .
\vskip0.1in
By { \it polyhedron}
we mean a convex bounded polyhedron in $ \RRR^2$ with integer vertices
having positive coordinates, in other words
a Newton polyhedron of a real polynomial in two variables of finite degree.
\vskip0.1in
The interior of a polyhedron $\D$ is denoted $\D^0$.
\vskip0.1in

A polyhedron is {\it non-degenerate} if its dimension is two,
in other words its interior is non-empty.
\vskip0.1in

A {\it subdivision} of a non-degenerate polyhedron  $\D$ is a set of
non-degenerate polyhedra $\D_1,...\D_r$
\sut~$\D=\cup_i\D_i$,
and any intersection $\D_i\cap \D_j$ is either empty or is a common
face of
both $\D_i, \D_j$.
We call {\it proper face} of a polyhedron $\D$ a face of codimension 1
and denote by $\GG'(\D)$ the set of proper faces of a polyhedron $\D$.
The set of faces of a polyhedron $\D$ is denoted by $\GG(\D)$.
\vskip0.1in
We denote $C(\D)$ the  vector subspace
of $\RRR^2$ which corresponds to the minimal affine subspace
containing $\D$.

\vskip0.1in

\subsection{Toric varieties}
\mlb{susu:torvar}
To each
polyhedron one can associate a variety called\\ $\KKK$-toric variety.
The following definition is deduced from (\cite{GKZ}
p.166-168).
\bede
\mlb{d:torvar}
\vskip0.1in
\cite{GKZ}
{\it Let $\D$  be a polyhedron.

Consider the finite set $A(\D)= \ZZZ^2 \cap \D$ and choose an ordering
of this set so that\\ $A(\D)=\{(i_1,j_1),...,(i_n,j_n) \}$.
Consider the subset
$$\KKK \D ^0= \{(x^{i_1}y^{j_1}:...: x^{i_n}y^{j_n}) \mid
(x,y) \in (\KKK^*)^2 \}$$ of $\KKK P^{n-1}$.

The closure  of $\KKK \D ^0$ is
called the $\KKK$-projective toric variety associated to $\D$ and the
ordering of $A$.}
\end{defi}
\bere
\mlb{r:rem1}
\vskip0.1in
\been
\item Another choice of the ordering of $A$ leads to an
isomorphic $\KKK$-projective toric variety, so we shall use
the term $\KKK$-projective toric variety associated to $\D$ whenever
no confusion is possible.
\item Moreover, let $B \subset \ZZZ^2$ and $T : \ZZZ^2 \to \ZZZ^m$ be an
affine integer injective transformation such that  $T(A) =B$,
the $\KKK$ projective variety associated to $\D$ and the set $B$ is
isomorphic to the $\KKK$-projective toric variety associated to $\D$ and an
ordering set $A(\D)$.
\enen
\enre

\vskip0.1in
\vskip0.1in
\vskip0.1in

Hence, denote by $\KKK \D$  the
$\KKK$-projective toric variety associated to a polyhedron $\D$.

\vskip0.1in
\bede
\mlb{d:torvar2}
\cite{GKZ}
\vskip0.1in
{\it Consider the action of the torus $(\KKK^*)^2$ given by formula
$(x,y). (z_1:...:z_n)=(x^{i_1}y^{j_1} z_1:...: x^{i_n}y^{j_n} z_n)$,
the variety $\KKK \D$ is the closure  of the orbit of
the point $(1:...:1)$ under this action.}
\end{defi}
\vskip0.1in
Let us recall some properties of toric varieties.
We refer to \cite{Vi},\cite{GKZ}
for more details about toric varieties.
The following Lemma is implicit in
\cite{GKZ}, Prop 1.9,p.171.
\bele \cite{GKZ}
{\it Let $\D$ be a polyhedron.
Let $\G$ be a face of $\D$.
   Let $A(\D)= \ZZZ^2 \cap \D$ and choose an ordering
of this set so that $A(\D)=\{(i_1,j_1),...,(i_n,j_n) \}$.
Consider the subset $A(\G)=\ZZZ^2 \cap \G$.\\
Let $\pr :A(\D) \to {0,1}$ be the map such that $\pr (A(\G)) =1$
$\pr(A(\D) \bk A(\G))=0$.
Consider the subset
$$\KKK' \G ^0= \{(x^{i_1}y^{j_1} \pr((i_1,j_1)):...:(x^{i_n}y^{j_n}
\pr((i_n,j_n)) \mid
(x,y) \in (\KKK^*)^2 \}$$ of $\KKK P^{n-1}$.
The closure  of $\KKK' \G ^0$
is a subvariety of  $\KKK \D$ which corresponds to the closure of the orbit
of the point $(\pr(i_1,j_1):...:\pr(i_n,j_n))$
under the action of the torus $(\KKK^*)^2$ given by formula
$(x,y). (z_1:...:z_n)=(x^{i_1}y^{j_1} z_1:...: x^{i_n}y^{j_n} z_n)$.
Denote it by $\KKK' \G$.
The variety $\KKK' \G$ is isomorphic to $\KKK \G$
the $\KKK$-projective toric variety associated to $\G$.
In the cases when it does not lead to confusion we shall identify  $\KKK' \G$
with $\KKK \G$ and $\KKK' \G^0$ with $\KKK \G^0$.}
\enle
\bepr \cite{GKZ}
{ \it
Let $\D$ a polyhedron.
Let $\KKK \D$ be the $\KKK$-projective
toric variety associated to $\D$.
Translations $qh_{(i,j),t}$ and $S_{(x',y')}$  of $(\KKK ^*)^2$
define naturally translations $qh_{(i,j),t}$ and $S_{(x',y')}$ on
$\KKK \D^0$ which can be extended to transformations on $\KKK \D$.}
In what follows,
we keep the same notations for translations
of $\KKK \D^0$ and their extensions to transformations of $\KKK \D$.
\enpr
{\bf proof:}
\vskip0.1in
On the assumptions of  definition \mrf{d:torvar}, let $\KKK \D$ be the
$\KKK$-projective toric variety associated to $\D$.
Since $\KKK' \G^0$ and their closure $\KKK'\G$,with $\G$ a face of $\D$,
are contained in the affine parts of $\KKK P^{n-1}$,
extensions of
transformations $qh_{(i,j),t}$ and $S_{(x',y')}$  on varieties
$\KKK' \G$ give transformations
on the affine parts of $\KKK P^{n-1}$ and thus define transformations
of $\KKK \D$. Q.E.D
\vskip0.1in

\bede \cite{GKZ}
\mlb{d:mo}
{\it Let $\D$ be a polyhedron.
On the same assumptions as in Definition \mrf{d:torvar},
let $\KKK \D$  be the closure  of $\KKK \D^0$
$$\KKK \D ^0= \{(x^{i_1}y^{j_1}:...: x^{i_n}y^{j_n}) \mid
(x,y) \in (\KKK^*)^2 \}$$ subset of $\KKK P^{n-1}$.
Let
$\RRR_+ \D$  be the closure (in the classic topology) of
$$\RRR_+ \D ^0= \{(x^{i_1}y^{j_1}:...: x^{i_n}y^{j_n}) \mid
(x,y) \in (\RRR_+^*)^2 \}$$ subset of $\KKK \D$.
The subvariety $\RRR _+ \D$ of $\KKK \D$ is
homeomorphic,as stratified space,to the polyhedron $\D$ stratified by its
faces.
An explicit homeomorphism is given by
the moment map $\mu:\KKK \D \to \D,\mu(x,y)= \frac {\sum_{(i,j) \in
A(\D)} \vert x^iy^j
\vert .(
i,j)} { \sum_{ (i,j) \in A(\D)} \vert x^iy^j \vert } $ }
\end{defi}

{\bf proof:}
\vskip0.1in
It follows easily from the properties  of
the moment map $\mu$(\cite{At}; see also \cite{GKZ}, p.198 Theorem1.11)
$\mu:\KKK \D \to \D $.Q.E.D
\vskip0.1in
\vskip0.1in

Consider the action of the torus $U_{\KKK}^2$
on $\KKK \D$ given by formula
$(x,y). (z_1:...:z_n)=(x^{i_1}y^{j_1} z_1:...: x^{i_n}y^{j_n} z_n)$.
\vskip0.1in
For a face $\G$ of $\D$,
denote by $U_{\G}$ the subgroup of $U_{\KKK}^2$ consisting of elements
$(e^{i \pi k}, e^{i \pi l})$ with $(k,l)$ a vector orthogonal
to $C(\G)$.
\vskip0.1in
\bele  \cite{GKZ}
\mlb{l:mo}
{\it
Let $\D$ be a polyhedron.
The map $\rho : \RRR_+\D \ti U_{\KKK}^2 \to \KKK \D$ defined by formula
$\rho((x,y),(x',y'))= S_{(x',y')}(x,y)$ is a proper surjection.
The variety $\KKK \D$ is homeomorphic to the
quotient space of $\RRR_+\D \ti U_{\KKK}^2$
with respect to the partition into sets
$x \ti y U_{\G}$, $x \in \RRR_+ \D \cap \KKK \G^0, y \in U_{\KKK}^2 $.}
\enle
{\bf proof:}\\
Let $\D$ be a polyhedron.  Let  $\RRR_+ \D$  be the closure of
$$\RRR_+ \D ^0= \{(x^{i_1}y^{j_1}:...: x^{i_n}y^{j_n}) \mid
(x,y) \in (\RRR_+^*)^2 \}$$ subset of $\KKK \D$.
Translations $S_{(x',y')}, (x',y')\in U_{\KKK}^2$ define an action of
$U_{\KKK}^2$ in
$\KKK \D$ such that intersection of $\RRR_+ \D$  with each
orbit under this action consists of one
point.
Furthermore, for $x$ in the interior of $\G$ where $\G$ is a face of
$\D$, the stationary subgroup of action of $U_{\KKK}^2$ consists of
transformations $S_{(e^{i \pi k}, e^{i \pi l})}$ where $(k,l)$ is a vector
orthogonal to $C(\G)$.  Since $\KKK \D$
is locally compact and Hausdorff,it follows the Lemma.
Q.E.D
\vskip0.1in
\bele
{\it The $\KKK$-projective toric variety associated to a
 non-degenerate triangle $\D$ is $\KKK P^2$.
The $\KKK$-projective toric variety associated to a degenerate triangle
$\D$ is $\KKK P^1$.}
\enle
{\bf  proof:}\\
It can be easily deduced from the Lemma \mrf{l:mo} above.
Let $\KKK =\RRR$.
Place $\D \ti U_{R}^2$ in $\RRR^2$ identifying
$(x,y) (x',y') \in \D \ti U_R^2$ with
$S_{(x',y')}(x,y)$.
The surface $\RRR \D$ can be obtained by an appropriate gluing of four
copies of $\D$. \\
Let $\KKK =\CCC$.
Place $\D \ti U_{C}^2$ in $\RRR^2$ identifying
$(x,y) (x',y') \in \D \ti U_C^2$ with
$S_{(x',y')}(x,y)$.
Let $\D$ be a non-degenerate triangle,
it is not difficult to get the usual handlebody decomposition of $\CCC P^2$
as the union of three 4-balls.
In the same way,
let $\D$ be a degenerate triangle, it is not difficult to get $\CCC P^1$
as the union of two 2-discs.Q.E.D

\subsection{Hypersurfaces of toric varieties}
\mlb{susu:HypTorvar}
Let $f$ be a real polynomial and $\D(f)$ its Newton polyhedron.
Let $\D$ be a polyhedron such that  $C(\D(f)) \subset C(\D)$.
The equation $f=0$
defines in $\KKK \D^0$ an hypersurface $V_{\KKK \D^0}(f)$ of $\KKK\D^0$.
Denote by $V_{\KKK\D}(f)$ its closure (in the Zarisky topology) in
$\KKK \D$.
(In the case $\KKK = \CCC$ the classic topology gives the same result,
but in the case $\KKK = \RRR$ the usual closure may be a non-algebraic set)
\vskip0.1in
\bele
(\cite{Vi}, p.175)
{\it
Let $f$ a real polynomial and $\D(f)$ its
Newton polyhedron.
Let $\D$ a non-degenerate polyhedron ( $C(\D(f)) \subset C(\D)$).
For any vector $(k,l) \in C(\D)$
orthogonal to $C(\D(f)) $, the hypersurface
$V_{\KKK_{\D}}(f)$ is invariant under transformations
$S_{(e^{i \pi k}, e^{i \pi l})}$ and $qh_{(k,l),t}$ of $\KKK \D$ .}
\enle
\vskip0.1in
Let $\D$ be a polyhedron.
An homeomorphism  $\phi: \D \to \D$
is said {\it admissible}
if it maps any face $\G$  of $\D$ to itself.
(Thus, an admissible  homeomorphism $\phi : \D \to \D$
extends to an homeomorphism $\tilde {\phi}:\D \ti U_{\KKK}^2 \to
 \D \ti U_{\KKK}^2$
by equivariance).

Hence, the following definition is naturally introduced.

\vskip0.1in
\bede \cite{Vi}
\mlb{d:chart}
{\it
\been
\item Let $f$  be a real polynomial in two variables,
$\D(f)$ its Newton polyhedron
and  $\KKK \D(f)$ the toric variety associated to $\D(f)$.
Denote by $h= \mu \vert_{\RRR_+ \D(f)}^{-1}$ the homeomorphism
$h : \D(f) \to \RRR_+ \D(f)$.
Let $\rho : \RRR_+ \D(f) \ti U_{\KKK}^2 \to \KKK \D(f)$ be the surjection
defined  by formula $\rho ((x,y),(x',y'))= S_{(x',y')}(x,y) = (x.x',y.y')$.
$$ \D(f) \ti U_{\KKK}^2 \to^{h \ti id} \RRR_+ \D(f)
\ti U_{\KKK}^2 \to^{\rho} \KKK\D(f)$$
A pair consisting of $(\D(f)) \ti U_{\KKK}^2$ and its subset $l$ which is
the pre-image of $V_{\KKK \D(f)}$ under $\rho \circ (h \ti id) $ is a
{\it canonical $\KKK $-chart} of $f$.
 \item
Let $ (\D(f) \ti U_{\KKK}^2, l)$
be the canonical $\KKK$-chart of $f$,
and\\ $\phi~:\D(f) \ti U_{\KKK}^2 \to \D(f) \ti U_{\KKK}^2$
the extended homeomorphism of an admissible homeomorphism
$\phi: \D(f) \to \D(f)$, the pair $( \phi~( \D(f) \ti U_{\KKK}^2, \phi~(l))$
is a {\it $\KKK $-chart} of $f$.
\item
Let $(\RRR)^2_+$ be set of positive integers
$\{ (x,y) \in \RRR^2 \mid x \ge 0, y \ge 0 \}$
Consider the map: $(\RRR^2)_+\ti U_{\KKK}^2 \to \CCC^2: ((x,y)(x',y')) \to
S_{(x',y')}(x,y)$.
Call {\it reduced $\KKK$-chart} of $f$ the image of a $\KKK$-chart of
$f$ under this map.
\enen}
\end{defi}
\vskip0.1in
{\it Singular  Hypersurfaces}
\vskip0.1in
Let
$f=f(x,y)=\sum_{(i,j) \in \NNN^2}a_{i,j}x^iy^j$
be a real polynomial.
\vskip0.1in
The set
 $$SV_{(\CCC ^*) ^2}(f)=
 \{(x,y) \in (\CCC^*)^2 \mid f(x,y)= \frac {\pr f}{\pr x}(x,y)
= \frac {\pr f}{\pr y}(x,y) =0 \} $$
defines the set of {\it singular points} of the variety defined
in $(\CCC^*)^2 $ by the equation $f=0$.
\vskip0.1in
For a set
$\G\sbs\NNN^2$ consider the polynomial
$ \sum_{(i,j)\in\G}a_{i,j}x^iy^j$.
It is called {\it $\G$-truncation} of $f$ and is denoted by $f^{\G}(x,y)$.

The polynomial $f$ is {\it completely non-degenerate}
if for any face $\G$ (including $\D(f)$)
of its Newton \pol~$\D(f)$
the variety defined in
$(\CCC^*)^2 $
by the equation $f^{\G}=0$ is non-singular.
\vskip0.1in
\bele\cite{Vi}
\mlb{l:le11}
{\it  Completely non-degenerate real polynomials form a
Zarisky open subset of the space of real polynomials
with a given Newton polyhedron.}
\enle

\section{Patchworking polynomials and T-curves}
\mlb{su:Patch}

In this section, we shall explain the procedure of patchworking
polynomials and define the concept of $T$-curve.

\subsection{Patchworking polynomials}
\mlb{susu:patch}

A reduced $\RRR$-chart of a non-degenerate polynomial
$a(x,y)=a_{i_1,j_1}x^{i_1} y^{j_1} +
a_{i_2,j_2}x^{i_2} y^{j_2} +
a_{i_3,j_3}x^{i_3} y^{j_3}$ may be constructed from the following method.

\vskip.1in
Let $a=a(x,y)=a_{i_1,j_1}x^{i_1} y^{j_1} +
a_{i_2,j_2}x^{i_2} y^{j_2} +
a_{i_3,j_3}x^{i_3} y^{j_3}$
be a real polynomial
with non-degenerate Newton polyhedron $\D$.
Assign to each vertex $a_{i,j}$ of $\D$
the sign
$$
\e_{i,j}=\frac {a_{i,j}}{\vert a_{i,j}\vert}
$$
Consider the union of $\D$ and its symmetric copies
$\D_x= s_x(\D)$
$\D_y= s_y(\D)$
$\D_{xy}=s(\D)$
where $s_x,s_y$ are reflections with respect to the coordinate axes
and $s = s_x \circ s_y$.
For each vector $(w_1,w_2)$  orthogonal
to $\D$ with integer relatively prime coordinates
glue the points $(x,y)$ and $((-1)^{w_1}x,(-1)^{w_2}y)$ of
the union
$\D \cup \D_{x} \cup \D_{y} \cup \D_{xy}$.
\vskip0.1in
Denote $\D_*$ the resulting space.
Extend the distribution of signs to $\D_*$ so that
$g^*(\e_{i,j}x^iy^j)=
\e_{g(i,j)}x^iy^j$
for
$g=s_x, s_y,s$.
Let $\bigtriangledown$ be one of the  four triangles
$\D, \D_{x}, \D_{y}, \D_{xy}$. If $\bigtriangledown$ has vertices of
different signs, consider the midline of $\bigtriangledown$
separating them.
Denote by $L$ the union of such midlines.

We shall say  that
{\it the pair $( \D_*, L)$ is obtained from $a$
 by combinatorial patchworking}.

\bele \cite{Vi}
\mlb{l:smallpatch}
{\bf "The smallest patch" }
\vskip0.1in
{\it Let $a(x,y)=a_{i_1,j_1}x^{i_1} y^{j_1} +
a_{i_2,j_2}x^{i_2} y^{j_2} +
a_{i_3,j_3}x^{i_3} y^{j_3}$
be a real polynomial
with non-degenerate Newton polyhedron $\D$.
Let $(\D_*, L)$, obtained from $a(x,y)$ by combinatorial patchworking.
The pair $(\D_*, L)$ is a reduced $\RRR$-chart of $a$.}
\enle
\vskip0.1in
It is easy to deduce the following characterization of
the set of real points $\{ (x,y) \in \RRR ^2, a(x,y)=0 \}$
and $ \{ (x_0:x_1:x_2) \in \RRR P^2, A(x_0,x_1,x_2)=0 \}$ where
$A$ is the homogenization of $a$.
\been
\item
Remove from $\D_*$ the sides of
$\D \cup \D_{x} \cup \D_{y} \cup \D_{xy}$
which are not glued in the construction of $\D_*$.
It turns the polyhedron $\D_*$ to the polyhedron $\D'$ homeomorphic to
$\RRR^2$ and the set $L$ to a set $L'$ such that the pair $(\D',v')$
is homeomorphic to $(\RRR^2, \{ (x,y) \in \RRR^2, a(x,y)=0 \})$.
\vskip0.1in
Glue by $s$ the opposite sides of $\D_*$.
The resulting space $\bar \D_*$ is homeomorphic
to the projective plane $\RRR P^2$.
Denote $\bar L$ the image of $L$ in $\bar\D_*$.
Let $A(x_0,x_1,x_2)$ the homogenization of $a(x,y)$.
It defines a curve $\mathcal A$.
Then there exists an homeomorphism
$(\bar \D_*,\bar L)\to(\RRR P^2, \RRR {\mathcal A})$.
\enen

\vskip0.1in
One can generalize this description to more complicated \poms~.
\vskip0.1in
Let $a_1,..., a_r$ be completely non-degenerate real polynomials
in two variables with Newton \pols~ $\D_1,..., \D_r$ in such a way
that   $\D_1,..., \D_r$
form a subdivision of a non-degenerate
\pol~$\D$ and
$a_i^{\D_i \cap \D_j} = a_j^{\D_i \cap \D_j}$

\bele \cite{Vi}
\mlb{l:patchpoly}
{\bf Patchworking polynomials}

{ Assume that the subdivision is {\it regular}, that is,
there exists a convex non-negative function $\nu:\D\to\RRR$, satisfying
the following conditions:
\begin{enumerate}
\item
all the restrictions $\nu\vert\D_i$ are linear.
\item
if the restriction of $\nu$
to an open set is linear, then this set is contained
in one of
$\D_i$
\item
$\nu(\D\cap\NNN ^2)\sbs\NNN$.
\end{enumerate}

Such a function is said {\it convexifying} the subdivision
 $\D_1,..., \D_r$ of $\D$.

There exists a unique \pom~ $a$ with $\D(a)=\D, a^{\D_i}=a_i$ for
$i=1,...,r$.  Let it be $a(x,y)=\sum_{(i,j)\in\NNN ^2}a_{i,j}x^iy^j $.}

\enle
Then introduce the one-parameter family of polynomials
$$
b_t=b_t(x,y)=
\sum_{(i,j)\in\NNN ^2}a_{i,j}x^iy^jt^{\nu(i,j)}
$$

We say that polynomials $b_t$
are obtained by
{\it Patchworking}
the polynomials
$a_1,...,a_r$ by $\nu$.

\bede
\mlb{d:patchchart} \cite{Vi}
{\bf Patchworking Charts }
{\it A pair $(\D \ti U_{\KKK}^2, L)$ is said to be
obtained by patchworking from $\KKK$-charts of polynomials
$a_1,..., a_r$
if $\D=\cup_{i=1}^r \D_i$ (where $\D_i$ denotes the Newton polyhedron
of $a_i$)
and one can choose $\KKK$-charts
$(\D_{i}\ti U_{\KKK}^2, l_i)$ of polynomials
$a_i$ such that $L=\cup_{i=1}^r (l_i)$.}
\end{defi}
\vskip0.1in
\beth
\mlb{t:patchtheo}\cite{Vi}
{\bf Patchwork Theorem}\footnote{This theorem is a simplified
version of the Patchworking Theorem for $L$-polynomials \cite{Vi}}
{\it
Let $a_1,..., a_r$ be completely
non-degenerate polynomial in two variables with Newton
polyhedra $\D_1$,...,$\D_r$.
Assume furthermore that $\D_1$,...,$\D_r$ form a regular subdivision of a
non-degenerate polyhedron $\D$.
Then there exists $t_0>0$ \sut~for any
$t\in ]0,t_0]$  the polynomial
$b_t$ is completely non-degenerate and
its $\KKK$-chart is obtained by
patchworking $\KKK$-charts of the
polynomials
$a_1,...,a_r$. }
\enth

Further, we say that such polynomials $b_t$ are obtained by
patchworking {\it process}.

We shall give the main ideas of the proof of the patchworking theorem
in the section \mrf{susu:Proof}.

\subsection{T-curves}

Bringing together the Lemma \mrf{l:smallpatch}
 ("the smallest patch ") and the patchworking
theorem \mrf{t:patchtheo}
 we deduce a combinatorial method of construction of
curves. These curves are called {\it T-curves}.

 Let $m$ be a positive integer.
 Let $\D$ be the triangle
 $$
\{(x,y)\in\RRR^2\mid x\geq 0, y\geq 0, x+y\leq m\}
$$
(Up to  linear change of coordinates $(x_0:x_1:x_2)$
of $\CCC P^2$, the convex hull of $\D$ may be the Newton polyhedron
of the affine polynomial $f(x,y)= F(1,x,y)$
associated to a homogeneous polynomial $F(x_0,x_1,x_2)$ of degree $m$.)

 Let $\tau$ be a regular triangulation of $\D$ whose vertices have integer
 coordinates. Suppose that some distribution of signs $\chi$  at the
 vertices of the triangulation is given.
 Denote the sign $\pm$ at the vertex with coordinates $(i,j)$ by
 $\epsilon_{(i,j)} $.
 Take the square $\D_*$ made of $\D$
and its symmetric copies
$\D_x=s_x(\D),
\D_y=s_y(\D),
\D_{xy}=s_{xy}(\D)$
where
$s_x,s_y, s=s_x\circ s_y$ are reflections with respect to
the coordinate axes. The resulting space is homeomorphic to $\RRR^2$.
Extend the triangulation $\tau$ of $\D$ to a symmetric triangulation
$\tau_*$ of $\D_*$.
Extend the distribution of signs to a distribution at the vertices
of $\D_*$ which verifies the modular properties:
$g^*(\e_{i,j})x^iy^j=
\e_{g(i,j)}x^iy^j$
for
$g=s_x, s_y,s$.

If a triangle of $\D_*$
has vertices of different signs, consider a midline separating them.
Denote
by $L$ the union of such midlines.
Glue by $s$ the opposite sides of $\D_*$.
The resulting space $\bar\D_*$ is homeomorphic
to the projective plane $\RRR P^2$.
Denote $\bar L$ the image of $L$ in $\bar\D_*$.
Call $\D_*$, $L$,
$\bar \D_*$,$\bar L$
 obtained from $\D,\tau,\chi$ by combinatorial patchworking.

\vskip0.1in

\beth
{\bf Patchwork Theorem from the real view point\\
 and T-curves.}\cite{Vi}
{\it Define the one-parameter family of polynomials
$$b_t=b_t(x,y)= \sum_{ (i,j) vertices of \tau} \e_{i,j}x^iy^jt^{\nu(i,j)}$$
where $\nu$ is a function convexifying the triangulation $\tau$
of $\D$.
Let,
$\bar \D_*$,$\bar L$
 obtained from $\D,\tau,\chi$ by combinatorial patchworking.

Denote by $B_t=B_t(x_0,x_1,x_2)$ the corresponding homogeneous polynomials:
$$B_t(x_0,x_1,x_2)=x_0^mb_t(x_1/x_0,x_2/x_0)$$

\vskip0.1in
Then there exists $t_0 > 0$ such that for any $t \in ]0,t_0]$
the equation $B_t(x_0,x_1,x_2)=0$ defines in $\RRR P^2$
the set of real points of an algebraic curve $C_t$ such that the pair
$(\RRR P^2, \RRR C_t)$
is homeomorphic to the pair $(\bar\D_*,\bar L)$.\\
Such a curve is called a {\it $T$-curve }.}
\enth
\vskip0.1in
{\bf proof :}
\vskip0.1in
It is an immediate consequence of the Lemma \mrf{l:patchpoly} and the
theorem \mrf{t:patchtheo}.
Let $ \D_*$,$ L$
obtained from $\D,\tau,\chi$ by combinatorial patchworking.
There exists $t_0$ such that for any $t \in ]0,t_0]$
the pair $(\D_*, L)$ is a reduced  $\RRR$-chart of $b_t$.
Q.E.D
\vskip0.1in
\bere
\been
\item
Remove the sides $x >0,y>0, x+y=m$ and its symmetric copies to $\D_*$.
It turns the polyhedron $\D_*$ to the polyhedron $\D'$ homeomorphic to
$\RRR^2$ and the set $L$ to a set $L'$ such that
for any $t \in ]0,t_0]$
the pair $(\D',v')$
is homeomorphic to $(\RRR^2, \{ (x,y) \in \RRR^2, b_t(x,y)=0 \})$.
\item
It is possible to recover in some cases whether
the set of real points of the
$T$-curve $B$ divides the set of its complex points.
Assume the triangulation $\tau$ of $\D$
sufficiently fine such that each triangle
$\D_i$ in the triangulation $\tau$ is the Newton
polyhedron of a polynomial $a_i$
which  defines a curve with orientable real set of points.
Denote $L_i$ the union
of midlines in $(\D_i)_*$ homeomorphic to this set.
If one can (resp, not)
choose an orientation on each $L_i$ compatible with an
orientation of the union of $\cup_i(L_i)$ then
the set of real points of the $T$-curve $B$  divides (resp does not divide)
the set of its complex points.
\enen
\enre
\subsection{Main Ideas of the Proof of the Patchworking Theorem}
\mlb{susu:Proof}
\vskip0.1in

In what follows,
we shall give definitions and statements
related to the patchworking method.
Afterwards, we shall give a sketch of proof of the Patchwork Theorem.
\vskip0.1in

{\it Preliminaries}\\
Let $a_1,..., a_r$ be completely non-degenerate polynomials
in two variables with Newton
polyhedra $\D_1$,...,$\D_r$
such that $\D_1$,...,$\D_r$ form  a regular subdivision of a non-degenerate
triangle $\D$.
Let $b_t$ be a polynomial obtained by patchworking process.
Let $B_t$ be the homogenization of $b_t$.
It defines a real algebraic curve ${\mathcal C}_t$.

In this section, we  explain the construction of $b_t$ and
give definitions we  shall need in the next sections.

\vskip0.1in

Recall that real polynomials in  two variables
belong to the ring
of $L$-polynomials, which is isomorphic to the ring of regular function
of the variety $(\KKK^*)^2$.

As already introduced, set $la$  the diffeomorphism
$$la=(l,a): (\KKK^*)^2 \to \RRR^2 \times U_{\KKK}^2$$

\vskip0.1in
Call {\it $\e$-tubular neighborhood} $N$ of a smooth submanifold $M$
of $(\KKK ^*) ^2$ the  normal tubular neighborhood of the smooth
submanifold $la(M)$ of
$\RRR^2 \ti U_{\KKK}^2$ whose fibers lie in  fibers
$\RRR \ti t \ti U_{\KKK} \ti s $, consist of segments of geodesics
 which are orthogonal to intersection of $la(M)$ with these
$\RRR \ti t \ti U_{\KKK} \ti s $, and
 are contained in balls of radius $\e$
 centered in the points of intersection of
 $la(M)$ with these
$\RRR \ti t \ti U_{\KKK} \ti s $.
The intersection of such tubular neighborhood  of $M$
with the fiber
$\RRR \ti t \ti U_{\KKK} \ti s $
is a normal tubular neighborhood of
$M \cap \RRR \ti t \ti U_{\KKK} \ti s $
in $R \ti t \ti U_{\KKK}\ti s$.
(Obviously, such definition of {\it $\e$ tubular neighborhood} of $M$
in $(\KKK ^*)^2$ requires
$la(M)$ transversal to
$R \ti t \ti U_{\KKK} \ti s $.)
\vskip0.1in
The following statement describes cases for which
such tubular neighborhood exists.

We refer to (\cite{Kh} p.59)
and also to(\cite{Vi} p.175, p178) for its proof.
\newpage
\vskip0.1in
\bele
\cite{Kh}
\mlb{l:Le1}
 {\it Let $f$ be a  polynomial and $\D(f)$ its
Newton polyhedron.
 Let $\G$ be a face of $\D(f)$ such that $f^{\G}$ is completely non-
 degenerate. Let $\rho: \RRR_+ \D(f) \ti U_{\KKK}^2 \to \KKK \D(f)$
 the  natural surjection .
Then the hypersurface $\rho^{-1}( V_{\KKK \D(f)}(f))$ is transversal
to $\RRR_+\G \ti U_{\KKK}^2$.}
\enle
\vskip0.1in
\vskip0.1in
Let $M$ be a smooth manifold of
$(\KKK^*)^2$ and $\D$ be a non-degenerate polyhedron.
Denote
$\rho : \RRR_+ \D \ti U_{\KKK}^2 \to \KKK \D$
the natural surjection.

Using the moment map
$\mu : \KKK \D \to \D$ which maps the torus
$(\KKK^*)^2 \approx \KKK \D^0$ onto the interior $\D^0$ of $\D$,
one can identify points of $\RRR_+ \D^0$ with points of $\D^0$.
In such a way, using this identification,
we call $\e$-tubular neighborhood of $M$ defined from
$p$ in  $\D^0$,
the $\e$-tubular neighborhood of $M$ in
$\rho(D(p,\e) \ti U_{\KKK}^2) \subset (\KKK^*)^2$
where $D(p,\e)$ is an open (euclidian) disc around $p$ of radius $\e$
inside $\D^0$.
Thereby, given a triangulation of $\D$, one can choose
$D(p,\e)$ sufficiently small such that $D(p,\e)$ intersects only
one proper $\G$ face of the triangulation.

\vskip0.1in
\bede \cite{Vi}
\mlb{d:def6}
{\it
Let $M$ a smooth manifold of
$(\KKK^*)^2$.
Let $\D$ be a non-degenerate polyhedron
and
$\rho : \RRR_+ \D \ti U_{\KKK}^2 \to \KKK \D$
be the natural surjection.
Given $\e >0$.
We call $\e$-tubular neighborhood of
$M \subset (\KKK^*)^2$
defined from a point
$p$ in  $\RRR_+ \D^0$,
the $\e$-tubular neighborhood of  $M$ in
$U(p)=\rho(D(p,\e) \ti U_{\KKK}^2) \subset (\KKK^*)^2$
where $D(p,\e)$ is an open (euclidian) disc around $p$ of radius $\e$
contained in $\RRR_+ \D^0$.}
\end{defi}
\vskip0.1in
{\it
Let  $D(p,\e) \in \RRR_+ T_m ^0$
be an open (euclidian) disc around $p$ of radius $\e$ such that
$\mu(D(p,\e))$ intersects only the face $\G$ of the triangulation
of $T_m$ and contains $\G$.
We call  $\e$-tubular neighborhood of $M$
defined from points of $\G^0$,
the $\e$-tubular neighborhood of $M$ in
$U(p)=\rho(D(p,\e) \ti U_{\KKK}^2) \subset (\KKK^*)^2$.
We call
$U(p)=\rho(D(p,\e) \ti U_{\KKK}^2) \subset (\KKK^*)^2 \subset (\KKK^*)^2$
the $\e$-neighborhood of $M$
defined from $\G^0$.}

\vskip0.1in

\bede  \cite{Vi}
{\it Let the norm in vector spaces of L-polynomials be:
$$ \vert \vert \sum a_wx^w \vert \vert =max \{ \vert a_w \vert | w\in
 \ZZZ^2 \}$$

Let $a$ a L-polynomial over $\KKK$ and $U$ a subset of $(\KKK^*)^2$
we say that {\it in $U$ the truncation $a^{\G}$ is $\e$-sufficient for $a$ }
if for any $L$-polynomial $b$ over $K$
such that $\D(b) \subset \D(a), b^{\G}= a^{\G},
\vert \vert b- b^{\G} \vert \vert \le \vert \vert a - a^{\G} \vert \vert$
the following conditions are verified:
\been
\item
$U \cap SV_{(\KKK^*)^2} (b) =\emptyset$
\item
the set $la (U \cap V_{(\KKK^*)^2}(b))$ lies in
a tubular $\e$-neighborhood  of \\
$la(V_{(\KKK^*)^2}(a^{\G}) \bk SV_{(\KKK^*)^2}(a^{\G}))$
\item
$la(U \cap V_{(\KKK^*)^2}(b))$ can be extended to the image
of a smooth section of the tubular fibration
 $ N \to la(V_{(\KKK^*)^2}(a^{\G}) \bk SV_{(\KKK^*)^2}(a^{\G}))$
\enen }
\end{defi}
\vskip0.1in

Let $\D$ a polyhedron.
Denote by $C_{\D}(\G)$ the cone
$C_{\D}(\G) =\cup_{r \in \RRR_+}r.(\D - y)$,
where  y  is  a  point of $\G \bk \pr \G$.
It contains the minimal affine subspace $C(\G)$
of $\RRR ^2$ containing $\G$, $C(\G) \subset C_{\D}(\G)$.
Denote by $DC_{\D}^-(\G) $ the set
$\{ x \in \RRR^2, a \in C_{\D}(\G) a.x \le 0 \}$.
For $A \subset \RRR^2$ and $\rho >0$, set
${\mathcal N}_{\rho}(A)
=\{ x \in \RRR ^2 \bk dist(x,A) < \rho \}$.
Let $\phi:
 \GG(\D) \to \RRR$ be a positive function,
denote by
$DC_{\D,\phi} (\G) $ the set ${\mathcal N}_{\phi (\G) } (D C_{\D}^- (\G)) \bk
\cup_{ \Sigma \in \GG (\D), \G \in \GG (\Sigma)}
{\mathcal N}_{ \phi (\Sigma) } (DC_{\D}^- (\Sigma))$.
\vskip0.1in

\bede \cite{Vi}
{\it Let $f$ be a Laurent polynomial over $\KKK$ and $\D$ its Newton
polyhedron.
A positive function $\phi: \GG(\D) \to \RRR$
 {\it describes  domain of $\e$-sufficiency for $f$} if
for any proper face $\G \in \G(\D)$,
for which truncation is completely non-degenerate and the hypersurface
$la(V_{\KKK \RRR^2}(f^{\G}))$ has an
$\e$-tubular neighborhood, the truncation $f^{\G}$ is
$\e$-sufficient
for $f$ in some neighborhood of $l^{-1}(DC_{\D,\phi}(\G))$.}
\end{defi}
\newpage
{\it Main Ideas of the proof of the Patchworking Theorem}
\vskip0.1in
We shall just give the main ideas of a proof of
the Patchworking theorem. We refer the reader to \cite{Vi}
for the complete proof.

\vskip0.1in

  First notice that from the Lemma \mrf{l:Le1}
  since polynomials $a_1$,...,$a_r$ are completely
  non-degenerate one can consider $\e$-tubular
  neighborhood around points of any hypersurface $a_i=0$ in $\KKK \D_i$
  which belong to $\KKK \G$ with $\G \in \GG'(\D_i)$.

  Besides, recall that to each non-degenerate
  polyhedron $\D$ one can associate projective toric variety
  $\KKK \D$.
  The variety $\KKK \D$ can be seen as  the completion of
  $(\KKK^*)^2$ in such a way that
  the toric varieties associated to the proper faces of $\D$
  cover $\KKK \D \bk (\KKK^*)^2$.

\vskip0.1in

Recall the following statement extracted from \cite{Vi}.

\vskip0.1in
{\bf Lemma 0 }\cite{Vi}
\vskip0.1in
{\it For any polynomial $f$ and $\e >0$ there exists a function \\
$\phi : \GG(\D(f)) \to \RRR $ describing domain of $\e$-sufficiency for
$f$.}
\vskip0.1in
Denote $\GG=\cup_{i=1}^r \GG(\D_i)$.
Define $b(x,y,t)=b_t(x,y)$ the polynomial in the three variables $(x,y,t)$.
Denote $\D"$ the Newton polyhedron of $b$.
$\D"$ is the convex hull of the graph of $\nu$.
For $\G \in \GG$ denote by $\G"$
the face of $\D"$ which is
the graph of $\nu\vert_{\G}$.
For $t>0$, let $j_t$ be the embedding  $\RRR ^2 \to \RRR ^3$ given by
formula $j_t (x,y) = (x,y,lnt)$.
 Let $\psi : \GG \to \RRR $ be a positive function.
For $\G \in \GG$ denote  by ${\mathcal E}_{t,\psi}$ the subset :
$${\mathcal N}_{\psi (\G) }j_t^{-1} (D C_{\D"}^- (\G" )) \bk
\cup_{ \Sigma \in \GG (\D), \G \in \GG (\Sigma)}
{\mathcal N}_{ \phi (\Sigma) } j_t^{-1} (DC_{\D"}^- (\Sigma"))$$
\vskip0.1in
The following Lemma is extracted from \cite{Vi}.
\vskip0.1in
{\bf Lemma 1} \cite{Vi}
{\it Let $a_1$,..$a_r$  be completely non-degenerate polynomials.
For any $\e > 0$, there exists a function
$\psi :\GG \to \RRR $ such that
for any $ t \in ]0,t_0]$ and any face $\G \in \GG$, there exists a
neighborhood   ${\mathcal E}_{t,\psi} (\G)$
of $\G$ such that
$b_t^{\G}$ is $\e$-sufficient for $b_t$ in $l^{-1}{\mathcal E}_{t,\psi} (\G)$.}
\vskip0.1in
{sketch of proof:}\\
It follows from the existence  of a positive
function $\psi : \GG \to \RRR$ such that
for each  proper face $\G"$ of $\D"$,
$b^{\G"}$ is $\e$-sufficient for $b$ in a neighborhood
which contains ${\mathcal E}_{t,\psi} (\G)$ for any $t$ sufficiently small.
Q.E.D.
\vskip0.1in
Let $b_t$ be a polynomial obtained by patchworking the polynomials
$a_1$,...$a_r$.
The Newton polyhedron $\D(b_t)$ of the polynomial
$b_t$ is the polyhedron $\D$.
The following lemma follows immediately from Lemma 0.
\vskip0.1in
{\bf Lemma 2}\cite{Vi}
\vskip0.1in
{\it For any polynomial $b_t$
obtained by patchworking the polynomials
$a_1$,...$a_r$
and $\e >0$ there exists a function
$\phi : \GG(\D) \to \RRR $ describing domain of $\e$-sufficiency for
$b_t$.}
\vskip0.1in

Besides (see \cite{Vi}), for a well chosen function
$\phi : \GG(\D) \to \RRR $,
the $\e$-sufficiency of
$b_t^{\G}$ for $b_t$ in $l^{-1}{\mathcal E}_{t,\psi} (\G)$
for $\G \in \GG(\D)\bk \D$
implies $\e$-sufficiency for $b_t$
in some neighborhood of $l^{-1}(DC_{\D,\phi}(\G))$
and of $l^{-1}(DC_{\D,\phi}(\D))$ in such a way that
for any $t\in ]0,t_0]$
a $\KKK$-chart of the polynomial $b_t$ is obtained by
patchworking $\KKK$-chart of the
polynomials $a_1,...,a_r$.
Hence, it follows the Patchworking Theorem.
\vskip0.1in
Furthermore, one can make the following remark.\\
Since polynomials $a_{i, i \in \{1,...,r\}}$ and $b_t^{\D_i}$
for $t$ sufficiently small
are completely non-degenerate,
from the Lemma 1 and the existence of a well chosen function
$\phi : \GG(\D) \to \RRR $ as above,
the following approximations of $V_{\KKK \D}(b_t)$ can
be easily deduced.
\vskip0.1in
Denote the gradient of the restriction of $\nu$ on $\G$ by
$\nabla (\nu \vert_{\G})$.
The truncation $b_t^{\G}$ equals $a^{\G} \circ
qh_{\nabla (\nu \vert_{\G}),t}$
where $a$ is the unique \pom~  with $\D(a)=\D, a^{\D_i}=a_i$ for
$i=1,...,r$.
From the Lemma 1, we get that
for $t \in ]0,t_0 ]$ the space $\KKK \RRR^2$  is covered by regions in which
$V_{\KKK \RRR^2}(b_t)$ is approximated by
$qh_{\nabla (\nu \vert_{\D_i}),t}^{-1}(V_{\KKK \RRR^2}(a_i))$.
Extending translations
$qh_{\nabla (\nu \vert_{\G}),t}$ to $\G \in \GG(\D_i)$
and $S_{(x',y')}(x,y)$
to the whole space $\KKK \D$,
it follows the description of $V_{\KKK \D} (b_t)$.

\section{Metric on $\CCC P^2$ and $\RRR P^2$}

The projective space $\CCC P^2$
as any differentiable manifold admits a Riemannian Metric.
In what follows, we shall present the Riemannian Metric of $\CCC P^2$
usually called Fubini-Study metric.
\subsection*{Local Euclidian metric}

The topological space $\CCC P^2$,
as complex manifold  looks like locally as the $2$-dimensional affine
space $\CCC^2$.
In other words, $\CCC P^2$ admits a covering by open set
$\CCC U_i=\{ (z_0:z_1:z_2) \in  \CCC P^2 | z_i \not= 0 \}$,
$i \in \{0,1,2 \}$
provided with
charts $\phi_i: \CCC U_i \to \CCC ^2$, $i \in \{0,1,2\}$
($\phi_1 (z_0:z_1:z_2)
= (\frac{z_0}{z_1}, \frac{z_2}{z_1})$;
$\phi_0$ and $\phi_2$ are symmetrically defined), where all the
$\phi_i \circ \phi_j^{-1}$ are complex analytic.

With their help, one can use the  cartesian coordinates in $\CCC^2$
as local coordinates in
$\CCC U_i \subset \CCC P^2$ and carry out complex analysis on $\CCC P^2$.\\

In such a way, we may  locally define
$4$-ball of $\CCC P^2$ in an open $\CCC U_i$ as
$\phi_i^{-1}(B^4) \supset  \CCC U_i$
where $B^4$ is a usual $4$-ball of $\CCC^2$.
Passing from complex to real set of points,
any
$2$-disc of $\RRR P^2 \subset \CCC P^2$  is defined in an open
$\RRR U_i \subset  \CCC U_i$,
as $\phi_i^{-1}(D^2)  \supset  \RRR U_i$
where $D^2$ is a usual $2$-disc of $\RRR^2$.\\
\vskip0.1in

\subsection*{Riemannian metric}

One can also provide $\CCC P^2$ with another
structure by giving a euclidian metric on the tangent space of each point
$p \in \CCC P^2$ that is
by introducing a Riemannian metric
on the complex projective space $\CCC P^2$.
This metric on $\CCC P^2$ is called the
Fubini-Study metric.\\

Let $S^5= \{ (z_0,z_1,z_2) \in \CCC ^3 |
 z_0.\bar{z_0} +
 z_1.\bar{z_1} +
 z_2.\bar{z_2}=1 \}$
and $\pi_{\CCC}: S^5 \to \CCC P^2$  be the natural projection.
For any point $p \in \CCC P^2$ choose a representative $(p_0:p_1:p_2)$
 with
$p_0^2+p_1^2+p_2^2=1$,
$\pi_{\CCC}^{-1} (p)$
is the circle  $\{ p' | p'= \lambda.(p_0,p_1,p_2)~with~ \lambda \in \CCC,
 |\lambda|=1 \}$.

Given two points
$p=(p_0:p_1:p_2) \in \CCC P^2$, $q=(q_0:q_1:q_2) \in \CCC P^2$
,( $p_0^2+p_1^2+p_2^2=1$, $q_0^2+q_1^2+q_2^2=1$),
we define the distance
$\delta$ between $p$ and $q$ in
the Fubini-Study metric as
the length of the geodesic between the circles
$\pi_{\CCC}^{-1}(p)$
$\pi_{\CCC}^{-1}(q)$ of $S^5$.

It is easy to verify that
the Fubini-Study metric
is invariant under the action of complex conjugation.
Besides, when consider
the restriction of $S^5$ to its real set of points  $S^2$,
one define the metric on $\RRR P^2$.\\
Denote
$\pi_{\RRR}: S^2 \to \RRR P^2$ the natural projection.

Given two points
$p=(p_0:p_1:p_2)\in \RRR P^2$,  $q=(q_0:q_1:q_2) \in \RRR P^2$,
($p_0^2+p_1^2+p_2^2=1$, $q_0^2+q_1^2+q_2^2=1$),
$\pi_{\RRR}^{-1}(p)=\{
(p_0,p_1,p_2),
(-p_0,-p_1,-p_2)\}$\\
$\pi_{\RRR}^{-1}(q)=\{ (q_0,q_1,q_2),
(-q_0,-q_1,-q_2) \}$, noticing that the
antipodal mapping $A: S^2 \to S^2$
is an isometry, we define the distance
$\delta$ between $p$ and $q$ in
the Fubini-Study metric as
the minimal arc length of the arcs of $S^2$
between
$(p_0,p_1,p_2) \in \pi_{\RRR}^{-1}(p)$,
and
$\pi_{\RRR}^{-1}(q)=\{ (q_0,q_1,q_2),
(-q_0,-q_1,-q_2) \}$.\\
Given two real points of $S^2 \subset S^5$,
it is easy to see (since geodesic minimizes arc length
between two of its points)
that any point of the geodesic
from these two points is real.\\
In such a way,
in the Fubini-Study metric of $\CCC P^2$,
we define
any $4$-ball of $\CCC P^2$ by
 $\pi_{\CCC} ({\mathcal E})$
 where ${\mathcal E}$ is an
an ellipsoid of $S^5$; and
any $2$-disc of $\RRR P^2 \subset \CCC P^2$ by
$\pi_{\CCC}(E)=\pi_{\RRR} (E)$ where $E$ is an ellipse of
$S^2 \subset S^5$.\\

\part{\bf Arnold Surfaces of Harnack Curves}

The maximal number \cite{Har} of connected components of
the real point set of curves of degree $m$ is $ \frac {(m-1)(m-2)}{2} +1$.
Curves with this maximal number are called $M$-curves.
In this section, we study the construction of
some $M$-curves called Harnack curves and prove that Arnold surfaces
of the so-called Harnack curves of even degree are standard.
\vskip0.1in
\chapter{Combinatorial patchworking construction for Harnack curves}
\lb{ch:PatchHar}

{\bf Preliminaries }
Recall that the pair $(\RRR P^2,\RRR \AA)$ where $\AA$
is a non-singular real plane curve
is determined up to homeomorphism by the real components of the curve and
their relative location.
\vskip0.1in
In the case of even degree curve $\AA$,
each connected  component of the
real point set $\RRR \AA$ is situated in $\RRR P^2$
as the boundary of an embedded disc and is called
an {\it oval}.
\vskip0.1in
 In the case of odd degree curve $\AA$, the
 real point set $\RRR \AA$  has besides oval  one connected component
 situated in $\RRR P^2$
 as an embedded projective line. It is called the
 {\it one side component} of $\RRR \AA$.
An oval divides $\RRR P^2$ into two components. The orientable component
(i.e the component homeomorphic to a disc)
is called the inside of the oval.
The non-orientable component (i.e the component
homeomorphic to a M\"{o}bius strip)
is called the outside of the oval.
\vskip0.1in
Given $\AA$ a non-singular real plane curve
and $F(x_0,x_1,x_2)$ its polynomial, denote by $\RRR P^2_+$
the subset of $\RRR P^2$ $\{ x \in \RRR P^2~| F(x_0,x_1,x_2) \ge 0 \}$.
An oval is said {\it outer}  (resp {\it inner})
if it bounds a component
of $\RRR P^2_+$ from the outside (resp, the inside).
In the case of an even degree curve $\AA$,
we can always assume (changing the
sign of the polynomial $F(x_0,x_1,x_2)$ giving $\AA$ if necessary) that
\break $\RRR P^2_+ =\{ x \in \RRR P^2 ~| F(x_0,x_1,x_2) \ge 0 \}$
is the orientable component of $\RRR P^2$.
In such away, in case of an even degree curve,
 since ovals lying in an even number of consecutive ovals are
outer while
the ovals lying in an odd number of consecutive ovals are inner;
one also calls {\it even} ovals the
 outer ovals  and {\it odd} ovals the inner ovals.

In the case of an odd degree curve $\AA$, the definition of outer and
inner oval apply only to oval which does not intersect the line at
infinity. Call {\it zero oval} an oval of a curve of odd degree intersecting
the line at infinity.

\vskip0.1in
The pair $(\RRR P^2,\RRR \AA)$ where$\AA$
is a  non singular curve of degree $m$
is defined by the scheme of disposition of the real component
$\RRR \AA$.
This scheme is called the {\it real scheme of the curve} $\AA$.
In what follows, we use the following usual system of notations
for real scheme (see \cite{Vi}).
A single oval is denoted by $\big \langle 1 \big \rangle$.
The one side component is denoted by  $\big \langle J \big \rangle$.
The empty real scheme is denoted by $\big \langle 0 \big \rangle$.
If one set of ovals is denoted by $\big \langle A \big \rangle$,
then the set of ovals obtained by addition of an oval which contained
$\big \langle A \big \rangle$ in its inside component
is denoted by
   $\big \langle 1\big \langle A \big \rangle \big \rangle$.
If the real scheme of a curve consists of two disjoint sets of ovals
denoted by $\big \langle A \big \rangle$ and$\big \langle B \big \rangle$
in such a way that no oval of one set contains an oval of the other set in
its inside component, then the real scheme of this curve is denoted by :
 $\big \langle A \sqcup B \big \rangle$.
If one set of ovals is denoted by $A$, then the set $A \cup ...\cup A$
where $A$ occurs n times, is denoted by $\big \langle n \ti A \big \rangle$;
a set $\big \langle n \ti 1 \big \rangle $ is denoted
by $\big \langle n \big \rangle $.

\vskip0.1in
Furthermore,
if $\AA$ is a curve of which the set of real points
$\RRR \AA$
divides the set of complex points $\CCC \AA$
in two connected components
( which induces two opposite orientations on $\RRR \AA$),
then the curve $\AA$
is said of type $I$. Otherwise, $\AA$ is said of type $II$.
Hence, the real scheme of a curve $\AA$ of degree $m$ is of type $I$
(resp, of type $II$), if any curve of degree $m$ having this real scheme is
of type $I$(resp, of type $II$). Otherwise (i.e if there exists both curves
of type $I$ and curves of type $II$ with the given real scheme), we say that
the real scheme is of {\it indeterminate type}.
\vskip0.in

Thus, isotopy of $\CCC P^2$ which
connects complex points of curves
and commutes with
the complex conjugation of $\CCC P^2$
, called {\it conj-equivariant} isotopy of $\CCC P^2$,
provides a convenient way to classify  pairs $(\CCC P^2,\CCC \AA)$.
We can already notice that
since $M$-curves
are curves of type $I$,
their real scheme is sufficient to classify
their complex set of points in $\CCC P^2$ up to
conj-equivariant isotopy.

\section{ Constructing Harnack curves with Harnack's initial method}
\mlb{su:init}
\vskip0.1in
\vskip0.1in

In 1876, Harnack proposed a method for constructing
$M$-curves.
Recall briefly Harnack's initial
construction \cite{Har} of Harnack's $M$-curves.
(The detailed method can be found in \cite{Gu})

Start with a line as  $M$-curve of degree $1$.
Then, consider a line $L$ which intersects it in one point.
Assume an $M$-curve of degree $m$
$\AA_m$
has been constructed at the step $m$ of the construction.
At the step $m+1$ of the construction, the  $M$-curve
of degree $m+1$
is obtained from classical deformation
(see \cite{Vi} Classical Small Perturbation Theorem
 for the definition of classical deformation)
of $\AA_m \cup L$.
The resulting $M$-curve of
degree $m+1$ intersects the line $L$ in $m+1$ real points.
\vskip0.1in

Passing from curves to polynomials giving these curves,
an $M$-curve $\tilde{X}_{m+1}=0$ is
constructed from an $M$-curve
$\tilde{X}_{m}=0$ of degree $m$
by the formula
$\tilde{X}_{m+1}=x_0. \tilde{X}_m +t.C_{m+1}$
where  $x_0$ is a line,
$C_{m+1}=0$ is $m+1$ parallels lines which intersect $x_0=0$
in $m+1$ points.
It is essential in the construction that the curves
$\tilde{X}_{m}=0$, $x_0=0$ and $C_{m+1}=0$
do not have common points.

Thus, one can choose projective coordinates $(x_0:x_1:x_2)$
of $\CCC P^2$, in such a way that
$\tilde{X}_{m}=0$
is an homogeneous polynomial of degree $m$
in the variables $x_0$, $x_1$,$x_2$ and
$C_{m+1}=0$ is
an homogeneous polynomial of degree $m+1$
in the variables $x_2$,$x_1$.
For sufficiently small $t>0$,
$\tilde{X}_{m+1}=0$ is an Harnack curve of degree $m+1$.

In particular,
this method provides curves with real scheme:\\
for even $m=2k$
\begin{equation}
\lb{e:defhar2}
\big\langle1\langle \frac {(k-1)(k-2)}2 \rangle \sqcup \frac {3k(k-1)}2
\big\rangle
\end{equation}

for odd $m=2k+1$
\begin{equation}
\lb{e:defhar1}
\big\langle J \sqcup k(2k-1)
\big\rangle
\end{equation}

Let us denote $\HH_m$ and call {\it  Harnack curve of degree $m$}, any curve
of degree $m$ with the real scheme above.
Besides, we shall
call {\it Harnack polynomial} of degree $m$ any polynomial
giving the Harnack curve of degree $m$.
\vskip0.1in
We shall call curve
of type $\HH$  a curve $\HH_m$
which intersects the line at infinity of $\CCC P^2$
in $m$ real distinct points.
Without loss of generality, one can always assume that the line $L$  is the
line at infinity.
Therefore, curves $\HH_m$ constructed by the Harnack's method
are curves of type $\HH$.

\section{ Patchworking construction for Harnack curves}
\mlb{su:patch}
\vskip0.1in
\vskip0.1in
In what follows, we shall give a construction of
Harnack curves provided by the patchworking construction of curves .
\vskip0.1in

Recall briefly the patchworking construction procedure
of curves due to Viro \cite{Vi}. (We refer the reader to \cite{Vi}
and also to the Preliminary section for details.)

{\bf Initial data}
\been
\item
 Let $m$ be a positive integer.
 Let $T_m$ be the triangle
 $$
\{(x,y)\in\RRR^2\mid x\geq 0, y\geq 0, x+y\leq m\}
$$
(Up to  linear change of coordinates $(x_0:x_1:x_2)$
of $\CCC P^2$,  the convex hull of $T_m$ may be the Newton polyhedron
of the affine polynomial $f(x,y)= F(1,x,y)$
associated to a homogeneous polynomial $F(x_0,x_1,x_2)$ of degree $m$.)

\item
Let $\tau$ be a triangulation of $T_m$ whose vertices have integer
coordinates.
Call regular triangulation of $T_m$
a triangulation of $T_m$
such that there exists with a convexifying function $\nu : T_m \to \RRR$
(that is a
 piecewise-linear function  $\nu: T_m \to \RRR$ which is linear on
 each triangle
 of the triangulation $\tau$ and not linear on the union of two triangles.)
Assume $\tau$ regular.
\item
 Suppose that some distribution of signs $\chi$  at the
vertices of the triangulation
is given.
 Denote the sign $\pm$ at the vertex with coordinates $(i,j)$ by
 $\epsilon_{(i,j)} $.

\enen

{\bf Combinatorial procedure}

Take the square $T_m^*$ made of $T_m$
and its symmetric copies
$T_m^x=s_x(T_m),
T_m^y=s_y(T_m),
T_m^{xy}=s_{xy}(T_m)$
where
$s_x,s_y, s=s_x\circ s_y$ are reflections with respect to
the coordinate axes. The resulting space is homeomorphic to $\RRR^2$.
Extend the triangulation $\tau$ of $T_m$ to a symmetric triangulations$\tau_*$
of $T_m^*$.
Extend the distribution of signs to a distribution at the vertices
of $T_m^*$ which verifies the modular properties:

$g^*(\e_{i,j})x^iy^j=
\e_{g(i,j)}x^iy^j$
for
$g=s_x, s_y,s$.

If a triangle of $T_m^*$
has vertices of different signs, consider a midline separating them.
Denote
by $L$ the union of such midlines.
Glue by $s$ the opposite sides of $T_m^*$.
The resulting space $\bar T_m^*$ is homeomorphic
to the projective plane $\RRR P^2$.
Denote $\bar L$ the image of $L$ in $\bar T_m^*$.
We shall say that
$T_m^*$, $L$,
$\bar T_m^*$,$\bar L$
are obtained from $T_m,\tau,\chi$ by combinatorial patchworking.

\beth
\mlb{t:Harpatch}
{\bf Polynomial Patchworking}
{\it Define the one-parameter family of polynomials
$$b_t=b_t(x,y)= \sum_{ (i,j) vertices of \tau} \e_{i,j}x^iy^jt^{\nu(i,j)}$$
where $\nu$ is a function convexifying the triangulation $\tau$
of $T_m$.
Let
$\bar T_m^*$,$\bar L$
obtained from $T_m,\tau,\chi$ by combinatorial patchworking.

Denote by $B_t=B_t(x_0,x_1,x_2)$ the corresponding homogeneous polynomials:
$$B_t(x_0,x_1,x_2)=x_0^mb_t(x_1/x_0,x_2/x_0)$$

\vskip0.1in
Then there exists $t_0 > 0$ such that for any $t \in ]0,t_0]$
the equation $B_t(x_0,x_1,x_2)=0$ defines in $\RRR P^2$
the set of real points of an algebraic curve $C_t$ such that the pair
$(\RRR P^2, \RRR C_t)$
is homeomorphic to the pair $(\bar T_m^*,\bar L)$.
Such a curve is called a {\it $T$-curve }.}
\enth

\newpage
Let $m$ be a positive integer.

 Let $T_m$ be the
 triangle
$$
\{(x,y)\in\RRR^2\mid x\geq 0, y\geq 0, x+y\leq m\}
$$
\vskip0.1in

 Define the regular triangulation $\tau$ of $T_m$ as follows:
 \been
 \item
 All integer points of $T_m$ are vertices of the triangulation $\tau$.
 \item
 Each proper face of the triangulation is contained in one segment of
 the set:
$$
\{(x+y=i)_{i=1,...,m}, (x=i, y \le m-i)_{i=1,...,m},
(y=i, x\le m-i)_{i=1,...,m} \}
$$
  \enen

Choose the following distribution of signs called " Harnack distribution"
at each integer points $(i,j)$ of $T_m$ :
\been\item
If i,j are both even, the integer point $(i,j)$ gets the sign
$\e_{(i,j)}= -$
\item
\vskip0.1in
If $i$ or $j$ is odd, the integer point $(i,j)$ gets the sign
$\e_{(i,j)}=+$.
\enen

\vskip0.1in
The following propositions are deduced from a more general statement
due to I.Itenberg \cite{It}.
\vskip0.1in
\bepr
\cite{It}
{\it Let $m=2k$ be a positive even integer.
The patchworking process applied to
$T_m$ with regular triangulation $\tau$
and the Harnack distribution of signs
at the vertices produces a $T$-curve which is
the Harnack curve $\HH_{2k}$ of degree
$2k$. }
\enpr

\vskip0.1in

\bepr\cite{It}
\mlb{p:patchodd}

{\it Let $m=2k+1$ be a positive even integer.
The patchworking process applied to $T_m$
with regular triangulation $\tau$ and the Harnack distribution of signs
at the vertices produces a $T$-curve which is
the Harnack curve
 $\HH_{2k+1}$ of degree $2k+1$ .}
\enpr

\vskip0.1in

Now we shall work out the polynomial entering in this patchworking
construction (following essentially Viro's initial method
with some modifications necessary for the sequel).
We want in particular to stress the recursive character of the patchworking
construction of Harnack' s curves.

Fix $m>0$ integer
\been
\item
An essential element in
the construction is the so-called
convexifying function.
We choose the
convexifying function
 $\nu_m :T_m \to \RRR$
in such a way that it would be the restriction of $\nu_{m+1}$.
Namely,
we extend the convexifying function
$\nu_m :T_m \to \RRR$
to a convexifying function  $\nu_{m+1}$
of the triangulation of $T_{m+1}= T_{m} \cup D_{m+1,m}$ with
$$
D_{m+1,m}=
\{(x,y)\in\RRR^2\mid 0 \leq x\leq m, 0 \leq y\leq m, x+y\leq m+1\}
$$

We shall assume the same notation for $\nu_{m}$ and its extension
$\nu_{m+1}$,
and for any $m>1$ denote by $\nu$ a convexifying function of the triangulation
of $T_m$.

\item
Now, we construct a polynomial in two variables

$$x_{m,t} (x,y)= \sum_{(i,j) vertices of T_m} \e_{i,j}a_{i,j}
x^iy^jt^{\nu(i,j)}$$

where:
$\epsilon_{(i,j)} $ is the sign at the vertex $(i,j)$
given by the Harnack distribution of sign,
$\nu(i,j)$ is the value of the convexifying function
$\nu:T_m \to \RRR$ on $(i,j)$.
\vskip0.1in
The numbers $t$ and $a_{i,j}$ are parameters which shall be chosen.
We shall denote the polynomial $x_{m,t} (x,y)$
by $x_{m(x,y;\vec{a}_m,t)}$ or $x_{m;\vec{a}_m,t}$.
Let $\tilde{X}_{m,t}(x_0,x_1,x_2)$ be the homogenization of $x_{m,t}$
obtained by remaining $x$ as $x_1$, $y$ as $x_2$, and adding the
third variables $x_0$  so as to obtain :

$$\tilde{X}_{m,t} (x_0,x_1,x_2)= \sum_{(i,j)~ vertices~ of~ T_m}
\e_{i,j}a_{i,j}
x_1^ix_2^jx_0^{m-i-j} t^{\nu(i,j)}$$

We shall denote the polynomial $\tilde{X}_{m,t} (x_0,x_1,x_2)$
by $\tilde{X}_{m(x_0,x_1,x_2;\vec{a}_m,t)}$
or $\tilde{X}_{m;\vec{a}_m,t}$.

Note that
 $$\tilde{X}_{m,t}= x_0.\tilde{X}_{m-1,t} + C_{m,t}$$
where the polynomial $C_{m,t}$ is an homogeneous
polynomial of degree in only $x_1$ and $x_2$.
\enen

\vskip0.1in
Using the patchworking method,
one can give an inductive construction, we shall call
{\it T-inductive construction}
of some Harnack polynomials.

Assume $\vec{a}_{m},t_{m}$ given in such a way
that for any $t \in ]0, t_{m}[$,
$\tilde{X}_{m;\vec{a}_{m},t}$ is a Harnack polynomial of degree $m$.
Then,
the Harnack polynomial of degree $m+1$ is constructed as follows:

{\it
\been
\item
{\bf Property 1 }
Choose $\vec{c}_{m+1}$ in such a way that
for any $t \in ]0,t_{m}[$,
singularities of $x_0. \tilde{X}_{m;\vec{a}_m,t}$
do not belong to the curve ${\mathcal C}_ {m+1;\vec{c}_{m+1},t}$.
\item
{\bf Property 2 }
Set $\vec{a}_{m+1}=(\vec{a}_{m},\vec{c}_{m+1})$.
Then,
$t_{m+1}$ is chosen as follows:
\been
\item
$t_{m+1} \le t_{m}$

\item
$t_{m+1}$ is the biggest $\tau \le t_{m}$
 such that for any $t \in ]0,\tau[$
 the polynomial
 $\tilde{X}_{m+1;\vec{a}_{m+1},\tau}=
  x_0.\tilde{X}_{m;\vec{a}_m, \tau}+ C_{m+1;\vec{c}_{m+1},\tau}$  is
 a Harnack polynomial of degree $m+1$.
\enen
\enen}

Following the preliminary section, consider
$\rho^{m} : \RRR_+ T_{m} \ti U_{\CCC}^2 \to \CCC T_{m}$ the
natural surjection.
Denote $l_m=\{ x \ge 0, y \ge 0, x+y=m \}$
the hypotenuse  of $T_m$.
The variety $\CCC T_m$
is homeomorphic to $\CCC P^2$
with line at infinity $L$ of $\CCC P^2$
such that $\CCC l_m \approx L$.

(The real line $\RRR L \subset \RRR P^2$
is obtained from the square
$\{(x,y)\mid \vert x\vert +\vert y \vert =m\}$
by gluing of the sides by $s$ the central symmetry
with center $0=(0,0)$).
\vskip0.1in

From the Patchworking construction of Harnack curves
it follows that:
\begin{coro}
\mlb{c:patch}
\vskip0.1in
{\it
 The Harnack curve $\HH _m$ of degree $m$
 obtained by patchworking process applied to $T_m$ with triangulation
 $\tau$ and Harnack distribution of signs at the vertices intersects
 the line at infinity
 $L \approx \CCC l_m$ of $ \CCC P^2 \approx \CCC T_m$
 in $m$ real points.
 Thus, the resulting curve $\HH_m$
  is a curve of type $\HH$.}
\vskip0.1in
 Let us denote $a_1,..., a_m$ the real points of the intersection
 $$\RRR \HH _m \cap \RRR L =\{ a_1,...,a_m\} $$
\end{coro}

\vskip0.1in
{\bf proof:}
\vskip0.1in
First of all, recall that in the patchworking construction
of $\HH_m$, the real line
at infinity  $\RRR L$ of $\RRR P^2$ is obtained from the
square $\{ (x,y) \in \RRR^2 |  |x| + |y| =m \} $ by gluing
of the sides by $s$  the central symmetry with center $0=(0,0)$.
Furthermore, according to the patchworking construction,
the real point set of the Harnack curve separates any two consecutive
integer points of the
square $\{ (x,y) \in \RRR^2 |  |x| + |y| =m \} $ with different signs.
Thus, the corollary follows immediately from the
Harnack distribution of signs on $T_m$ and its extension to its symmetric
copies.Q.E.D
\vskip0.1in

Moreover, one can make the following remark:

\bere
\mlb{r:ori}
  Harnack curves are curves of type $I$, in other words
  their real set of points are orientable.
  Let $\HH _m$ be the Harnack curve of degree $m$.
  Assume $\RRR \HH_m$ given with an orientation.
  There is only one orientation of $\RRR L$ compatible
  with the deformation of
  $\HH_m \cup L$ to $\HH_{m+1}$ and an orientation
  of $\RRR \HH _{m+1}$.
  In such a way, we obtain the relative orientation
  of the connected components of the real point set of the Harnack curves.
\enre
\vskip0.1in\chapter{Morse-Petrovskii's Theory of Harnack curves}
\lb{ch:Mope}
In this chapter, at first we shall  study critical points of
Harnack polynomials from a view-point
initiated by Petrovskii. This investigation is analogous to those of
Morse on critical points of  functions.\\
Then,
using rigid isotopy classification of Harnack curves,
we shall construct deformation of Harnack polynomial and
deduce a description up to conj-equivariant
isotopy of $\CCC P^2$ of the complex
set of points of Harnack curve (see proposition \mrf{p:prop8}
and theorem \mrf{t:theo1}).\\
We shall divide this chapter into three sections.\\
In the first section,
we start with proposition  \mrf{p:crit} in which we prove
that up to slightly modify coefficients of a Harnack polynomial,
the number of its critical points depends only on its degree.
Then, in proposition \mrf{p:crit2}, we
precise, up to change the
system of projective coordinates if necessary,
the number of critical points of index $1$ with positive
and respectively negative critical value of such modified Harnack polynomials.
Such results may be obtained from the Viro's patchworking method for
$T$-Harnack polynomials,
it is explicitely explained in Shustin's paper \cite{Shu}.\\
In the second section, we consider only real points set of Harnack curves.
In Theorem \ref{t:rigiso}, we prove that isotopy also implies rigid isotopy
for Harnack curves $\HH_m$ of degree $m$.
In the third section, we define
in Proposition \mrf{p:prop7},
a deformation of $\HH_m$,
to a singular irreducible curve
of which  singular points are  critical points of index 1
with positive critical value of $\HH_m$.\\
The main result of this chapter is given
in the third and last section.
It provides link between properties of $L$-curves (\cite{Fi})
and the Harnack curves.
In Theorem \mrf{t:theo1}, we present any Harnack curve $\HH_m$
of degree $m$  up to conj-equivariant isotopy of  $\CCC P^2$
as follows.
Outside
a finite number of $4$-balls $B(a_i)$ globally invariant
by complex conjugation, $\HH_m$ splits in
$m$ non-intersecting projective lines minus their intersection
with the $B(a_i)$;
inside any $4$-ball $B(a_i)$
it is the perturbation of a crossing.

\section{ Critical points of Harnack polynomials}
\mlb{s:top}
Recall  a Harnack polynomial  of degree $m$
is a non-singular homogeneous polynomial in three variables
such that its set of zeros has real scheme:
\vskip0.1in
for even $m=2k$:
$$\big\langle1\langle \frac {(k-1)(k-2)}2 \rangle \sqcup \frac {3k(k-1)}2
\big\rangle$$
for odd $m=2k+1$:
$$\big\langle J \sqcup k(2k-1)
\big\rangle$$
In this section, we study the critical points
and the level surfaces of Harnack polynomials. \quad

Given $x_0:=0$ the line at infinity and
$R(x_0,x_1,x_2)$ an homogeneous polynomial of degree $m$,
we call affine restriction  $r(x,y)$ of $R(x_0,x_1,x_2)$ the unique
polynomial such that $R(x_0,x_1,x_2)=x_0^m.r(x_1/x_0,x_2/x_0)$.
In such a way, we consider $\CCC P^2$ as the completion
of $\CCC^2$ with $x_0:=0$ the line at infinity.
We call critical point of $r(x,y)$,(and by misuse critical point
of $R(x_0,x_1,x_2)$),  any point $(x_0,y_0)$
(finite or not) such that
$r_{x}(x_0,y_0)=0$, $r_{y}(x_0,y_0)=0$.
\vskip0.1in

We shall consider generic polynomials.
To be more precise, let us call an homogeneous non-singular real polynomial
$R(x_0,x_1,x_2)$ with affine  restriction $r(x,y)$
{\it regular} if :
\been
\item
none of critical points of $r(x,y)$
lies on the line at infinity
\item
any two different real critical points $r(x,y)$ have different critical
values.
\enen
From (\cite{Petr}, p.359 Lemma 2),
given a polynomial of a smooth algebraic curve,
one can always perturb its coefficients
so as to have $1)$ and $2)$ as above.

Bringing together  Lemma \mrf{l:le11} and Petrovskii's Lemma
(\cite{Petr}, p.359 Lemma 2),
it follows that the set of regular polynomials is open and dense in the
set of all polynomials.
\vskip0.1in

We shall call {\it polynomial of type $\HH$}
a regular Harnack polynomial of degree $m$  giving  a Harnack curve $\HH_m$
 which intersects the
line at infinity of $\RRR P^2$ in $m$  real distincts points.
We shall call {\it polynomial of type $\HH^0$} a polynomial of type $\HH$
such that the line at infinity of $\RRR P^2$ intersects the
Harnack curve $\HH_m$ of degree $m$
in $m$ real distinct points which
belong to the same connected component of $\RRR \HH_m$: in case
$m=2k+1$  the real connected component of $\HH_{2k+1}$
homeomorphic to the projective line;
in case $m=2k$ the non-empty oval of $\HH_{2k}$.

\vskip0.1in

Let  us start by the proposition \mrf{p:crit}
in which we find the number of critical points of all indices of any regular
Harnack polynomial of degree $m$.

For a regular polynomial of homogeneous degree $m$,
$R(x_0,x_1,x_2)$ with affine restriction $r(x,y)$
,$R(x_0,x_1,x_2)=x_0^m.r(x_1/x_0,x_2/x_0)$,
let us denote by $c_i(R)$ the number
critical points of index $i$ of $r(x,y)$ and set
 $c(R)=(c_0(R),c_1(R),c_2(R))$.

\vskip0.1in
\bepr
\mlb{p:crit}
{\it
Let $B_m(x_0,x_1,x_2)$ be a  Harnack polynomial
of degree $m$ of type $\HH^0$.\\
Then, up to change the sign of $B_m$:
\been
\item
$c_0(B_m)+c_1(B_m)+c_2(B_m)=(m-1)^2$
\item
\been
\item
For even  $m=2k$,
$$c(B_{2k})=\big( \frac{(k-1)(k-2)}{2},k(2k-1),\frac {3k(k-1)}{2} \big)$$
\item
For odd $m=2k+1$,
$$c(B_{2k+1})=\big( \frac {k(k-1)}{2},k(2k+1), \frac{k(3k-1)}{2} \big)$$
\enen
\enen}
\enpr
\vskip0.1in
{\bf proof:}
\vskip0.1in
Let $b_m(x,y)$ be the affine restriction of the polynomial $B_m(x_0,x_1,x_2)$.
For any oval which does
not intersect the line at infinity,
consider the disc with boundary the oval.
In case of even oval, the gradient of the affine polynomial $b_m(x,y)$
points inward.
Therefore, the maximum of the polynomial $b_m$ on the disc is in the interior
of the oval and is necessary a critical point of index $2$ of $b_m$.

In case of odd oval, the gradient of the affine polynomial $b_m(x,y)$
points outward.
Therefore, the minimum of the polynomial $b_m$ on the disc is in the interior
of the oval and is necessarily a critical point of index $0$ of $b_m$.
This implies the following inequalities (up to change the sign of $B_m$
and in case of even degree curve according to the convention.
-ovals not lying within other ovals or lying inside an even number
number of consecutive ovals are outer, while ovals lying within an odd
number of ovals are inner-)

\vskip0.1in

for even $m=2k$:

$$c_0 \ge  \frac{(k-1)(k-2)}{2}, \quad
c_2 \ge \frac {3k(k-1)}{2}$$

for odd $m=2k+1$:

$$c_0 \ge\frac {k(k-1)}{2}, \quad
c_2 \ge  \frac{k(3k-1)}{2}$$

Furthermore, since $\HH_m$  intersects the line at infinity
in $m$ real points we have $c_0 -c_1 + c_2 = 1 -m$ (see \cite {Du}).
Thus, it follows  from the inequalities above that
for even $m=2k$ :
\begin{equation}
\lb{e:c2}
c_1 \ge k(2k-1)
\end{equation}
for odd $m=2k+1$ :
\begin{equation}
\lb{e:c1}
c_1 \ge k(2k+1)
\end{equation}

and thus
$c_0+c_1+c_2 \ge (m -1)^2$.

From the Bezout's theorem, the number of non-degenerate critical
points of a homogeneous polynomial of degree $m$ does not
exceed $(m-1)^2$.
Therefore $c_0 +c_1 +c_2 = (m-1)^2$ and the inequalities \ref{e:c2}
and \ref{e:c1} are equalities.

Q.E.D.
\vskip0.1in
\vskip0.1in

Our next aim is to
 find the number of  critical points of indice 1 with
negative (resp, positive) critical value of a any polynomial
of type $\HH^0$. This will be done in Proposition \mrf{p:crit2}.\\
For a regular polynomial, let us denote by $c_1^-(R)$  (resp, $c_1^+(R)$)
the number of critical points of index 1 of $R$
of negative (resp, positive) critical value.
Besides,
denote by $c'_1(R)$ the number of critical points of
index 1 of $R$ with positive critical value $c_0$
such that as $c$ increases from
$c_0-\e$ to $c_0+\e$ the number of  connected components of
$M_c=
\{(x,y)\in\RRR^2\mid
r(x,y)>c\}$ of which boundary contains
the line at infinity increases by $1$.\quad
We shall denote $S'_{2k}$ the set constituted by
the $c_1'(B_{2k})$ critical points of
index 1 of a polynomial $B_{2k}$ of
type $\HH^0$.

\vskip0.1in
Our proof uses Petrovskii's theory introduced in \cite{Petr}.
It is based  on the consideration of the lines $b_m(x,y)=c$ when
$c$ crosses the critical value of a polynomial $B_m$ of type $\HH^0$.

These last investigations are analogous to Morse theory.

Let us recall Petrovskii's Lemmas
implicit in \cite{Petr}\\
(p.360-361).
\vskip0.1in
First, introduce definitions.

Define $\RRR P^2$ as a M\"{o}bius strip and a disc $D^2$ glued along their
boundaries.
  The core of the M\"{o}bius strip intersects the line at infinity in
  a finite number of points.
  Call "odd component  of the affine plane $\RRR ^2$" the core of
  the M\"{o}bius strip minus its intersection with the line
  at infinity.

\vskip0.1in

\bele \cite{Petr}
\mlb{l:lp2}
Let $A$ a curve of degree $m$
given by a regular polynomial $R(x_0,x_1,x_2)=x_0^mr(x_1/x_0,x_2/x_0)$
with $(m-1)^2$ critical points.
  Let $(x_0,y_0)$ be a finite real critical point of $r(x,y)$.
  Setting $r(x_0,y_0)=c_0$, assume  that $r^{-1}[c_0- \e,c_0 + \e]$
  contains no critical point other than $(x_0,y_0)$.

   Assume the Hessian $H(x_0,y_0) < 0$, (i.e $(x_0,y_0)$ is a
   critical point of index $1$).
 \vskip0.1in
 (*) Assume $r(x,y)$ is of even degree,
 and $c$ decreases from $c_0+ \e $ to $c_0 -\e$.
 Then there are three possibilities:

 \been
 \item
     When $c$ takes value
    $c_0$, one outer oval touches another oval
     (outer or inner),
     then one outer oval disappears.
 \item
     When $c$ takes value  $c_0$,
     one oval (outer or inner) touches itself.

 \been\item
 \lb{e:ep2}

          If one of the components
 $r(x,y) >r(x_0,y_0) + \e $ ;
 $r(x,y) <r(x_0,y_0) - \e $
 contains the odd component of $\RRR ^2$,
 then one inner oval appears.
 \item
  Otherwise,
  one outer oval becomes inner.
  (This case may present itself at most once
   as $c$ varies from $+ \infty$
  to $- \infty$.)
\enen
\enen
 (*) Assume r(x,y) is of odd degree and $c$ decreases from $c_0+ \e $ to
 $c_0 -\e$.
 Then there are three possibilities:
\been
\item
  When $c$ takes value $c_0$,
  one outer oval touches an other oval
 or the one-side component,
  then the outer oval disappears.

 \item
 When $c$ takes value  $c_0$,
 one oval (outer or inner) touches itself
  then one inner oval appears.
 \item
 \lb{e:ep3}
 When $c$ takes value  $c_0$,
 one zero oval or the odd-component of the curve
 touches itself or an other zero oval,
  then there are three possibilities:

   Consider the regions $G_i$ $i=1,...,r$ where $r(x,y) >c_0$ which contain
  segment of the line at infinity of $\RRR P^2$ on their boundary.
  \been
  \item
  the boundary
 of one region $G_i$ touches itself
  (in one point  which does not belong to the line at infinity),
  then one inner oval appears.

  \item
  one region $G_i$ touches the boundary
    (not contained in the line at infinity)  of
  another such regions $G_j$
    then the regions coalesce and
   a zero-oval disappears or a zero oval appears.
\enen
 \enen
\enle

The next Lemma
is implicit in(\cite{Petr},
Lemma 4, p.361)\\

\bele
\mlb{l:lp3}
(line at infinity (\cite{Petr})\label{lininf}
{\it
Let $A$ a curve of degree $m$
given by a regular polynomial $R(x_0,x_1,x_2)=x_0^m.r(x_1/x_0,x_2/x_0)$.
with $(m-1)^2$ critical points.

Assume furthermore that $A$ intersects the real line at infinity
in $k$ $(k>1)$ distinct points.
Let
$M_c=
\{(x,y)\in\RRR^2\mid
r(x,y)>c\}$
Let
$c_M$ be the maximal critical value, and
$c_m$ be the minimal critical value.
\been
\item
For $c>c_M$, $M_c$ has $k$ connected components.
\item
For $c<c_m$, $M_c$
is simply connected.
Two of the connected components above
can coalesce only when
$r(x,y)$ passes a critical value $c$ of index 1.
\enen}
\enle
\vskip0.1in

\bepr
\mlb{p:crit2}
{\it
Let $m \ge 0$ and $B_m(x_0,x_1,x_2)$ be a Harnack polynomial
of type $\HH^0$.
Then, up to change the sign of $B_m$:

\been
\item
for even $m=2k$;

$$c_1^-(B_{2k})  = \frac {k (3k-1)}{2}$$
$$c_1^+(B_{2k})=  \frac{ k(k-1)}{2},\quad c'_1(B_{2k})=k-1$$

\item
for odd $m=2k+1$
$$c_1^- (B_{2k+1}) = \frac{ k(3k+3)}{2}$$
$$c_1^+(B_{2k+1}) = \frac{ k(k-1)}{2}, \quad c'_1(B_{2k+1})=0$$

\enen}
\enpr

\vskip0.1in

{\bf proof :}

The proof of Proposition \mrf{p:crit2} is based on
the Lemma \mrf{l:lp2} and Lemma \mrf{l:lp3}.

Let  $\HH _m$ be the Harnack curve of degree $m$ given by
the polynomial $B_{m} (x_0,x_1,x_2)= x_0^mb_{m}(x_1/x_0,x_2/x_0)$
of type $\HH^0$.
Consider the pencil
of curves given by polynomials
$x_0^m(b_m(x_1/x_0,x_2/x_0)-c)$ with $c \in \RRR$.

\been
\item
Let $c$ decrease from $0$ to $c_m -\e$, $\e>0$, the inner ovals
shrink and disappear, the outer disappear in
a unique oval. It follows from the Lemma \mrf{l:lp3} that the set
$M_{c_m - \e}$ is simply connected.
Furthermore, from Lemma \mrf{l:lp2} :
one  outer oval can touch another outer oval
and then disappears as
$c$ decreases from $c_0+ \e$ to $c_0 -\e$, if and only if
$c_0$ is a critical value of index $1$ of the polynomial
$b_m$.
\vskip0.1in
(*)Let $\HH_{2k}$ be the Harnack curve of even degree $ 2k$

When $c$ decreases from $0$ to $c_m - \e$, outer ovals expand, then
coalesce and finally disappear in a unique oval.
Moreover, from Lemma \mrf{l:lp2},
to each  coalescence of a outer oval ${\mathcal O}$
with an other outer oval
is associated
a critical point $(x_0,y_0)$ of index 1  with
$b_{2k}(x_0,y_0)=c_0$; $c_0 \in ]c_m, 0[$ such that as $c$ decreases from
$c_0+\e$ to $c_0 -\e$
the oval ${\mathcal O}$ disappears.
Denote $S_{2k}^-$ the set of these critical points.
The curve $\HH_{2k}$ has $\frac {3k (k-1)} {2}$ outer
ovals and $(k+1)$ connected components of
$M_0= \{ (x,y) \in \RRR^2 | b_{2k}(x,y) > 0 \}$ contain a
segment of the line at infinity of $\RRR P^2$ on their boundary.
Therefore, the number of critical points of $S_{2k}^-$ is
\begin{equation}
\lb{e:e1}
c_1^-(B_{2k}) \ge \frac {3k (k-1)} {2} -1 +(k+1) = \frac {k (3k-1)}{2}
\end{equation}

\vskip0.1in
(*)Let $\HH _{2k+1}$ be the Harnack curve of odd degree $2k+1 $

When $c$ decreases from $0$ to $c_m - \e$, outer ovals and
regions $G_i$ where $b_ {2k+1} > c_0$
which intersect the line at infinity of $\RRR P^2$ on their boundary,
expand,
then coalesce, and finally disappear in a simply connected region.

Moreover, from Lemma \mrf{l:lp2}, one can associate
to each coalescence of one outer oval ${\mathcal O}$ with
 an other outer oval or with one of the regions $G_i$
 a critical point $(x_0,y_0)$  of index 1 with
$b_{2k}(x_0,y_0)=c_0$; $c_0 \in ]c_m, 0[$
with the property that ${\mathcal O}$ disappears as $c$ decreases
from $c_0+\e$  to $c_0-\e$.

Denote $S_{2k+1} ^-$ the set of these critical points.
The curve $\HH_{2k+1}$ has $\frac {k(3k-1)} {2}$ outer
ovals and $(2k+1)$ connected components\\
$M_0= \{ (x,y) \in \RRR^2 | b_{2k+1}(x,y) > 0 \}$ which contain
a segment of the line at infinity of $\RRR P^2$ on their boundary.
Therefore, the number of critical points of $S_{2k+1}^-$ is:
\begin{equation}
\lb{e:e2}
c_1^-(B_{2k+1}) \ge \frac {k (3k-1)} {2} + 2k= \frac{ k(3k+3)}{2}
\end{equation}

\item
Let $c$ increase from $0$ to $c_M + \e$, then the outer oval
shrink and disappear.
Moreover, from the Lemma \mrf{l:lp3}, the set  $M_{c_M+\e}$
has $m$ components.
\vskip0.1in
 (*)Let $\HH _{2k}$ be the Harnack curve of even degree $2k$

The set $M_0$ has $(2k -(k-1))$ components intersecting the line at infinity;
the curve  $\HH _{2k}$ has $\frac {(k-1)(k-2)} {2}$ inner ovals.
It is easy to deduce from the Lemma \mrf{l:lp2}  that
when $c$ increases from $0$ to $c_M + \e$, the inner ovals expand,
touch the non-empty outer oval
or another inner oval
and  then disappear.
In other words, when $c$ decreases from $c_M + \e$
to $0$, $\frac {(k-1)(k-2)}{2}$ inner ovals appear.
Hence, from Lemma \mrf{l:lp2}(\ref{e:ep2}), one can associate
to each inner oval  ${\mathcal O}$ a critical point
  $(x_0,y_0)$ of index 1
with
$b_{2k}(x_0,y_0)=c_0$; $c_0 \in ]0, c_M[$
with the property that ${\mathcal O}$
disappears as $c$ increases from  $c_0 - \e$ to $c_0 + \e$.

\vskip0.1in
Furthermore, the set  $M_{c_M+\e}$
has $2k$ components.
Thus, from the Lemma \mrf{l:lp3},
$c$ has passed at least $(k-1)$ critical values of index $1$ in his way
from $0$ to $c_M + \e$.
We denote by $S'_{2k}$ the set of critical points associated to
these critical values.

\vskip0.1in
Denote $S_{2k}^+ $ the set of
critical points of index 1 $(x_0,y_0)$ such that
$b_{2k}(x_0,y_0)=c_0$; $c_0 \in ]0, c_M[$.
Therefore, the number of critical points of $S_{2k}^+$ is:
\begin{equation}
\lb{e:e3}
c_1^+(B_{2k}) \ge \frac{ (k-1)(k-2)}{2} +(k-1)=\frac{ k(k-1)}{2}
\end{equation}

\vskip0.1in
 (*)Let $\HH _{2k+1}$ be the Harnack curve of odd degree $2k+1$.

The set $M_0$
has $2k+1$ components intersecting the line at infinity;
the curve $\HH _{2k+1}$ has $\frac {k(k-1)} {2}$ inner ovals.
When $c$ increases from $0$ to $c_M + \e$;
inner ovals expand, then touch
the one-side component or another inner oval,
and finally disappear.

Moreover, the set  $M_{c_M+ \e}$ has
$2k+1$ components.

From Lemma \mrf{l:lp2} (\ref{e:ep3}), one can associate
to each inner oval  ${\mathcal O}$ of
$\HH_{2k+1}$ a critical point
$(x_0,y_0)$ of index 1 with
$b_{2k+1}(x_0,y_0)=c_0$; $c_0 \in ]0, c_M[$ with the property that
${\mathcal O}$ disappears as $c$ increases from $c_0-\e$ to $c_0 + \e$.

Denote $S_{2k+1}^-$ the set of these critical points.
The curve $\HH_{2k+1}$ has $\frac {k(k-1)} {2}$ inner
ovals.
Hence, the number of critical points of $S_{2k+1}^+$ is:
\begin{equation}
\lb{e:e4}
c_1^+(B_{2k+1}) \ge \frac{ k(k-1)}{2}
\end{equation}

\enen

\vskip0.1in
Hence, according to Proposition \mrf{p:crit},
inequalities \ref{e:e1}, \ref{e:e2}, \ref{e:e3}, \ref{e:e4}
are equalities.
This implies the Proposition \mrf{p:crit2}.
 Q.E.D
In particular,
this method provides curves with real scheme:\\
for even $m=2k$\section{Harnack Curves from a real viewpoint-Rigid Isotopy Classification-}
\mlb{su:comp}
Passing from polynomials to real set of points, it follows
that real algebraic curves form a real projective space of dimension
$\frac {m (m+3)} 2$.
We shall denote this space by the symbol $\RRR {\mathcal C}_m$
and by $\RRR {\mathcal D}_m$ the subset of
$\RRR {\mathcal C}_m$ corresponding
to real singular curves.
We call a  path in the complement
$\RRR {\mathcal C}_m \bk \RRR {\mathcal D}_m$
of the discriminant hypersurface in $\RRR {\mathcal C}_m$
a {\it rigid isotopy} of real point set
of nonsingular curves of degree $m$.
The classification of real point set of
curves of degree $\le 4$ up to rigid isotopy
is known since the $19^{th}$ century. It was completed  for curves
of degree $\le 5$ and $\le 6$ at the end of the seventies.
Up to rigid isotopy a curve
of degree $\le 4$ is determined by its real scheme; up to rigid
isotopy a curve of degree $5$ or $6$, is determined by its real scheme
and its type.
\vskip0.1in
Let us recall that we denote by  $\HH_m$ and call Harnack curve
any curve
with real scheme:\\
- for even $m=2k$
$$
\big\langle1\langle \frac {(k-1)(k-2)}2 \rangle \sqcup \frac {3k(k-1)}2
\big\rangle$$
- for odd $m=2k+1$
$$
\big\langle J \sqcup k(2k-1)
\big\rangle $$
\vskip0.1in

The main result of this section is
given in the Theorem \mrf{t:rigiso}
where we establish
the rigid isotopy classification of real point set of Harnack curves
$\RRR \HH_m$.
Harnack curves of degree $i \le 6$, as any $M$-curve
of degree $\le 6$, are rigidly isotopic.
In Theorem \mrf{t:rigiso}, we
extend this property to Harnack curves of arbitrary degree.
Precisely, we prove that
isotopy also implies
rigid isotopy for real point set of Harnack curves.
This section until its end is devoted to the proof of this result.
\beth
\mlb{t:rigiso}
\centerline{Rigid Isotopy Classification Theorem}
Harnack curve $\HH_m$ of degree $m$
are rigidly isotopic.
\enth

The proof of Theorem \mrf{t:rigiso} is based on
a modification of Harnack polynomials.
Let us call {\it regular modification}
of a polynomial any modification on its coefficients
with the property that
real point set of curves of the modified polynomial and
the initial polynomial  are  rigidly isotopic.\\
(In particular, given a polynomial
of a smooth algebraic curve, the modification of its coefficients
such that the modified polynomial is a regular polynomial (i.e
none of its critical points lies on the line at infinity,
any two different critical points have distinct critical values)
is obviously regular.)\\
Let ${\mathcal A}$ and ${\mathcal B}$,
be two smooth curves
such that the union ${\mathcal A} \cup {\mathcal B}$
is a singular curve all of whose singular points are crossings.
Denote by ${\mathcal H}$ the set of curves which result from
the classical deformation of ${\mathcal A} \cup {\mathcal B}$.
We shall say that a curve is deduced from
{\it deformation} of ${\mathcal A} \cup {\mathcal B}$
if it is rigidly isotopic to a curve ${\mathcal C}$ of
the set ${\mathcal H}$.\\
Any
two Harnack curves $\HH_m$ constructed from the Harnack's method
are rigidly isotopic.(We recall Harnack's method and
 propose a proof of this
statement in the Appendix)
The proof of the Theorem \mrf{t:rigiso} is based on
a regular modification of Harnack polynomials.
From properties of this regular modification,
we shall deduce that any Harnack curve $\HH_m$
is, up to rigid isotopy, constructed from the Harnack's method.
In this way,
we get the rigid isotopy
Theorem \mrf{t:rigiso}.

\subsection{ Harnack curve of type $\HH^0$}
\mlb{su:HH0}
\vskip0.1in
Let us recall that we call {\it polynomial of type $\HH$}
a regular Harnack polynomial of degree $m$  giving  a Harnack curve $\HH_m$
which intersects the
line at infinity $\RRR L$ of $\RRR P^2$ in $m$  real distinct points.
We shall call {\it polynomial of type $\HH^0$}
(relatively to $L$) a polynomial of type $\HH$
such that the line at infinity $L$ of $\RRR P^2$ intersects the
Harnack curve $\HH_m$ of degree $m$
in $m$ real distinct points which
belong to the same connected component of $\RRR \HH_m$:
in case $m=2k$, the non-empty oval of $\HH_{2k}$;
in case $m=2k+1$, the odd component of $\HH_{2k+1}$.
The odd component of $\HH_{2k+1}$ is divided
into $2k+1$ arcs which delimit a region of which boundary
contains a segment of the
line at infinity. The line at infinity is chosen such that
only one of this region contains ovals of
$\HH_{2k+1}$.
We say that a curve is of type $\HH^0$ if its polynomial is of type $\HH^0$.
In particular,
Harnack curves constructed from the
Harnack's method (see Appendix) are
curves of type $\HH^0$.

\vskip0.1in
The main result of this subsection is the Theorem \mrf{t:rig}
where we give the rigid isotopy classification
of curves of type $\HH^0$.

\beth
\mlb{t:rig}
Harnack curve $\HH_m$ of degree $m$ and type $\HH^0$
are rigidly isotopic.
\enth
Let $\HH_m$  be a  Harnack curve with polynomial
$B_{m}(x_0,x_1,x_2)$ of type $\HH^0$.
Denote by
$C_i(x_1,x_2)$ the unique
homogeneous polynomial of degree $i$ in the variables $x_1,x_2$
such that:
$B_m(x_0,x_1,x_2) = x_0^{m-1}.B_1(x_0,x_1,x_2)+
 \Sigma_{i=2}^m x_0^{m-i}.C_i(x_1,x_2)$. \\
Let $b_m(x,y)$ be the affine polynomial associated to $B_m(x_0,x_1,x_2)$,
$$B_m(x_0,x_1,x_2)=x_0^m.b_m(x_1/x_0,x_2/x_0)$$
Let $b_i(x,y)$ be the truncation of $b_m(x,y)$ on the monomials
$x^{\a}.y^{\b}$ with $0 \le \a + \b \le i $ and
$B_i(x_0,x_1,x_2)= x_0^i.b_i(x_1/x_0,x_2/x_0)$  be
the homogeneous polynomial associated to $b_i$.
We shall denote by ${\mathcal B}_i$ the curve
with polynomial $B_i(x_0,x_1,x_2)$.
\vskip0.1in
Consider the norm in the vector space of polynomials
$$ \vert \vert \sum a_{i,j}x_1^ix_2^j \vert \vert =max \{ \vert a_{i,j}
 \vert~ | (i,j) \in \NNN^2 \}$$
(Given $A_m(x_0,x_1,x_2)$ an homogeneous polynomial,
$$ \vert \vert A_m(x_0,x_1,x_2) \vert \vert =
 \vert \vert A_m(1,x_1,x_2) \vert \vert $$)

The Theorem \mrf{t:rig} may be also formulated as follows:

\bepr
\mlb{p:ll}
On the assumption that $B_m(x_0,x_1,x_2)$ is of type $\HH^0$
relatively to the line at infinity $x_0=0$,
up to regular modification of $B_m(x_0,x_1,x_2)$, for any
$B_{i}(x_{0},x_1,x_2)$, $i \ge 1$, we have:
\been
\item
\lb{i:a41}
$B_{i}(x_{0},x_1,x_2)$
 is smooth  of type $\HH^0$
\item
\lb{i:a4}
none of the critical points of $B_{i}(x_{0},x_1,x_2)$
belongs to the line at infinity;\\
for any critical point $(1,x_{0,1},x_{0,2}) \in \RRR P^2$
of $B_i(1,x_1,x_2) = b_i(x_1,x_2)$
its representative in $S^2$
$\frac {1} {(1+x_{0,1}^2+x_{0,2}^2)^{1/2}} (1,x_{0,1},x_{0,2})$
is such that:\\
$B_i(
\frac {1} {(1+x_{0,1}^2+x_{0,2}^2)^{1/2}} (1,x_{0,1},x_{0,2}))
\notin [-
\Sigma_{j=i+1}^m
\vert \vert C_j \vert  \vert ,
\Sigma_{j=i+1}^m
\vert \vert C_j \vert  \vert ]$
\enen
\enpr

On the assumption that $B_{m}$ is of type $\HH^0$, according to
proposition \mrf{p:crit}, the following
equalities are verified:\\
for even $m=2k$:
$c(B_{2k})=(\frac {(k-1) (k-2)}{2},k.(2k-1),\frac {3k.(k-1)} {2})$\\
and for odd $m=2k+1$:
$c(B_{2k+1})=(\frac { k.(k-1)} 2, k.(2k+1), \frac {k.(3k-1)} {2})$\\

{\bf proof:}\\
Our proof is based on Morse Lemma and Petrovskii's theory.
We shall proceed by descending induction on the degree $i$ of $B_i$.
Let us assume that assumptions (\ref{i:a41}) and (\ref{i:a4})
of Proposition \mrf{p:ll} are satisfied for
$B_n$ with $i+1 \le n \le m$,
and prove that $B_i$ also satisfies (\ref{i:a41}) and (\ref{i:a4}).\\
On these assumptions, we shall prove that,
the Harnack
curve $\HH_{i+1}$ with polynomial $B_{i+1}$
is deduced from deformation of
${\mathcal B_{i}} \cup L$
,$B_{i+1}(x_{0},x_1,x_2)=x_{0}.B_{i}(x_{0},x_1,x_2)+ C_{i+1}(x_1,x_2)$,
where $L:=\{ x_{0}=0 \}$ and $B_i$ is of type $\HH^0$ relatively to $L$.
\vskip0.1in
\vskip0.1in
{\bf I)}
Let us in a first part
study  curves ${\mathcal B}_i$ of degree $\ge 4$.
Up to regular modification of the Harnack polynomial $B_{i+1}$
one can always assume that $B_{i}$ is smooth.\\
\vskip0.1in
{\bf Introduction}
\vskip0.1in
Consider a real projective line $\RRR L \subset \RRR P^2$;
its tubular neighborhood in $\RRR P^2$
is homeomorphic to a M\"{o}bius band.
The core of the M\"{o}bius band, i.e the real projective line $\RRR L$,
is the circle
with framing $\pm 1$.
In such a way, one can identify points of $\RRR L$ with points
of two halves (oriented) circles.
Halves intersect each other in two points
we shall call {\it extremities}.
We shall denote $\RRR L^{+}$ and $\RRR L^-$ the halves of $\RRR L$.
On each half $\RRR L^{\pm}$
of $\RRR L$ one can consider an oriented tubular fibration.
We shall denote by ${\mathcal M}^{\pm}$
the half of the M\"{o}bius band which is the
tubular neighborhood of $L^{\pm}$.
In such a way, the boundary of the half of M\"{o}bius band
${\mathcal M}^{\pm}$ contains
the union of two real projective
lines $\RRR L_1^{\pm}$, $\RRR L_2^{\pm}$.
The  half  $\RRR L_i^{\pm}$ of $\RRR L_i$, $i \in \{1,2\}$,
is the image of a smooth
section of the tubular oriented fibration of
a tubular neighborhood of  $\RRR L^{\pm}$;
${\mathcal M}^- \cup {\mathcal M}^+={\mathcal M}$,
$\pr {\mathcal M} \supset \RRR L_1^{\pm} \cup \RRR L_2^{\pm}$.
\vskip0.1in
Let $\HH_m$ be
the Harnack curve of degree $m$ with polynomial $B_m$ of type $\HH^0$
relatively to $L$.
We shall study $\RRR \HH_m$  in a tubular neighborhood ${\mathcal M}$
of the line at infinity $\RRR L$.\\
Let ${\mathcal A}$ be the set of
arcs of $\RRR \HH_m$ bounding regions
which contain a segment of the line at infinity on  their boundary
and do not contain ovals of $\RRR \HH_m$.
We shall say that a critical point of index $1$
is {\it associated} to an arc $\gamma  \in {\mathcal A}$ of $\RRR \HH_{m}$
if there exists $p$ a critical point of $b_{m}$
with critical value $c_0$
such that as $c$ varies from $0$ to $c_0$,
the region which contains $\gamma$ on its boundary
varies in such a way that for $c=c_0$ it touches an other arc $\gamma'$.
The arc $\gamma'$ is said {\it associated} to the arc $\gamma$
and the pair $(\gamma, \gamma')$ is said associated to the critical
point $p$.
We shall denote
${\mathcal A}'$ the set of arcs ${\gamma}'$.
\vskip0.1in
From the study of the Petrovskii's pencil
$x_0^{m}.(b_{m}(x_1/x_0,x_2/x_0)-c)$, $c \in \RRR$
over $\HH_{m}$ with $L:=x_0=0$, we shall define a set ${\mathcal P}$
of critical points of $b_m$ and a set of arcs of $\RRR \HH_m$
associated to ${\mathcal A}$.\\
The set ${\mathcal A}$ is the image of a smooth
section $s$ of a tubular fibration of $\RRR L^+$ minus 2 points.
The proof of Theorem \mrf{t:rig}
is based on the characterization (up to regular deformation
of $B_m$) of a subset  of $\RRR \HH_{m}$
which is the image of the
extension of $s$ to a smooth section of a tubular fibration
of $\RRR L$ minus a finite number of points.\\
We shall proceed as follows.
Using Morse Lemma and Petrovskii's theory,
we shall define a set of arcs
${\mathcal D}$ of $\RRR \HH_m$
with the property that up to regular modification of $B_m$
the line at infinity $\RRR L$ divides  any arc
 $\xi \in {\mathcal D}$  into two halves
which belong respectively to $\RRR {\mathcal B}_{m-1}$ and
$\RRR  {\mathcal  B}_{m} \bk \RRR {\mathcal B}_{m-1}$
where ${\mathcal B}_{m-1}$
is the Harnack curve of degree $m-1$.
In such a way, we shall deduce
that the Harnack curve $\HH_m$ is deduced from deformation of
$\HH_{m-1} \cup L$
$B_{m}(x_{0},x_1,x_2)=x_{0}.B_{m-1}(x_0,x_1,x_2)+ C_{m}(x_1,x_2)$,
where $L:=\{ x_{0}=0 \}$ and $B_{m-1}$ is of type $\HH^0$.

\vskip0.1in
Let us distinguish in parts {\bf I.i)} and {\bf I.ii)}
curves of even and odd degree.\\
{\bf I.i) -Harnack curves of even degree-}\\

Let  $\HH_{2k}$ be a Harnack curve of degree $2k$ with
polynomial $B_{2k}(x_0,x_1,x_2)$
of type $\HH^0$ relatively to the line at infinity
$L:=x_0=0$.
We shall prove the following Proposition \mrf{p:stv}.
\vskip0.1in
\bepr
\mlb{p:stv}
Let  $\HH_{2k}$ be a Harnack curve of degree $2k$ with
polynomial $B_{2k}(x_0,x_1,x_2)$
of type $\HH^0$ relatively to the line at infinity
$L:=x_0=0$.
Up to regular modification of $B_{2k}$,
$B_{2k-1}(x_0,x_1,x_2)$ is also of type $\HH^0$ relatively to $L$.
\enpr
{\bf proof:}\\
The part I.i) until its end is devoted to the proof of
Proposition \mrf{p:stv}.\\
The non-empty oval of $\HH_{2k}$ is divided
into $2k$ arcs which  delimit a region which contains
a segment of the line at infinity on its boundary.
Denote ${\mathcal A}$  the set of
the $2k-1$ arcs which delimit a region
which does not contain oval of $\RRR \HH_{2k}$;
$k-1$ (resp, $k$) of them delimit
regions $\{ x \in \RRR^2 | b_{2k}(x) < 0 \}$
(resp, $\{ x \in \RRR^2 | b_{2k} (x) > 0\}$).
The union of arcs of ${\mathcal A}$
is a connected orientable part of the non-empty oval
which belongs to a tubular neighborhood of $L$.
The connected surface (with boundary the non-empty oval)
obtained from removing the interior of the inner ovals to
the interior of the non-empty oval is also orientable.
Thus, the orientation of the connected
union of the $2k-1$ arcs of ${\mathcal A}$
is induced by a half of $\RRR L$. We shall denote $\RRR L^+$ this half. \\
The tubular neighborhood of the line $\RRR L$
is a M\"{o}bius band ${\mathcal M}={\mathcal M}^- \cup {\mathcal M}^+$
$\pr {\mathcal M}^{\pm} \supset \RRR L_1^{\pm} \cup \RRR L_2^{\pm}$.
The set ${\mathcal A}$ is the image of a smooth
section of a tubular fibration of $\RRR L^+$
(-precisely, $\RRR L^+$
minus its extremities, these 2 points belong to the
arc of $\RRR \HH_{2k}$ bounding the non-empty region with
a segment of the line at infinity on its boundary-)
\vskip0.1in

It is easy to see that
there exists an isotopy  of $\RRR P^2$
which pushes the $k$ (oriented) arcs of ${\mathcal A} \cap \RRR^2$
bounding positive regions
$\{ x \in \RRR^2 | b_{2k}(x) > 0 \}$  to
one positive line $\RRR L_i^+$, $i \in \{1,2 \}$.
Let it be $\RRR L_1^+$.
In this way,
there exists an isotopy
 of $\RRR P^2$
which pushes
the $k-1$ (oriented) arcs of ${\mathcal A} \cap \RRR^2$
bounding  negative regions
$\{ x \in \RRR^2 | b_{2k}(x) < 0 \}$ to  $\RRR L_2^+$.
\vskip0.1in
Let us prove that one can associate
to the set of $2k-1$ arcs ${\mathcal A}$
a set of critical points ${\mathcal P}$
of index $1$ of $b_{2k}$ and a set of arcs
${\mathcal A}'$ of $\RRR \HH_{2k}$.\\

To this end,  we shall study the pencil of curves
$x_0^{2k}.(b_{2k}(x_1/x_0,x_2/x_0)-c)$, $c \in \RRR$,
over $\HH_{2k}$.\\
Let us recall that (see Petrovskii's Lemmas \mrf{l:lp2} \mrf{l:lp3})
as $c$ decreases from $0$
positive regions and positive ovals expand.
As $c$ decreases from $0$
the number of regions $G_i$ of $\{ x \in \RRR^2 | b_{2k}(x) > 0 \}$
which contain a segment of the line at infinity of $\RRR P^2$ on their
boundary decreases from $k+1$ to $1$.\\
As $c$ increases from $0$
the number of regions $G_i$ of $\{ x \in \RRR^2 | b_{2k}(x) > 0 \}$
which contain a segment of the line at infinity of $\RRR P^2$ on their
boundary increases from $k+1$ to $2k$.
Hence, to each arc of ${\mathcal A}$ is associated a critical of index $1$.\\

-Let us study the set ${\mathcal P}$ of critical point $1$ associated to
the set  ${\mathcal A}$ of arcs of $\RRR \HH_{2k}$.\\
Recall that
we consider on $\RRR P^2$ the Fubini-Study  metric induced
by the projection $\pi_{\RRR} :S^2 \to \RRR P^2$.\\

Without loss of generality, one can assume
that critical points of $B_{2k}$
and points of the intersection of $\HH_{2k}$ with $L$
do not belong to the line $x_2=0$.
On such assumption,
consider the function
$b_{2k}(0,x_1/x_2,1)=b_{2k}(y)$
$$\frac {\pr} {\pr x_1} b_{2k}(0,x_1/x_2,1)=
\frac {\pr} {\pr y} b_{2k}(y) \frac{1} {x_2}$$
$$\frac {\pr} {\pr x_2} b_{2k}(0,x_1/x_2,1)=
\frac {\pr} {\pr y} b_{2k}(y) \frac{x_1} {-(x_2)^2}$$
Critical points of the
function
$b_{2k}(0,x_1/x_2,1)$
 may be defined from
critical points $b_{2k}(y)$.
(Obviously, for any such critical point
of $b_{2k}(0,x_1/x_2,1)$
the following equality is verified
$(0:x_1/x_2:1)=(0:x_1:x_2)=(0:-x_1:-x_2)$.)
The function $b_{2k}(y)$ has exactly $2k$ zeroes
which coincide with intersection points of  $\HH_{2k}$ with
the line $L$.
Hence, by Rolle's Theorem $2k-1$ extrema at which
$b_{2k}'(y)$ must change from positive to negative.
By continuity of $b_{2k}$ and $b'_{2k}$, it follows the
alternation of
sign of $b'_{2k}(x_0,x_1,x_2)$
in an $\e$-tubular neighborhood  ${\mathcal M}_{\e}$
of the line $x_0=0$ in $\CCC P^2$ and
thereby the existence of a set of
$2k-1$ critical points of $B_{2k}$
(i.e critical points of the affine polynomial $b_{2k}(x_1/x_0,x_2/x_0)$)
in ${\mathcal M}_{\e}$.
According to Petrovskii's theory,
to describe critical points of $b_{2k}(x_1/x_0,x_2/x_0)$
by means of the pencil  $x_0^{2k}.(b_{2k}(x_1/x_0,x_2/x_0)-c)$
one can assume that $B_{2k}$ is a regular polynomial.
(i.e none of its critical points belongs to $x_0=0$,
 any two of them have different critical values.)
Up to slightly modify coefficients of $b_{2k}$,
any critical point  $y=(0,x_1/x_2,1)$
of $b_{2k}(y)$
gives rise to one critical point
which may be chosen
among two points
$(\e:\pm x_{0,1}: \pm x_{0,2})
 \in \RRR P^2$.\\

\vskip0.1in
Let us prove in Lemma \mrf{l:alp} that
without loss of generality
one can assume that the $2k-1$ points of ${\mathcal P}$
belong to a line $L_{\e_1}$ of
${\mathcal M}_{\e}$ the $\e$-tubular neighborhood of $L:=x_0=0$
where $\e>0$ is arbitrarily small.\\

\bele
\mlb{l:alp}
Let $L$ be the real projective line at infinity.
Denote by ${\mathcal M}_{\e}$ the $\e$-tubular
neighborhood of $L$ in $\RRR P^2$.\\
Given $B_{2k}$ a Harnack polynomial of type $\HH^0$
(relatively to $L$).
There exists a regular modification
of $B_{2k}$ such that:
\been
\item
the modified polynomial
is of type $\HH^0$ (relatively to $L$).
\item
the set ${\mathcal P}$ is a set of
$(2k-1)$ points of a  line $L_{\e_1} \subset {\mathcal M}_{\e}$
and $\e$ is arbitrarily small.
\enen
\enle
{\bf proof:}

Recall that any point  $p \in {\mathcal P}$
is a critical point of index $1$ associated  to
one of the $2k-1$ arcs of the set ${\mathcal A}$
(-the set of arcs of $\RRR \HH_{2k}$
which  delimit a  region which contains
a segment of the line at infinity on its boundary
and does not contain an oval of $\RRR \HH_{2k}$-)

-It is easy to get that we can regularly modify
$B_{2k}$ such that any $p \in {\mathcal P}$ belongs to ${\mathcal M}_{\e}$
where $\e >0$ is arbitrarily small.\\

-Assume that any point $p \in {\mathcal P}$
belongs to ${\mathcal M}_{\e}$ the $\e$-tubular neighborhood
of $L$ with $\e$ is arbitrarily small.
The polynomial $B_{2k}(x_0,x_1,x_2)$  is of type $\HH^0$
relatively to the line $L$.\\
Let us now prove that  one can
regularly modify $B_{2k}$ such that any $p \in {\mathcal P}$
belongs to a real projective line $L_{\e_1} \subset {\mathcal M}_{\e}$.\\
Let $L_{\e_1}$ be a real projective line of ${\mathcal M}_{\e}$.
Denote by
$x'_0= x_0+ \a.x_1+\b.x_2$ the polynomial of $L_{\e_1}$.
Then,
consider the linear change of projective coordinates
mapping $(x_0:x_1:x_2)$ to $(x'_0:x_1:x_2)$.
Such transformation carries
$B_{2k}(x_0,x_1,x_2)=x_0.B_{2k-1} (x_0,x_1,x_2)+ C_{2k}(x_1,x_2)$ to\\
$B'_{2k}(x'_0,x_1,x_2) =x'_0. B'_{2k-1}(x'_0,x_1,x_2) +C'_{2k}(x_1,x_2)$
where\\
$B'_{2k-1}(x'_0,x_1,x_2) =x'_0. B'_{2k-2}(x'_0,x_1,x_2) +C'_{2k-1}(x_1,x_2)$.
We shall  prove that, up to regular modification of $B'_{2k}$ and
thus of $B_{2k}$,
 one can choose
$C'_{2k-1}(x_1,x_2)$
in such a way that
the $2k-1$ points of ${\mathcal P}$ belong to $L_{\e_1}$.
Using the linear change of projective coordinates
mapping $(x'_0:x_1:x_2)$ to $(x_0:x_1:x_2)$,
we shall get the polynomial
$B_{2k}(x_0,x_1,x_2)$.\\
Let us detail this construction.\\
Given $B_{2k}(x_0,x_1,x_2)$  of type $\HH^0$,
the  polynomial $C_{2k}(x_1,x_2)$ such that
$B_{2k}(x_0,x_1,x_2)=x_0.B_{2k-1} (x_0,x_1,x_2)+ C_{2k}(x_1,x_2)$
has $2k$ distinct roots
on the line at infinity $L:=x_0=0$.
Thus, $C_{2k}(x_1,x_2)=
\Sigma_{i=0}^{2k} a_ix_1^{2k-i}x_2^i$ with
$a_i \not=0$.

It follows from the previous local study around $p \in {\mathcal P}$
that, up to regular modification,  one can assume the
polynomial  $B_{2k}$
(of type $\HH^0$ relatively to $L$)
such that
any point $p$ of ${\mathcal P}$
belongs to the $\e$-tubular neighborhood  ${\mathcal M}_{\e}$ of $L$
where $\e>0 $ is arbitrarily small.
Denote by $B_{2k,0}$ such a polynomial and
by ${\mathcal P}_0$ its respective set ${\mathcal P}$ of critical points.

Set $x_0=\e>0$.\\
$B_{2k,0}(\e,x_1,x_2)=\e.B_{2k-1,0}(\e,x_1,x_2)+C_{2k,0}(x_1,x_2)$\\
$=\e.(\e.B_{2k-2,0}(\e,x_1,x_2)+C_{2k-1,0}(x_1,x_2))
+C_{2k,0}(x_1,x_2)$
with $B_i=x_0^i.b_i$ where $b_i$ is the truncation
of $b_{2k,0}=B_{2k,0}(1,x_1,x_2)$ monomials of
 $x_1^a.x_2^b$ of  degree $a+b \le i$.
Local coordinates defined
in a neighborhood $U(p)$ of a point $p$
depend principally on the first derivative and the second derivative
of the function $B_{2k}$ around $p$.
Hence, to describe the curve in neighborhood
of points  ${\mathcal P}$, it is sufficient to consider
the truncation $D_{2k,0}$ of $B_{2k,0}$ on monomials
$x_1^a.x_2^b$ of  degree  $2k-2 \le a + b \le 2k$.
Therefore,
to define  a regular modification
$B_{2k,t}$, $t \in [0,1]$,
of $B_{2k,0}$ of type $\HH^0$ with ${\mathcal P}_0  \subset {\mathcal M}_{\e}$
to  $B_{2k,1}$ of type $\HH^0$
with ${\mathcal P} \in L_{\e_1}$ it sufficient to consider the truncation
$D_{2k,0}$ of $B_{2k,0}$.\\
Without loss of generality, one can assume
that the $2k-1$ points $p \in {\mathcal P}_{0}$
do not belong to $x_2=0$.
In this way, for any $p \in {\mathcal P}_0$,
$\frac{\pr B_{2k,0}} {\pr x_0} (p) =0$
 and
$\frac {\pr {B}_{2k,0}} {\pr x_1} (p)=0$.
Moreover, $B_{2k,0}$ has $2k$ roots on the line $x_0=0$
and none of these points is a point of ${\mathcal P}$.
Therefore, it follows that
the set ${\mathcal P}_{0} \subset {\mathcal M}_{\e}$
may be chosen such
that the truncation
$D_{2k,0}(x_0,x_1,x_2)$ of $B_{2k,0}(x_0,x_1,x_2)$
is of the form
$D_{2k,0}=\Sigma_{2k-2\le i+j \le 2k} a_{i,j}x_0^{2k-(i+j)}x_1^{i}x_2^j$ with
$a_{i,j} \not = 0$.\\
Set $B'_{2k}=B'_{2k,0}=x'_0.B'_{2k-1}+C'_{2k}$. We shall
prove that
there exists a regular
modification of $B'_{2k,0}$ (that is of $B_{2k,0}$)
which pushes any point $p \in {\mathcal P}_0$
to a point
of the set
$\tilde{\mathcal P}_0 :=B'_{2k-1,1}=0 \cap L_{\e_1}$.

In such a way,
any point $p \in \tilde{\mathcal P}_0$
verifies
\begin{equation}
\lb{s:x1}
\frac {\pr B'_{2k,1}} {\pr x'_0} (p) =0
\end{equation}

Without loss of generality,
on can assume that any point
$p \in \tilde{\mathcal P}_0$
does not belong to $x_2=0$.\\

Thus, any point $p \in \tilde{\mathcal P}_0$
also verifies:
\begin{equation}
\lb{s:x3}
\frac {\pr B'_{2k,1}} {\pr x_1} (p) =0
\end{equation}

\vskip0.1in

Let us prove that the modification
$B_{2k,t}$ $t \in [0,1]$, $B_{2k,0}=B_{2k}$
such that $B_{2k,1}$ verifies properties (\ref{s:x1})
and (\ref{s:x3}) is regular.
We shall proceed by induction.\\
Denote $p_i,~1 \le i \le 2k -1$ the set of $2k-1$ points of ${\mathcal P}_0$.
Denote $L_0^{i,j}$ the unique line such that $p_i,p_j \in L_0^{i,j}$
Setting $L_{0}^{i,j}:=x_0'=0$,
properties (\ref{s:x1}) and (\ref{s:x3}) are verified for $p_i,p_j$.
Choose $p_1,p_2$
two points of ${\mathcal P}_0$ associated to respectively $\gamma_1$
$\gamma_2$ arcs  of ${\mathcal A}$  which
intersect each other in one point of the line at infinity $L$.

Let $p_3 \in {\mathcal P}_0$ associated to $\gamma_3 \in {\mathcal A}$,
such that $\gamma_3$ and  $\gamma_2$ intersect each other in one point.

Setting $L_{\e_1}=L_{0}^{1,2}:=x'_0=0$,
$L_0^{2,3}:=x"_0=x'_0+\a'.x_1+\b'.x_2=0$.
The linear change of projective coordinates
mapping $(x'_0:x_1:x_2)$ to $(x"_0:x_1:x_2)$ carries
$B'_{2k}(x'_0,x_1,x_2)=x'_0.B_{2k-1} (x'_0,x_1,x_2)+ C_{2k}(x_1,x_2)$ to\\
$B"_{2k}(x"_0,x_1,x_2) =x"_0. B"_{2k-1}(x"_0,x_1,x_2) +C"_{2k}(x_1,x_2)$\\
Moreover, $\frac{\pr B"_{2k,0}} {\pr x"_0} (p_3) =0$,
 $\frac{\pr B"_{2k,0}} {\pr x_1} (p_3) =0$\\
One can modify regularly coefficients of $B_{2k,0}$
such that $p_3 \in L_{\e_1}= L_0^{1,2}$.
(i.e there exists regular  modification
on gradient trajectories
such that $p_3$ is pushed to a point of the line $L_{\e_1}=L_0^{1,2}$)
This modification may be expressed from the path
$x"_{0,t}=x'_0+(\a'-t.\a').x_1 + (\b'-t.\b').x_2$, $t \in [0,1]$
from $L_0^{2,3}$
to $L_{\e_1}$
and polynomials $B"_{2k-1,t}(x"_{0,t},x_1,x_2)$ and $C"_{2k,t}(x_1,x_2)$
$B"_{2k,t}=x"_{0,t}. B"_{2k-1,t}(x"_{0,t},x_1,x_2) +C"_{2k,t}(x_1,x_2)$.
The deformation
moves $p_{3,t} \in L_t$  in such a way that
$\frac {\pr B"_{2k,t}} {\pr x"_{0,t}} (p_{3,t})$ and
$\frac {\pr B"_{2k,t}} {\pr x_1} (p_{3,t})$.
Note that $B"_{2k,0}=B"_{2k}$,$p_{3,0}=p_3 \in L_0^{2,3}$,
$p_{3,t} \in L_t$ where $x"_{0,t}$ is the polynomial of $L_t$
and $L_1=L_{\e_1}$)
($B"_{2k,1}=x"_{0,1}. B"_{2k-1,1}(x"_{0,1},x_1,x_2) +C"_{2k,1}(x_1,x_2)$;
$B"_{2k,1}=x'_{0}. B"_{2k-1,1}(x'_{0},x_1,x_2) +C'_{2k,1}(x_1,x_2)
\not =x'_0.B_{2k-1} (x'_0,x_1,x_2)+ C_{2k}(x_1,x_2)$.)
Since $p_1,p_2,p_3 \in {\mathcal M}_{\e}$, where $\e>0$ is arbitrarily
small; $|\a'|$ and $|\b'|$ are also arbitrarily small.
In this way,
(see Lemma 1 and Lemma 2 of \cite{Petr}), the deformation
$B"_{2k,t}=x"_{0,t}. B"_{2k-1,t}(x"_{0,t},x_1,x_2) +C"_{2k,t}(x_1,x_2)$
such that $p_{3,t}$ is a root of $B"_{2k-1,t}(0,x_1,x_2)$
(i.e of $\frac {\pr B"_{2k,t}} {\pr x"_{0,t}}$)
and of $\frac {\pr B"_{2k,t}} {\pr x_1}$, $t \in [0,1]$ is regular.

Iterating this argument for $p_i$, $4 \le i \le 2k-1$, we get the announced
regular deformation $B_{2k,t}$, $t \in [0,1]$, of $B_{2k}$.
Any $p_i \in {\mathcal P}_0$ is
pushed to a root of $C'_{2k-1,1}(x_1,1)$
(that is a root of $\frac {\pr B'_{2k,1}} {\pr x'_0}$)
which is also a root of
$\frac {\pr {C}'_{2k,1}} {\pr x_1} (x_1,1)$
(that is a root of $\frac {\pr B'_{2k,1}} {\pr x_1}$)
\vskip0.1in
For any
$p_i$, $1 \le i \le 2k-1$, of
$\tilde{\mathcal P}_0:=B'_{2k-1,1}=0 \cap L_{\e_1}$
consider its representative  $(0:x_i:1)$.

Let  $s_j$, $j \in \{1,...,2k-1 \}$
be the elementary roots symmetric polynomial
of degree $j$ on the  $p_i$.

The polynomial $C'_{2k,1}(x_1,x_2)$ of
$B'_{2k,1}(x'_0,x_1,x_2)=x'_0.B'_{2k-1,1}(x_0',x_1,x_2)+C'_{2k,1}(x_1,x_2)$
is (up to multiplication by a constant) such that:
$$\frac {\pr {B}'_{2k,1}} {\pr x_1} (0,x_1,x_2)
= x_1^{2k-1} -s_1x_1^{2k-2}x_2 + ...(-1)^{2k-1}s_{2k-1}x_2^{2k-1}$$

Bringing together
the resulting  equalities, it follows the polynomial
$B'_{2k,1}(x'_0,x_1,x_2)$
and thus the definition of $B_{2k,1}(x_0,x_1,x_2)$
with ${\mathcal P}= \tilde{\mathcal P}_0 \subset L_{\e_1}$.
Since
$L_{\e_1} \subset {\mathcal M}_{\e}$,
with $\e>0$
arbitrarily small, the modification
$B_{2k,t}$, $t \in [0,1]$,
of $B_{2k,0}$ of type $\HH^0$ with ${\mathcal P}_0  \subset {\mathcal M}_{\e}$
is such that  $B_{2k,1}$ with
${\mathcal P}_1 \in L_{\e_1}$ is also of type $\HH^0$.
\vskip0.1in
\bere
Note that the rigid proposed above requires a balance of critical value.
Hence, it does not push
arcs ${\mathcal A}$ and ${\mathcal A}'$ to ${\mathcal M}_{\e}$.
Indeed,  consider the Fubini-Study metric on $S^2$
and in such a way consider $\RRR \HH_{2k} \subset S^2$.
Any point $p \in {\mathcal P}$ has representative in $S^2$
$(x_0,x_1,x_2)= (\e,p_1,p_2)$ with $\e << 1 $.
$$B_{2k}(x_0,x_1,x_2)
=x_0^{2k-2}B_{2k-2}(x_0,x_1,x_2)
+D_{2k}$$
Hence,  since $\e^{2k-2}.||B_{2k-2}|| << 1$ for $\e << 1$;
In a neighborhood of $p$,
$$B_{2k}(\e,p_1,p_2)
=\e^{2k-2}B_{2k-2}(\e,p_1,p_2)
+D_{2k}(\e,p_1,p_2)
\approx
D_{2k}(\e,p_1,p_2)$$
It follows that
the "depth of waves"
${\mathcal A} \cup {\mathcal A}'$
depends principally  on coefficients of
monomials $x_1^{2k-i}x_2^i$ of $C_{2k}$; $B_{2k}=x_0.B_{2k-1}+C_{2k}$)
\enre
It concludes the proof of Lemma \mrf{l:alp}.
Q.E.D
\vskip0.1in

-Let us describe the set ${\mathcal A}'$
associated to ${\mathcal P}$ and ${\mathcal A}$.
To this end, we shall study the pencil $x_0^{2k}.(b_{2k}(x_0,x_1,x_2)-c)$
$c \in \RRR$.\\

Let us recall, for sake of clarity, some properties of this pencil
(see Lemma \mrf{l:lp3}).
\vskip0.1in
Let $c_M$ be the maximal critical value,
let $c_m$ be the minimal critical value.
Denote by $M_c = \{ x \in \RRR^2 | b_{2k} >0 \}$
For $c > c_M$, $M_c$ has $2k$ connected components.
For $c < c_m$, $M_c$ is  connected.\\
Two components of the connected components can coalescence
only when $b_{2k}$ passes a critical value $1$.
Bringing together this statement with Lemma \mrf{l:lp2},
and the fact that any curve of the pencil over $\HH_{2k}$
intersects  the line $x_0=0$ in $2k$ points,
it follows easily
that an arc of ${\mathcal A}$ bounding  a positive region
may not touch an arc of ${\mathcal A}$
bounding  a  negative region as $c$ varies from $0$.\\

\vskip0.1in
Let us prove the following Lemma.

\bele
\mlb{l:iny}
As $c$ varies from $0$ to $+ \infty$ or $- \infty$
any two arcs $\gamma, \gamma'$ of ${\mathcal A}$
may not be deformed in such a way
that for $c=c_0$ they
touch each other.
\enle
{\bf proof:}
Let $\HH_{2k}$ be a curve of type $\HH^0$ relatively to $L$
and $B_{2k}$ its polynomial.
The proof of Lemma \mrf{l:iny}
is based on a study of $\RRR \HH_{2k}$ on  a tubular neighborhood
${\mathcal M}_{\e}$  of $\RRR L$ and in particular the
two following  facts.\\
According to
Lemma \mrf{l:alp}, given a curve $\HH_{2k}$
of type $\HH^0$
relatively to $L$, one can push by a rigid isotopy
the set ${\mathcal P}$ associated to ${\mathcal A}$ to a line
of the tubular neighborhood ${\mathcal M}_{\e}$
of $\RRR L$ where $\e$ positive is arbitrarily small.\\
According to the Petrovskii's theory, \cite{Petr},
up to slightly modify coefficients of $B_{2k}$
without changing topological structure of $\HH_{2k}$
one can assume that none of its critical point belongs
to the line at infinity.
\vskip0.1in
Choose ${\mathcal M}_{\e}$ the tubular neighborhood of $\RRR L$
such that  any line of${\mathcal M}_{\e}$ intersects the non-empty
oval of $\HH_{2k}$ in $2k$ points.
According to proof of Lemma \mrf{l:alp}, one can
assume that critical points of $B_{2k}$  of the set ${\mathcal P}$
belong to $L_{\e_1}$ where $\RRR L_{\e_1} \subset {\mathcal M}_{\e}$
and $\e >0$ is arbitrarily small.
In this way,
as $c$ increases two arcs $\gamma, \gamma'$ of ${\mathcal A}$
bounding positive regions
are deformed and  for $c=c_0$ they
touch each other in a point of the line
$\RRR L_{\e_1} \subset {\mathcal M}_{\e}$.
The connected set ${\mathcal A}$ is such that
these two arcs are separated from each other
by an arc of ${\mathcal A}$ which bounds
a negative region. As $c$ increases from $c_0$ to
$c_0 + \e$, a negative oval ${\mathcal O}$  appears.
The oval ${\mathcal O}$ does not intersect $\RRR L_{\e_1}$,
it intersects the line at infinity $\RRR L$ in two points.
Up to slightly modify coefficients of $B_{2k}$
without changing topological structure of $\HH_{2k}$
one can consider  $\RRR L_{\e_1}$ as the line at infinity of $\RRR P^2$.
The creation of this negative oval leads to contradiction
with the number of critical points of $B_{2k}$ and also with
Petrovskii's Lemma \mrf{l:lp3}.
Indeed, as $c$ grows, this negative oval
shrinks in point and then disappears;
the number of intersection points of $\HH_{2k}$
with the line $L_{\e_1}$ decreases from $2k$ to $2k-2$.\\
An argumentation analogous to the previous one may be used
to prove that  as $c$ increases from $0$ to $+\infty$,
two arcs bounding a negative region may not be deformed
and touch each other.\\
Moreover,  when consider the Petrovskii's pencil over
$b_{2k}$, it follows easily that one arc bounding  a positive region
may not touch an arc bounding  a  negative region
as $c$ varies from $0$.\\
This proves Lemma \mrf{l:iny}.\\
Q.E.D

\vskip0.1in
We may describe, according to
Petrovskii's Lemmas \mrf{l:lp2} and  \mrf{l:lp3},
critical points of index $1$
associated  to arcs of ${\mathcal A}$ as follows.
There exist $2k-3$ critical points with critical value $c_0 > 0$
such that for $c=c_0$, an arc $\gamma$ of ${\mathcal A}$
touches an oval.-(Each of the $(k-1)$ arcs bounding
negative region touches a positive oval,
the $(k-2)$ arcs which bound positive region
touch negative ovals)-\\
There exist $2$ arcs of ${\mathcal A}$ bounding a positive
region $\{x \in \RRR^2 | b_{2k}(x)>0 \}$
for which  there exists $c_0 < 0$
with the following property.
As $c$ decreases from $0$ to $c_0$, the region expands
and for $c=c_0$ the arc touches a positive oval.
(Note that, at this time of our proof, any positive oval
under consideration above  may be the non-empty one).
In such a way,
according to Morse Lemma
since any point $p \in {\mathcal P}$ is a critical point
of index $1$ of $B_{2k}(x_0,x_1,x_2)$,
any arc $\gamma$ of ${\mathcal A}$
is associated to an  arc $\gamma'$
of $\RRR \HH_{2k}$.
We shall denote by ${\mathcal A}'$ the set of arcs $\gamma'$.
\vskip0.1in
Denote by ${\mathcal M}$ a tubular neighborhood of the line at infinity.
As already noticed, arcs ${\mathcal A} \cap \RRR^2$
may be seen (up to isotopy)
as arcs of $\RRR L_1^+$ and $\RRR L_2^+$.
In the same way,
arcs $\RRR \HH_{2k}$ of ${\mathcal A}'$
may be pushed on the boundary and the inside of ${\mathcal M}$.
Any oriented arc of the $k-1$ arcs of
$\RRR L_2^+$ is associated (to a critical point with positive critical value)
to an arc of $\RRR L_1^+$, any oriented arc of the $k-3$ arcs
of $\RRR L_1^+$ associated
to a critical point with positive critical value
is associated to an arc of $\RRR L_2^+$.
(According to Petrovskii's theory,
as $c$ increases, these arcs recede from $\RRR L$. )
These arcs are associated to a critical point
with positive critical value.
The two other arcs of $\RRR L_1^+$ are
associated to a critical point with negative critical value
and to an arc lying in the inside of ${\mathcal M}^-$.
These two arcs contain an extremity of $\RRR L^+$.
(Let $(\gamma, \gamma') \in ({\mathcal A}, {\mathcal A}')$
be such a pair of associated arcs.
According to Petrovskii's theory,
as $c$ decreases, the arc ${\gamma}$
comes closer to
$\RRR L \subset {\mathcal M}$,
the arc ${\gamma}'$
comes closer to an extremity of $\RRR L^+$.)
\vskip0.1in

-Let us extend this description  of $\RRR \HH_{2k}$
to a description in the whole ${\mathcal M}$.

Let ${\mathcal A}$  be a curve and $A$ its polynomial.
We shall say that the truncation $D$  of $A$
is  {\it sufficient} for ${\mathcal C} \subset {\mathcal A}$
in $\RRR P^2$ if the curve ${\mathcal D}$ with polynomial $D$
is such that :
$\RRR {\mathcal D}$ is embedded in $\RRR {\mathcal A}$
and
$\RRR {\mathcal D}$ is homeomorphic to
a subset  of $\RRR {\mathcal A}$
which contains ${\mathcal C}$.
\vskip0.1in

Let us prove the following Lemma:

\bele
\mlb{l:MorsTr}
Let $B_{2k}(x_0,x_1,x_2)$ be a Harnack polynomial of degree $2k$
and type $\HH^0$.
Up to regular modification of $B_{2k}(x_0,x_1,x_2)$,
\been
\item
\lb{i:tru1}
the truncation of $b_{2k}^{\D}(x,y)$ on
monomials of homogeneous degree $2k-2 \le i \le 2k$ is
sufficient for
${\mathcal A} \cup {\mathcal A}' $.
\item
\lb{i:tru2}
there exists a truncation  (on four monomials
$x^cy^d, x^{c+1}y^d, x^cy^{d+1},
x^{c+1}y^{d+1}$ with $c+d=2k-2$)
 $B_{2k}^{S}$ of $B_{2k}^{\D}$ for which
 $B_{2k}^{S}$ is sufficient for
  $\gamma \cup \gamma'$
 where
  $\gamma  \in {\mathcal A} $
  $\gamma' \in  {\mathcal A}' $.
\enen
\enle

Denote by $\HH_{2k}^{\D}$ the curve with polynomial $B_{2k}^{\D}$.
-Roughly speaking, it means that monomials which are not in $b_{2k}^{\D}$
(resp, $b_{2k}^{S}$)
have a small influence on
${\mathcal A} \cup {\mathcal A}'
\subset  \RRR \HH_{2k}$
(resp,
$\gamma \cup \gamma' \subset  \RRR \HH_{2k}$)-
\vskip0.1in

{\bf proof:}
The proof of Lemma \mrf{l:MorsTr} is based on Lemma \mrf{l:alp}.

According to the proof of Morse Lemma, local coordinates defined
in a neighborhood $U(p)$ of a non-degenerate critical point $p$
of a function $f$ depend principally on the first
derivative and the second derivative
of the function around this point.
Let $L_{\e_1}$ be  a real projective line in the
 $\e$-tubular neighborhood of $x_0=0$, where $\e>0$ is arbitrarily small.
On the assumption that any critical point $p$ belongs to the line $L_{\e_1}$,
monomials which are not in $b_{2k}^{\D}$
have a small influence on
$\{(x_0,x_1,x_2) \in \RRR P^2 | B_{2k}(x_0,x_1,x_2)=0 \} \cap U(p)$.
Without loss of generality, we may assume that
points at infinity of $B_{2k}$ and points of ${\mathcal P}$
do not belong $x_2=0$.
In such a way, any  point at infinity of an arc
of ${\mathcal A}$ is a root of the polynomial in one variable
$b_{2k}(0,x_1/x_2,1)=b_{2k}(y)$.
Any  arc of ${\mathcal A}$ intersects the line $x_0=0$
in two points. Up to regular modification
the $2k-1$ points of ${\mathcal P}$  belong to $L_{\e_1}$.
According to Rolle's Theorem,
$b_{2k} (y)$ has exactly $2k-1$ extrema at which
$b_{2k}' (y)$ must change from positive to negative.
Thus, $b_{2k}"(y)$ has  $2k-2$
zeroes.
By continuity, it gives the sign of the function $b_{2k}(x_0,x_1,x_2)$
in a neighborhood
${\mathcal M}=\{ (x_0:x_1:x_2) \in \RRR P^2 | |x_0|<\e \}$
of $x_0=0$.
In a neighborhood $U(p)$ of an extremum,
$\{ x \in \RRR P^2 | b_{2k}=0  \} \cap U(p)$
is described as a desingularized crossing.
From the previous study, it follows that up to regular modification
of $B_{2k}$, properties \ref{i:tru1}.
and \ref{i:tru2}. are verified.
One can also notice that the alternation of sign of $B_{2k}$
in the inside of ${\mathcal M}$ is equivalent to the modular property of the
distribution of sign of $b_{2k}^{\D}$
in the patchworking construction.
Set ${\mathcal N}$  tubular neighborhood of $L$ in $\CCC P^2$
with real part ${\mathcal M}$.
It is easy to verify that
replacing ${\mathcal N}$ in place of
${\mathcal M}$
one get  a general  formulation of
the Lemma \mrf{l:MorsTr} in $\CCC P^2$
This last remark concludes and completes our proof.
Q.E.D

\vskip0.1in
\vskip0.1in
According to the Lemma \mrf{l:MorsTr} the previous description
of $\RRR \HH_{2k}$ in $U(p)$ enlarges to a description of
$\RRR \HH_{2k}$
in $U \supset U(p)$ where
$$U=\{ z= <u,p>=(u_0.p_0:u_1.p_1:u_2.p_2) \in \CCC P^2 |
u=(u_0:u_1:u_2) \in U_{\CCC}^3, p=(p_0:p_1:p_2) \in U(p) \}$$
Local coordinates $y_1,y_2$ in $U(p)$ extend to local
coordinates in $U$ as follows.
Given $z=<u,p> \in U$ with
$u=(1:u_1:u_2) \in U_{\CCC}^3$, $p=(p_0:p_1:p_2) \in U(p) \}$,
we may set
$y_1(z)=u_1.y_1(p)$, $y_2(z)=u_2.y_2(p)$.
In $U(p)$,
the truncation $b_{2k}^{S}$ is $\e$-sufficient for $b_{2k}$.\\
$b_{2k}^S(x,y)=l(x,y)+\e.k(x,y)$\\
with
$l(x,y)=
 a_{c,d}x^cy^d  +
 a_{c+1,d}x^{c+1}y^d $,
$k(x,y)=
 a_{c,d+1}x^cy^{d+1}  +
-a_{c+1,d+1}x^{c+1}y^{d+1}$,
(with $a_{c,d}>0, a_{c+1,d}>0, a_{c,d+1}>0, a_{c+1,d+1}>0$
and $c+d=2k-2$ and $\e>0$)
\vskip0.1in

Note that up to modify the coefficients
$a_{c,d},a_{c,d+1},a_{c+1,d},a_{c+1,d+1}$ if necessary,
the point $p= (x_0,y_0)$ is
(up to homeomorphism)
a critical point of the function
$\frac {l(x,y)} {k(x,y)}$ with positive critical value.
Hence, it follows from the equalities \\
$l(x,-y)=-l(x,y)$, $k(x,-y)=k(x,y)$,\\
$\frac {\pr l} {\pr x}(x,-y) =-\frac {\pr l} {\pr x}(x,y)$,
$\frac {\pr l} {\pr y} (x,-y) =\frac {\pr l} {\pr y} (x,y)$,\\
$\frac {\pr k} {\pr x} (x,-y)= \frac {\pr k} {\pr x} (x,y)$,
$\frac {\pr k} {\pr y} (x,-y) = -\frac {\pr k} {\pr y} (x,y)$\\
that $(x_0,- y_0)$ is also a critical  point of
the function $\frac {l(x,y)} {k(x,y)}$ with negative critical value.
\vskip0.1in

The tubular neighborhood of $\RRR L^+$
is ${\mathcal M}^+= \cup_{p \in {\mathcal P}} U(p)$.
The transformation $(x,y) \to (-x,y)$
defined locally inside any $U(p)$, $p \in {\mathcal P}$
maps the set of
arcs ${\mathcal A}\cup{\mathcal A}'$ to a set of
arcs ${\mathcal B}$ and the set of points ${\mathcal P}$
to a set of points ${\mathcal S}$.
Any open $U(p)$ is mapped to an open $U(s)$
in such a way that
the tubular neighborhood of $L^-$ is
${\mathcal M}^-= \cup_{s \in {\mathcal S}} U(s)$
and $L^+$ is mapped to $L^-$.
\vskip0.1in
Precisely, given $x,y$ local coordinates  in
a neighborhood $U(p)$,
the pair $(\gamma,\gamma') \in ({\mathcal A},{\mathcal A}')$
of $\{ (x_0:x_1:x_2) \in \RRR P^2 | B_{2k}(x_0,x_1,x_2)=0 \} \cap U(p)$
is up to homeomorphism
the hyperbole $x.y=1/2$.
By means  of the transformation
$(x,y) \to (x,-y)$
one maps  $U(p)$ onto $U(s)$
and the pair of arcs
$(\gamma,\gamma') \in ({\mathcal A},{\mathcal A}')$
defined up to homeomorphism by $x.y=1/2$
to the  pair of arcs
$(\xi,\xi')$
of $\{ (x_0:x_1:x_2) \in \RRR P^2 | B_{2k}(x_0,x_1,x_2)=0 \} \cap U(p)$
is up homeomorphism the hyperbole $x.y=-1/2$.
(In $U(s)$,
one can say that
arcs  $\xi$,$\xi'$ are {\it associated}
and also associated to $s$.)
\vskip0.1in

Let us in Lemma \mrf{l:cntb}
describe the set
${\mathcal D}=
{\mathcal M}^- \cap \RRR \HH_{2k}^{\D}$.

\bele
\mlb{l:cntb}
Let $\HH_{2k}$ be a Harnack curve of type $\HH^0$
with polynomial $B_{2k}$.
Up to regular modification of the Harnack polynomial $B_{2k}$,
$\RRR \HH_{2k}^{\D}$ intersects ${\mathcal M}^-$ in
$2k-2$ positive ovals.
\enle

{\bf proof:}
\vskip0.1in
The set ${\mathcal A}$ is  a connected set of arcs
which intersects $\RRR L$ the line at infinity of $\RRR P^2$ in $2k$ points
which are intersection points of two arcs of ${\mathcal A}$.
Using Morse Theory, one can consider these points as
limit points of an half hyperbole.
Up to isotopy of $\RRR P^2$,
one can consider
arcs of ${\mathcal A} \cup {\mathcal A'}$
as arcs lying on $\RRR L_1^+$
and $\RRR L_2^+$.

Except the two arcs pushed by an isotopy of $\RRR P^2$
to a segment of $\RRR L_1^+$ which contains an extremity of $\RRR L_1^+$,
each arc of ${\mathcal A}$ is associated to an arc of ${\mathcal A'}$
in such a way that
the set ${\mathcal A} \cup {\mathcal A}'$
is the union of two connected
parts lying respectively on $\RRR L_1^+$ and $\RRR L_2^+$.
\vskip0.1in
Each of the two arcs which contains an extremity
of $\RRR L_1^+$ is associated
to an arc  ${\mathcal A'}$ which belongs to the part  ${\mathcal M}^-$
of the M\"{o}bius band lying between $\RRR L_1^-$ and $\RRR L_2^-$.
In our identification of arcs with
connected real parts lying in the inside of ${\mathcal M}$,
one can assume that these two pairs
remain the same under the action $(x,y) \to (x,-y)$
which maps arcs ${\mathcal A} \cup {\mathcal A}'$ to
arcs of ${\mathcal B}$.
(Indeed, ${\mathcal P}$ intersects ${\mathcal S}$
in two points which are extremities of $\RRR L_{\e_1}^+ \approx
\RRR L$.)

It follows that the union of
$({\mathcal A} \cup {\mathcal A}')$
is mapped to $2(k-1)$ ovals constituted
by the union of arcs ${\mathcal B}$ union two arcs of ${\mathcal M}^+$
which contains an extremity of $\RRR L_1^+$.

These ovals lie in the inside
of the M\"{o}bius band
${\mathcal M}^-$
and intersects $\RRR L^-$ in two points.
This last remark concludes our proof.
Q.E.D

\vskip0.1in
\bere
It follows from the Lemma \mrf{l:cntb}
a description, up to isotopy of $\RRR P^2$,
of $\RRR \HH_{2k}$ in ${\mathcal M}$
which extends our previous description of
${\mathcal A} \cup {\mathcal A}' \subset \RRR \HH_{2k}$
Note that this description is consistent
with the study of the one variable
function  $b_{2k}(y)=b_{2k}(0,x_1/x_2,1)$.
The function $b_{2k}"(y)$ has $2k-2$ zeroes.
Any  point at infinity of an arc
of ${\mathcal A}$ is a root of the polynomial in one variable
$b_{2k}(0,x_1/x_2,1)=b_{2k}(y)$.
According to Rolle's Theorem,
$b_{2k} (y)$ has exactly $2k-1$ extrema at which
$b_{2k}' (y)$ must change from positive to negative.
The zeroes of $b_{2k}"(y)$
are simple zeroes.
These zeroes are the intersection-points of  the set
${\mathcal B}$ with $L$.
\enre

\vskip0.1in

The set ${\mathcal A}$ is the image of a smooth
section of a tubular fibration of $\RRR L^+$ minus 2 points.
In Lemma \mrf{l:ist},
we shall prove that
there exists a subset ${\mathcal D}^1$
of
 ${\mathcal D} = \RRR \HH_{2k} \cap {\mathcal M}^-$
with the property that
the set
${\mathcal A} \cup {\mathcal D}^1$ is the image of a smooth
section of a tubular fibration of $\RRR L$ minus a finite  number of points.
\vskip0.1in

By cutting
${\mathcal M}^{-}$ along
$L^-$, one get two surfaces.
Denote by ${\mathcal M}^{1,-}$
the one which contains $L_1^-$.
Consider the intersection
$\RRR \HH_{2k} \cap {\mathcal M}^{1,-}$

Any arc $\xi$ of
${\mathcal D} =\RRR \HH_{2k} \cap {\mathcal M}$
intersects the line at infinity $L$.
Hence, it
is divided into two halves
with common point $L \cap \xi$.
Denote by
$\xi^1$ (resp, $\xi^2$)
the half of $\xi$ which belongs
to the inside of M\"{o}bius delimited by $\RRR L$ and $\RRR L_1$
(resp, $\RRR L$ and $\RRR L_2$)
minus its intersection with $L$.
The set
$\RRR \HH_{2k} \cap {\mathcal M}^{1,-}$
is the set of arcs $\xi^1$.

In Lemma \mrf{l:ist} we prove that
${\mathcal D}^1=\RRR \HH_{2k} \cap {\mathcal M}^{1,-}$

\bele
\mlb{l:ist}
Let $\HH_{2k}$ be a Harnack curve of type $\HH^0$
with polynomial $B_{2k}$.\\
Up to regular modification of the Harnack polynomial $B_{2k}$,
the set
${\mathcal A} \cup {\mathcal D}^1$,
where
${\mathcal D}^1 =\RRR \HH_{2k}^{\D} \cap {\mathcal M}^{1,-}$,
is the image of a smooth
section of a tubular fibration of $\RRR L$ minus $2k-1$ points.
\enle

{\bf proof:}
The set ${\mathcal A}$ is connected.
There exists an isotopy of $\RRR P^2$
which pushes the union of the $k$ arcs of ${\mathcal A} \cap \RRR^2$
bounding positive regions $\{ x \in \RRR^2 | b_{2k}(x) > 0 \}$
with $k-1$ arcs of ${\mathcal A}'$
(associated to the $k-1$  arcs of ${\mathcal A} \cap \RRR^2$
bounding negative regions$\{ x \in \RRR^2 | b_{2k}(x) > 0 \}$)
onto  $\RRR L_1^+$.
Let us denote by ${\mathcal A}^1$ this set of
$k$ arcs of ${\mathcal A}$ with $(k-1)$ arcs
of ${\mathcal A}'$ arcs.
The transformation $(x,y) \to (x,-y)$
maps  $U(p) \to U(s)$
and the pair of arcs
$(\gamma,\gamma') \in ({\mathcal A},{\mathcal A}')$
defined up to homeomorphism by $x.y=1/2$
to  the  pair of arcs
$(\xi,\tilde{\xi})$ of ${\mathcal B}$
defined up to homeomorphism by $x.y=-1/2$.
Hence, identifying each arc of ${\mathcal A}^1$,
with a segment of $\RRR L_1^+ \cap U(p)$,
it follows the next descriptions in $U(s)$.
Let us first consider an open $U(s)$ which does not contain
an extremity of $L$.
According to Morse Lemma
in $U(s)$, it follows that the pair
$(\xi^1,\tilde{\xi}^1)$ of halves of associated arcs in $U(s)$
is the image of a smooth
section of a tubular fibration of
$U(s) \cap \RRR L_1^-$ minus one
point.

If $U(s)$ contains an extremity of $L$,
then the pair of arcs
$(\xi^1,
 \tilde{\xi}^1)=(\xi^1 \cap {\mathcal M}^-,
 \tilde{\xi}^1 \cap {\mathcal M}^+)$
   of halves of associated arcs in $U(s)$ is such that
${\xi}^1$ is the image of a smooth
section of a tubular fibration of $U(s) \cap \RRR L_1^-$ minus the
extremity of  $\RRR L_1^- \cap U(s)$.\\
By means of the transformation $(x,y) \to (x,-y)$,
inside any open $U(p)$,
the tubular neighborhood
${\mathcal M}^+= \cup_{p \in {\mathcal P}} U(p)$
 of $L^+$ is mapped to
 the tubular neighborhood
${\mathcal M}^-= \cup_{s \in {\mathcal S}} U(s)$
 of $L^-$.
It follows  from Morse Lemma that
the set ${\mathcal D}^1$
 is the image of a smooth
section of tubular fibration of $\RRR L_1^-$ minus $2k-1$ points.

Hence,
the set ${\mathcal D}^1$
 is also the image of a smooth
section of tubular fibration of $\RRR L^-$ minus $2k-1$ points.

Since
${\mathcal P}$ intersects ${\mathcal S}$ in the extremities
of $\RRR L_{\e_1}^+$ (see proof of Lemma \mrf{l:cntb}),
two of these $2k-1$ points of $\RRR L^-$
are extremities.

Hence, ${\mathcal D}^1$
is the image of a smooth
section of tubular fibration of $\RRR L^-$ minus ${\mathcal S}$.
Therefore, since ${\mathcal A}$
is the image of a smooth
section of tubular fibration of $\RRR L_+$ minus its extremities,
${\mathcal A} \cup {\mathcal D}^1$
is the image of a smooth
section of tubular fibration of $\RRR L$ minus
$2k-1$ points.
Q.E.D
\vskip0.1in

According to Proposition \mrf{p:crit2},
given $B_{2k}$ and $B_{2k-1}$  Harnack polynomials of
type $\HH^0$ and respective degree $2k$ and $2k-1$
$$c(B_{2k})=(\frac {(k-1) (k-2)}{2},k.(2k-1),\frac {3k.(k-1)} {2})) $$
 $$c(B_{2k-1})=(\frac {(k-1).(k-2)} 2,
(k-1).(2k-1),\frac {(k-1).(3k-4)} {2})$$
The Harnack curve $\HH_{2k}$ of degree $2k$
has
$2(k-2)$ ( $:=c_2(B_{2k})-c_2(B_{2k-1})$)
more positive ovals
than
the Harnack curve $\HH_{2k-1}$ of degree $2k-1$.
It follows from Lemma \mrf{l:MorsTr} and Lemma \mrf{l:cntb}
that, up to regular to modification of the Harnack
polynomial $B_{2k}$,
these ovals are the ovals of $\RRR \HH_{2k}^{\D} \cap {\mathcal M}^-$.
\vskip0.1in
In such a way, it follows from Lemma \mrf{l:ist}
that, up to regular modification of $B_{2k}$,
$B_{2k}=x_0.B_{2k-1} + C_{2k}$
where $B_{2k-1}$ is a Harnack polynomial of type $\HH^0$
relatively to $L:=x_0=0$ the line at infinity.
(Indeed, as we have noticed, the truncation $B_{2k}^{\D}$
is sufficient for $B_{2k}$ and $L_{\e_1}$ may be chosen
arbitrarily close to $L$.)
It concludes our proof of Proposition \mrf{p:stv}. {\bf Q.E.D}

\vskip0.1in

{\bf I.ii)-Harnack curves of odd degree-}\\
Let  $\HH_{2k+1}$ be a Harnack curve of degree $2k+1 \ge 5$ with
polynomial $B_{2k+1}(x_0,x_1,x_2)$
of type $\HH^0$ relatively to the line at infinity
$L:=x_0=0$.

We shall prove the following Proposition \mrf{p:sto}.
\vskip0.1in
\bepr
\mlb{p:sto}
Let  $\HH_{2k+1}$ be a Harnack curve of degree $2k+1$ with
polynomial $B_{2k+1}(x_0,x_1,x_2)$
of type $\HH^0$ relatively to the line at infinity
$L:=x_0=0$.
Up to regular modification of $B_{2k+1}$,
$B_{2k}(x_0,x_1,x_2)$ is also of type $\HH^0$ relatively to $L$.
\enpr
{\bf proof:}\\
The part I.ii) until its end is devoted to the proof of
Proposition \mrf{p:sto}.\\
Our argumentation is a slightly modified version of the one given
in the proof of Proposition \mrf{p:sto}.
When it is possible, we refer to the first part I.i).\\

The odd component of $\HH_{2k+1}$ is divided
into $2k+1$ arcs which delimit a region bounding the
line at infinity. One of this region contains ovals of
$\HH_{2k+1}$.
Denote ${\mathcal A}$
the set of $2k$ arcs of $\RRR \HH_{2k+1}$ which delimit a positive region
$\{ x \in \RRR^2 | b_{2k+1} >0 \}$
(i.e region which does not contain oval of $\HH_{2k+1}$).
The union of arcs of ${\mathcal A}$ is a connected orientable part of
the odd component of $\RRR \HH_{2k+1}$.
The odd component of $\RRR \HH_{2k+1}$ is homeomorphic to a projective line.
Consider an orientation of the line at infinity $\RRR L$.
The surface delimited by
the $2k$ arcs of ${\mathcal A}$ and $\RRR L$ is orientable.
Thus, one can assume that
the orientation of the connected
union of the $2k$ arcs of ${\mathcal A}$ is induced by a half of $\RRR L$.
We shall denote by $\RRR L^+$ this half.
The set ${\mathcal A}$ is the image of a smooth
section of a tubular fibration of $\RRR L^+$ minus 2 points.

\vskip0.1in
Let us prove that one can associate
to the set ${\mathcal A}$ of the $2k$ arcs
a set of critical points ${\mathcal P}$
of index $1$ of $b_{2k+1}$ and a set of arcs ${\mathcal A}'$ of
$\RRR \HH_{2k+1}$. \\
To this end, we shall study the Petrovskii's pencil
$x_0^{2k+1}.(b_{2k+1}(x_1/x_0,x_2/x_0)-c)$, $c \in \RRR$,
over  $\HH_{2k+1}$.\\
As $c$ decreases from $0$
positive regions and positive ovals expand.
As $c$ decreases from $0$
the number of regions $G_i$ of $\{x \in \RRR^2 | b_{2k}(x) > 0 \}$
which contain a segment of the line at infinity of $\RRR P^2$ on their
boundary decreases from $k+1$ to $1$.\\
As $c$ increases from $0$  negative ovals disappear;
the number of regions $G_i$ of $\{ x \in \RRR^2 | b_{2k+1}(x) > 0 \}$
which contain a segment of the line at infinity of $\RRR P^2$ on their
boundary remains the same.
Hence, to each arc of ${\mathcal A}$ is associated a critical of index $1$.\\

-Let us study the set ${\mathcal P}$ of critical point $1$ associated to
the set  ${\mathcal A}$ of arcs of $\RRR \HH_{2k}$.
\vskip0.1in

Let us prove the following Lemma:
\bele
\mlb{l:inyo}
There exist at most two arcs of the set ${\mathcal A}$
which may be deformed
as $c$ varies from $0$ to $- \infty$ or $+ \infty$
in such a way that for $c=c_0$
they touch each other.
\enle
{\bf proof:}
The projective plane $\RRR P^2$  may be seen
as a M\"{o}bius band ${\mathcal M}$ and a disc $D^2$
glued along their common boundary where the core of the M\"{o}bius
band is the line at infinity $\RRR L$.
Assume that there exists $c_0$,
such that as $c$ varies from $0$ to $c_0$
an arc $\gamma$ of ${\mathcal A}$ is deformed
and for $c=c_0$  it touches an oval ${\mathcal O}$ of $\HH_{2k+1}$.
The arc $\gamma$ bounds positive region .
Hence, if $c_0$ exists it  is negative.
The inside of the oval ${\mathcal O}$ is an orientable part of $\RRR P^2$
homeomorphic to $D^2$.
The arc $\gamma$ may be pushed by an isotopy to the boundary $\RRR L_1^+$
of the M\"{o}bius band.
In this way, for $c=c_0$, $\gamma$ may be also pushed  to
an arc  $\tilde {\gamma}$
of the oval ${\mathcal O}$.
Since $c_0 <0$, according to Petrovskii's Lemmas,
for $c=c_0 + \e$ one can trace a non-orientable branch (i.e
an arc of $\RRR L$ which intersects $\RRR L^+$ and $\RRR L^-$)
of $\RRR P^2$
in the region $b_{2k} < c_0-\e$ containing $\gamma$ and $\tilde{\gamma}$.
The set ${\mathcal A}$ may be pushed by an isotopy
onto $\RRR L_1^+$.
In a neighborhood of $\RRR L_1^+$,
a non-oriented branch of $\RRR P^2$ may be traced
only near the extremities of $\RRR L_1^+$.
Hence, the only  two arcs of ${\mathcal A}$
which may glue with an oval are those which intersect
an extremity of $\RRR L_1^+$.
Q.E.D
\vskip0.1in
It follows from  Petrovskii's Lemmas \mrf{l:lp2}, \mrf{l:lp3}
and Lemma \mrf{l:inyo},
that there exist $2k-2$ critical points with critical value $c_0 < 0$
such that as $c$ varies from $0$ to $c_0$
two arcs are deformed in such a way that
for $c=c_0$ they touch each other.

As already done for Harnack curves of even degree
one can study critical points ${\mathcal P}$
of $B_{2k+1}$
using the function in one variable
$b_{2k+1}(y)=b_{2k+1}(0,x_1/x_2,1)$.
According to Rolle's Theorem
applied to
$b_{2k+1}(y)=b_{2k+1}(0,x_1/x_2,1)$
for any  pair ($\gamma, \gamma')$ of associated arcs  of ${\mathcal A}$
there exists $\gamma" \in {\mathcal A}$
such $\gamma \cap \gamma" \not= \emptyset$
$\gamma' \cap \gamma"\not= \emptyset$.

According to Lemma \mrf{l:inyo} and its proof,
for the $2$ arcs of ${\mathcal A}$ which intersect
$L$ in one of its extremities,
there exists $c_0 < 0$
such that as $c$ varies from $0$ to $c_0$, the arc is deformed
and for $c=c_0$ it touches a positive oval.

In such a way, any arc $\gamma$ of ${\mathcal A}$
is associated to a critical point $p$ and to an other arc $\gamma'$.
We denote by ${\mathcal P}$
the set of points associated to ${\mathcal A}$
and by ${\mathcal A}'$
the set of arcs associated to ${\mathcal A}$.
The sets ${\mathcal A}'$ and ${\mathcal A}$
consist of $2k$ arcs of $\HH_{2k+1}$; $2k-2$ of them
belong to ${\mathcal A}$ and ${\mathcal A}'$.

\vskip0.1in
-Let us give a description up to rigid isotopy
 of $\RRR \HH_{2k}^{\D}$ in the whole ${\mathcal M}$.
\vskip0.1in

Let $B_{2k+1}$ be a Harnack polynomial of degree $2k+1$
and type $\HH^0$.

From a version of Lemma \mrf{l:alp} for odd degree curves $\HH_{2k+1}$,
it follows that there exists a rigid isotopy  $B_{2k+1,t}$
,$t \in [0,1]$, $B_{2k+1,0}=B_{2k+1}$,
such that the modified polynomial $B_{2k+1,1}$
has the following the property:
the $2k-1$ critical points ${\mathcal P}$
associated to ${\mathcal A}$
belong to a line $L_{\e_1}$
of ${\mathcal M}_{\e}$ the $\e$-tubular
neighborhood of $L$ in $\RRR P^2$ with $\e >0$ arbitrarily small.\\
From a version of Lemma \mrf{l:MorsTr} for odd degree,
the truncation $b_{2k+1}^{\D}(x,y)$ of $b_{2k+1}(x,y)$
on monomials of homogeneous degree $2k-1 \le i \le 2k +1$ is
sufficient for
${\mathcal A} \cup {\mathcal A}'$.\\
Denote by $\HH_{2k+1}^{\D}$ the curve with polynomial
$B_{2k+1}^{\D}(x_0,x_1,x_2)=x_0^{2k+1}.b_{2k+1}^{\D}(x,y)$
the description
of $\RRR \HH_{2k+1}^{\D}$ in any neighborhood  $U(p)$ $p \in {\mathcal P}$
enlarges to a description of $\RRR \HH_{2k+1}^{\D}$
in $U \supset U(p)$ where
$$U=\{ z= <u,p>=(u_0.p_0:u_1.p_1:u_2.p_2) \in \CCC P^2 |
u=(u_0:u_1:u_2) \in U_{\CCC}^3, p=(p_0:p_1:p_2) \in U(p) \}$$\\
\vskip0.1in

Consider ${\mathcal M}^+$ the tubular neighborhood of $L^+$,
${\mathcal M}^+= \cup_{p \in {\mathcal P}} U(p)$
and
${\mathcal M}^-$ the tubular neighborhood of $L^-$.

From an argumentation analogous to the one given in the part I.i),
according to Lemma \mrf{l:MorsTr}
we define a set ${\mathcal B}$ of arcs of $\RRR \HH_{2k+1}^{\D}$
in ${\mathcal M}^-$.

Let $(x,y)$ be local coordinates in $U(p)$
, $p \in {\mathcal P}$.
The transformation, defined locally in any $U(p)$
by $(x,y) \to (-x,y)$ maps $U(p) \to U(s)$
and in this way  ${\mathcal M}^+=
 \cup_{p \in {\mathcal P}} U(p)$
to ${\mathcal M}^- =
 \cup_{s \in {\mathcal S}} U(s)$.
It maps the set of
arcs ${\mathcal A}\cup{\mathcal A}'$ to a set of
arcs ${\mathcal B}$ and the set of points ${\mathcal P}$
to a set of points ${\mathcal S}$.
\vskip0.1in

Let us in Lemma \mrf{l:cntbo}
describe the set ${\mathcal D}=
{\mathcal M}^- \cap \RRR \HH_{2k+1}^{\D}$.

\bele
\mlb{l:cntbo}
Up to regular modification of the Harnack polynomial $B_{2k+1}$,
$\RRR \HH_{2k+1}^{\D}$ intersects ${\mathcal M}^-$
in $k-1$ negative ovals and
$k$ positive ovals.
\enle

{\bf proof:}
\vskip0.1in
Our argumentation is
similar to the one given in the proof of Lemma \mrf{l:cntb}.\\
The set ${\mathcal A}$ is a connected set of arcs which intersects
$L$ the line at infinity of $\RRR P^2$ in $2k$ points.
Using Morse Theory, one can consider intersection points
of the set ${\mathcal A}$ with the line $L$ as
limit points of an half hyperbole.
Except the two arcs of ${\mathcal A}$
which contain an extremity of $L^+$,
any arc of ${\mathcal A}$
is associated to another arc of ${\mathcal A}$.
Given $U(p)$ with local coordinates $(x,y)$,
the transformation $(x,y) \to (x,-y)$ maps
any pair of arcs of ${\mathcal A}$
to another pair of arcs $\RRR \HH_{2k+1}$.
It follows from the alternation of signs of $b_{2k+1}$
in any $U(p)$ of ${\mathcal M}^+$,
the alternation of sign $b_{2k}$
in any $U(s)$ of ${\mathcal M}^-$
and thereby around branches of $\RRR \HH_{2k+1}$.
It follows that
any two arcs of ${\mathcal A}$ associated to each other in $U(p)$
are mapped to an arc of a positive oval
and an arc of a negative oval in $U(s)$.
Each of the two arcs which contains an extremity
of $\RRR L^+$ is associated
to an arc  ${\mathcal A'}$ which belongs to
${\mathcal M}^-$.
One can assume that these two pairs
remain the same under the action $(x,y) \to (x,-y)$
which maps arcs ${\mathcal A} \cup {\mathcal A}'$ to
arcs of ${\mathcal B}$.
It follows that the union of
$({\mathcal A} \cup {\mathcal A}')$
is mapped to the union ${\mathcal B}$ of
two arcs of the odd component of $\HH_{2k+1}$
which contain extremities  of $\RRR L_1^+$
(these arcs lie in ${\mathcal M}^+$)
and
$k-1$ negative ovals and
$k$ positive ovals of $\HH_{2k+1}$
constituted
by the union of arcs ${\mathcal B}$.
These ovals lie in the inside ${\mathcal M}^-$
of the M\"{o}bius band
and intersect $\RRR L^-$ in two points.
This last remark concludes our proof.
Q.E.D
\vskip0.1in

\vskip0.1in

The set ${\mathcal A}$ is the image of a smooth
section of a tubular fibration of $\RRR L^+$ minus 2 points.
In Lemma \mrf{l:isto},
we shall prove that
there exists a subset ${\mathcal D}^1$
of
${\mathcal D} = \RRR \HH_{2k+1}^{\D} \cap {\mathcal M}^-$
with the property that
the set
${\mathcal A} \cup {\mathcal D}^1$ is the image of a smooth
section of a tubular fibration of $\RRR L$ minus a finite  number of points.
\vskip0.1in

By cutting
${\mathcal M}^{-}$ along
$L^-$, one get two surfaces.
Denote by ${\mathcal M}^{1,-}$
the one which contains $L_1^-$.
Consider the intersection
$\RRR \HH_{2k+1}^{\D} \cap {\mathcal M}^{1,-}$

Any arc $\xi$ of
${\mathcal D} =\RRR \HH_{2k+1}^{\D} \cap {\mathcal M}$
intersects the line at infinity $L$.
Hence, it
is divided into two halves
with common point $L \cap \xi$.
Denote by
$\xi^1$ (resp, $\xi^2$)
the half of $\xi$ which belongs
to the inside of M\"{o}bius delimited by $\RRR L$ and $\RRR L_1$
(resp, $\RRR L$ and $\RRR L_2$)
minus its intersection with $L$.
The set
$\RRR \HH_{2k}^{\D} \cap {\mathcal M}^{1,-}$
is the set of arcs $\xi^1$.

In Lemma \mrf{l:ist} we prove that
${\mathcal D}^1=\RRR \HH_{2k+1}^{\D} \cap {\mathcal M}^{1,-}$

\bele
\mlb{l:isto}
Up to regular modification of the Harnack polynomial $B_{2k+1}$,
the set
${\mathcal A} \cup {\mathcal D}^1$ is the image of a smooth
section of a tubular fibration of $\RRR L$ minus $2k$ points.
\enle

{\bf proof:}
Our proof is analogous to the one of Lemma \mrf{l:ist}.
The set ${\mathcal A}$  is the image of a smooth
section of a tubular fibration of $\RRR L^+$ minus 2 points
The transformation $(x,y) \to (x,-y)$
maps  $U(p) \to U(s)$
and the pair of arcs
$(\gamma,\gamma') \in ({\mathcal A},{\mathcal A}')$
defined up to homeomorphism by $x.y=1/2$
to  the  pair of arcs
$(\xi,\tilde{\xi})$ of ${\mathcal B}$
defined up to homeomorphism by $x.y=-1/2$.

Hence, identifying each of arcs of ${\mathcal A}$,
with a segment of $\RRR L^+ \cap U(p)$,
it follows the next descriptions in $U(s)$.

Let us first consider open $U(s)$ which does not contain
an extremity of $L$.
According to Morse Lemma
in $U(s)$, it follows that the pair
$(\xi^1,\tilde{\xi}^1)$ of halves of associated arcs in $U(s)$
is the image of a smooth
section of a tubular fibration of $U(s) \cap \RRR L_1^-$ minus one point.\\
If $U(s)$ contains an extremity of $L$,
then the pair of arcs
$(\xi^1,
 \tilde{\xi}^1)=(\xi^1 \cap {\mathcal M}^-,
 \tilde{\xi}^1 \cap {\mathcal M}^+)$
   of halves of associated arcs in $U(s)$ is such that
${\xi}^1$ is the image of a smooth
section of a tubular fibration of $U(s) \cap \RRR L_1^-$ minus the
extremity of $\RRR L_1^- \cap U(s)$.\\

By means of the transformation $(x,y) \to (x,-y)$,
inside any open $U(p)$,
the tubular neighborhood
${\mathcal M}^+= \cup_{p \in {\mathcal P}} U(p)$
 of $L^+$ is mapped to
 the tubular neighborhood
${\mathcal M}^-= \cup_{s \in {\mathcal S}} U(s)$
 of $L^-$.
It follows  from Morse Lemma that
the set ${\mathcal D}^1$
 is the image of a smooth
section of tubular fibration of $\RRR L^-$ minus $2k$ points.

Two of these $2k$ points of $\RRR L^-$ are  extremities
of $\RRR L^+$

Hence, ${\mathcal D}^1$
is the image of a smooth
section of tubular fibration of $\RRR L_-$ minus ${\mathcal S}$.
Therefore,
${\mathcal A} \cup {\mathcal D}^1$
is the image of a smooth
section of tubular fibration of $\RRR L$ minus
the $2k$ points.
Q.E.D
\vskip0.1in

According to Proposition \mrf{p:crit2},
given $B_{2k+1}$ and $B_{2k}$  Harnack polynomials of
type $\HH^0$ and respective degree $2k$ and $2k-1$
 $$c(B_{2k+1})=(\frac {(k).(k-1)} 2,
k.(2k+1),\frac {(k).(3k-1)} {2})$$
$$c(B_{2k})=(\frac {(k-1) (k-2)}{2},k.(2k-1),\frac {3k.(k-1)} {2})$$

The Harnack curve $\HH_{2k+1}$ of degree $2k+1$
has
$k-1$ negative ovals and
$k$ positive ovals more than $\HH_{2k}$.

It follows from the version of Lemma \mrf{l:MorsTr}
for odd degree curves and Lemma \mrf{l:cntbo}
that up to regular to modification of the Harnack
polynomial $B_{2k+1}$,
these ovals are the ovals of $\RRR \HH_{2k+1}^{\D} \cap {\mathcal M}^-$.
\vskip0.1in
In such a way, it follows from Lemma \mrf{l:isto}
that, up to regular modification of $B_{2k+1}$,
$B_{2k+1}=x_0.B_{2k} + C_{2k+1}$
where $B_{2k}$ is a Harnack polynomial of type $\HH^0$
relatively to $x_0=0$ the line at infinity.
It concludes our proof of Proposition \mrf{p:sto}.
{\bf Q.E.D}

{\bf II)}
For curves of degree $i \le 4$ the Proposition \mrf{p:ll} is immediate.
It is sufficient to use the fact that
Harnack curves of degree $i \le 4$
are rigidly  isotopy and are  $L$-curves.\\
(When consider the curve ${\mathcal B}_3$,
one can also give the following argument.
Assume ${\mathcal B}_4$ is a curve of type $\HH^0$
which results from deformation of ${\mathcal B}_3 \cup L$.
The curve ${\mathcal B}_3$ has necessarily two real connected components.
Noticing that none of the polynomials $x_0^3+a_1.x_1^3+a_2.x_2^3$
with $a_1, a_2 \in \RRR$
is a polynomial of a curve of degree $3$ with two real connected components,
we get the Proposition \mrf{p:ll} for ${\mathcal B}_3$.)
Q.E.D
\vskip0.1in
{\bf Conclusion: Proof of Theorem \mrf{t:rig}:}

Bringing together
Theorem \mrf{t:ty0} of the Appendix,
Proposition \mrf{p:stv}
Proposition \mrf{p:sto}, and this last result for curves of degree $\le 4$
we get by induction the Theorem \mrf{t:rig}.
Q.E.D
\vskip0.1in

The following remark may  easily deduced
from Theorem \mrf{t:rig}.

\bere
\mlb{r:rrig}
Let $L$ be a real projective line and $S$ a set of $m$ real
points lying on $L$.
Denote by $\HH_{m}'$ a curve of degree $m$
of type $\HH^0$ relatively to a real projective line $L'$.
There exists a rigid isotopy of $\RRR P^2$
which maps $\HH_{m}'$ of type $\HH^0$
to a curve $\tilde{\HH_{m}}$
of type $\HH^0$ relatively to $L$
such that $\tilde{\HH_{m}} \cap L = S$.
\enre
\subsection{ A generalization of Harnack's Method}
\vskip0.1in
Harnack curves $\HH_i$ of degree $i \le 4$ are $L$-curves.
Hence, up to rigid isotopy,
they result from the successive classical
deformation of the union $\HH_{j} \cup L_{j+1}$, $1 \le j \le 4$,
where $L_{j+1}$ is a real projective line which intersects $\HH_{j}$ in
$j$ real points.
In Theorem \mrf{t:bli}, we extend this  property
to Harnack curves $\HH_m$ of arbitrary degree.
Precisely, we prove that up to regular
modification of its polynomial
any Harnack curve results from the Harnack's method (i.e $L_{j+1}=L$
for any $j \ge 0$).
In this way,
the Rigid Isotopy Classification Theorem \mrf{t:rigiso}.
follows from
the Theorem \mrf{t:bli}.
\vskip0.1in

Given
$\HH_m$  a  Harnack curve with polynomial
$B_{m}(x_0,x_1,x_2)$,
we denote by ${\mathcal B}_i$ the curve
with polynomial
$B_i(x_0,x_1,x_2)= x_0^i.b_i(x_1/x_0,x_2/x_0)$
where $B_i(x_0,x_1,x_2)$ is
the homogeneous polynomial associated to
the truncation $b_i(x,y)$ of $b_m(x,y)$ on the monomials
$x^{\a}.y^{\b}$ with $0 \le \a + \b \le i $.\\
Let us prove the following Proposition:
\bepr
\mlb{p:l07}
Up to rigid isotopy, one can always assume that
${\mathcal B}_i$, $ 1 \le i \le 3$, is of type $\HH^0$
relatively to the line at infinity $L:=x_0=0$.
\enpr
{\bf proof :}\\
In the proof of the
Rigid Isotopy Classification Theorem \mrf{t:rigiso},
we shall use only the fact that
up to rigid isotopy, one can always assume that
${\mathcal B}_i$, $ 1 \le i \le 3$, is of type $\HH^0$
relatively to the line at infinity $L:=x_0=0$.
Our argumentation is based on the fact that
up to rigid isotopy there exists only one $M$-curve
of degree $3$.\\
The proof is based on two lemmas.
\vskip0.1in
{\bf I)}
Noticing that, up to regular modification of $B_m(x_0,x_1,x_2)$,
the truncation $B_3(x_0,x_1,x_2)$
is the polynomial of a smooth curve,
we shall prove that, up to regular modification of $B_m(x_0,x_1,x_2)$,
$B_3(x_0,x_1,x_2)$ is of type $\HH^0$.
Our argumentation is based on the following statement:\\
{\it A nonsingular cubic has exactly three real inflection points.
These inflection points are collinear.
\cite{Briesk}}\\

Denote by $s_i$, $0 \le i \le 2$, the real inflection points
of $B_3$ and by $L_i$, $0 \le i \le 2$, the tangent to ${\mathcal B}_3$
at $s_i$.

Let us prove the following Lemma:
\bele
\mlb{l:cub3}
For suitable projective coordinates,
and up to rigid isotopy of $B_m(x_0,x_1,x_2)$,
${\mathcal B}_3$
results from the classical deformation:\\
$B_3(x_0,x_1,x_2)=x_0.B_2(x_0,x_1,x_2)
+\e.L_0.L_1.L_2$
where $B_2(x_0,x_1,x_2)$ is of type $\HH^0$
relatively to $x_0=0$.\\
($L_i \not= x_0=0$,
$B_3(0,1,1)=\e.L_0.L_1.L_2(0,1,1)$)
\enle

{\bf proof:}
Consider the Taylor expansion for $B_3(x_0,x_1,x_2)$ at each
inflection point $s_i$, $0 \le i \le 2$.

One can assume $s_i \notin \{ (x_0:x_1:x_2) \in \RRR P^2 | x_{i}=0 \}$
and consider  $b_{3,i}(x,y)= B_3(x_0,x_1,x_2)|_{ x_{i}=1 }$.
In a neighborhood $U(s_i)$ of $s_i=(a_i,b_i)$,\\
$b_{3,i}=  \Sigma_{n=0}^3
\frac {1} {n!}
\Sigma_{k=0}^n (\frac {\pr^n B_3}  {\pr^k x \pr^{n-k} y}(a_i,b_i)
(x-a)^k (y-b)^{n-k})$
Each tangent $L_i$ meets ${\mathcal B}_3$ with multiplicity $3$.
Since $\frac {\pr^2 B_3} {\pr^k x \pr^{2-k} y} (a_i,b_i)=0$
for any $s_i$, $0 \le i \le 2$,
when consider the Taylor expansion of $B_3$,
in any neighborhood $U(s_i) \subset \RRR P^2$
of the inflection point $s_i$
$B_2(s)=L_i(s)$ for any $s \in U(s_i)$;
and thus
$B_3(s) - L_i(s) =C_3(s)$  for any $s \in U(s_i)$.\\
Consider the linear change of coordinates which maps $(x_0,x_1,x_2)$
to $(x_0'=L_0,x_1'= L_1, x_2'= L_2)$.
It maps $B_3(x_0,x_1,x_2)$ to $B'_3
(x_0',x_1',x_2')$ .
Let us prove that it also maps $C_3$ to $\e.x_0'.x_1'.x_2'$.
It maps $s_i$ to $s'_i$
and one can assume $s'_i \notin \{
(x_0':x_1':x_2') \in \RRR P^2 |x'_{i}=0 \}$.
Consider the Taylor expansion of $b'_{3,i}$ in
a neighborhood
$U(s'_i)=\{ (x_0':x_1':x_2') \in \CCC P^2 | x'_i \not=0 \}$
of $s'_i$.
For any
 $i \in \{0,1,2 \}$,
 $j, k \not=i$
 $j,k \in \{0,1,2 \}$,
$b'_{3,i}(s) - x'_j(s).x'_k(s)=C'_3(s)$ for any $s \in U(s_i)$.
Hence,
it  follows  that $C'_3=\e.L_0.L_1.L_2$.
since $L_0.L_1.L_2(s)-\e.C'_3 (s) =0$
for any $s \in U(s_i)$.
It concludes the proof of Lemma \mrf{l:cub3}.
Q.E.D
\vskip0.1in
{\bf II)}
\bede
\mlb{d:l07}
Let ${\mathcal A}$ and ${\mathcal B}$ be
smooth algebraic curves of $\CCC P^2$
with respective order $i$ and $j$, $i \le j$.
Let $\RRR {\mathcal A}$
,resp. $\RRR {\mathcal B}$,
be the real point set
of ${\mathcal A}$, resp. ${\mathcal B}$.
We shall say ${\mathcal A}$ is ${\it immersed}$ in
${\mathcal B}$ if, up to rigid isotopy  of $\RRR {\mathcal B}$,
$\RRR {\mathcal A}$ is embedded in $\RRR {\mathcal B}$.
\end{defi}
\bede
Given $A(x_0,x_1,x_2)=x_0^i.a_(x_1/x_0,x_2/x_0)$ the polynomial of a curve
${\mathcal A}$ immersed in ${\mathcal B}$.
Assume $A$ regular.
Denote  by
 ${\mathcal P}_{\mathcal A}$
  the pencil of curves  over ${\mathcal A}$.
(i.e curves with polynomial
$x_0^i(a_(x_1/x_0,x_2/x_0)-c)$, $c \in \RRR$.)
\vskip0.1in
We shall say that ${\mathcal C} \in {\mathcal P}_{\mathcal A}$
is $M-immersed$  (over ${\mathcal A}$) in ${\mathcal B}$ if:
\been
\item
${\mathcal C}$ is immersed
in ${\mathcal B}$
\item
${\mathcal C}$ has the maximal number of real components
a curve of ${\mathcal P}_{\mathcal A}$ immersed in ${\mathcal B}$
may have.
\enen
\end{defi}
\vskip0.1in

We shall prove in  Lemma \mrf{l:l071} that a curve $\HH_3$
is immersed in $\HH_m$.

\bele
\mlb{l:l071}
Let $\HH_m$ be the curve with polynomial $B_m$.
Up to regular modification of $B_m$, the curve
${\mathcal B_3}$ is the curve $\HH_3$ where $\HH_3$
is immersed in $\HH_m$
as the classical deformation of $\HH_2 \cup L$ where $L:=x_0=0$.
\enle

\vskip0.1in

The  curve ${\mathcal B}_3$
is determined up to rigid  isotopy
by its real scheme.
The two possible  real schemes
for ${\mathcal B}_3$ are
$\big\langle J \sqcup 1 \big\rangle$
and
$\big\langle J \big\rangle$.

Hence,
according to Lemma \mrf{l:cub3},
up to rigid isotopy of $B_m(x_0,x_1,x_2)$,
which is also a rigid isotopy of $B_3(x_0,x_1,x_2)$,
one can assume
that $B_3(x_0,x_1,x_2)$ results from the classical deformation
of the union of a line $L$
with a curve ${\mathcal B}_2$  of degree $2$  with real scheme
$\big\langle 1 \big\rangle$;
the deformation is directed to the union of the
tangents at  real inflection points
of ${\mathcal B}_3$.
Real inflection points
of ${\mathcal B}_3$ are points of $L$.
It is not hard to see that if none of these points
belongs to the inner of ${\mathcal B}_2$,
then  ${\mathcal B}_3$ is of type $\HH^0$ relatively $L$.
Otherwise, if at least one of this point
belongs to the inner of ${\mathcal B}_2$
,then  ${\mathcal B}_3$ has real scheme
$\big\langle J \big\rangle$.\\
Denote $L$ the line infinity $x_0=0$,
and let $b_3(x,y)=B_3(1,x_1,x_2)$
be the  affine polynomial associated to
$B_3(x_0,x_1,x_2)$.
Consider the pencil $x_0^{m}.
(b_3(x_1/x_0,x_2/x_0)-c)$.
For any $i \ge 4$,
let $C_i(x_1,x_2)$ be
the polynomial of degree $i$ in the variables $x_1,x_2$
such that:
$B_m(x_0,x_1,x_2) = x_0^{m-3}.B_3(x_0,x_1,x_2)+
 \Sigma_{i=4}^m x_0^{m-i}.C_i(x_1,x_2)$. \\
Consider the pencil $x_0^{m-3}.
(b_3(x_1/x_0,x_2/x_0)-c)$.
Assume $B_3$ regular.
For any critical point
$(1,x_{0,1},x_{0,2})$
of $B_3(1,x_1,x_2)=b_3(x_1,x_2)$,
choose its representative in $S^2$
$\frac {1} {(1+x_{0,1}^2+x_{0,2})^{1/2}} (1,x_{0,1},x_{0,2})$
In such a way,  it is easy to see that a curve of degree $3$
is immersed in ${\mathcal B}_m$.
Since ${\mathcal B}_m$ has the maximal number of real components
a curve of degree $m$ may have,
an $M$-curve of degree $3$  is $M$-immersed in ${\mathcal B}_m$.
Indeed,
if the curve ${\mathcal B}_3$ is not an $M$-curve,
it is an easy consequence of the Petrovskii's theory that
as $c$ varies from
$[-
\Sigma_{j=4}^m
\vert \vert C_j \vert  \vert ,
\Sigma_{j=4}^m
\vert \vert C_j \vert  \vert ]$
the real point set of the
curve ${\mathcal B}_3$ undergoes at least one Morse modification
(i.e
at least one critical point passes through its critical value.)
According to Lemma \mrf{l:cub3},
this is equivalent to say that any union of lines $L_i \cup L$
$1 \le i \le 3$ is perturbed in such a way
none of the resulting real branches intersects the oval
${\mathcal B}_2$.
It follows that, up to regular deformation of
 $B_m(x_0,x_1,x_2)$,
$B_3(x_0,x_1,x_2)$ is of type $\HH^0$.
Since curves of degree $m \le 3$ are defined up to rigid isotopy
by their real scheme, up to rigid isotopy of
 $B_m(x_0,x_1,x_2)$,
$B_2(x_0,x_1,x_2)$, and in this way also $B_1(x_0,x_1,x_2)$
are  of type $\HH^0$.
It concludes the proof of Lemma \mrf{l:l071}.
Q.E.D
\vskip0.1in
The Proposition \mrf{p:l07} is a straightforward consequence of
the Lemma \mrf{l:l071}.
Q.E.D

\vskip0.1in
Denote $\Omega_j$ the set  of $M$-curves with real scheme:\\
-for even $j=2k$
\begin{eqnarray}
\big\langle1\langle  {\alpha} \rangle \sqcup \beta
\big\rangle  &\nonumber
\end{eqnarray}
with $\alpha + \beta=\frac {(j-1)(j-2)} 2$\\
-for odd $j=2k+1$
\begin{eqnarray}
\big\langle J \sqcup {\gamma}
\big\rangle &\nonumber
\end{eqnarray}
with $\gamma=\frac {(j-1)(j-2)} 2 $\\

\bepr
\mlb{p:l08}
Let $\HH_m$ be a Harnack curve of degree $m \ge 4$.
Denote by $B_m$ its polynomial.
Up to regular modification of $B_m$,
the curve $\HH_m$ results from the deformation of
the union  ${\mathcal A}_{m-i} \cup \HH_i$, $i \le 3$,
where ${\mathcal A}_{m-i} \in \Omega_{m-i}$
and $\HH_i$ is a Harnack curve of degree $i \le 3$.
\enpr

{\bf proof:}\\

Set $L:=x_0=0$.
It follows from Proposition \mrf{p:l07} that one can set
$$B_{m}(x_0,x_1,x_2) =
A_{m-i}(x_0,x_1,x_2)
.B_i(x_0,x_1,x_2) +
C_m(x_0,x_1,x_2)$$
where $B_i(x_0,x_1,x_2)$ is a Harnack of degree $i \le 3$
and type $\HH^0$ relatively to $L$ and
$A_{m-i}(x_0,x_1,x_2)=x_0^{m-i}$.\\
Let us prove that there exists a regular modification of
$B_{m,t}(x_0,x_1,x_2)$, $t \in [0,1]$,
$$B_{m,0}(x_0,x_1,x_2)=B_{m}(x_0,x_1,x_2)$$

$$B_{m,t}(x_0,x_1,x_2)=A_{m-i,t}(x_0,x_1,x_2).B_{i,t}(x_0,x_1,x_2) +
C_{m,t}(x_0,x_1,x_2)$$
such that:
\been
\item
$B_{i,t}(x_0,x_1,x_2)$
is a regular modification of $B_i(x_0,x_1,x_2)$ the Harnack polynomial
of a Harnack curve $\HH_i$,
\item
$A_{m-i,1}$
is the polynomial of a smooth curve of degree $m-i$
($A_{m-i,0}=x_0^{m-i}$)
\item
The polynomial $B_{m,1}$ is the polynomial of a
curve $\HH_m$ which results from the classical
deformation of the union of ${\mathcal A}_{m-i,1} \cup \HH_i$
with polynomial $A_{m-i,1} \cup B_{i,1}$
\enen
(It is obvious that
the regular modification  $B_{m,t}$ of $B_m$
is not regular on $A_{m-i}$)
Our proof is based on the fact that
curves of degree $i \le 4$
are determined up to  rigid isotopy by their real schemes
and may be realized as $L$-curves
(i.e may realized by classical small perturbation
of $i$ lines in general position) \cite{Fi}.
We shall consider immersion
(see the  definition \mrf{d:l07} of
the proof of Proposition \mrf{p:l07}),
of such curves
in the curve $\HH_m$ with polynomial $B_m$.\\

We shall consider curves
$\HH_m$ of degree $m \le 7$ and degree $m\ge 8$
separately.
\vskip0.1in

{\bf 1)}\\
Let us consider curves $\HH_m$ of degree $m \le 7$.
According to the Proposition \mrf{p:l07}, up to rigid isotopy,
the polynomial $B_m$ of $\HH_{m}$ is of the form
$$B_{m}(x_0,x_1,x_2) =
x_0^{m-3}.B_3(x_0,x_1,x_2) +
C_m(x_0,x_1,x_2)$$
with $m \le 7$ where $B_3$ is the polynomial of a Harnack curve of degree
$3$.\\
{\bf -}
One can deform the polynomial
$A_{m-3,0}= x_0^{m-3}$
into
the polynomial
$\Pi_{j=1}^{m-3}(x_{0,j})$
of $m-3$ lines in general position
in such a way that
the path
 $B_{m,t}$, $t \in [0,t_0]$,
is rigid isotopy from
$B_{m,0}=B_{m}$ to $B_{m,t_0}$.

$$B_{m,t_0}(x_0,x_1,x_2)=
\Pi_{j=1}^{m-3}(x_{0,j})
.B_{3}(x_0,x_1,x_2)+
C_{m}(x_0,x_1,x_2)$$
{\bf -}
According to the Lemma \mrf{l:l071},
one can assume that the curve $\HH_3$ with polynomial $B_3$ is
immersed in $\HH_m$ as the classical deformation
of $\HH_2 \cup L$ where $L:=x_0=0$.\\
Therefore, one can consider a rigid isotopy
$B_{m,t}$, $t \in [t_0,t_1]$,
$$B_{m,t_1}=B_{3}.B_{m-3} + C_{m}$$
where $B_{m-3}$ is  a regular polynomial
of an $L$-curve which results from the deformation of
$\Pi_{j=1}^{m-3}(x_{0,j})$.

Let $||B_{3}||$ be the norm of $B_{3}$
and let
$b_{m-3}.||B_{3}||$
be the affine polynomial associated to
$B_{m-3}.||B_{3}||$.
Consider the pencil of curves
$x_0^{m}(b_{m-3}.||B_{3}|| - c)$.\\
For any critical point $(x_0,x_1,x_2)$ of
$B_{m-3}.||B_{3}||$
choose its representative in $S^2$
$\frac {1} {(1+x_{0,1}^2+x_{0,2}^2)^{1/2}} (1,x_{0,1},x_{0,2})$.
As $c$ varies from $[-||C_{m}||,||C_{m}||]$,
for any  critical point of $b_{m-3}$ which
goes through its
critical value; the
real point set of a curve of the pencil
undergoes a Morse modification.
Since $\HH_m$ has the maximal number of real connected components
a curve of degree $m$ may have,
it is not hard to see that
an $M$-curve of degree $m-3$ is $M$-immersed
in $\HH_m$.
Up to rigid isotopy of $B_{m}$, it is realized as the $L$-curve
which results from the perturbation of the union of $m-3$ lines
with polynomial $\Pi_{j=1}^{m-3}(x_{r,j})$.\\
{\bf -}
It follows that there exists
a rigid isotopy
$B_{m,t}$, $t \in [0,1]$
from
$B_{m,0}$
to
$$B_{m,1}=B_{3}.B_{m-3} + C_{m,1}$$
where $B_{m-3}$ is the polynomial of $M$-curve of degree $m-3 \le 7$;
and the curve the curve $\HH_m$
results from classical small deformation
of $\HH_{m-3} \cup \HH_3$ with polynomial $B_{m-3}.B_{3}$.
Since
any curve $\HH_{m-3}$ of degree $m-3 \le 4$
belongs to $\Omega_{m-3}$, it proves the Proposition \mrf{p:l08}
for curves $\HH_m$ of degree $m \le 7$.
\vskip0.1in
{\bf 2)}\\
Let us consider curves $\HH_m$ of degree $m \ge 8$.
According to the Proposition \mrf{p:l07}, up to rigid isotopy,
the polynomial $B_m$ of $\HH_{m}$ is of the form
$$B_{m}(x_0,x_1,x_2) =
x_0^{m-3}.B_3(x_0,x_1,x_2) +
C_m(x_0,x_1,x_2)$$
where $B_3$ is the polynomial of a Harnack of degree
$3$ and
$m-3=4.r+s$ with $r,s \in \NNN, s \le 3$.\\

Set $A_{m-3,0}=x_0^{m-3} =x_0^{4r+s} =(x_0^4)^r. x_0^s=
(x_0^4+...x_0^4).x_0^s$.
For any of the $n^{th}$ ($1 \le n \le r$)
product $x_0^4$ of the polynomial
$(x_0^4)^r=x_0^4+...x_0^4)$
consider the deformation
of $x_0^4$
into the polynomial
$\Pi_{j=1}^4(x_{n,j})$
of the union
of four
lines in generic position.
In the same way, consider the deformation of $x_0^s$ into
the polynomial
$\Pi_{j=1}^s(x,{r+1,j})$
of $s$ lines
in generic position.\\
{\bf -} For any of these union of $4$ (resp, $s \le 4$),
lines in generic position,
consider the  rigid isotopy
$B_{m,t}$, $t \in [0,t_0]$
from
$B_{m,0}=B_{m}$ to $B_{m,t_0}$.\\
For $1 \le n \le r$,
$$B_{m,t_0}(x_0,x_1,x_2)=
x_0^{4(r-1)+s}.\Pi_{j=1}^{4}(x_{n,j})
.B_{3}(x_0,x_1,x_2)+
C_{m}(x_0,x_1,x_2)$$
(resp,
$$B_{m,t_0}(x_0,x_1,x_2)=
x_0^{4r}.
\Pi_{j=1}^{s}(x_{r+1,j})
.B_{3}(x_0,x_1,x_2)+
C_{m}(x_0,x_1,x_2))$$
From an argumentation analogous to the one given in ${\bf 1)}$,
it follows that, up to regular deformation of $B_m$,
any  of the union of ${4}$ (resp, $s$) lines in generic position
may be deformed into an $L$-curve of degree $4$ (resp, $s$)
which is an $M$-curve immersed in $\HH_m$.
Denote by  $\HH_{4,n}$, $1 \le n \le r$,
(resp, $\HH_s$)
the corresponding curves.
Denote by $B_{4,n}$
(resp, $B_{s}$)
the polynomial of $\HH_{4,n}$ (resp, $\HH_s$).
{\bf -}
It follows that
there exists a  rigid isotopy
$B_{m,t}$, $t \in [0,t_1]$,
from
$B_{m,0}=B_{m}$ to $B_{m,t_1}$.\\
$$B_{m,t_1}(x_0,x_1,x_2)=
\Pi_{n=1}^{r}(B_{4,n})
.B_{s}.B_3+
C_{m,t_1}(x_0,x_1,x_2)$$
\vskip0.1in

\bede
Call {\it succesive} small perturbation of a finite  union of
$\cup_{i=1}^n{\mathcal A}_i$ of curve ${\mathcal A}_i$
the result of the following recursive classical small perturbation
of the union of two curves:
\been
\item
the union ${\mathcal A}_{1} \cup {\mathcal A}_{2}$
is a singular curve all of whose singular points are crossings;
the classical small perturbation
${\mathcal A}_1 \cup {\mathcal A}_2$
leads to a curve ${\mathcal B}_2$
\item
For $2 \le i \le n-1$,
the classical small perturbation
${\mathcal B}_i \cup {\mathcal A}_{i+1}$
leads ${\mathcal B}_{i+1}$ such that
the union ${\mathcal B}_{i+1} \cup {\mathcal A}_{i+2}$
is a singular curve all of whose singular points are crossings.
\enen
\end{defi}
Let us prove that a
curve ${\mathcal A}_{m-3} \in \Omega_{m-3}$
of degree $m-3$
which results from a
successive classical small deformation of the union
$\cup_{n=1}^{r} \HH_{4,r} \cup \HH_s$
is $M$-immersed in $\HH_m$.\\

Let us consider the union of curves
$\cup_{n=1}^{r}\HH_{4,n} \cup \HH_s$
with polynomial $\Pi_{n=1}^{r}B_{4,n}.B_s$.
Up to rigid isotopy of $B_m$,
one can assume $\HH_{4,n}=\HH_{4,1}$, $2 \le n \le r$,
where $\HH_{4,1}$ is the $L$-curve of degree $4$ which results from
the deformation of the union
of four lines
in generic position with polynomial
$\Pi_{j=1}^4 x_{1,j}$.
One can also assume that
$\HH_s$, $s \le 4$,is the $L$-curve of degree $s$ which results
from the classical deformation of the union of $s$ lines
with polynomial $\Pi_{j=1}^s x_{1,j}$.\\
Since $\HH_m$ is an $M$-curve, for any $i~, 1 \le i \le r$,
there exists an $M$-curve of degree $4i+s$
immersed in $\HH_m$
as the  result of
the successive classical perturbation of
$\cup_{k=1}^i \HH_{4,k} \cup \HH_s$.
Denote by ${\mathcal B}_{s+4i}$
the $M$-curve of degree $4+si$ immersed in $\HH_m$
as the  classical perturbation of
$\cup_{k=1}^i \HH_{4,k} \cup \HH_s$.

The curve ${\mathcal B}_{4+s}$ immersed in $\HH_m$
belongs to  $\Omega_{4+s}$.
Indeed, none of the pair
of ovals $({\mathcal O},{\mathcal O}')
\in (\RRR {\mathcal \HH}_{4,1}, \RRR \HH_{s})$
is injective.
In the same way,
since $\HH_{4,n}=\HH_{4,1}$ for $2 \le n \le r$
none of the pair of ovals
$({\mathcal O},{\mathcal O}')
\in (\RRR {\mathcal B}_{s+4i}, \RRR \HH_{4,{i+1}})$
is an injective pair of ovals.
Hence,
at any step of the successive classical small deformation
we get an $M$-curve ${\mathcal B}_{s+4i} \in \Omega_{s+4i}$.
\vskip0.1in

For clarity,
let us verify that, up to rigid isotopy of $B_m$,
the successive classical deformation of
$\cup_{k=1}^r \HH_{4,k} \cup \HH_s$ as curve immersed in $\HH_m$
does not depend on the order of the deformations.

Let us notice that:
\bere
\mlb{r:clm}
The classical small deformation of the union of two $M$-curves
${\mathcal A} \cup {\mathcal B}$ leads to an $M$ curve
only if
common points of ${\mathcal A}$ and
${\mathcal B}$
belong  to one connected real connected of ${\mathcal A}$
and one connected real connected of ${\mathcal B}$.
\enre
{\bf proof:}
Otherwise, it would lead to contradiction with
the number of real connected components of $\HH_m$.
Q.E.D
\vskip0.1in

According to remark \mrf{r:clm},
there exists ${\mathcal O} \in \RRR \HH_{4,1}$ such that
any two curves immersed in $\HH_m$
as the successive classical deformation
$\cup_{n=1}^l \HH_{4,n}$  and respectively
$\cup_{n=1}^{l'} \HH_{4,n}$, $ l\not l', 2 \le l,l' \le r$
(where $\HH_{4,2} \not= \HH_{4,1}$)
intersects $\HH_{4,1}$ in ${\mathcal O}$.
In the same way,
there exists ${\mathcal O} \in \RRR \HH_{s}$
(if $s$ is odd, ${\mathcal O}=J$ the odd component of $\HH_s$)
 such that
any two curves immersed in $\HH_m$
as the successive classical deformation
$\cup_{n=1}^l \HH_{4,n}$  and
$\cup_{n=1}^{l'} \HH_{4,n}$, $ l \not=  l',~ 2 \le l,l' \le r$
(where $\HH_{4,2} \not= \HH_{4,1}$)
intersect $\HH_{s}$ in ${\mathcal O}$.
\vskip0.1in

Hence, the $M$-curve  ${\mathcal A}_{m-3}$ immersed in $\HH_m$
as the result of the successive classical small deformation
of $\cup_{k=1}^r \HH_{4,k} \cup \HH_s$ is well defined.
Besides, it belongs to  $\Omega_{m-3}$.

{\bf -}
It follows that there exists
a rigid isotopy
$B_{m,t}$, $t \in [0,1]$
from
$B_{m,0}$
to
$$B_{m,1}=B_{3}.B_{m-3} + C_{m,1}$$
where
$B_{m-3}$ is the polynomial of
${\mathcal A}_{m-3} \in \Omega_{m-3}$
and $\HH_m$
results from classical small deformation
of the curve ${\mathcal A}_{m-3}  \cup \HH_3$
with polynomial $B_{m-3}.B_{3}$.\\
It concludes the proof of Proposition \mrf{p:l08}.
Q.E.D

\vskip0.1in
In the subsection \mrf{su:HH0},
we have proven
that if  $\HH_n$ is of type $\HH^0$, then,  up to regular modification
of its polynomial $B_n$, the polynomial $B_{n-1}$
is also of type $\HH^0$.
Let us now in Proposition \mrf{p:HH0} prove the converse.\\

\bepr
\mlb{p:HH0}
Let $\HH_n$ be a Harnack curve of degree $n$, $n \ge 4$,
and $B_n(x_0,x_1,x_2)=x_0^{n}.b_{n}(x_1/x_0,x_2/x_0)$
its polynomial.
Denote by $B_{n-1}(x_0,x_1,x_2)=
 x_0^{n-1}.b_{n-1}(x_1/x_0,x_2/x_0)$
the homogeneous polynomial associated to
the truncation $b_{n-1}(x,y)$ of $b_n(x,y)$ on the monomials
$x^{\a}.y^{\b}$ with $0 \le \a + \b \le n-1 $.\\
If $B_{n-1}(x_0,x_1,x_2)$
is of type $\HH^0$
relatively to the line $L_n:=x_0=0$, then
up to regular modification of $B_n$,
there exists a line $L_{n+1}$  with the property that
$\HH_{n}$ is of type $\HH^0$ relatively to $L_{n+1}$.
\enpr
{\bf proof :}\\
Our proof uses an argumentation similar to
the one given in the subsection \mrf{su:HH0} to prove
that if  $\HH_n$ is of type $\HH^0$ then up to regular modification
of its polynomial $B_n$, the curve ${\mathcal B}_{n-1}$
is also of type $\HH^0$.\\
Let $L_n:=x_0=0$ and
${\mathcal M}$ be a neighborhood of the line $L_n$.
The curve $\HH_n$ with polynomial $B_n$ results from deformation
of $\HH_{n-1} \cup  L_n$
$$B_n=B_{n-1}(x_0,x_1,x_2).x_0 +C_n(x_0,x_1,x_2)$$
where
$C_n(x_0,x_1,x_2)$ is homogeneous of degree $n$.
Let  $U(s) \subset \RRR P^2$ be
a  neighborhood
of a singular point  of
$\RRR \HH_{n-1} \cup \RRR L_{n}$
 such that
$\RRR \HH_n  \cap U(s) \not = \emptyset$.
According to the topology of $\HH_n$,
for each singular point $s$ there exists
a homeomorphism $h: U(s) \to D^1 \ti D^1$ such that
$h(\RRR \HH_{n-1} \cup L_{n})=D^1 \ti 0 \cup  0 \ti D^1$
and $h(\RRR \HH_{n} \cap U(s))=\{ (x,y) \in D^1 \ti D^1 |xy=1/2 \}$.
Since the polynomial of $\HH_{n}$
is smooth, it follows by continuity,
that real connected components of $\RRR \HH_{n}$ define
entirely and uniquely
the homomorphism $h$.\\
We shall  enlarge
the previous description  of $\{ x \in \RRR P^2| B_{n}(x)=0 \}$
in $U(s)$ to a  description of
in $U_s$ where
$$U_s=\{ z= <u,p>=(u_0.p_0:u_1.p_1:u_2.p_2) \in \CCC P^2 |
u=(u_0:u_1:u_2) \in U_{\CCC}^3, p=(p_0:p_1:p_2) \in U(s) \}$$
To this end, let us prove the following Lemma \mrf{l:MorsTr2}.

\bele
\mlb{l:MorsTr2}
Let
$B_n(x_0,x_1,x_2)=x_0^{n}.b_{n}(x_1/x_0,x_2/x_0)$
be a Harnack polynomial of degree $n$.
Denote by $B_{n-1}(x_0,x_1,x_2)=
 x_0^{n-1}.b_{n-1}(x_1/x_0,x_2/x_0)$
the homogeneous polynomial associated to
the truncation $b_{n-1}(x,y)$ of $b_n(x,y)$ on the monomials
$x^{\a}.y^{\b}$ with $0 \le \a + \b \le n-1 $.\\
Assume that $B_{n-1}$ is  of type $\HH^0$
relatively to $L:=x_0=0$.
\been
\item
\lb{i:tru11}
there exists ${\mathcal M}$ neighborhood
of $L$ such that
the polynomial $B_{n}^{\D}(x_0,x_1,x_2)=(x_0^n.b_n^{\D}(x_1/x_0,x_2/x_0)$
where $b_n^{\D}$ is the truncation of $b_n$
on monomials of homogeneous degree $n-2 \le i \le n$ is
$\e$-sufficient for $B_{n}(x_0,x_1,x_2)$.
\item
\lb{i:tru22}
Up to regular modification of $B_n$,
there exist ${\mathcal M}_{\e}  \subset {\mathcal M}$,
${\mathcal M}_{\e}= \cup_{s \in \RRR \HH_{m-1} \cap L} U_s$
$\e$-tubular neighborhood of $L$ and a polynomial $\tilde{B}_n$
with the following properties:\\
-$\tilde{B}_n$ is $\e$-sufficient for $B_n$ in ${\mathcal M}_{\e}$\\
-for any $U_s$, $U_s \subset {\mathcal M}_{\e}$,
the truncation $\tilde{B}_{n}^{S}$ of $\tilde{B}_n^{S}$
(on four monomials
$x^cy^dz^2, x^{c+1}y^dz, x^cy^{d+1}z,x^{c+1}y^{d+1}$
 with $c+d+2=n$)
is $\e$-sufficient for $\tilde{B}_{n}$ (and thus for $B_n$) in $U_s$.
(For any $s \not = s'$,
such that the truncation of $\tilde{B}^{n}$
on $x^cy^dz^2, x^{c+1}y^dz, x^cy^{d+1}z,x^{c+1}y^{d+1}$
(resp,on
$x^{c'}y^{d'}z^2, x^{c'+1}y^{d'}z, x^{c'}y^{d'+1}z,x^{c'+1}y^{d'+1}$)
is
$\e$-sufficient for $\tilde{B}_{n}$ (and thus for $B_n$) in $U_s$
(resp, in  $U(s')$,
 $(c,d) \not = (c',d')$~)
\\
\enen
\enle
\bere
Note that if the curve $\HH_{n}$ is not of type $\HH^0$ (relatively to $L$),
then
${\mathcal M}_{\e} \subset {\mathcal M}$,
${\mathcal M}_{\e} \not = {\mathcal M}$.
\enre
\vskip0.1in
{\bf proof:}
Assume $B_{n-1}$ of type $\HH^0$ relatively to $L_n :=x_0=0$.
According to the proof of Morse Lemma, local coordinates defined
in a neighborhood $U(s)$ of a point $s$
depend principally on the first
derivative and the second derivative
of the function $b_n$ around this point.
Hence, the truncation of $b_{n}^{\D}(x,y)$ on
monomials of homogeneous degree $n-2 \le i \le n$ is
$\e$-sufficient for $b_{n}(x,y)$.\\
According to Morse Lemma \cite{Mil},
in neighborhood $U(s)$ of any non-degenerate singular point
$s$ of $x_0.B_{n-1}$
one can choose local coordinates system  $z_1, z_2$ with $z_1(s)=0,z_2(s)=0$
and $x_0.B_{n-1}=z_1.z_2$.
Real connected components of $\RRR \HH_{n}$ define
entirely and uniquely
the homomorphism $h$ defined in $U(s)$
$h: U(s) \to D^1 \ti D^1$ such that
$h(\RRR \HH_{n-1} \cup L_{n})=D^1 \ti 0 \cup  0 \ti D^1$
and $h(\RRR \HH_{n} \cap U(s))=\{ (x,y) \in D^1 \ti D^1 |xy=1/2 \}$.
In such a way, in a neighborhood $U(s)$ of a point $s$,
one can choose local coordinates system $y_1, y_2$
with $y_1(s)=0,y_2(s)=0$; $B_n =y_1.y_2 + B_n(s)$.
In $U(s)$, any point $s$ is a local extremum of $b_n$.
Any singular point $s$ of $\RRR B_{n-1} \cup \RRR L_n$
belongs to the line $x_0=0$,
hence it is a local  extremum of
the function $b_n(y)=b_{n}(0,x_1/x_2,1)$.
One can assume
without loss of generality,
that any point $s$
does not belong to the line $x_2=0$.
On such assumption gradient trajectories
of $b_{n}(x_0,x_1,1)$ in $U(s)$ give the direction
of the perturbation of the crossing $s$.
Thus, up to regular modification of $B_n$
and for sufficiently small $\e >0$,
there exist
$U(s) \subset {\mathcal M}_{\e}$
and a polynomial $\tilde{b}_{n}^{S}$ on four monomials
($x^cy^d, x^{c+1}y^d, x^cy^{d+1},x^{c+1}y^{d+1}$
with $c+d+2=n$)  which is $\e$-sufficient for
$B_{n}(1,x_1,x_2)$
in $U(s) \subset {\mathcal M}_{\e}$;
$\tilde{b}_{n}^S(x,y)=l(x,y)+k(x,y)$
with
$l(x,y)=
 a_{c,d}x^cy^d  +
 a_{c+1,d}x^{c+1}y^d $\\
$k(x,y)=
 a_{c,d+1}x^cy^{d+1}  +
a_{c+1,d+1}x^{c+1}y^{d+1}$
with $c+d+2=n$.
From this last observation, we deduce the definition
of $\tilde{b}_n$.
Q.E.D
\vskip0.1in
Let us consider the polynomial $\tilde{b}_n$
with truncation $\tilde{b}_{n}^{S}$  $\e$-sufficient for
$B_{n}(1,x_1,x_2)$ in $U_s$.\\
$\tilde{b}_{n}^S(x,y)=l(x,y)+k(x,y)$
with
$l(x,y)=
 a_{c,d}x^cy^d  +
 a_{c+1,d}x^{c+1}y^d $\\
$k(x,y)=
 a_{c,d+1}x^cy^{d+1}  +
a_{c+1,d+1}x^{c+1}y^{d+1}$
with $c+d+2=n$.

Up to modify the coefficients
$a_{c,d},a_{c,d+1},a_{c+1,d},a_{c+1,d+1}$ if necessary,
the point $s= (x_0,y_0)$ is, up to homeomorphism,
a critical point of the function
$\frac {l(x,y)} {k(x,y)}$.
Hence, it follows from the equalities
$l(x,-y)=-l(x,y)$, $k(x,-y)=k(x,y)$,\\
$\frac {\pr l} {\pr x}(x,-y) =-\frac {\pr l} {\pr x}(x,y)$,
$\frac {\pr l} {\pr y} (x,-y) =\frac {\pr l} {\pr y} (x,y)$,\\
$\frac {\pr k} {\pr x} (x,-y)= \frac {\pr k} {\pr x} (x,y)$,
$\frac {\pr k} {\pr y} (x,-y) = -\frac {\pr k} {\pr y} (x,y)$\\
that $(x_0,- y_0)$ is also a critical  point of
the function $\frac {l(x,y)} {k(x,y)}$.
In such a way, we  define a map
which maps the set ${\mathcal S}$ of
singular points  of $\HH_{n-1} \cup L_n$
to a set ${\mathcal P}$ of $n-1$ real points.
In a neighborhood $U(p)$ of $p \in {\mathcal P}$,
$\RRR \HH_m$ is desingularized crossing i.e
consists of two real arcs $\gamma,\gamma'$.
The set ${\mathcal S}$ is a set of aligned points.
Hence, up to regular modification, the  set
${\mathcal P}$ is also a set of aligned points which belong
to a line $L$.
Let ${\mathcal M}_{\e}= {\mathcal M}^+ \cup {\mathcal M}^-$
,$\pr{\mathcal M}^{\pm}=L_1^{\pm} \cup L_2^{\pm}$,
be a neighborhood of the line $L$;
$\cup_{p \in {\mathcal P}} U(p) \subset {\mathcal M}$.\\
In the case $n=2k$, any pair of arcs of $\RRR \HH_{2k} \cap U(p)$,
$p \in {\mathcal P}$, is such that
one arc belongs to the non-empty oval of $\RRR \HH_{2k}$.
The non-empty oval is connected. Hence, it
intersects a projective line $L_1$
in $2k$ points.\\
In the case $n=2k+1$, any pair of arcs of $\RRR \HH_{2k+1} \cap U(p)$,
$p \in {\mathcal P}$, is such that
at least one arc belongs to the odd component of $\RRR \HH_{2k+1}$.
The odd component is connected. Hence, it
intersects a projective line $L_1$
in $2k+1$ points.\\
In such a way, we get that
if $\HH_{n-1}$ is of type $\HH^0$
relatively to $L_{n}$, there exists $L_1:=L_{n+1}$ such
that the curve $\HH_{n}$ is of type $\HH^0$ relatively
to $L_1=L_{n+1}$.
This concludes our proof
of Proposition \mrf{p:HH0}.
Q.E.D
\vskip0.1in
\vskip0.1in

Let us prove in Theorem \mrf{t:seqli}
, that up to regular modification of the Harnack polynomial
$B_m(x_0,x_1,x_2)$,
there exists a sequence of lines $L_{i+1}$, $ 1 \le i \le m$,
such that $B_i$ is of type $\HH^0$
relatively to the line $L_{i+1}$
and describe relative position of these lines.
We shall call rotation of $\RRR P^2$ the transformation of
$\RRR P^2$ which extends continuously a usual rotation of $\RRR^2$.
Relative position
of lines $L_i$, $ 1 \le i \le  m+1$ is described by means of
rotations of $\RRR P^2$.

\vskip0.1in

\beth
\mlb{t:seqli}
Let $\HH_m$ be a Harnack curve of degree $m$.
Up to rigid isotopy, there exists a sequence of lines $\HH_1=L_1$,
$L_2$, ..., $L_{m+1}$ such that $\HH_i$ is of type $\HH^0$
relatively to $L_{i+1}$ and $\HH_{i+1}$ results from the deformation
of $\HH_{i} \cup L_{i+1}$.\\
The sequence of lines $L_i$ may be ordered in three different ways:
\been
\item
\lb{i:klo}
 There exist $p \in L_2 \cap L_j$ , $3 \le j \le m+1$,
 and a  rotation of center $p$ and angle $\theta$
 (where $\theta$ may be chosen arbitrarily small)
 which maps $\RRR L_i$ to  $\RRR L_{i+1}$.
\item
\lb{i:kli}
 There exists $p \in L_2 \cap L_j$ , $3 \le j \le m+1$,
 such that:\\
 for $i$, $2 \le  i \le 2.[(m+1)/2]-2$
 the rotation of center $p$ and angle $\theta$
 maps $\RRR L_i$ to $\RRR L_{i+2}$.
 for  $i$, $1 \le  i \le 2.[(m+1)/2]-2$
 the rotation of center $p$ and angle $-\theta$
 maps $\RRR L_{i+1}$ to $\RRR L_{i+3}$
($\theta$ may be chosen arbitrarily small)\\
-~[m/2] denotes the integer $n$,  $ (m-1)/2 \le  n \le m/2 $ -
\item
\lb{i:klio}
 There exists $p \in L_2 \cap L_j$ , $3 \le j \le m+1$,
 such that:\\
 for $i$, $2 \le i \le l$,
 arrangement of lines $L_{i+1}$  is given by  \ref{i:kli}\\
 for $i$, $l \le i \le m+1$
 arrangement of lines $L_{i+1}$  is given by  \ref{i:klo}
\enen
\enth

\vskip0.1in

{\bf proof:}
In a first part, we shall prove Theorem \mrf{t:seqli}
for curves $\HH_m$ of even degree $m$.
Then, we shall extend our argumentation to Harnack curves of odd degree.
\vskip0.1in

We shall use the rigid isotopy classification of
Harnack curves of type $\HH^0$ given in Theorem \mrf{t:rig}.
The following remark may be easily  deduced from the
Harnack's construction of curves.
(see \cite{Har}, \cite{Vi}, tables of isotopy types of Harnack $M$-curves)
\bere
\mlb{r:Hcr}
The Harnack curve $\HH_{n+1}$
is the only $M$-curve of degree $n+1$ with:\\
-for even $n+1$ only one non-empty  oval\\
- for odd $n+1$ only empty ovals\\
which results from
the  classical deformation of the union
$\HH_n \cup L_{n+1}$ of a  Harnack $\HH_n$
of type $\HH^0$ relatively to a line $L_{n+1}$.
\enre
\vskip0.1in
{\bf I) Harnack Curves of even degree $\HH_{2k}$}\\
Our proof is based on the construction
of a sequence of curves
${\mathcal A}_{2i} \in \Omega_{2i}$, $1\le i \le k-1$ such that
the deformation of
${\mathcal A}_{2i} \cup \HH_{2k-2i}$ leads to the curve $\HH_{2k}$.\\
According to Proposition \mrf{p:l08},
one can assume that $\HH_{2k}$ results from the deformation
of  ${\mathcal A}_{2k-2} \cup \HH_2$
where ${\mathcal A}_{2k-2} \in \Omega_{2k-2}$.
It follows from the Bezout's theorem
Morse Lemma and the real components of $\RRR \HH_{2k}$
that the $(2k-2).2$ common points of  ${\mathcal A}_{2k-2}$ and $\HH_2$
belong to their non-empty ovals.
We shall consider
connected arcs $I_{j}$, $3 \le j \le {2k}$
of the non-empty oval of ${\mathcal A}_{2k-2}$;
$\cup I_{j} \supset \RRR \HH_2 \cap {\mathcal A}_{2k-2}$
and
identify  these arcs $I_j$ with lines $\RRR L_j$.
In this way,
we shall construct a sequence of curves $\HH_j$
of type $\HH^0$ relatively to $L_{j+1}$
with the property that
curves  $\HH_{2k-2i}$ are such that
${\mathcal A}_{2i} \cup \HH_{2k-2i}$ leads to the curve $\HH_{2k}$.\\
We shall call {\it extreme point} a point of the intersection
$\RRR {\mathcal A}_{2k-2} \cap \RRR \HH_2$
the two points connected by an arc of  $\RRR {\mathcal A}_{2k-2} $
which does not contain
other points $\RRR {\mathcal A}_{2k-2} \cap \RRR \HH_2$.
We shall call  {\it extreme arc} an arc which contains an extreme point.
\vskip0.1in
Let us detail the construction of the sequences of curves
${\mathcal A}_{2i} \in \Omega_{2i}$, $1\le i \le k-1$.
Any Harnack curves $\HH_2$ are rigidly isotopic.
Therefore,
without changing relative position of real connected components
of  ${\mathcal A}_{2k-2}$ and $\HH_2$,
we may assume $\HH_2$ of type $\HH^0$ relatively
to an (extreme) arc $I_3 \approx \RRR L_3$
of ${\mathcal A}_{2k-2}$.
(Choose $I_3$ of small lenght.
up to glue its extremities, we may identify
$I_3$ with a projective line  $\RRR L_3$.)
In such a way, according to the remark \mrf{r:Hcr},
the deformation of $I_3 \cup \RRR \HH_2$
leads necessarily to the Harnack curve $\HH_3$.
(Otherwise, it would lead to contradiction with the real connected
components of $\HH_{2k}$.)
Moreover, it follows from Proposition \mrf{p:HH0},
that the curve $\HH_3$
is of type $\HH^0$ relatively to one line.
According to  Theorem \mrf{t:rig} and remark \mrf{r:rrig},
we may, up to regular modification,
identify
(without changing relative position of real connected
components of
$\RRR {\mathcal A}_{2k-2}$ and $\RRR \HH_2$)
this line with an other arc $I_4$
of the non-empty oval of the curve ${\mathcal A}_{2k-2}$.
The arc  $I_4$  contains,
besides a set of $2$ consecutive points
which  are up to rigid isotopy  points of
$\RRR {\mathcal A}_{2k-2} \cap \RRR \HH_2$, an other point.
Identifying arcs  with lines
$I_3 \approx \RRR L_3$,
$I_4 \approx \RRR L_4$, this last point is provided by the intersection
$L_3 \cap L_4$.
From an argument analogous to the previous one,
we get a curve $\HH_4$ of type $\HH^0$
relatively to an arc $I_5 \approx \RRR L_5$.
By induction, it follows that $\HH_{2k}$ results from the deformation
of $\HH_4 \cup {\mathcal  A}_{2k-4}$ where
${\mathcal A}_{2k-4} \in \Omega_{2k-4}$.\\
Iterating our argumentation, we get a sequence of curves
${\mathcal A}_{2i} \in \Omega_{2i}$, $1\le i \le k-1$, such that
the deformation of
${\mathcal A}_{2i} \cup \HH_{2k-2i}$ leads to the curve $\HH_{2k}$.
We also define a set of arcs $I_{j+1}$,$2 \le j \le (2k-1)$,
with the property that the curve $\HH_j$ is of type $\HH^0$
relatively to $I_{j+1} \approx \RRR L_{j+1}$.
Arcs $I_{j}$, $3 \le j \le  2k$,
belong to the non-empty oval of ${\mathcal A}_{2k-2}$;
the set of arcs  $\cup_{j} I_{j}$, $2k-2i +1 \le j \le 2k$ is mapped
to a set of arcs of ${\mathcal A}_{2i}$, $1 \le i \le k-1$.

At the final step, we get $\HH_{2k}$
from the deformation of ${\mathcal A}_2 \cup \HH_{2k-2}$.\\

In the construction of the sequence of curves ${\mathcal A}_{2i}$,
we have to take into account the way $\HH_{2k-2i}$ intersects
${\mathcal A}_{2i}$.
Curves ${\mathcal A}_{2i}$ may intersect curves $\HH_{2k-2i}$
in only two different ways:
the non-empty oval of one curve contains all inner ovals
of the other curve or it contains none of them.
Let us first consider this last case.
\vskip0.1in
{\bf 1)}
Assume that the non-empty oval of
${\mathcal A}_{2i}$ does not contain ovals of the non-empty
oval of $\HH_{2k-2i}$.
The curve $\HH_{2k}$ results from the deformation
of the union $\HH_{2k-2i} \cup {\mathcal A}_{2i}$.
In this case,  any arc $I_{j}$, $2k-2i+1 \le j \le 2k$,
may be identified with the same line $\RRR L_{2k-2i+1}$.
Thus, from step of the construction of $\HH_{2k-2i+1}$
the method proposed for the construction of $\HH_{2k}$ is analogous to the
Harnack's construction.
Since  arcs  $I_j$,
 $2k-2i+1 \le j \le 2k$,
are distinct, one can say equivalently
-identifying  $I_j$ with lines $\RRR L_j$-
that there exists $p \in L_{2k-2i+1} \cap L_j$ , $2k-2i+1 \le j \le m+1$,
and a rotation of center $p$ and angle $\theta$
which maps $\RRR L_i$ to $\RRR L_{i+1}$.
For simplicity, we choose the sequence of arcs
$I_j$, $2k-2i+1 \le j \le 2k$ as follows.
The arc $I_{j+1}$ intersects
the arcs $I_{j}$ and $I_{j+2}$; $I_{2k}$ is extreme.
\vskip0.1in
{\bf 2)}Let us now assume that
the non-empty oval of the curve
$\HH_{2}$ contains all inner ovals ${\mathcal A}_{2k-2}$.
In the case {\bf 1)},
we have chosen arcs $I_j$, $ 3 \le j \le 2k$
such that an arc intersects at most an other arc in one point.
We need here one more assumption on the choice of arcs  $I_j$
of the non-empty oval of ${\mathcal A}_{2k-2}$.\\
Assume that for any $i ,1 \le i < k$, the non-empty oval of the
curve ${\mathcal A}_{2k-2i}$ does not contain inner ovals
of $\HH_{2i}$.
Let $I_3$ and $I_4$ be the two extreme arcs of ${\mathcal A}_{2k-2}$.
Then, we get $\HH_4$ and a curve ${\mathcal A}_{2k-4} \in \Omega_{2k-4}$
with the property  that $\HH_{2k}$ results from the deformation of
$\HH_4 \cup {\mathcal A}_{2k-4}$.
We iterate this process.
Consider the union $I_{j}$, $ 3 \le j \le {2k}$,
of arcs of ${\mathcal A}_{2k-2}$ as follows:
$I_3$ and $I_4$  are extreme arcs of ${\mathcal A}_{2k-2}$
Up to rigid isotopy, arcs $I_{j}$, $j >5$,
are arcs of the non-empty oval empty of ${\mathcal A}_{2k-4}$.
We choose $I_5$ and $I_6$ such that
$I_5$ intersects $I_3$ in one point,
$I_6$ intersects $I_4$ in one point.
In such a way, they are mapped to
extreme arcs of ${\mathcal A}_{2k-4}$.
By induction, we choose
$I_{3+2l}$, $I_{4+2l}$, $1 \le l \le (k-1)$ such that
$I_{3+2l}$ intersects $I_{3+2(l-1)}$ in one point,
$I_{4+2l}$ intersects $I_{4+2(l-1)}$ in one point;
and are
extreme arcs of ${\mathcal A}_{2k-2(l+1)}$.
At the end of the process,
we get $\HH_{2k}$ from the deformation of
${\mathcal A}_2 \cup \HH_{2k-2}$
where ${\mathcal A}_{2}= \HH_{2}$.
(Identifying arcs $I_j$ with lines,
 $\RRR L_{2k-1} \approx I_{2k-1}$,
 $\RRR L_{2k} \approx I_{2k}$
the curve
${\mathcal A}_{2}= \HH_{2}$, as $L$-curve,
results from the perturbation of the union $L_{2k-1} \cup L_{2k}$.)

It is an easy property of the proposed
construction that
the oval of ${\mathcal A}_{2}= \HH_{2}$
contains besides inner ovals of $\HH_{2k-2}$
inner ovals of $\HH_{2k}$ in its inner component.
Iterating this process, we get $\HH_{2k-2i}$, $i < k $,
from $\HH_{2k-2i -2} \cup {\mathcal A}_2$
where
the oval of ${\mathcal A}_{2}= \HH_{2}$
contains,
besides inner ovals of $\HH_{2k-2i-2}$
inner ovals of $\HH_{2k-2i}$ in its inner component.
(Identifying arcs $I_j$ with lines,
the curve
${\mathcal A}_{2}= \HH_{2}$, as $L$-curve,
results from the perturbation of
two lines
$\RRR L_{2k-2i-1} \approx I_{2k-2i-1}$,
$\RRR L_{2k-2i} \approx I_{2k-2i}$ in general position.)
\vskip0.1in

If there exists $j \le k $ such that the non-empty oval of
${\mathcal A}_{2k-2j}$ does not contain inner ovals of
$\HH_{2j}$, then
from the step when
the deformation $\HH_{2j} \cup {\mathcal A}_{2k-2j}$
leads to $\HH_{2k}$ to the final step
- the deformation of
 ${\mathcal A}_2 \cup \HH_{2k-2}$ leads to  $\HH_{2k}$~-
, we are in the case {\bf 1)}.\\
\vskip0.1in
In the previous argumentation, we have assumed that
the non-empty oval of
$\RRR {\mathcal A}_{2k-2}$
is the union
of an arc which does not contain points
of $\RRR {\mathcal A}_{2k-2} \cap  \RRR \HH_2$.
with a set of arcs $I_j$,
$3 \le j  \le 2k$ which contains points of
$\RRR {\mathcal A}_{2k-2} \cap  \RRR \HH_2$
From identification of these $2k-2$ arcs $I_j$
with lines $\RRR L_j$, $3 \le j  \le 2k$,
we have described components of $\RRR {\mathcal A}_{2k-2}$
as an arrangement of
pieces of these $2k-2$ lines.
The degree of ${\mathcal A}_{2k-2}$ is $2k-2$.
Hence, its real components may not result from more than $2k-2$
lines. Moreover, since ${\mathcal A}_{2k-2}$ is an $M$-curve,
its real components may not result from less than $2k-2$ lines.
Our initial assumption is therefore always possible.
Moreover,
one can assume that $\HH_{2k}$ results from the deformation
of  ${\mathcal A}_{2k-2} \cup \HH_2$
where ${\mathcal A}_{2k-2} \in \Omega_{2k-2}$.
and $\RRR {\mathcal A}_{2k-2}$ is entirely covered
by the set of arcs $\cup_{j=3}^{2k}I_j$.
Hence, since $\HH_2$ is, up to rigid isotopy, an $L$-curve;
we get Theorem \mrf{t:seqli} for Harnack curves $\HH_{2k}$
of even degree.
\vskip0.1in
{\bf II) Harnack Curves of odd degree $\HH_{2k+1}$}
\vskip0.1in
According to Proposition \mrf{p:l07},
up to regular modification of its Harnack polynomial $B_{2k+1}$,
one can assume
that $\HH_{2k+1}$ results from the deformation of
the union  ${\mathcal A}_{2k-2} \cup \HH_3$,
where ${\mathcal A}_{2k-2} \in \Omega_{2k-2}$
and $\HH_3$ is the Harnack curve of degree $3$.
Curves $\HH_3$ are rigidly isotopic.
One can consider $\HH_3$ as an $L$-curve i.e  as the result
of the perturbation of three real lines  $L_0$, $L_1$,$L_2$
in general position.
According to remark \mrf{r:rrig}, we shall assume
lines $L_0$,$L_1$, $L_2$ chosen as follows.
Outside a neighborhood ${\mathcal B}$
of singular points $L_0 \cup L_1 \cup L_2$,
there exists a rigid isotopy $j_t$ of $\RRR P^2$
which pushes $\cup_{i=0}^2 \RRR L_i \bk {\mathcal B}$
to $\RRR \HH_3 \bk {\mathcal B}$
such that
$\RRR \HH_{3}$ intersects $\RRR {\mathcal A}_{2k-2}$
only its part
$j_1 (\cup_{i=0}^2 \RRR L_i \bk {\mathcal B}) \bk j_1 (\RRR L_0)$.
Consider the deformation
$\HH_2 \cup {\mathcal A}_{2k-2}$
where $\HH_2$
is an $L$-curve which results from
the perturbation of the two real lines $L_1$,$L_2$
in general position.
Combining this construction of curves $\HH_{2k+1}$
and our previous study of curves of even degree,
we obtain  the Theorem \mrf{t:seqli}
for curves $\HH_{2k+1}$.
Q.E.D
\vskip0.1in
To prove
the rigid isotopy Theorem \mrf{t:rigiso}
it is sufficient to note that
according to Proposition \mrf{p:HH0} and
the Theorem \mrf{t:seqli} and its proof,
any curve $\HH_m$ is of type $\HH^0$.
Thus, the rigid isotopy Theorem \mrf{t:rigiso}
is  a straightforward consequence of the Theorem \mrf{t:rig}.

However, for clarity,
we propose to prove the Theorem \mrf{t:bli}.\\

\beth
\mlb{t:bli}
Up to rigid isotopy, any curve $\HH_{m}$, $m \ge 1$,
results from the successive classical small deformation
of the union $L_{i+1} \cup \HH_{i}$ , $ 1 \le i \le n-1$,
where
$\HH_{i}$ is a Harnack curve of degree $i$,
$L_{i+1}$ is a projective line and  the curve $\HH_{i}$ is of type
$\HH^0$ relatively to $L_{i+1}$.
\enth

Let us give in Proposition \mrf{p:Lemma A'}
an other formulation of the Theorem \mrf{t:bli}.
\bepr
\mlb{p:Lemma A'}
Let $\HH_m$ be a Harnack curve of degree $m$ and
$B_{m}(x_0,x_1,x_2)$ its polynomial.\\
Denote by $b_m(x,y)$ the affine polynomial associated to $B_m(x_0,x_1,x_2)$,
$$B_m(x_0,x_1,x_2)=x_0^m.b_m(x_1/x_0,x_2/x_0)$$
Denote by $b_i(x,y)$ the truncation of $b_m(x,y)$ on the monomials
$x^{\a}.y^{\b}$ with $0 \le \a + \b \le i $ and by
$B_i(x_0,x_1,x_2)= x_0^i.b_i(x_1/x_0,x_2/x_0)$
the homogeneous polynomial associated to $b_i$.
Then, up to linear change of projective coordinates, and
up to slightly perturb coefficients of $B_m(x_0,x_1,x_2)$
without changing either order or topological structure of $\HH_m$,
we may assume that any polynomial $B_i$ is of type $\HH^0$
relatively to $x_0=0$.
\enpr

\vskip0.1in

{\bf proof:}
According to Theorem \mrf{t:seqli},
up to rigid isotopy, there exists a sequence of lines $\HH_1=L_1$,
$L_2$, ..., $L_{m+1}$ such that $\HH_i$ is of type $\HH^0$
relatively to $L_{i+1}$ and $\HH_{i+1}$ results
from the deformation of $\HH_{i} \cup L_{i+1}$.
The sequence of lines $L_i$ may be ordered in three different ways
described in  (\ref{i:klo}), (\ref{i:kli}) and (\ref{i:klio}).
To prove to Theorem \ref{t:bli}, it is sufficient to prove
that the description (\ref{i:klo}) is rigidly isotopic
to  the description (\ref{i:kli}).\\
Without loss of generality, one can
assume that points $\HH_i \cap L_{i+1}$
do not belong to the line at infinity of $\RRR P^2$.
Any affine line $l_j=\RRR^2\cap L_j$ is divided into two
halves $l_j^+$ $l_j^-$ with common point $p$.
(The rotation of center $p$ maps a positive (resp, negative) half
of one line to a positive (resp, negative) half of an other line.)
Consider the half plane $\RRR^+$ (resp, $\RRR^-$)
which contains positive (resp, negative) halves lines.
Let $H$ be the common line of $\RRR^+$ and $\RRR^-$
(where $H \cap L_j=p$, for any $L_j$, $1 \le j \le m$).
The symmetry $s_p$ of center $p$ exchanges halves of lines.
Let us denote $s_H$ the orthogonal reflection with respect to $H$.
Assume, as described in (\ref{i:kli}),
that
there exists $p \in L_2 \cap L_j$ , $3 \le j \le m+1$,
 such that:\\
- for $i$, $2 \le  i \le 2.[(m+1)/2]-2$
 the rotation of center $p$ and angle $\theta$
 maps $\RRR L_i$ to $\RRR L_{i+2}$.\\
- for  $i$, $1 \le  i \le 2.[(m+1)/2]-2$
 the rotation of center $p$ and angle $-\theta$
  maps $\RRR L_{i+1}$ to $\RRR L_{i+3}$
It follows from the remark \ref{r:rrig}
that  without loss of generality,
one can assume that $\HH_{i}$ intersects $L_{i+1}$ in its half $l_{i+1}^+$.
Choose the lines $L_j$ such that the symmetry $s_H$ exchanges
lines $L_{2l-1}$ and $L_{2l}$,
$s_H(L_{2l-1})=L_{2l}$,$s_H(L_{2l})=L_{2l-1}$
According to remark \mrf{r:rrig}, the
rigid isotopy is trivially preserved under the action of $s_H$.\\
Let us now choose lines $L_{j}$
such that the rotation of center $p$ and angle $\pm \frac {\theta} 2$
maps  $s_H(L_{2l})$ to $L_{2l \pm 1}$.
Obviously, the description of $\HH_m$ which results from
this arrangement of lines $L_j$ and the one where
$s_H(L_{2l-1})=L_{2l}$,$s_H(L_{2l})=L_{2l-1}$ are rigidly isotopic.
The rigid isotopy is preserved by the symmetry $s_p$
which maps the intersection of $\HH_i$ with $L_{i+1}$ (namely,
with $l_{i+1}^+$)
to a set of the half $l_{i+1}^-$.
The rigid isotopy is
also preserved when consider $s_H(L_{2j})=L'_{2j}$
and $L_{2j-1}$ i.e $l_{2j}^{'+}$ and $l_{2j-1}^+$
instead of halves lines $l_{2j-1}^-$, $l_{2j}^-$
In this way, we get that
the description (\ref{i:kli})
is rigidly isotopic to the description  (\ref{i:klo}).
This proves the Theorem \mrf{t:bli}.
Q.E.D

\vskip0.2in
Harnack curves $\HH_m$ constructed from the Harnack's method
are rigidly isotopic.
It follows that the Theorem \mrf{t:bli} implies
the Rigid Isotopy Classification Theorem \mrf{t:rigiso}.

{\bf Appendix}
\vskip0.1in

Curves $\HH_m$ of type $\HH^0$ may result from
the Harnack's inductive construction (\cite{Har}, \cite{Vi})
as follows.\\
Given a projective line $L$, each curve $\HH_{i+1}$,
$1 \le  i \le m-1$, results
from the classical small perturbation of the union
$\HH_i \cup L$ where $\HH_i$ is of type $\HH^0$ relatively
to $L$; $\HH_1$
is a projective line which intersects $L$ in one point.\\
We take the auxiliary curve ${\mathcal C}_{i+1}$ which perturb the union
$\HH_i \cup L$ to be a union of $i+1$ lines.
\vskip0.1in
Consider $L$ as the line at infinity of $\RRR P^2$.
From the Harnack's method \cite{Har} initiated with $L$,
one can obtain
isotopy of curves different from $\HH_m$
and curves $\HH_m$ which are not of type $\HH^0$
relatively to $L$.
The isotopy type of curves (projective and affine)
which result from  the Harnack's method
depends on the choice of auxiliary curves.\\
Let us prove in Theorem \mrf{t:ty0}
that isotopy implies also rigid isotopy for
Harnack curves of type $\HH^0$ obtained from this method.
For curves of degree $m \le 6$,
the Theorem is obvious.\\

\beth
\mlb{t:ty0}
Any two Harnack curves $\HH_m$ of type $\HH^0$
constructed from the Harnack's method
are rigidly isotopic.
\enth

{\bf proof:}
We shall proceed by ascending induction on the degree.\\
The  Theorem \mrf{t:ty0} is trivial
for Harnack curves of degree $\HH_i$, $i \le 6$.
Let us start our induction
with Harnack curves of degree $1$.
Let $L_1$ and $L_2$ be two projective lines,
and $\HH_1^1$, $\HH_1^2$ be Harnack curves of degree $1$.
Consider the classical small perturbation of
$\HH_1^1 \cup L_1$ and $\HH_1^2 \cup L_2$.
Denote by
$\HH_2^1$,(resp, $\HH_2^2$), the Harnack curve of degree $2$
which results from the classical deformation of
$\HH_1^1 \cup L_1$ (resp, $\HH_1^2 \cup L_2$).

Any two real projective lines are rigidly isotopic.
Let $h_t$ be the rigid isotopy
$[0,1] \to  \RRR {\mathcal C}_1 \bk \RRR {\mathcal D}_1$
,$h_{0}(\RRR L_{1})= (\RRR L_{1})$, $h_{1}(\RRR L_{1})
=\RRR L_{2}$.
Denote by $\RRR  L_{t+1}$ the line $h_t(\RRR L_{1})$.
Any two projective lines intersect each other
in one point. Hence, one can assume, without loss of generality
that $\HH_1=\HH_2=L_0$.
Along the rigid isotopy $h_t$, any line
$\RRR L_{t+1}=h_t(\RRR L_{1})$ intersects $\RRR \HH_0$ in one point.
Noticing that the classical perturbation of the union of two lines
which intersect each other in one point leads to the Harnack curve $\HH_2$,
the path $h_t * id$:
$[0,1] \to \RRR {\mathcal D}_{2}$
extends to a path $g_t:$
$[0,1] \to \RRR {\mathcal C}_{2} \bk \RRR {\mathcal D}_2$,
$g(0) = \RRR \HH_{2}^1$,
$g(1)= \RRR \HH_{2}^2$ which is a rigid isotopy.\\

Given $\HH_{m}^1$, (resp, $\HH_{m}^2$), $i \ge 2$
the Harnack curve $m$ deduced from the Harnack's method
initiated with the line $L_1$, (resp, $L_2$).
We shall prove that curves $\HH_{m}^1$
and $\HH_{m}^2$ are rigidly isotopic for arbitrary degree $m$.
\vskip0.1in

Let us start with the following remark.
\bele
\mlb{l:ii}
Let $L$ be a real projective line.
Consider the Harnack's inductive construction of curves
initiated with the line $L$.
Any two
Harnack curves $\HH_m$ of type $\HH^0$
obtained from the Harnack's method initiated with the line $L$
are rigidly isotopic.
\enle
{\bf proof:}
Indeed, consider
the Harnack's construction  of the curve $\HH_m$
initiated with the line $L$.
In this construction,
any curve $\HH_i$, $1 \le i \le m-1$,
is of type $\HH^0$ relatively to $L$.
Consider $L$ as the line  at infinity of $\RRR P^2$.
Given a Harnack curve $\HH_i$, denote by $h_i$ the affine
corresponding curve.
Denote by ${\mathcal C}_{i}$ the curve which deforms the union of $L$
with the curve of degree $i-1$, $i \le m$.
Let
the auxiliary curve of degree $j$
${\mathcal C}_j$
be the union of $j$ lines of a pencil through a point.
Noticing that the isotopy type of the affine curve $h_i$ ($1 \le i \le m$)
depends exclusively on the position of ${\mathcal C}_{j} \cap L$
,$1 \le j \le i $, we shall prove the Lemma \mrf{l:ii}.
Let us detail our argumentation.
\vskip0.1in

Let us prove that  any two curves $\HH_{i+1}$ and $\HH_{i+1}'$
resulting from  the classical perturbation of  $\HH_i \cup L$
directed to a union of $i+1$ lines
${\mathcal C}_{i+1}$, resp ${\mathcal C}_{i+1}'$
 are rigidly isotopic.\\
Denote by $x_0$, $B_i$,
 $C_{i+1}$,(resp, $C_{i+1}'$), $B_{i+1}$
,(resp,$B_{i+1}'$)
the polynomial of $L$, $\HH_i$,
and
${\mathcal C}_{i+1}$, (resp ${\mathcal C}_{i+1}'$), $\HH_{i+1}$,
(resp, $\HH'_{i+1}$).\\
In the construction of $\HH_{i+1}$ of type $\HH^0$ from $\HH_{i} \cup L$,
the curve ${\mathcal C}_{i+1}$ intersects $L$ in
$i+1$ distinct points lying\\
-for even $i+1$
(i.e construction of $\HH_{2k}$ from $\HH_{2k-1} \cup L$)
in the component $S_i$ of $\RRR L$ which is a boundary of
the unique  non-empty positive region $\{ x \in \RRR^2 | b_{i}(x)>0 \}$
the non-empty region of $\RRR ^2$ with boundary a part of the
odd component of $\RRR \HH_i$).\\
-for odd $i+1$,
(i.e construction of $\HH_{2k+1}$ from $\HH_{2k} \cup L$)
in the component $S_i$ of $\RRR L \bk \RRR \HH_{i}$
containing $\RRR L \cap \RRR \HH_{i-1}$.\\
(It is also necessary that ${\mathcal C}_{i+1}$ does not intersect
$\HH_{i} \cup L$ in its singular points).\\
{\bf i.)}
If ${\mathcal C}_{i+1} \cap L ={\mathcal C}'_{i+1} \cap L$,
one can vary continuously the direction of lines
of ${\mathcal C}_{i+1}$ and in this way define a one parameter family
${\mathcal C}_{i+1,t;~t \in [0,1] }$,
${\mathcal C}_{i+1,0}={\mathcal C}_{i+1}$,
${\mathcal C}_{i+1,1}={\mathcal C}'_{i+1}$, and thus a rigid isotopy
$x_0.B_i+C_{i+1,t}$ from $\HH_{i+1}$ to
$\HH_{i+1}'$.\\
{\bf ii.)}
Otherwise ${\mathcal C}_{i+1} \cap L \not= {\mathcal C}'_{i+1} \cap L$.
In the construction of $\HH_{i+1}$ of type $\HH^0$ from
classical deformation of $\HH_i \cup L$, any  auxiliary
curve which deforms  the union $\HH_i \cup L$
intersects $\RRR L$ in $i+1$ points lying
in  a connected part of $S_i$, $S_i \subset\RRR L$.
Given ${\mathcal C}_{i+1}$,
denote ${\mathcal C}_{i+1} \cap L$ by $I_{i+1}$
and let $I'_{i+1}$ be a set of $i+1$ real points of $S_i$.
Let us distinguish the cases $i+1$ even and odd.\\
-In case $i+1=2k$ even,
(i.e construction of $\HH_{2k}$ from $\HH_{2k-1} \cup L$)
 $S_i=S_{2k-1}$ is connected.
According to the previous study, one can assume that lines of
${\mathcal C}_{2k}$ and ${\mathcal C}'_{2k}$ have the same direction.
It is not hard to see that there exists a rigid isotopy
$B_{2k,t}$ with $t \in [0,1]$
of $B_{2k,0}=x_0.B_{2k-1} +C_{2k}$, $B_{2k,1}=x_0.B_{2k-1}+C'_{2k}$ such that
${\mathcal C}'_{2k}(x_0,x_1,x_2)  \cap L=I'_{2k}$.\\
-The case $i+1=2k+1$ odd
(i.e construction of $\HH_{2k+1}$
from $\HH_{2k} \cup L$) differs from the case $i$ odd
in the sense
that the set $S_i=S_{2k}$ has two connected parts.
These two connected parts  may  be defined using a projection
to $L$ in a direction perpendicular to $L$
as follows.
The non-empty oval of $\RRR \HH_{2k}$ intersects the line $L$
in $2k$ points.
They belong to a segment of $L$
with extremities $2$ points of $\RRR \HH_{2k} \cap L$.
Using a pencil of lines in a direction (for example
in a direction perpendicular) to $L$,
project the non-empty oval of $\RRR \HH_{2k}$ to a segment $S$ of $L$.
The set $S_{2k}$ consists of the two  segments  $S_{2k}^1=[a_1,b_1[$
and $S_{2k}^2=[a_2,b_2[$ with extremities an
extremity of $S$ and a point of $\RRR \HH_{2k} \cap L$
(this last one is the open extremity of the segment).
Changing the direction of the projection, one defines
two other segments.
(The idea in our proof is to find a way to join
this two segments. It may be done
by changing the direction of the projection.)\\
Let $\HH_{2k+1}$, such that $\HH_{2k+1} \cap S_{2k}^2 = \emptyset$
and $\HH'_{2k+1} \cap S_{2k}^2= \not \emptyset$.
To prove that the curves $\HH_{2k+1}$ and $\HH'_{2k+1}$ are
rigidly isotopic, it is sufficient to prove that
there exists a continuous path  $\HH_{2k+1,t}$
$t \in [0,1]$, $\HH_{2k+1,0}=\HH_{2k+1}$, $\HH_{2k+1,1}$
with the following properties:\\
-for any $t \in [0,1]$
$\HH_{2k+1,t}$ is a Harnack curve of degree $2k+1$\\
-as $t$ varies from $0$ to $1$
a line of $\HH_{2k+1}$ (i.e of ${\mathcal C}_{2k}$)
intersecting
 $S_{2k+1}^1$
 is deformed to a line of ($\HH'_{2k+1,1}$)
intersecting
 $S_{2k}^2$.

(Recall that
if ${\mathcal C}_{2k+1} \cap L ={\mathcal C}_{2k+1}' \cap L$
there exists a rigid of isotopy
$x_0.B_{2k}+C_{2k+1,t}$, $t \in [0,1]$, from $\HH_{2k+1}$ to
$\HH'_{2k+1}$.
Hence, we shall consider
${\mathcal C}_{2k+1}$  and ${\mathcal C}'_{2k+1}$
up to such rigid isotopy.)
\vskip0.1in
Let $\gamma$ be the non-empty arc (of the non-empty oval)
of $\HH_{2k}$
which intersects $L$ in two points $a_1 \in S_{2k}^1=[a_1,b_1[$
$a_2 \in S_{2k}^2=[a_2,b_2[$.
Let $p \in {\gamma}$, consider a path ${\gamma}_t$,
${\gamma}_0=p$ which moves $p$ along $\gamma$.

Let $p_1 \in  S_{2k}^1$, and $L_1=(pp_1)$ be a line
through  $p$.
Consider the pencil of lines through the point $p$
which intersect $\gamma$.
Move continuously the point $p$ along
the path $\gamma_t \subset \gamma$ and
in this way
lines $(p_1 \gamma_t)$ of the pencil of lines through $p_1$.
The projection of $p$ (in a direction perpendicular to $L$)
to $L$ belongs to $S_{2k}^1$.
It is not hard to see that
one can choose the path ${\gamma}_t$ such that,
as $p$ moves along $\gamma_t$,
the projection ${\gamma}_t$  to $L$
moves from $S_{2k}^1$ to $p_2 \in S_{2k}^2$.

In this way,
considering together pencil of lines through any point $p_t :={\gamma}_t$
and pencils through the point $p_1 \in S_{2k}^1$
and pencils through the point $p_2  \in S_{2k}^2$
it follows that there exists a rigid isotopy
of $\HH_{2k+1}$  from $\HH'_{2k+1}$
which maps a line of ${\mathcal C}_{2k+1}$ which intersects $S_{2k}^1$
to a line of a curve ${\mathcal C}'_{2k+1}$ which intersects $S_{2k}^2$.

It concludes the proof of the Lemma \mrf{l:ii}.
Q.E.D
\vskip0.1in

It follows from the Lemma \mrf{l:ii} that
when one considers Harnack's construction of curves $\HH_m$
up to rigid isotopy,
we may restrict our study to
classical deformation of $\HH_i \cup L$, $1 \le i \le (m-1)$,
directed to a union of $i+1$ lines of a pencil
through a point $p$ chosen outside $L$ and $\HH_i$.

Denote ${\mathcal C}_{i+1}$ the union of lines
which deform the union of $L$
with the curve of degree $i$, $i \le m$.
The isotopy type of curves of degree $m$
obtained from the Harnack's  recursive
method depends on the relative position of
intersection points of ${\mathcal C}_{i+1}$ with
the line $L$.
\vskip0.1in

Consider the Harnack's construction of curves $\HH_i$,
$1 \le i  \le m$.
Denote by $B_i$ the polynomial of $\HH_i$.
Consider $L$ as the line at infinity, and denote by
$b_i$ the affine polynomial
associated to $B_i$.
Curves $\HH_i$ constructed from the Harnack's method
are of type $\HH^0$ relatively to $L$.\\
Let $\HH_5$ be the Harnack curve of degree $5$.
There  exists a unique  positive region $\{ x \in \RRR^2 | b_5(x)>0 \}$
with a segment $S_5$ of the line $\RRR L$
on its boundary
which contains
an oval of $\HH_5$.
To get $\HH_6$  from $\HH_5 \cup L$, it is necessary that the intersection
of ${\mathcal C}_6$ with $L$ consists of $6$ points lying in $S_5$.
The non-empty oval of $\HH_6$
intersects $L$ in $6$ points.
Denote  $S_6$  the connected part of the line $\RRR L$
which contains the  intersection points of $\HH_6$
with $L$.
To get $\HH_7$ from $\HH_6 \cup L$, it is necessary that the intersection
of ${\mathcal C}_7$ with $L$ consists of $7$ points lying in
$L \bk S_6$ and that it intersects $\HH_6$ in its non-empty real component.
Iterating this process, (set $2k+1$ (resp, $2k$)
instead of $5$  (resp, $6$))
we define a sequence $S_i$ of connected components
of $\RRR L$ such that $S_i \subset {\mathcal C}_{i+1} \cap L$.\\
In other words, we define relative position of lines ${\mathcal C}_{i+1}$
as follows.
We take ${\mathcal C}_m$ to be a union of $m$ lines
which intersect $L$ in $m$ distinct points lying,
for even $m$
in the component of $\RRR L$ which is a boundary of
the unique  non-empty positive region $\{ x \in \RRR^2 | b_{m-1}(x)>0 \}$,
for odd $m$ in the component of $\RRR L \bk \RRR \HH_{m-1}$
containing $\RRR L \cap \RRR \HH_{m-2}$.
(It is also necessary to choose ${\mathcal C}_{m}$ such that
it  does not intersect
$\HH_{m-1} \cup L$ in its singular points).\\

Denote by $B_i^1$ (resp, $B_i^2$)
the polynomial of $\HH_i^1$ (resp, $\HH_i^2$)
and by $C_{i+1}^1$
(resp, $C_{i+1}^2$)
the union of $i+1$ lines
which deform
$\HH_{i}^1 \cup L_1$
(resp,
$\HH_{i}^2 \cup L_1$)
 to $\HH_{i+1}^1$
(resp, $\HH_{i+1}^2$).

To prove  Theorem \mrf{t:ty0},
is sufficient to construct
for $i,~1 \le i \le (m-1)$
a continuous one parameter family
with parameter $t \in [0,1]$
of curves ${\mathcal C}_{i+1}^{1+t}$
such that as $t$ varies from $0$ to $1$
the relative position of lines ${\mathcal C}_{i+1}^{1+t}$
remains the same.

The rigid isotopy $h_t$ is continuous.
Therefore, the function which associates to any line $h_t(\RRR L_1)$
its normal vector $\vec{n}_t$ and its tangent vector $\vec{v}_t$
is also continuous.
It follows from the continuity of $h_t$, that one can choose lines
${\mathcal C}_{i+1}^{1+t}$ such that
the relative position of intersections
${\mathcal C}_{i+1}^{1+t} \cap L^{1+t}$,
and also the direction of lines of ${\mathcal C}_{i+1}^{1+t}$
relatively to $(\vec{v}_t,\vec{n}_t)$
remains the same as $t$ varies from $0$ to $1$.
In such a way,
for any $i \le (m-2)$
we construct by induction
curves $\HH_{i+2}^t$, $t \in [0,1]$,
of degree $i+2$ rigidly isotopic.
Hence, for any $i \le (m-2)$,
curves $\HH_{i+2}^1$ and $\HH_{i+2}^2$ are rigidly isotopic.
In particular,
curves $\HH_{m}^1$ and $\HH_{m}^2$
are rigidly isotopic.
Q.E.D\section{Harnack Curves from a Complex viewpoint}
\mlb{s:comp}

In this section, we shall first  construct
particular deformations of Harnack polynomials.
Then, we shall deduce  from the properties of these deformations
a characterization of the complex set point of Harnack curves
in $\CCC P^2$.\\
We shall divide this section into two subsections.
Denote $\HH_m$ the Harnack curve of degree $m$  defined
up to isotopy  of  real points set.
In the first section, we define
(in proposition \mrf{p:prop7b})
deformation of any Harnack curve $\HH_m$ to a singular curve of
which singular points are critical points index $1$
of the Harnack curve $\HH_m$.
 At first, we  characterize in proposition \mrf{p:prop7}
such deformation for Harnack curves obtained via the pathchworking method
and then generalize it in proposition \mrf{p:prop7b}
to any Harnack curve.
Considering patchworking method seems
at first glance of less interest
at this time of our proof of the Rokhlin's Conjecture.
Nonetheless,
it will be useful in the second part (Perestroika Theory
on Harnack curves).
In the next section, we deduce from the proposition \mrf{p:prop7b}
a description of $(\CCC P^2, \CCC \HH_m)$
up to conj-equivariant isotopy and give the main result
(Theorem \mrf{t:theo1}) of this chapter.

\subsection{Deformation of Harnack Curves}
\vskip0.1in
Given a  non-singular curve $\PP_{m}$ in $\CCC P^2$,
call {\it simple deformation} $\PP_{m;t}$
of the curve $\PP_{m}$
a path
\begin{eqnarray}
& [0,1] \to \RRR {\mathcal C}_m &\nonumber\\
& t \to \PP_{m;t} &\nonumber
\end{eqnarray}
with the following properties:
\been
\item
For any $ 0 \le t < 1$, $\PP_{m;t}$ is a smooth curve.
If $\PP_{m;1}$ is singular, then it is irreducible
and any of its singular points is a
real crossing which is critical point with positive critical value
of the affine  Harnack polynomial.
Denote $S$ the set of singular points of$\PP_{m;1}$.
\item
For $\epsilon>0$,
let ${\mathcal D}_{\epsilon}= \cup_{a \in S}  D_(a,\e)$ be
the union of disc $D(a,\e)$ (in the Fubini-Study metric)
of center $a$ and radius $\epsilon$ taken over
crossings $a$ of $\PP_{m;1}$.\\
\enen

There exists $\epsilon_{0}>0$,
such that:
for any curve $\PP_{m,t}$, $t \in ]0,1[$,
$\CCC P_{m,t}$ lies in an $\e_0$-tubular $N$ neighborhood
of $\CCC P_{m,1} \bk \{ a \in S\}$.\\
$\CCC P_{m,t}$ can be extented to the image of a smooth
section of the tubular fibration
$N \to \CCC P_{m,1} \bk \{ a \in S\}$.\\

Under two distinct simple deformations $\PP_{m;t}$ and $\tilde{\PP_{m;t}}$
of $\PP_{m}$,
the number of real singular points of
 $\PP_{m;1}$ and
 $\tilde{\PP_{m;1}}$ may be distinct.

Among simple deformations we shall distinguish
deformations which lead to a curve $\PP_{m;1}$ with
the maximal number  of real singular points
a curve $\PP_{m;1}$ may have.
Denote by $\alpha$  this number.
We shall call
{\it maximal simple deformation} $\PP_{m;t}$, $t \in [0,1]$,
of the curve $\PP_{m}$ a simple deformation $\PP_{m;t}$
such that the curve  $\PP_{m;1}$  has
$\alpha$ real singular points.

In proposition \mrf{p:prop7b}, we define
some deformation  of $\RRR \HH_m$,
\begin{eqnarray}
& [0,1] \to \RRR {\mathcal C}_m &\nonumber\\
& t \to \RRR \HH_{m;t} &\nonumber
\end{eqnarray}
to the real points set $\RRR {\mathcal A}$
of a singular curve ${\mathcal A}$ of degree $m$
of which singular points are critical points of index $1$
of the Harnack curve $\HH_m$.
At first, we  characterize in proposition \mrf{p:prop7}
such deformation for Harnack curves obtained via the pathchworking method
and then generalize it in proposition \mrf{p:prop7b}
to any Harnack curve.
\vskip0.1in
Recall that given a Harnack polynomial $B_{2k}$ of degree
$2k$ type $\HH^0$  and associated affine polynomial $b_{2k}$,
$S'_{2k}$ denotes the  subset of critical points
$(x_0,y_0)$
of positive critical value $c_0$ with the property that as
$c$ increases from $c_0 - \e$ to $c_0 + \e$ the number of real connected
components of
$M_c =\{((x,y) \in \RRR ^2 |  b_{2k} > c \}$
bounding the line at infinity increases.\\
\subsubsection{Maximal Deformation of $T$-Harnack curves}
\vskip0.1in
Let us in Proposition \mrf{p:prop7}
characterize maximal simple deformation of $T$-Harnack curves.
It is obious that the Proposition \mrf{p:prop7} may be
generalized to any Harnack curve.
The general statement \mrf{p:prop7b}
may be proven independently of Proposition \mrf{p:prop7b}.
For sake of clarity, we present here deformation of $T$-Harnack
curves.

\bepr
\mlb{p:prop7}
{\it
\vskip0.2in
\been
\item
Along any simple deformation of a $T$-Harnack curve of odd degree $2k+1$
curves are smooth.
\item
Given $\HH_{2k}$ a T-Harnack curve  of even degree $2k$,
any maximal simple deformation of $\HH_{2k}$
\begin{eqnarray}
& [0,1] \to \CCC {\mathcal C}_m &\nonumber\\
& t \to \HH_{m;t} &\nonumber
\end{eqnarray}
is such that $\HH_{2k;1}$
has $2k-3$ crossings.\\
Let $S$ be the set of singular points of
$\HH_{2k;1}$.
Then, there exists a Harnack polynomial $B_{2k}$ of degree $2k$
and type $\HH^0$
of which the set of critical points contains the set $S$.\\
These points are critical points with positive critical value.
Besides, $k-1$ of these points are points of $S'_{2k}$.
\enen
\vskip0.1in}
\enpr
As explained in section \mrf{su:Patch},
one can give via the Patchworking method
an inductive construction of Harnack curves
called $T$-inductive contruction.\\

Recall that in the $T$-inductive construction of Harnack polynomials
$\vec{a}_{m},t_{m}$ denotes the pair constituted by a vector
$\vec{a}_m$ and a real $t_m$
such that for any $t \in ]0, t_{m}[$,
$\tilde{X}_{m;\vec{a}_{m},t}=
\sum_{(i,j)~ vertices~ of~ T_m} \e_{i,j}a_{i,j}
x_1^ix_2^jx_0^{m-i-j} t^{\nu(i,j)}$
is a Harnack polynomial of degree $m$.\\
Given $\vec{a}_{m},t_{m}$,
$\vec{a}_{m+1}=(\vec{a}_{m},\vec{c}_{m+1})$.

\bere
\been
\item
There exists $\vec{c}_{2k+1}$ such that
$t_{2k+1}=t_{2k}$.
\item
There exists $\vec{c}_{2k}$ and $\tau <t_{2k-1}$
such that for $t=\tau$
the curve given by $\tilde{X}_{2k;\vec{a}_{2k},\tau}$
has $2k-3$ singular points
which are collinear crossings.
\enen
\enre

{\bf proof:}
\vskip0.1in
Consider a $T$-Harnack curve of degree $m$.\\
For $m=1$ and $m=2$,
the proposition \mrf{p:prop7} is trivially verified.
For $m >2$, it is based on
the previous proposition \mrf{p:crit2}.
\vskip0.1in
We shall  use the terminology introduced in the
chapter \ref{ch:PatchHar} and denote
$\tilde{X}_{m;t}$, $t \in ]0, t_{m}[$ a T-Harnack polynomial of degree
$m$.
It is obvious that any $T$-Harnack polynomial is of type $\HH^0$.
In what follows, we shall deal with one-parameter polynomials
$\tilde{X}_{m;t}$, $t \in [t_{m}, t_{m-1}[$.

When no confusion is possible,
we shall denote $B_m$ a regular Harnack polynomial
of degree $m$ of type $\HH^0$ and $S^-_m$, respectively $S^+_m$,
the set of its critical points of negative, respectively positive,
critical value.

\vskip0.1in
\vskip0.1in
Let us recall a result of the patchworking theory.

\vskip0.1in
\bele
\mlb{l:Lemma A}
(\cite {Vi}, p.190)
{\it
Let $a$ be  a polynomial such that
$a=0$ admits an  $\e$-tubular neighborhood.
If a set $U \subset \CCC P^2$  is compact and contains no
singular points of $a=0$,
then for any $\e>0$ and any  polyhedron $\D \supset \D(a)$
there exists $\a>0$ such that for
any polynomial $b$ with $\D(b) \subset \D$, $ ||b-a|| < \a$
and  $b^{\D(a)}=a$ the truncation
$b^{\D(a)}$ is $\e$-sufficient in $U$.}
\enle
\vskip0.1in
We shall consider polynomials
$\tilde{X}_{m+1;\tau}$ with the following properties.
For any elementary triangle $\D$
of the triangulation of $T_{m+1}$:
\been
\item
the truncation
$\tilde{X}_{m+1;\tau}^{\D}$
is completely non-degenerate.
\item
$\tilde{X}_{m+1;\tau}^{\D}$ is $\e$-sufficient for
$\tilde{X}_{m+1;\tau}$ in $\rho^{m+1} (\RRR_+ \D^0 \ti U_{\CCC^2})$
\enen

We shall say that such polynomial $\tilde{X}_{m+1;\tau}$
satisfies the {\it "good truncation properties".}
It is easy to see that any polynomial $\tilde{X}_{m+1;\tau}$ with
good truncation properties
is such that for any elementary triangle $\D$,
the truncation
$\RRR \tilde{X}_{m+1;\tau}^{\D}$ is isotopic
to $\RRR \HH_{m+1}$ in the open
$\rho^m(\RRR_+ \D^0 \ti U_{\CCC^2}) \subset (\CCC^*)^2$.
(Obviously, any Harnack polynomial has good truncation properties.)
\vskip0.1in
Let us give the main ideas which motivate our study and
 explain the method of our proof.
\vskip0.1in
{\it Motivation}\\
The $T$-inductive construction of Harnack polynomials
can be considered as a slightly modified version of Harnack's initial one.

Fix $t_0 \in ]0,t_{m}[$, and denote by $B_m:=\tilde{X}_{m;t_0}=
\sum_{(i,j)~ vertices~ of~ T_m} \e_{i,j}a_{i,j}
x_1^ix_2^jx_0^{m-i-j} t_0^{\nu(i,j)}$
the corresponding T-Harnack
polynomial.

The T-Harnack polynomial $\tilde{X}_{m+1}
=\sum_{(i,j)~ vertices~ of~ T_{m+1}} \e_{i,j}a_{i,j}
x_1^ix_2^jx_0^{m+1-i-j} t_0^{\nu(i,j)}$
of degree $m+1$ may be deduced
by the formula:
$\tilde{X}_{m+1;t}=x_0.B_m +t_0.C_{m+1}$
where  $x_0$ is a line,  the curve given by
$C_{m+1}$ is the union of $m+1$ parallels lines which
intersect $x_0=0$
and do not pass through the singular points of
the curve $\AA_{m+1}$
of degree $m+1$ with polynomial
$x_0.B_m$.
For $t$ sufficiently small, in particular for $t_0 < t_{m+1}$,
$\tilde{X}_{m+1;t}=x_0.B_m +t.C_{m+1}$ is a Harnack polynomial $B_{m+1}$.
\vskip0.1in
Therefore, outside $\RRR C_{m+1}$ the curve
$\RRR \tilde{X}_{m+1;t_0}$
is  a level curve of the function
$\frac {x_0.B_m} {C_{m+1}}$.
On $\RRR \AA_{m+1} \bk \RRR C_{m+1}$,
this level curve has critical points only at the
singular points of $\RRR \AA_{m+1}$.
These points are non-degenerate singular points.
Hence, the behavior of
$\RRR \tilde{X}_{m+1;t_0}$
outside $\RRR C_{m+1}$ is described by the implicit function theorem
and Morse Lemma.
In particular, we have the following description:
\been
\item
Let $\{a_1,...,a_m\}$ be the set of crossings of
$x_0.B_m=0$.
Denote $D(a_i,\e)$ a small neighborhood of radius $\e$ around
$a_i$ in $\RRR P^2$.

Denote ${\mathcal D}_{\e}= \cup_{i=1}^{m} D(a_i,\e)$
the neighborhood of the set
of singular points of $\RRR \AA_{m+1}$
in $\RRR P^2$.
Let $N$ be a  tubular neighborhood of
$\RRR \AA_{m+1} \bk {\mathcal D}_{\e}$ in $\RRR P^2 \bk {\mathcal D}_{\e}$.
Then, there exists $\e_0$ such that  for
any $0 <\e \le \e_0$ and any $D(a_i,\e)$ of ${\mathcal D}_{\e}$
there exists a homeomorphism
$ h :D(a_i,\e) \to D^1 \ti D^1$ (where $D^1$ is a one-disc of $\RRR$)
such that
$$h(\RRR \HH_{m+1} \cap D(a_i,\e))=
\{ (x,y) \in D^1 \ti D^1 | x.y= \frac {1} {2} \}$$
\item
Moreover, for any $0 < \e \le \e_0$,
$\RRR \HH_{m+1} \bk {\mathcal D}_{\e}$ is a section of the tubular fibration \\
$N \to \RRR \AA_{m+1} \bk {\mathcal D}_{\e}$.
\enen
\vskip0.1in

{\it Method}\\
Assume that there exists a singular polynomial
$\tilde{X}_{m+1;\tau}$
with good truncation properties
and denote by $S$ the set of its singular points.
On this assumption, the Harnack curve $\HH_{m+1}$
is the image of a smooth section of the tubular fibration
$N \to \{ (x_0:x_1:x_2) \in \CCC P^2 | \tilde{X}_{m+1;\tau}=0  \}\bk {S}$.

Moreover,
outside $x_0.B_{m}=0$ and outside $C_{m+1}=0$,
curves $\tilde{X}_{m+1;\tau}=0$ and $\tilde{X}_{m+1;t}=0$, $0 < t< t_{m+1}$,
are isotopic.
Thus, $\HH_{m+1}$ and
$\tilde{X}_{m+1;\tau}=0$ may be not isotopic
only in $\e$-neighborhood $U(p)$ of $\HH_{m+1}$ defined from points
of faces $\G \subset l_m \subset T_{m+1}$
(see definition \mrf{d:def6}).
\vskip0.1in

Denote $U(p_i)$, $1 \le i \le m$, the $\e$-neighborhood defined from points
of a face  $\G_i=\{ x \ge (m-i), y \le i, x+y=m \}$,
$\G_i \subset l_m \subset T_{m+1}$.
The union ${\mathcal  B}= \cup_{i=1}^{m} U(p_i)$
is a neighborhood of the set of singular points of $\CCC \AA_{m+1}$
in $\CCC P^2$.
Let $N$ be the $\e$-tubular neighborhood of
$\CCC \AA_{m+1} \bk {\mathcal  B} $ in $\CCC P^2 \bk {\mathcal  B}$.
It follows immediately  from the Lemma \mrf{l:Lemma A} that
given $U \subset \CCC P^2$ compact which contains no singular points
of $\CCC \AA_{m+1}$, any polynomial $x_0.\tilde{X}_{m,t}$ with $t \in ]0,t_m[$,
is $\e$-sufficient for $\tilde{X}_{m+1,t}$ in $U$.
In other words, for $t \in ]0,t_{m+1}[$,
the intersection $U \cap \CCC \HH_{m+1}$ is contained in $N$
and can be extended to the image of a smooth section of  a tubular
fibration $N \to \CCC \AA_{m+1} \bk {\mathcal  B}$.\\
According to the corollary \mrf{c:patch} of chapter \ref{ch:PatchHar} and
its proof, in the patchworking scheme
crossings of $\HH_{m} \cup L$ are in bijective correspondence with
faces $\G_i=\{ x \ge (m-i), y \le i, x+y=m \}$
of the triangulation  of $T_{m+1}$.
Hence, we may assume that any crossing $a_i \in \HH_m \cup L $
belongs to the $\e$-neighborhood $U(p_i)$ of $\HH_{m+1}$
defined from points $\G_i^0$,
and consider the truncation
of the affine restriction of a Harnack polynomial $b_{m+1}$
on the monomials
$x^ {m-i}y^ {i-1}, x^{m-i}y^{i}, x^{m-i+1}y^{i-1}, x^{m -i+1}y^{i}$
which is $\e$-sufficient for $b_{m+1}$ in $U(p_i)$.

Using these truncations,
and the fact that for any neighborhood $B(p)$ of a singular point $p$
of $\tilde{X}_{m+1;\tau}=0$,
$\CCC \HH_{m+1} \cap B(p)$ is the image of a smooth section
of the tubular fibration \\
$N \cap B(p) \to \{(x_0:x_1:x_2) \in \CCC P^2 | \tilde{X}_{m+1;\tau}=0 \}
\cap B(p) \bk {p}$
we shall characterize singular points of curves $\tilde{X}_{m+1;\tau}=0$.
Then, according to the definition of the vector $\vec{a}_{m+1}$,
we shall construct singular polynomials  $\tilde{X}_{m+1;\tau}$
with singular points on ${\mathcal  B}$ and thus characterize
maximal simple deformation of $T$-Harnack curves.
\vskip0.1in
The exposition of the main ideas of our proof is now finished and
we shall now proceed to precise arguments.
\vskip0.1in

In order to construct singular polynomials
$\tilde{X}_{m;t}$, $t \in [t_{m}, t_{m-1}[$ with good truncation
properties,
we shall distinguish the cases of even  and odd $m$.

We shall work with
notations and definitions
related to the patchworking theory given in the chapter \ref{ch:PatchHar}.

{\bf i)} Consider the case $m=2k$.

{\bf i.1)}
Let us prove the following Lemma.
\bele
\mlb{l:even}
Any polynomial $\tilde{X}_{2k;\vec{a}_{2k},t}$ which
satisfies good truncation properties has at most $2k-3$
singular points. These points are crossings.
\enle

{\bf proof:}
\vskip0.1in
The proof is based on the Patchworking theory.
From a local study
on the $e$-neighborhoods $U(p_i)$ of $\HH_{2k}$
defined from points $\G_i^0$ of faces
 $\G_i=\{ x \ge (m-i), y \le i, x+y=m \}$,
we shall deduce where singular points
of a polynomial $\tilde{X}_{2k;t}$ may appear.
\vskip0.1in
{\bf i.a)}
Let $S$ be a square defined by vertices\\
$(c,d), (c+1,d), (c,d+1),(c+1,d+1)$ with $c,d$ odd and\\ $c+d =2k-2$.
It is contained in the interior of the triangle $T_{2k}$.
(Obviously, there are $(k-1)$ such squares contained in the interior
of $T_{2k}$.)

Consider polynomials
$\tilde{X}_{2k;t}$, $t \in ]0,t_{2k-1}[$\\
$\tilde{X}_{2k;t}= x_0.\tilde{X}_{2k-1;\vec{a}_{2k-1},t}+
C_{2k;\vec{c}_{2k},t}$
satisfying good truncation properties.
(Given  the vector $\vec{a}_{2k-1}$,
the vector $\vec{c}_{2k} \in ((\RRR ^*)^+)^{2k+1}$ is the vector which
may be chosen.)

Denote $x_{2k;t}$ the affine restriction of a polynomial
$\tilde{X}_{2k;t}$, $t \in ]0,t_{2k-1}[$.

Denote
$x_{2k;t}^S (x,y)$ its truncation on the monomials\\
$x^cy^d, x^cy^{d+1}, x^{c+1}y^d, x^{c+1}y^{d+1}$.
Namely,
$$x_{2k;t}^S (x,y)=
a_{c,d}t^{\nu(c,d)}x^cy^d
+ a_{c,d+1}t^{\nu(c,d+1)}x^cy^{d+1}
+a_{c+1,d} t^{\nu(c+1,d)}x^{c+1}y^d$$
$$-a_{c+1,d+1}t^{\nu(c+1,d+1)}x^{c+1}y^{d+1}$$
with $a_{c,d}>0, a_{c+1,d}>0, a_{c,d+1}>0, a_{c+1,d+1}>0$.

Let $\G$ be the face in the triangulation  of $T_{2k}$
given by coordinates $(c,d+1), (c+1,d)$.
According to the patchworking construction, it follows that
for $t \in ]0,t_{2k}[$
the truncation $x_{2k;t}^S$ is $\e$-sufficient for $x_{2k;t}$
in an $\e$-neighborhood $U(p)$ of $x_{2k;t}^S=0$
defined from  points of $\G^0$.)
\vskip0.1in

Fix $t \in ]0,t_{2k-1}[$,
and let $a_t$ be the crossing of $\HH_{2k-1} \cup L$ contained in $U(p)$.

According to the patchworking construction,
there exists an homeomorphism
$\tilde{h} :\CCC \HH_{2k} \cap U(p) \to
 \{(x,y) \in (\CCC^*)^2
| x_{2k;t}^S =0 \}\cap U(p)$
such that $\tilde{h}(a_t)=(x_1,y_1)_t$
is a critical point of $x_{2k;t}^S(x,y)$.
Obviously, (see corollary \mrf{c:patch}),
$(x_1,y_1)_t$ is such that $x_1.y_1 <0$.

Set
$x_{2k;t}^S(x,y)=l_t(x,y)+t.k_t(x,y)$
with\\

$l_t(x,y)=
 a_{c,d}  t^{\nu(c,d)}x^cy^d  +
 a_{c+1,d} t^{\nu(c,d+1)}x^{c+1}y^d $\\
$k_t(x,y)=
 a_{c,d+1} t^{\nu(c+1,d) -1}x^cy^{d+1}  +
-a_{c+1,d+1} t^{\nu(c+1,d+1) -1}x^{c+1}y^{d+1}$.
\vskip0.1in
Note that up to modify the coefficients
$a_{c,d},a_{c,d+1},a_{c+1,d},a_{c+1,d+1}$ if necessary, the point
$\tilde{h}(a_t)=(x_1,y_1)_t$ is
a critical point of the function
$\frac {l_t(x,y)} {k_t(x,y)}$ with positive critical value
$\le t< t_{2k-1}$.

On this assumption, it follows from the equalities \\
$l_t(x,-y)=-l_t(x,y)$, $k_t(x,-y)=k_t(x,y)$,\\
$\frac {\pr l_t} {\pr x}(x,-y) =-\frac {\pr l_t} {\pr x}(x,y)$,
$\frac {\pr l_t} {\pr y} (x,-y) =\frac {\pr l_t} {\pr y} (x,y)$,\\
$\frac {\pr k_t} {\pr x} (x,-y)= \frac {\pr k_t} {\pr x} (x,y)$,
$\frac {\pr k_t} {\pr y} (x,-y) = -\frac {\pr k_t} {\pr y} (x,y)$\\
that $(x_1,- y_1)_t$ is a critical  point of
the function $\frac {l_t(x,y)} {k_t(x,y)}$ with negative critical value
$\ge - t>- t_{2k-1}$.

Therefore,
for fixed $t \in [t_{2k}, t_{2k-1}[$,
one can choose $C_{2k;t}$
and thus
$\tilde{X}_{2k;t}$, in such a way that
$\tilde{h}^{-1}(x_1,- y_1)_t$   with
$x_1.y_1<0$  is a singular point
of
$\tilde{X}_{2k;t} =0$, and one-parameter curves
$\tilde{X}_{2k;t}$ are Harnack curves for $t>0$ sufficiently
small.
\vskip0.1in
Moreover,
the singular situation, we shall denote by
{\bf ${\mathcal S}'$}, is as follows:
\vskip0.1in

{\it
{\bf ${\mathcal S}'$}:
 one outer oval of the part
$j_1(\RRR \HH_{2k-1}\bk {\mathcal D}_{\e})$
of $\RRR \HH_{2k}$  touches the non-empty outer
oval in the part
$\RRR \HH_{2k} \bk j_1(\RRR \HH_{2k}\bk {\mathcal D}_{\e})$.}\\
(Namely, one branch of a outer oval
 contained in the patchworking scheme in the subset
$\rho^{2k}(D_{2k-1,2k-2} \ti U_{\RRR}^2) \subset \RRR P^2$
touches a branch of the non-empty outer oval.)
\vskip0.1in
Thus, according to the
Proposition \mrf{p:crit2} and Petrovskii's theory
(Petrovskii's Lemmas (Lemma.2 and Lemma.3)), it follows that
the point $\tilde{h}^{-1}(x_1,-y_1)_t$
is a critical point of a Harnack polynomial $B_{2k}$ of type $\HH^0$
and belongs to $S'_{2k}$.
\vskip0.1in
{\bf i.b)}
Likewise, consider squares $S$ defined by vertices\\
$(c+1,d), (c,d), (c,d+1), (c+1,d+1)$ with $c+d =2k-2$
and\\ $c >0, d >0$ even.
(Obviously, there are $(k-2)$ such squares contained into $T_{2k}$.)
Denote
$x_{2k;t}^S (x,y)$ the truncation
of a polynomial  $x_{2k;t}$ on the monomials
$x^cy^d, x^cy^{d+1}, x^{c+1}y^d, x^{c+1}y^{d+1}$.
Let $\G$ be the face in the triangulation  of $T_{2k}$
given by coordinates $(c,d+1), (c+1,d)$.

(From the patchworking construction, it follows that  for $t \in ]0,t_{2k}[$
the truncation $x_{2k;t}^S$ is $\e$-sufficient for $x_{2k;t}$
in an $\e$-neighborhood $U(p)$ of $x_{2k;t}^S=0$
defined from  points of $\G^0$.)

Fix $t \in ]0,t_{2k-1}[$,
and let $a_t$ be a crossing of $\HH_{2k-1} \cup L$ contained in $U(p)$.

According to the patchworking construction,
there exists an homeomorphism
$\tilde{h} :\CCC \HH_{2k} \cap U(p) \to
 \{(x,y) \in (\CCC^*)^2 | x_{2k;t}^S =0 \}\cap U(p)$
such that $\tilde{h}(a_t)=(x_1,y_1)_t$
is a critical point of $x_{2k;t}^S(x,y)$.
Obviously, (see corollary \mrf{c:patch}),
$(x_1,y_1)_t$ is such that $x_1.y_1 <0$.

From an argumentation similar to the one above, keeping the
same notations, it follows that, up to modify the coefficients
$a_{c,d},a_{c,d+1},a_{c+1,d},a_{c+1,d+1}$ if necessary,
one can get a critical point $\tilde{h}(a_t)$
$(x_1,y_1)_t$, $x_1.y_1 < 0$
of the function
$\frac {l_t(x,y)} {k_t(x,y)}$
with positive critical value
$\le t< t_{2k-1}$.
Therefore,
for fixed
$t \in ]0,t_{2k-1}[$,
one can choose $C_{2k;t}$ and thus
$\tilde{X}_{2k;t}$,
in such a way that
$\tilde{h}^{-1}(x_1,-y_1)_t$ with $x_1.y_1<0$
is a singular point
of
$\tilde{X}_{2k;t} =0$,
and one-parameter curves
$\tilde{X}_{2k;t}$ are Harnack curves for $t$ sufficiently
small.
\vskip0.1in
Moreover, the singular situation, we shall
denote by {\bf ${\mathcal S}"$}, is as follows:
\vskip0.1in

{\it
{\bf ${\mathcal S}"$}:
 one branch of a inner oval of the part
$j_1(\RRR \HH_{2k-1}\bk {\mathcal D}_{\e})$
of $\RRR \HH_{2k}$ not isotopic to a part of
$\RRR \HH_{2k-2}$\\
(namely, in the patchworking scheme,
 one inner oval  contained in the subset
$\rho^{2k}(D_{2k-1,2k-2} \ti U_{\RRR}^2) \subset \RRR P^2$)
touches a branch of the non-empty outer oval.}
\vskip0.1in
Thus, according to the Proposition \mrf{p:crit2} and
the Petrovskii's theory,
it follows that the point $\tilde{h}^{-1}((x_1,-y_1)_t)$
is a critical point of positive critical value
of a Harnack polynomial $B_{2k}$.
(Obviously, $\tilde{h}^{-1}((x_1,-y_1)_t \in S^+_{2k} \bk S'_{2k}$.)
\vskip0.1in
From an argumentation similar to the previous one,
it is easy to deduce that there are no other points
which may be singular points of a polynomial
$\tilde{X}_{2k;t}$.
Otherwise, it would contradict the strict positiveness of the coefficients
of the vector $\vec{a}_{2k}$.
(Indeed, consider the squares $S$
 given by  vertices
$(c+1,d), (c,d), (c,d+1), (c+1,d+1)$ with $c+d =2k-2$ and $c=0$ or $d=0$.
Keeping the previous notations,
from a crossing $a_t$ of $\HH_{2k-1} \cup L$,
one can get
a critical point
$\tilde{h}(a_t)=(x_1,y_1)_t$
 of the function
$\frac {l_t(x,y)} {k_t(x,y)}$
 with negative critical value.)
\vskip0.1in

Hence, we have obtained the Lemma \mrf{l:even}.
\vskip0.1in

{\bf i.2)}\\
Let us  now construct a polynomial
$\tilde{X}_{2k;\tau}$, $\tau \in [t_{2k}, t_{2k-1}[$
with $2k-3$ singular points: $k-1$ points
described locally by the singular situation {\bf ${\mathcal S}'$},
and $k-2$ points described locally
by the singular situation {\bf ${\mathcal S}"$}.

Let us fix $\tau < t_{2k-1}$.
Crossings of the union $\HH_{2k-1} \cup L$ belong to the
line $L$.

Let $p_i=\tilde{h}^{-1}(x,-y)_{\tau}$
with $(x,y)_{\tau}$  defined up to homeomorphism
from  a  crossing of the curve $\HH_{2k-1} \cup L$
with polynomial $x_0.\tilde{X}_{2k-1;\vec{a}_{2k-1},\tau}$.
We shall prove that given the $2k-3$ points
$p_1$,...,$p_{2k-3}$
the polynomial $\tilde{X}_{2k;\tau}$ may be entirely defined.

According to the previous description
of points $p_1,...,p_{2k-3}$, we shall
assume that the $2k-3$ singular points
$\tilde{X}_{2k,\tau}$ belong to the line
$x'_0=0 $ with $x'_0= x_0+ \a.x_1+\b.x_2$.
Then,
consider the linear change of complex projective coordinates
mapping $(x_0:x_1:x_2)$ to $(x'_0:x_1:x_2)$.
Such transformation carries
$x_0.\tilde{X}_{2k-1,\tau} (x_0,x_1,x_2)+ C_{2k,\tau}$ to
$\tilde{X}'_{2k,\tau}=
x'_0. \tilde{X}'_{2k-1,\tau}(x'_0,x_1,x_2) +C'_{2k, \tau}$.
We shall  prove that one can
choose $C'_{2k}(x_1,x_2)=0$
a union of $m$ parallel lines
in such a way that
$C'_{2k}(x_1,x_2)=0$, $\tilde{X}'_{2k-1} (x'_0,x_1,x_2)=0$,
and $x'_0=0$ have
$2k-3$ common points.

In such a way, using the linear change of projective coordinates
mapping $(x'_0:x_1:x_2)$ to $(x_0:x_1:x_2)$,
we shall get the polynomial
$\tilde{X}_{2k,\tau}(x_0,x_1,x_2)
=x_0.\tilde{X}_{2k-1,\tau}(x_0,x_1,x_2) + C_{2k}(x_1,x_2)$
from  the polynomial
$\tilde{X}'_{2k,\tau}(x_0,x_1,x_2)
=x'_0.\tilde{X}'_{2k-1,\tau}(x'_0,x_1,x_2) + C'_{2k}(x_1,x_2)$
and in that way the vector $\vec{c}_{2k}$.
\vskip0.1in
Let us detail this construction.
\vskip0.1in
Assume singular points $p=(0:p_1:p_2)$ of  $\tilde{X'}_{2k}=0$
such that $\frac {\pr \tilde{X'}_{2k}} {\pr x'_0} (p) =0$
(i.e belong to the curve
$C'_{2k}=0$.)\\
Assume
$(\frac {\pr \tilde{X'}_{2k-1}} {\pr_{x'_0}}=0)
 \cap (x'_0=0)
 \cap (\tilde{X}'_{2k-1}=0)
 = \emptyset$,

Set
$$ D_{2k} (x)=
x_1. \frac {\pr \tilde{X}'_{2k}} {\pr x_1}
+ x_2. \frac {\pr \tilde{X}'_{2k}}  {\pr x_2} $$
$$ D_{2k-1}(x) = x_1.\frac {\pr^2 \tilde{X}'_{2k}}  {\pr x'_0 \pr x_1}
+ x_2.\frac {\pr^2 \tilde{X}'_{2k}} {\pr x'_0 \pr x_2 }$$
where $D_{2k-1} (0,x_1,x_2)$
is an homogeneous polynomial of degree $(2k-1)$ in the variables
$x_1,x_2$.
\vskip0.1in
Since
$(\frac {\pr \tilde{X'}_{2k-1}} {\pr_{x'_0}}=0)
 \cap (x'_0=0)
 \cap (\tilde{X'}_{2k-1}=0)= \emptyset$,
$\frac {\pr^2 \tilde{X'}_{2k}} {\pr^2 x'_0} (p) \not=0$
  for any singular point $p$ of $\tilde{X'}_{2k}$.

Thus, for $2k>2$,
$D_{2k}(x'_0,x_1,x_2)$ is of degree at least $1$ in
$x'_0$ and any point $p$ belongs to
$\frac {\pr D_{2k}} {\pr_{x'_0}}=D_{2k-1}=0$.

According to the Euler formula,
for any singular point $p$ of $\tilde{X}'_{2k}$
the following equalities are verified:

\begin{eqnarray}
\lb{e:(i)}
 D_{2k} (p)=
p_1. \frac {\pr \tilde{X}'_{2k}} {\pr x_1} (p)
+ p_2. \frac {\pr \tilde{X}'_{2k}}  {\pr x_2} (p) =0\\
\lb{e:(ii)}
D_{2k-1} (p) = p_1.\frac {\pr^2 \tilde{X}'_{2k}}  {\pr x'_0 \pr x_1}(p)
+ p_2.\frac {\pr^2 \tilde{X}'_{2k}} {\pr x'_0 \pr x_2 }(p) =0\\
\lb{e:(iii)}
\frac {\pr \tilde{X}'_{2k}} {\pr x_1} (p)=0\\
or\nonumber\\
\frac {\pr \tilde{X}'_{2k}} {\pr x_2} (p)=0\nonumber
\end{eqnarray}

It is easy to deduce from the Newton's binomial formula
(applied to $(x_0')^i=(x_0+ \a.x_1+\b.x_2)^i$.
that coefficients of
the polynomial $\tilde{X}'_{2k}= x'_0.\tilde{X}'_{2k-1}+ C'_{2k}$
which depend on the vector $\vec{c}_{2k}$
(and do not depend on $\vec{c}_{2k-1}$)
are coefficients of monomials $x_1^{2k-i}x_2^i$,
(namely coefficients of $C'_{2k}$)
and coefficients of  monomials  $x'_0.x_1^{2k-1-i}x_2^i$
and in that way coefficients of $D_{2k-1}(0,x_1,x_2)$.

Therefore, any  singular point of $\tilde{X}'_{2k}=0$
belongs to the intersection of $x'_0=0$
with $C'_{2k}=0$ and $\tilde{X}'_{2k-1}=0$.
Let us detail the construction of $\tilde{X}'_{2k}$.
\vskip0.1in
Since any singular point $p$ of $\tilde{X}'_{2k}$
belongs to
$O_1=\{ (0:x_1:x_2) \in \CCC P^2 | x_1 \not=0   \}$
or
$O_2=\{ (0:x_1:x_2) \in \CCC P^2  | x_2 \not= 0  \}$,
it is sufficient to consider singular points
of $\tilde{X}'_{2k}$ in the affine charts associated to $O_1$ and $O_2$.

It follows from the form of
$\tilde{X}'_{2k} (x_0,x_1,x_2)$
that points
$(0:0:1)$  and $(0:1:0)$ are not singular points of
the polynomial
$\tilde{X}'_{2k} (x_0,x_1,x_2)$.
Consequently, any singular point $p$ of $\tilde{X}'_{2k}$
is such that
$$ \frac {\pr \tilde{X}'_{2k}} {\pr x_1} (p)=0$$
or
$$ \frac {\pr \tilde{X}'_{2k}} {\pr x_2} (p)=0$$
\vskip0.1in
Up to constant factor,
the following equalities are verified:
\begin{equation}
\lb{e:x1}
\frac {\pr \tilde{X}'_{2k}} {\pr x_2} (0,x_1,x_2)
= x_2^{2k-1}-c_1x_2^{2k-2}x_1+ ...(-1)^{2k-1}c_{2k-1}x_1^{2k-1}
\end{equation}
where $c_i$, $i \in \{1,...,2k-1 \}$
is the elementary roots symmetric polynomial
of degree $i$
of $\frac {\pr \tilde{X}'_{2k}} {\pr x_2} (0,1,x_2)$

\begin{equation}
\lb{e:x2}
\frac {\pr \tilde{X}'_{2k}} {\pr x_1} (0,x_1,x_2)
= x_1^{2k-1} -s_1x_1^{2k-2}x_2 + ...(-1)^{2k-1}s_{2k-1}x_2^{2k-1}
\end{equation}

where $s_i$, $i \in \{1,...,2k-1 \}$
is the elementary roots symmetric polynomial
of degree $i$
of $\frac {\pr \tilde{X}'_{2k}} {\pr x_1} (0,x_1,1)$.
\vskip0.1in

Set $x_1=1$ (respectively, $x_2=1$)
in the equality (\ref{e:x1}) (respectively, (\ref{e:x2})).
Bringing together
the resulting  equalities, it follows that
$\tilde{X}'_{2k}(0,x_1,1)$
has at most $[\frac {2k-1} 2]$ singular points
(Indeed,
any singular point $(0,p_1,1)$
is of multiplicity $1$ of $D_{2k-1}(0,x_1,1)$
and of multiplicity at least $2$ of
$\tilde{X}'_{2k,\tau}(0,x_1,1)$).
Besides, $\tilde{X}'_{2k}(0:1:x_2)$
has at most $[\frac {2k-1} 2]$ singular points
(since
any singular point $(0,1,p_2)$
is of multiplicity $1$ of $D_{2k-1}(0,1,x_2)$
and  of multiplicity at least $2$ of
$\tilde{X}'_{2k,\tau}(0,1,x_2)$)).

Since $\tilde{X}'_{2k}=0$ is singular, it follows immediately
from the equalities (\ref{e:x1}), (\ref{e:x2}) and
 (\ref{e:(ii)}), (\ref{e:(iii)})
that at least one point
$p=(0:1:-1)=(0:-1:1)$
is a singular point of $\tilde{X}'_{2k}$.

Therefore, the construction of the
polynomial $\tilde{X}'_{2k;\tau}(x_0,x_1,x_2)$
with $2k-3$ singular points
$p_1$,...,$p_{2k-3}$ on the intersection of $x'_0=0$
with $C'_{2k}=0$ and $\tilde{X}'_{2k-1}=0$ follows immediately.

Moreover, according to the Lemma \mrf{l:Lemma A},
one can perturb coefficients
of one monomials $x_1^{2k}$ or $x_2^{2k}$
of the polynomial $\tilde{X}'_{2k}$
in such a way that the modified polynomial
 is a Harnack polynomial
of degree $2k$ with
 critical points $p_1$,..,$p_{2k-3}$.
\vskip0.1in
This concludes the part (2) of proposition \mrf{p:prop7}.

{\bf ii)} Consider the case $m=2k+1$.
\vskip0.1in
From an  argumentation similar to the previous one,
it is easy to deduce
that any polynomial $\tilde{X}_{2k+1;t}$, $t \in ]0,t_{2k}[$
which satisfies good truncation properties
is a Harnack polynomial.
Besides, it  is  also obvious that
the existence of a singular polynomial
$\tilde{X}_{2k+1;t}(x_0,x_1,x_2).x_0+ C_{2k+1;t}(x_1,x_2)$
with good truncation properties
is not compatible with the Harnack's distribution of signs
(i.e the strict positiveness of the coefficients
of the vector $\vec{a}_{2k+1}$.)

Indeed, consider the squares $S$ given by vertices\\
$(c+1,d), (c,d), (c,d+1), (c+1,d+1)$
with $c+d =2k-1$, $c \ge 0$ and $d \ge 0$.
Keeping the notations introduced in the case $m=2k$,
from a crossing $a_t$ of $\HH_{2k} \cup L$,
one can get a critical point $\tilde{h}(a_t)=(x_1,y_1)_t$
of the function $\frac {l_t(x,y)} {k_t(x,y)}$
with negative critical value.
Therefore,
from an argument similar to the one given in case $m=2k$,
it follows easily that any polynomial
$\tilde{X}_{2k+1;t}(x_0,x_1,x_2).x_0+ C_{2k+1;t}(x_1,x_2)$,
$t \in ]0,t_{2k}[$,
with good truncation properties is smooth.\\
Q.E.D
\vskip0.1in
\subsubsection{Maximal Deformation of any Harnack curves}

\vskip0.1in

We shall now
in  Proposition \mrf{p:prop7b}
describe  maximal simple deformation of  any
Harnack curve.\\
Let us start
with  Corollary \mrf{c:gpa} of Theorem \mrf{t:rigiso}.

\begin{coro}
\mlb{c:gpa}
{\it
 Let $\HH_m$ be a Harnack curve of degree $m$.
 Then, one can assume, without changing its topological properties,
that any Harnack  curve $\HH_m$ is obtained via the
 $T$-inductive construction of Harnack curves.}
\end{coro}
{\bf proof:}
It is straightforward consequence of the rigid isotopy Theorem
\mrf{t:rigiso}.
Q.E.D.

The generalization of Proposition \mrf{p:prop7}
to any Harnack curve
is a straightforward consequence of the Corollary \mrf{c:gpa}.\\

\vskip0.1in
\bepr
\mlb{p:prop7b}
{\it
Denote $\HH_m$ the Harnack curve of degree $m$ defined up to isotopy of
real points set.
\been
\item
Along any simple deformation of a Harnack curve of odd degree $2k+1$
curves are smooth.
\item
Given $\HH_{2k}$ a Harnack curve of even degree $2k$,
any maximal simple deformation of $\HH_{2k}$
\begin{eqnarray}
& [0,1] \to \CCC {\mathcal C}_m &\nonumber\\
& t \to \HH_{m;t} &\nonumber
\end{eqnarray}
is such that $\HH_{2k;1}$
has $2k-3$ crossings.\\
Let $S$ be the set of singular points of
the singular  curve $\HH_{2k;1}$.
Then, there exists a Harnack polynomial $B_{2k}$ of degree $2k$
and type $\HH^0$
of which the set of critical points contains the set $S$.\\
These points are critical points with positive critical value.
Besides, $k-1$ of these points are points of $S'_{2k}$.
\enen}
\enpr
\vskip0.1in
{\bf proof:}
The generalization of Proposition \mrf{p:prop7}
is a straightforward consequence of the Corollary \mrf{c:gpa}.\\
We propose here a generalization of the argumentation
given in the proof of proposition \mrf{p:prop7b}.
\vskip0.1in
Consider $B_{m}(x_0,x_1,x_2)$ a Harnack polynomial of degree $m$.

According to Proposition \mrf{p:Lemma A'}, it is possible to
choose projective coordinates such that
$B_{m}(x_0,x_1,x_2)=
 x_0.B_{m-1}(x_0,x_1,x_2)+C_{m}(x_1,x_2)$
 where
$B_{m-1}(x_0,x_1,x_2)$ is the Harnack polynomial of degree $m-1$.

According to Morse Lemma \cite{Mil},
around any non-degenerate singular point $s$ of $B_{m-1}.x_0$
one can choose a local coordinates system $y_1,y_2$
with $z_1(s)=0,z_2(s)=0$  and
$x_0.B_{m-1}(x_0,x_1,x_2) = y_1.y_2$.

Analogously, in a neighborhood of any critical point $p$
of $B_{m}$,
one can choose a local coordinates system $y_1,y_2$
with $y_1(p)=0,y_2(p)=0$  and
$B_{m}(x_0,x_1,x_2) = B_{m}(p) +  y_1.y_2$.
Therefore, given a singular crossing $s$ of $x_0.B_{m-1}$,
since $s$ is critical point of the curve $B_{m}=0$
considered as level curve of the function $\frac {x_0.B_{m-1}} {C_m}$
the crossing $s$ is perturbed in such a way that around $s$,
one can choose local coordinates in an open $U(s)$ with
$y_1(s)\not=0,y_2(s)\not=0$  and
$B_{m}(x_0,x_1,x_2)=\frac{y_1} {y_2} + t$.\\
Such description enlarges to a description of
$B_{m}(x_0,x_1,x_2) = 0$
inside $U \supset U(s)$
$U=\{ z= <u,p>=(u_0.p_0:u_1.p_1:u_2.p_2) \in \CCC P^2 |
u=(u_0:u_1:u_2) \in U_{\CCC}^3, p=(p_0:p_1:p_2) \in U(s) \}$
Local coordinates $y_1,y_2$ inside the open $U(s)$ extend to local
coordinates inside $U$ as follows.
Given $z=<u,p> \in U$ with
$u=(1:u_1:u_2) \in U_{\CCC}^3$, $p=(p_0:p_1:p_2) \in U(s) \}$,
we may set
$y_1(z)=u_1.y_1(p)$, $y_2(z)=u_2.y_2(p)$.

In such a way, given a singular crossing $s$ of $x_0.B_{m-1}$,
any point $p \in U$ with local coordinates $y_1(p),y_2(p)$
such that $|y_1(p)|=|y_1(s)|$, $|y_2(p)|=|y_2(s)|$ and
$y_1(p).y_2(p) = -y_1(s).y_2(s)$, is a critical point of
$B_{m}(x_0,x_1,x_2)$.
\vskip0.1in
Moreover, according to Rolle's Theorem, around any singular point $s$ of
$x_0.B_{m-1}$, in a neighborhood of the line $x_0=0$,
the sign of $B_{m}(0,x_1/x_2,1)$ alternates.
(One can refer also to the marking method (see \cite{Gu})
where a distribution of sign of $C_{m}$ and
$x_0.B_{m-1}$ is defined
around any point of intersection of $C_{m}$ with $x_0.B_{m-1}$.)
Therefore, when consider the critical points of
$B_{m}$ deduced locally
from singular points of $x_0.B_{m-1}$,
the sign of $B_{m}(0,x_1/x_2,1)$
alternates as well around these points.
By continuity of $B_{m}$, it follows the alternation of
sign of $B_{m}(x_0,x_1,x_2)$
in an $\e$-tubular neighborhood
of the line $x_0=0$ in $\CCC P^2$.
It is obvious that such alternation of sign is equivalent to the
modular property of the distribution of
sign in the patchworking construction of Harnack's curves.
Hence, our proposition is a straightforward consequence of
the version of proposition \mrf{p:prop7b} for $T$-Harnack curves.
For sake of clarity, we refer to proposition \mrf{p:prop7}
of this chapter where a proof of the version of
proposition \mrf{p:prop7b} for $T$-Harnack curves
is given.

Hence, it follows from  proposition \mrf{p:prop7},
that along any simple deformation of a Harnack curve of odd degree $2k+1$
curves are smooth.
From an argumentation
similar to the one given in the second part of proof of
it follows that in case $m=2k$
one can modify coefficients of the polynomial  in such a way
that it has at most  $2k-3$ singular points.
\vskip0.1in
Q.E.D

\subsection{Description of $(\CCC P^2, \CCC \HH_m)$
up to conj-equivariant isotopy}
\mlb{susu:cji}
\vskip0.1in
Let us  recall
in introduction some  properties
of the real point set of Harnack curves which may be easily deduced from
the Harnack curves construction and, according to proposition \mrf{p:prop7b},
may be generalized to any Harnack curve

\vskip0.1in

Given $\AA_{m+1}=\HH_{m} \cup L$.
Let ${\mathcal D}_{\e}=\cup_{a_i~ \in \HH_{m}\cap L}  D(a_i,\e) \in \RRR P^2$
be a neighborhood of the set
of singular points of $\AA_{m+1}$ in $\RRR P^2$
(where $D(a_i,\e)$ denotes a 2-disc of radius $\e$ around $a_i$
in the Fubini-Study metric).

Let $N$ be a tubular neighborhood of
$\RRR \AA_{m+1} \bk {\mathcal D}_{\e}$
in $\RRR P^2 \bk {\mathcal D}_{\e}$.
Then, the Harnack curve $\HH_{m+1}$ is a non-singular curve  of degree $m+1$
such that:
$\RRR \HH_{m+1} \bk {\mathcal D}_{\e}$
is a section of the tubular fibration
$N \to \RRR \AA_{m+1} \bk {\mathcal D}_{\e}$.

Therefore,
there exists $\e >0$ and  $\tilde{j}_t$, with $t \in [0,1]$,
an isotopy  of $\RRR P^2$
which pushes $\RRR \HH_{m} \bk  {\mathcal D}_{\e_0}$ onto a subset
of $\RRR \HH_{m+1}$.

Besides, there is the following biunivoque correspondence
between critical points of
index $0$ and $2$ of Harnack polynomials $B_m$ and ovals of
curves:
\been
\item
ovals of
the part
$\tilde{j}_1(\RRR \HH_{m}\bk {\mathcal D}_{\e})$
of $\RRR \HH_{m+1}$ with ovals of $\HH_{m}$ and thus with
critical points of index $0$ and $2$ of $B_{m+1}$.
\item
$m$ ovals of $\RRR \HH_{m+1} \bk \tilde{j}_1(\RRR \HH_{m}\bk
 {\mathcal D}_{\e})$
with
critical points of index $0$ and $2$ of
$B_{m+1}$.
\enen
Any oval of these $m$ ovals intersects
a $2$-disc $D(a_i) \subset \RRR P^2$
defined around a point $a_i \in \HH_{m} \cap L$.

We shall in Proposition \mrf{p:prop8} define
the minimal
 (minimal in the sense of the inclusion
of subsets of $\CCC P^2$),
subset ${\mathcal B}_{\e}$
of $\CCC P^2$
which consists of  union of $4$-ball
(in the Fubini-Study metric) of $\CCC P^2$ of radius $\e$,
 ${\mathcal D}_{\e} \subset {\mathcal B}_{\e}$,
outside of which the isotopy $\tilde{j}_t$, $t \in [0,1]$,
extends to a conj-equivariant isotopy  $\CCC P^2$
which pushes $\CCC \HH_{m} \bk  {\mathcal B}$ onto a subset
of $\CCC \HH_{m+1}$.

In such a way, we shall
enlarge the above description of Harnack curves to the complex
domain and then in theorem \mrf{t:theo1} describe the pair
$(\CCC \HH_m,\CCC P^2)$ up to conj-equivariant isotopy of $\CCC P^2$.
\vskip0.1in

In way of preparation, let us make some remarks and recall
generalities connected
with Harnack curves and complex topological
characteristics of real curves.
\vskip0.1in
A curve $\AA_m$ of type $I$ is such its real set of points
$\RRR \AA_m$ divides its complex set of points into two connected
pieces called halves $\CCC \AA_m^+$, $\CCC \AA_m^-$.
These two halves are exchanged by complex conjugation
$ conj (\CCC \AA_m^+) = \CCC \AA_m^-$.
The natural orientations of these two halves determine
two opposite orientations on $\RRR \AA_m$  (as their common boundary),
called the {\it complex orientations} of $\AA_m$.\\
It is easy to see that
the deformation turning $\HH_{m} \cup L$,
into $\HH_{m+1}$ brings the complex orientations of $\HH_{m}$ and $L$ to
the orientations of the corresponding pieces of
$\RRR \HH_{m+1}$ induced by a single orientation
of the whole $\RRR \HH_{m+1}$.\\
In such a way,
given an orientation of $\RRR \HH_{m}$,
there exists only one orientation of $\RRR L$ which induces
the orientation of $\RRR \HH_{m+1}$.
Moreover, it is easy to deduce from the recursive construction of Harnack
curves that orientations of $\RRR L$ alternate.
Namely, given an orientation of $\RRR \HH_{m-1}$, $m>1$,
if $\CCC L^+$ is the half of $\CCC L$ which
induces the orientation on $\RRR L$
carried on $\RRR \HH_{m}$,
then $\CCC L^-=conj(\CCC L^+)$ induces the orientation
on $\RRR L$ carried on $\RRR \HH_{m+1}$.
\vskip0.1in

Let us recall
the different ways to perturb one crossing $p$ of a singular
curve $\AA$.
In a neighborhood invariant by complex conjugation of the crossing $p$
, the curve $\AA$ can be considered
up to conj-equivariant isotopy
as the intersection of two lines in the point $p$.
Thus, let $L_1$ and $L_2$ be lines embedded in $\CCC P^2$
which intersect each other in a single point $p$.
Denote by $C$ the result of the perturbation of the union $L_1 \cup L_2$.
From the complex viewpoint, there are essentially two ways
to perturb the singular union $L_1\cup L_2$.
On the other hand, there are two ways to connect the halves
of their complexifications.
Indeed, the curve $C$ is a curve of type $I$ since it is a non-empty conic.
The halves of $\CCC L_i$ $i \in \{1,2\}$ connected
each other after perturbation correspond to the complex orientation
of $\RRR L_i$ which agrees with some perturbation of $\RRR C$.
One distinguishes two ways to  connect the halves  with each other.
We shall define them as {\it perturbation of type 1 and type 2} of the
crossing (see figure $1.1$ and figure $1.2$).
Let $B$ be a complex $4$-ball globally invariant by complex conjugation
around the crossing.
After a {\it perturbation of type $1$}, $\RRR C \cap B$ does not divide
$\CCC C \cap B$ into two connected pieces.
After a perturbation of type $2$, $\RRR C \cap B$ divides
$\CCC C \cap B$ into two connected pieces.
\vskip0.1in
According to Theorem \mrf{t:bli}, one can
assume that Harnack polynomials are deduced by induction
on the degree as follows:\\
{\bf i)}
The Harnack polynomial $B_1$ is the polynomial of a real projective line.\\
Given
$B_{2k-1}$ a Harnack polynomial of degree $2k-1 \ge 3$.
Then, (up to linear change of projective coordinates
and up to slightly modify its coefficients),
$B_{2k-1}$ is a Harnack polynomial of degree $2k-1$,
and type $\HH^0$ deduced from classical small deformation
$B_{2k-1}=x_0.B_{2k-2} + \epsilon_{2k} C_{2k-1}$  of $x_0.B_{2k-2}$
where $B_{2k-2}$ is Harnack polynomial of degree $2k-2$ and type $\HH^0$.
We shall denote $\AA_{2k-1}= \HH_{2k-2} \cup L$ the curve given by
$x_0.B_{2k-2}$.\\
For $k \ge 1$, we shall denote
$A_{2k-1}$
 and
call {\it set of  points perturbed
in a maximal simple deformation of $\HH_{2k-1}$}
 the set of points
$A_{2k-1}=\{a_1,...,a_{2k} \}$
where $a_1,...a_{2k}$ are
the crossings of $\HH_{2k-2} \cup L$.\\

{\bf ii)}
The Harnack polynomial $B_2$ is the polynomial of a curve of degree
2 of which  real part consists of an oval.
We shall denote $A_{2}$ the set $A_{2}=\{ a_1 \}$
where $a_1$ is the crossing of $\HH_{1}\cup L$.\\

Given
$B_{2k}$ a Harnack polynomial of degree $2k \le 4 $.
Then, (up to linear change of projective coordinates
and up to slightly modify its coefficients),
$B_{2k}$ is a Harnack polynomial of degree $2k$ and type $\HH^0$
deduced from classical small deformation
$B_{2k}=x_0.B_{2k-1} + C_{2k}$
 (with $||C_{2k}||$ arbitrarily  small)
 of $B_{2k-1}.x_0$,
where $B_{2k-1}$ is a Harnack polynomial of degree $2k-1$ and type $\HH^0$.
We shall denote $\AA_{2k}= \HH_{2k-1} \cup L$ the curve given by
$x_0.B_{2k-1}$.\\

Let
$\HH_{2k;t}$, $t \in [0,1]$, be a maximal simple deformation of
the Harnack curve $\HH_{2k}$ given by the polynomial $B_{2k}$.
For $k \ge 2$,
we shall denote $A_{2k}$
and call {\it set of  points perturbed
in a maximal simple deformation of $\HH_{2k}$}
the ordered set of points
$A_{2k}=\{a_1,...,a_{2k-1},...,a_{4k-2} \}$
where $a_1,...a_{2k-1}$ are
the crossings of $\HH_{2k-1}\cup L$
and $a_{2k},...,a_{4k-2}$ are the crossings
of $\HH_{2k;1}$.
\vskip0.1in

One can easily notice that
any point of $A_m$ is perturbed  by
a perturbation of type 1 to give $\HH_m$.
(i.e locally around any crossing of
$\HH_{m} \cup L$,
the deformation turning  $\HH_{m} \cup L$ into $\HH_{m+1}$ is
a perturbation of type 1 of the crossing.
In case of even $m=2k$, given a maximal simple deformation
$\HH_{2k,t}$, $0 \le t \le 1$, of $\HH_{2k}$;
locally around  any  crossing $a$ of $\HH_{2k;1}$,
the deformation turning $\HH_{2k;1}=0$ into $\HH_{2k}=0$ is
a perturbation of type 1 of $a$.)

\vskip0.1in
Let $\CCC\HH_{m}^+$ (resp, $\CCC L^+$)
be the half of $\CCC \HH_{m}$ (resp, $\CCC L$)
which induces orientation on $\RRR \HH_{m}$
(resp,$ \RRR L$).
Denote $\CCC {\mathcal  A}_{m}^+$ the union $\CCC\HH_{m}^+ \cup \CCC L^+$.

In proposition \mrf{p:prop8}, we construct a subset ${\mathcal B}_{\e}$ of
$\CCC P^2$, (minimal in the sense of the inclusion of subsets
of $\CCC P^2$), ${\mathcal D}_{\e} \subset {\mathcal B}_{\e}$,
with the property that outside ${\mathcal B}_{\e}$,
the isotopy $\tilde{j}_t$
which pushes $\RRR \HH_{m} \bk  {\mathcal D}_{\e}$ onto a subset
of $\RRR \HH_{m+1}$ extends to a
a conj-equivariant isotopy of $\CCC P^2$ which maps
$\CCC \HH_{m+1}^+ \bk {\mathcal  B}_{\e}$ to
$\CCC\HH_{m}^+ \cup \CCC L^+ \bk {\mathcal  B}_{\e}$.
The construction of ${\mathcal B}_{\e}$ is based
on the definition of a maximal simple deformation of $\HH_{m}$
and the way points of the set $A_m$
are perturbed to give $\HH_m$.

\bepr
\mlb{p:prop8}
{\it
\been
\item
\vskip0.1in
Let ${\mathcal  B}_{\e}=
 \cup_{a \in A_{2k+1}} B_(a,\e) \subset \CCC P^2$, the union of
conj-equivarant 4-ball $B_(a,\e)$ of radius
$\epsilon$ taken over
crossings $a$ of $\AA_{2k+1}=\HH_{2k} \cup L$,
be a neighborhood of the set of singular points of
 $\AA_{2k+1}=\HH_{2k} \cup L$  in $\CCC P^2$.
Let $N$ be an oriented tubular neighborhood of
$\CCC \AA_{2k+1}^+ \bk {\mathcal  B}_{\e}$
in $\CCC P^2 \bk {\mathcal  B}_{\e}$.\\
Then, there exists  $\e_0$ such that,
$\CCC \HH_{2k+1}^+ \bk {\mathcal  B}_{\e_0}$
is the image of a non-zero section of the oriented tubular fibration
$N \to \CCC \AA_{2k+1}^+ \bk {\mathcal  B}_{\e_0}$ and
there exists $j_t$ with $t \in [0,1]$
a conj-equivariant isotopy of $\CCC P^2$ which maps
$\CCC \HH_{2k+1}^+ \bk {\mathcal  B}_{\e_0}$ to
$\CCC\HH_{2k}^+ \cup \CCC L^+ \bk {\mathcal B}_{\e_0}$.
\item
\vskip0.1in
Let ${\mathcal  B}_{\e}= \cup_{a \in A_{2k}}  B(a,\e) \subset \CCC P^2$
the union of conj-equivarant 4-ball $B(a,\e)$ of radius
$\epsilon$ taken over
points $a$ of  $A_{2k}$,be a neighborhood
in $\CCC P^2$ of the set
of  points perturbed
in a maximal simple deformation of $\HH_{2k}$.
(${\mathcal  B}_{\e}$  contains
besides neighborhoods of singular point of
$\AA_{2k}=\HH_{2k-1} \cup L$,
the union of neighborhood of $2k-3$ crossings
which appear in a maximal simple deformation of $\HH_{2k}$.)
Let $N$ be an oriented tubular neighborhood of
$\CCC \AA_{2k}^+ \bk {\mathcal  B}_{\e} $
in $\CCC P^2 \bk {\mathcal  B}_{\e}$.\\
Then, there exists  $\e_0$ such that
$\CCC \HH_{2k}^+ \bk {\mathcal  B}_{\e_0}$
is the image of a non-zero section of the oriented tubular fibration
$N \to \CCC \AA_{2k}^+ \bk  {\mathcal  B}_{\e_0}$
and there exists $j_t$, with $t \in [0,1]$,
a conj-equivariant isotopy of $\CCC P^2$ which maps
$\CCC \HH_{2k}^+ \bk {\mathcal  B}_{\e_0}$  to
$\CCC\HH_{2k-1}^+ \cup \CCC L^+ \bk {\mathcal  B}_{\e_0}$.
\enen}
\enpr

\vskip0.1in
\vskip0.1in
{\bf proof:}
\vskip0.1in
Our proof is based on
results of proposition \mrf{p:prop7}  and the fact
that orientation of $\RRR L$ alternates in
the recursive construction of Harnack curves.
\vskip0.1in
Let us explain briefly the method of our proof.
\vskip0.1in
We shall consider the usual handlebody decomposition of
$\CCC P^2= B_0 \cup B_1 \cup B$  where $B_0$,$B_1$,$B$
are respectively 0,2 and 4 handles.
In such a way,
the canonical $\RRR P^2$ can be seen as the union of a M\"{o}bius
band $\mathcal  M$ and the disc $D^2 \subset B$ glued along their boundary.
The M\"{o}bius band $\mathcal  M$ lies in $B_0 \cap B_1 \approx S^1 \ti D^2$.
The complex conjugation switches $B_0$ and $B_1$ and
$\mathcal  M$ belongs to the set of its fixed points.

Let ${\mathcal D}_{\e}$ be  a neighborhood of the set
of singular points of $\AA_{m+1}$ in $\CCC P^2$
which is the union of
$2$-disc of radius $\epsilon$ around point of the set $\AA_m$.

We shall assume that there exists an isotopy $\tilde{j}_t$
of $\RRR P^2$ which maps $\RRR \HH_{m+1} \bk {\mathcal D}$
on the boundary of the M\"{o}bius band and
study whether or not $\tilde{j}_t$ extends to $j_t$ a
conj-equivariant isotopy of $\CCC P^2$
which maps
$\CCC \HH_{m+1}^+ \bk {\mathcal D}_{\e}$ to
$\CCC\HH_{m}^+ \cup \CCC L^+ \bk {\mathcal D}_{\e}$.

We shall notice that given a $4$-ball
$B(p) \subset (B_0 \cup B_1) \subset \CCC P^2$
globally invariant by complex conjugation and such that
$B(p)^-=B(p) \cap B_1$ and $B(p)^+=B(p) \cap B_0$,
if
$\CCC \HH_{m+1}^+ \cap B(p) \not = \emptyset$
and
$\CCC \HH_{m+1}^+ \cap B(p) \not =
\CCC \HH_{m+1}^+  \cap B(p)^{\pm}$ then
there is no conj-equivariant isotopy
$j_t$ of $\CCC P^2$ which maps
$\CCC \HH_{m+1}^+ \bk {\mathcal D}_{\e}  \cap B(p)$
to
$((\CCC \HH_{m}^+ \cup \CCC L^+) \bk {\mathcal D}_{\e}) \cap B(p)$.

Thereby, since orientation of $\RRR L$ alternates in the recursive
construction of Harnack curves
an  obstruction to construct a conj-equivariant isotopy in a $4$-ball
$B(p)$ will be provided by a crossing $a \in B(p)$ of a curve
$\HH_{m+1;1}$
where $\HH_{m+1;t}$ $t \in [0,1]$ is a simple deformation of the
Harnack curve $\HH_{m+1}$.
\vskip0.1in
We shall now proceed to precise arguments.\\
Assume $m+1=2k \ge 4$.
\vskip0.1in
First, we shall prove that the existence of a conj-equivariant
isotopy $j_t$ which maps
$\CCC \HH_{2k}^+ \bk {\mathcal D}_{\e}$  to
$\CCC\HH_{2k-1}^+ \cup \CCC L^+ \bk {\mathcal D}_{\e}$
is not compatible with
elementary topological properties of $\CCC \HH_{2k}$.
Then, using  proposition \mrf{p:prop7}, we shall define
a subset ${\mathcal  B}_{\e} \subset \CCC P^2$,
${\mathcal  B}_{\e} \supset {\mathcal D}_{\e}$,
such that there exists a conj-equivariant
isotopy $j_t$ which maps
$\CCC \HH_{2k}^+ \bk {\mathcal  B}_{\e}$  to
$\CCC\HH_{2k-1}^+ \cup \CCC L^+ \bk {\mathcal  B}_{\e}$.

Since $\HH_{2k}$ is of even degree,
its real point set $\RRR \HH_{2k}$
divides the projective plane $\RRR P^2$ in  two components
respectively orientable
and non-orientable whose common boundary is  $\RRR \HH_{2k}$.

Without loss of generality, one can assume that
$D^2 \subset \RRR P^2$
is a small real disc with boundary the line at infinity $L$ of
$\RRR P^2$ with
${\mathcal D}_{\e}
 \cap \RRR P^2
 \subset  D^2$.
The tubular neighborhood of the line $L$
is the M\"{o}bius band ${\mathcal  M}$ of $\RRR P^2$.
Consider $\tilde{j}_t$, with  $t \in [0,1]$
an isotopy of $\RRR P^2$ which
lets fix the interior of the disc $D^2$ of $\RRR P^2$ and
pushes the rest of $\RRR \HH_{2k} \bk {\mathcal D}_{\e}$
to the boundary of the
M\"{o}bius band. \\
Assume that the isotopy $\tilde{j}_t$ extends to a
conj-equivariant isotopy $j_t$ of $\CCC P^2$.
Then, one can set
$j_1(\CCC \HH_{2k}^+ \bk B) \subset B_0$.
Given $N$ a tubular neighborhood
of $j_1(\CCC \HH_{2k})$ in $\CCC P^2$,
$j_1(\CCC \HH_{2k})$ is the image of a non-zero section
of the oriented tubular fibration $N \to j_1(\CCC \HH_{2k})$.

Since $j_t$ is conj-equivariant,
there exists $N^+ \subset B_0$ a tubular
neighborhood of $\CCC \HH_{2k}^+ \bk B$ in
$\CCC P^2 \bk {\mathcal D}_{\e}$
such that \\ $N^+ \to \CCC \HH_{2k}^+ \bk B$
is an oriented tubular fibration.

Hence, since $\HH_{2k}$
is two-sided and has orientable real part,
there exists $N$ an oriented tubular
neighborhood of
$\tilde{j}_1(\RRR \HH_{2k} \bk {\mathcal D}_{\e})$
in $\RRR P^2 \bk {\mathcal D}_{\e}$
such that
$\tilde{j}_1(\RRR \HH_{2k} \bk {\mathcal D}_{\e}$
is the image of a non-zero smooth section
of the tubular fibration
$N \to \tilde{j}_1(\RRR \HH_{2k}\bk {\mathcal D}_{\e})$.

It leads to contradiction since
the  tubular neighborhood of
$j_1(\RRR \HH_{2k} \bk (D^2)^0)$ in
$\RRR P^2 \bk (D^2)^0$
is homeomorphic to the M\"{o}bius band and  therefore
has non-orientable normal fibration.

Therefore, on $\RRR \HH_{2k} \bk  {\mathcal D}_{\e}$,
the isotopy  $\tilde{j}_t$ does not extend to a
conj-equivariant isotopy $j_t$ of $\CCC P^2$.

Nonetheless,
using the proposition \mrf{p:prop7b},
we shall define a set
${\mathcal  B}_{\e} \supset {\mathcal D}_{\e}$  for which
on  $\RRR \HH_{2k} \bk {\mathcal  B}$
the isotopy  $\tilde{j_t}$ extends to a conj-equivariant isotopy
$j_t$ which maps
$\CCC \HH_{2k}^+ \bk {\mathcal  B}_{\e}$  to
$\CCC\HH_{2k-1}^+ \cup \CCC L^+ \bk {\mathcal B}_{\e}$.
\vskip0.1in

Given $\HH_{2k}$  the Harnack curve of degree $2k$,
consider $\HH_{2k;t}$, $t \in [0,1]$, a simple deformation of the Harnack
curve $\HH_{2k}$.
Let us prove the following Lemma.
\vskip0.1in
\bele
\mlb{l:Lemma B}
\vskip0.1in
{\it
There exists
$B(p) \subset B_0 \cup B_1$
such that \\
$\CCC \HH_{2k}^+ \cap B(p) \not =
\CCC \HH_{2k}^+  \cap B(p)^{\pm}$
if and only if
there exists a simple deformation
 $\HH_{2k;t}$ $t \in [0,1]$ of the Harnack
curve $\HH_{2k}$ such that a crossing of $\HH_{2k;1}$ belongs to $B(p)$.}
\enle
\vskip0.1in
{\bf proof:}
\vskip0.1in

For $\e > 0$ small,
the sets of complex points of curves
$\HH_{2k;t}$
 $t \in [1 -\e, 1]$
belong to a tubular neighborhood of the
set of complex points
of the Harnack curves $\HH_{2k;1}$.

Let $B(p) \subset \CCC P^2$ be a $4$-ball
centered in a singular point of $\HH_{2k;1}$
such that
$\CCC \HH_{2k;t~,
t \in [1 -\e, 1]} \cap B(p) \not=
\emptyset$.
Inside $B(p)$, consider the gradient trajectories of
the deformation turning
$\HH_{2k;1}$ into $\HH_{2k;1-\e}$.
Up to isotopy  $\tilde{j}_t$ of $\RRR P^2$,
one can assume
$B(p) \subset B_0 \cup B_1$ and set
$B(p)^-=B(p) \cap B_1 $ and $B(p)^+=B(p) \cap B_0$.

Since orientation of the line $L$
alternates in the recursive construction of Harnack curves,
the complex conjugation acts on
gradient trajectories as symmetry of center $p$.
Therefore, inside $B(p)$ we have
$\CCC \HH_{2k}^+ \cap B(p) \not =
\CCC \HH_{2k}^+  \cap B(p)^{\pm}$.
(i.e around $B(p)$,
the deformation turning $\HH_{2k;1}$  into $\HH_{2k;1 - \e}$
is  a deformation of type 1.)\\

Reciprocally, it is easy to see that  if there exists
$B(p) \subset B_0 \cup B_1$
such that \\
$\CCC \HH_{2k}^+ \cap B(p) \not =
\CCC \HH_{2k}^+  \cap B(p)^{\pm}$
then there exists a simple deformation
$\HH_{2k;t}$ $t \in [0,1]$ such that a crossing of $\HH_{2k;1}$
belongs to $B(p)$.\\
Q.E.D

\vskip0.1in
\vskip0.1in
The Lemma \mrf{l:Lemma B} implies that
outside  a neighborhood ${\mathcal  B}_{\e}$ of the set of singular points
of $\HH_{2k;1}$,
$\tilde{j}_t$ with $t \in [0,1]$
extends to  a conj-equivariant isotopy of $\CCC P^2$
which maps
$\CCC \HH_{2k}^+ \bk {\mathcal  B}_{\e}$  to
$\CCC\HH_{2k-1}^+ \cup \CCC L^+ \bk {\mathcal  B}_{\e}$.

Assume now $m+1=2k+1 \ge 3$,
\vskip0.1in
Without loss of generality one can assume that
$D^2$ is a small real disc  such that
$D^2 \cap \HH_{2k+1} \cup L= {\mathcal D}$.
Consider $\tilde{j}_t$, with  $t \in [0,1]$,
an isotopy of $\RRR P^2$ which
lets fix the disc $D^2$ of $\RRR P^2$ and
pushes $\RRR \HH_{2k+1} \bk {\mathcal D}_{\e}$ to the boundary of the
M\"{o}bius band.
As any curve of odd degree, $\HH_{2k+1}$ is one-sided.
Therefore, we may not use the argument
given in case of curve of even degree to refute the  existence of a
conj-equivariant isotopy $j_t$ of $\CCC P^2$
which maps
$\CCC \HH_{2k+1}^+ \bk {\mathcal D}_{\e}$ to
$\CCC\HH_{2k}^+ \cup \CCC L^+ \bk {\mathcal D}_{\e}$.

Besides, from an argument analogous to the one given in even degree,
using the odd version of the previous Lemma \mrf{l:Lemma B},
it follows that on $\RRR \HH_{2k} \bk  {\mathcal D}_{\e}$,
$\tilde{j}_t$ extends to a conj-equivariant isotopy $j_t$
 of $\CCC P^2$.

Otherwise, there would exist a simple deformation $\HH_{2k+1;t}$
$t \in [0, 1]$ of the  Harnack curve
 $\HH_{2k+1}$ of degree $2k+1$
which intersects the discriminant hypersurface
$\RRR {\mathcal D}_{2k+1}$.
According to proposition \mrf{p:prop7b},
it is impossible.
Q.E.D
\vskip0.1in

The following theorem is a straighforward consequence
of the Proposition \mrf{p:prop8}.

\beth
\mlb{t:theo1}
{\it
\vskip0.1in
Let $\HH_m$ be the Harnack curve of degree $m$.
There exists a finite number $I$
($I=1+2 ...m +  \Sigma_{k=2}^{k=[m/2]} 2k-3$)
of disjoint $4$-balls $B(a_i)$
invariant by complex conjugation centered in points $a_i$
of $\RRR P^2$ such that
up to conj-equivariant isotopy of $\CCC P^2$,
\been
\item $\HH_m \bk \cup_{i \in I} B(a_i)= \cup_{i=1}^m L_i \bk
  \cup_{i=1}^I B(a_i) $
 where $L_1,..., L_m$ are $m$ distinct projective lines
 with
 $$L_i \bk \cup_{i =1}^I B(a_i) \cap
  L_j \bk \cup_{i=1}^I  B(a_i) = \emptyset$$
  for any $i \not=j$, $1\le i, j \le m$.
\item
situation inside any $4$-ball $B(a_i)$
is described by the perturbation of type $1$ of the crossing $a_i$.
\enen
\vskip0.1in}
\enth

{\bf proof:}

For $m=1$, the theorem is trivially verified:
$\HH_1$ is a projective line.
For $m > 1$, it may be deduced by induction on $m$.
The induction is based on
the proposition \mrf{p:prop8} and
the inductive construction of Harnack curves.

As related above,
for each integer $m$ ($ m \ge 0$),
on can assume without loss of generality,
that
the curve $\HH_{m+1}$ results
from classical deformation  of
$\HH_m \cup L$.

According to the proposition \mrf{p:prop7b},
\been
\item
let $A_{2k}=\{a_1,...,a_{2k-1},...,a_{4k-2} \}$
be the set of points perturbed
in a maximal simple deformation of $\HH_{2k}$
\item
denote $A_{2k+1}=\{a_1,...,a_{2k} \}$
 the set of points perturbed
 in a maximal simple deformation of $\HH_{2k+1}$
\enen

\vskip0.1in
For any $a_i \in A_{m+1}$,
we may always choose a small open conj-symmetric $4$-ball $B(a_i)$ of
$a_i$, in such a way that
$B(a_i) \subset \rho^{m+1}(D(p_j)
\ti U_{\CCC}^2)$ with
$D(p_j)$ is  a small open real disc around
$p_j \in \G_j$.
\vskip0.1in
\been
\item
Given $\HH_{2k-1}$
the Harnack curve of degree $2k-1$,
from the proposition \mrf{p:prop8}, it results
the following characterization of $\HH_{2k}$:

There exists $j_t$, $t \in [0,1]$,
a conj-equivariant isotopy  of $\CCC P^2$,
which maps
$\HH_{2k}  \bk \cup_{i=1}^{4k-2} B(a_i)$ onto
\been
\item
 $\HH_{2k-1}  \bk
  \cup_{i=1}^{4k-2} B(a_i)$
(in the patchworking scheme, such part is contained in the restriction
$\rho ^{2k}(T_{2k-1} \ti U_{\CCC}^2)$ of
$\CCC P^2 \approx \CCC T_{2k}$.)
\item
union $L \bk \cup_{i=1}^{4k-2} B(a_i)$
(in the patchworking scheme, such part  is contained in the restriction
$\rho ^{2k}(D_{2k,2k-1} \ti U_{\CCC}^2 )$ of
$\CCC P^2 \approx \CCC T_{2k}$.)
\enen

\item

Given $\HH_{2k}$
the Harnack curve of degree $2k$.
From the proposition \mrf{p:prop8}, it results
the following characterization of $\HH_{2k+1}$:

There exists $j_t$, $t \in [0,1]$,
a conj-equivariant isotopy  of $\CCC P^2$,
which maps $\HH_{2k+1}
\bk  \cup_{i=1}^{2k} B(a_i)$ onto
\been
\item
$\HH_{2k} \bk \cup_{i=1}^{2k} B(a_i)$
(in the patchworking scheme, such part is contained in the restriction
$\rho ^{2k+1}(T_{2k} \ti U_{\CCC}^2)$ of
$\CCC P^2 \approx \CCC T_{2k+1}$.)
\item
union $L \bk  \cup_{i=1}^{2k} B(a_i)$
(in the patchworking scheme, such part  is contained in the restriction
$\rho ^{2k+1}(D_{2k+1,2k} \ti U_{\CCC}^2 )$ of
$\CCC P^2 \approx \CCC T_{2k+1}$.)
\enen

\enen

Therefore, it follows easily by induction on $m$
that outside a finite number of $4$-balls $B(a_i)$ the curve
$\HH_m$ is, up to conj-equivariant isotopy,
the union of $m$ (non-intersecting) projective lines
minus their intersections with $4$-balls $B(a_i)$.

\vskip0.1in

The $4$-balls $B(a_i)$ are centered in points $a_i$.
Inside any $4$-ball $B(a_i)$, one can get
the whole description of $\HH_{m+1}$
from the real set of points
$\RRR \HH_{m+1} \cap B(a_i)$ and its orientation
(i.e the type of deformation).

Situations inside $4$-balls centered in crossings  $a_i$
of  $\HH_m \cup L$ are easily deduced from the construction of
$\HH_{m+1}$:
the deformation turning $(\HH_{2k} \cup L) \cap B(a_i)$
into $\HH_{2k+1} \cap B(a_i)$ is of type $1$.

In case $m+1=2k$, one has  to consider
$4$-balls
centered in crossings $B_{2k;1}$.
The situation inside $4$-balls
centered in crossings $B_{2k;1}$ has been explicitly
described inside $D(a_i)$, and respectively inside $B(a_i)$,
in the proof of proposition \mrf{p:prop7},
and respectively in the proof of proposition \mrf{p:prop8}.
The deformation turning $\HH_{2k;1} \cap B(a_i)$
into $\HH_{2k} \cap B(a_i)$ is of type $1$.

Q.E.D

\chapter{ Arnold surfaces of Harnack curves}
\mlb{ch:ArnHar}
In this chapter, we valid in Theorem \mrf{t:theo2} the Rokhlin Conjecture
for the Harnack curves of even degree.
In other words, we prove in Theorem \mrf{t:theo2} that
Arnold surfaces of
Harnack curves  split into connected sum of
$\RRR P^2$ and $\overline{\RRR P^2}$.
This result is a straightforward consequence
of Theorem \mrf{t:theo1}.\\

\beth
\mlb{t:theo2}
{\it Arnold surfaces of Harnack curves of even degree are standard.
   Moreover,
    let $\HH_{2k}$ be the Harnack curve  of degree $2k$
    and ${\gA}_+$, ${\gA}_-$ its Arnold surfaces
    then, up to isotopy of $S^4$:
    for $k=1$
    $${\gA} _+ \approx D^2$$
    $${\gA} _- \approx \overline{\RRR P^2}$$
    for $k \ge 2$
    $${\gA} _+ \approx p \RRR P^2  \krest q_+\overline{\RRR P^2}$$
    $${\gA} _- \approx p \RRR P^2  \krest q_-\overline{\RRR P^2}$$

    where

$p=q_+= \frac {(k-1)(k-2)}{2}; q_-= \frac {k(5k-3)}{2} $ }

\enth
\vskip0.1in

{\bf proof:}

The proof,
is based on the theorem \mrf{t:theo1} stated previously,
and on the next Livingston's statement:
\vskip0.1in
{\it
{\bf Theorem \cite{Li} }: Let $ F \subset D^4 \subset S^4$ be a closed
surface which lies in $S^3 =\pr D^4$, except several two-dimensional discs
standardly embedded inside $D^4$.
Then $F$ is a standard surface in $S^4$.}
\vskip0.1in
Let us recall in details the handles decomposition of $\CCC P^2$
and \\ $\CCC P^2 /conj \approx S^4 $.
\vskip0.1in
Choose $b$ a point in $\RRR P^2$ and
$B$ a small neighborhood of $b$ in $\CCC P^2$.
Setting $b=(b_0:b_1:1) \in \RRR P^2$,
$B$ is defined, up to homeomorphism, in
$U_2 =\{ (x_0:x_1:x_2) \in \CCC P^2 | x_2 \not = 0 \} \subset \CCC P^2$
$U_2 \approx \CCC ^2$ as a usual $4$-ball of $\CCC^2$:
$B \approx \{ p =(p_0:p_1:1) \in \CCC ^2
 \mid \vert \vert
 (p_0-b_0, p_1-b_1) \vert \vert \le \e \}$.

The $4$-ball $B$ is globally invariant
by complex conjugation.
The real set of points of the ball $B$
is a two-dimensional disc $D^2$ which remains
fix by complex conjugation.
\vskip0.1in
Then,
consider the central projection $$p:\CCC P^2 \bk B \to \CCC L $$
from $b$ to some set of complex points of a real projective line $L$,
$b \not \in \CCC L$.
Denote $\CCC L_0$ and $\CCC L_1$ the two connected parts of
$\CCC L$ such that :
\been
\item
$\CCC L =\CCC L_0 \cup  \CCC L_1$
\item
 $\CCC L_0 \cap \CCC L_1 = \RRR L$
\enen
Since
  $\CCC L_0 \approx D^2 $,$ \CCC L_1 \approx D^2 $, $\RRR L \approx S^1$,
the fibrations:
 $$p^{-1} (\CCC L_0) \to
  \CCC L_0 \approx D^2 $$
 $$p^{-1} (\CCC L_1)
  \to \CCC L_1 \approx D^2 $$
 $$p^{-1}(\RRR L) \to
  \RRR L \approx S^1$$ are, up to homeomorphism, trivial fibrations.

Therefore, setting
$p^{-1} (\CCC L_0)=B_0$ and $p^{-1} (\CCC L_1)=B_1$
it follows :
$$B_0\approx D^2 \ti D^2,B_1 \approx D^2 \ti D^2,
B_0 \cap B_1 = p^{-1} (\RRR L) \approx S^1 \ti D^2$$

Thus, we shall
take the handlebody decomposition of $\CCC P^2$
$\CCC P^2= B_0 \cup B_1 \cup B$  where $B_0$,$B_1$,$B$
are respectively 0,2 and 4 handles.
\vskip0.1in
Moreover,
from the action of complex conjugation on \\
$\CCC P^2 =B_0 \cup B_1 \cup B$
it is easy to deduce
a decomposition of $\CCC P^2/ conj \approx S^4$.
\vskip0.1in
The canonical $\RRR P^2$ can be seen as the union of a M\"{o}bius
band $\mathcal M$ and the disc $D^2$ glued along their boundary.
In such a way, $\mathcal M$ lies in $B_0 \cap B_1 \approx S^1 \ti D^2$.
\vskip0.1in
The complex conjugation switches $B_0$ and $B_1$ and
$\mathcal M$ belongs to the set of its fixed points.
The quotient $(B_0 \cup B_1)/conj$
is a $4$-ball as well which
contains the quotient $(B_0 \cap B_1)/conj$.

\vskip0.1in

Consider now
$\HH_{2k}$ the Harnack curve of degree $2k$
and ${\gA}_+$
and ${\gA}_-$ its Arnold surfaces.
\vskip0.1in

Acording to  Theorem \mrf{t:theo1},
there  exists
a conj-equivariant isotopy of $\CCC P^2$,
$h_t$, $t \in [0,1]$,
$h(0)=\HH_{2k}$
and a finite number $I$ of disjoint $4$-balls $B(a_i)$
invariant by complex conjugation centered in points $a_i$
of $\RRR P^2$ such that:
\been
\item
$h_1(\HH_{2k} \bk \cup_{i =1}^I B(a_i))= \cup_{i=1}^{2k} L_i \bk
 \cup_{i=1}^I B(a_i) $
 where $L_1,..., L_{2k}$ are $2k$ distinct projective lines
 with
$$(L_i \bk \cup_{i =1}^I B(a_i)) \cap
  (L_j \bk \cup_{i=1}^I  B(a_i)) = \emptyset$$
 for any $i,j$, $1\le i \not =j \le 2k$,
 \item
inside any $4$-ball $B(a_i)$,
$h_1(\HH_{m} \cap B(a_i))$ is the perturbation
of a crossing as described in theorem \mrf{t:theo1}.

\enen

\vskip0.1in
Choose
the $4$-ball $B$
of the handlebody decomposition of $\CCC P^2$
centered in $b$ in such a way that:
\been
\item
 $b \in \RRR P^2_{\pm}$
(For example, choose $b=(0,0,1)$ inside $\RRR P^2_-$
and\\ $b=(\e,\e, 1)$, $0<\e<1$, inside $\RRR P^2_+$).
\item
$B$ does not intersect with
$\cup_{ i \in I} B(a_i)/conj $.
and lines $L_1,...,L_{2k}$ provided by the isotopy $h_t$.
\enen

Then, delete a small tubular neighborhood $U$ of
$(B_0 \cap B_1)/conj$ from $(B_0 \cup B_1) /conj$
and denote by $B'$ the resulting 4-ball.

From the decomposition of $\CCC P^2 /conj \approx S^4$
and the conj-equivariant isotopy $h_t$,
it follows an isotopy $h_t/conj $ of $S^4$
which pushes \quad \break
 $\RRR P^2_{\pm} \bk
 \cup_{i=1}^I B(a_i)/conj$
 on the boundary of $B'$
 and $(\HH_{2k}/conj \bk \cup_{i=1}^I B(a_i)/conj)$
 to $2k$ standards discs of $B'$.
(such discs are defined by \quad \break
quotient of lines $L_1/conj$, ...,$L_{2k}/conj$ minus their
intersection with $\cup_{ i \in I} B(a_i)/conj$.)
\vskip0.1in

Thus, since
${\gA}_{\pm} \bk
 \cup_{ i \in I} B(a_i)/conj  \cap Int U =\emptyset$,
the problem of construction of an isotopy which pushes
${\gA}_{\pm}$
to the boundary of a $4$-disc and to several $2$-discs standardly embedded
is reduced to the problem
of construction of such an isotopy inside
$ \cup_{ i \in I} B(a_i)/conj$.

Moreover, since the $4$-balls $B(a_i)$ are disjoints, it is sufficient
to study local questions inside $2$-discs $B(a_i)/conj$.

Inside any $4$-ball $B(a_i)$, $h_1(\HH_{2k} \cap B(a_i))$
is described by  the perturbation
of the singular crossing $a_i$.

Thus, inside any $2$-disc $B(a_i)/conj$,\\
$ h_1/conj (\HH_{2k}/conj \cup \RRR P^2_{\pm} \cap B(a_i)/conj)$
is determined by the relative positions
of the point $a_i$, the point $b$ and $\RRR P^2 _{\pm}$.
\vskip0.1in
Using an argumentation similar to the one
given in [\cite{Fi}, see p.6], in any case
it is always possible to push by an isotopy
the part ${\gA}_{\pm}$, contained inside
$B(a_i) / conj$  in  $\pr B(a_i) / conj$ leaving a $2$-disc inside
$B(a_i) / conj$ if necessary.

Hence, the part ${\gA}_{\pm}$, contained inside $U$,
can be pushed by an isotopy
(which
coincides with $h_t/conj$ inside $ S^4 \bk \cup_{i=1}^I B(a_i)/conj$
and inside any $B(a_i) / conj$ satisfies the requirement above)
into a $4$-ball $B"$  obtained by taking union of $B'$
with disc $B(a_i) / conj$ or excising $B(a_i) / conj$.
(The choice of union or excision of the disc $B(a_i) / conj$
depends on relative position of $b$ and $\RRR P^2_{\pm} \cap
B(a_i) / conj$
[see \cite{Fi}])
\vskip0.1in
 After such isotopy, ${\gA}_{\pm}$ lies in the boundary
of $B"$, except $2k$  discs left from
$\HH_{2k}/conj \bk \cup_{i=1}^I B(a_i) /conj$
and  several discs which lie all inside $B"$ and are unknotted.
\vskip0.1in
Thus from  Livingston's Theorem  it follows that :
\vskip0.1in

{\it Arnold surfaces
${\gA}_{\pm}$ of Harnack curve $\HH_{2k}$ are standard surfaces in $S^4$.}
\vskip0.1in
The decomposition announced in theorem \mrf{t:theo2}
is immediate when considers the two double coverings branched respectively over
 ${\gA}_{-}$ and ${\gA}_+$ (see \cite{ArS}),
However, we propose to recover this decomposition from
the study of local situations
inside $4$-balls $B(a_i)$ centered in  crossings points $a_i$
perturbed in maximal simple deformation of the Harnack curves.
\vskip0.2in
\bele
\mlb{l:lars1}
\vskip0.1in
{\it Let ${\mathcal C}_0$ be a singular curve of
degree $m$ with one singular  crossing.
Consider a variation ${\mathcal C}_t$, $t \in [-1,1]$  which crosses
transversally the discriminant $\D_m$ at ${\mathcal C}_0$
Let $p$ be the singular crossing of ${\mathcal C}_0$.
Let ${\gA}_t$ denote either the positive or negative
Arnold surface of ${\mathcal C}_t$.
Assume that inside a small neighborhood of $p$:
\been
\item
For $t \le 0$, ${\gA}_t$ is the union of ${\mathcal C}_t/conj$ and a disc
component.
\item
For $t>0$, ${\gA}_t$
is the union of ${\mathcal C}_t/conj$ and a M\"{o}bius
component.
Then the perturbation ${\mathcal C}_{0}$
to ${\mathcal C}_{1}$ corresponds
on Arnold surface to $\krest \overline{\RRR P^2}$
(The direction of twisting of the M\"{o}bius band
is  obviously standard.)
\enen}
\enle
{\bf proof:}
\vskip0.1in
It can be easily deduced from an algebraic model.
In a neighborhood of $p$, one can choose coordinates such that
${\mathcal C}_t$ is given by $x_1^2-x_2^2=t$, the projection
$q : \CCC P^2 \to \CCC P^2 \bk conj$ is given by the map $\CCC^2 \to \CCC^2$
$(x_1,x_2) \to (x_1^2,x_2)$.  Q.E.D

\vskip0.1in
\bele
\mlb{l:lars2}
{\it Let ${\mathcal C}_0$ be a singular curve of
degree $2k$ with two crossings $p_1$ and $p_2$ such that
each $p_i$ $i \in \{1,2\}$ belongs to a real branch and
an oval (with a given orientation);the two real branches
which contain $p_i$ have opposite orientation.
(see figure 3.1)
\vskip0.1in
Consider the deformation ${\mathcal C}_t$, $t \in [0,1]$  which crosses
transversally the discriminant $\RRR {\mathcal D}_m$ at ${\mathcal C}_0$.
Then the deformation ${\mathcal C}_{0}$ to ${\mathcal C}_{1}$
implies
${\gA}_1 = {\gA}_0 \krest \RRR P^2  \krest \overline{\RRR P^2}$.}
\enle
{\bf proof:}
\vskip0.1in

In case of Arnold surface ${\gA}_-$, it may be deduced
from the preceding Lemma \mrf{l:lars1}.
It is easy to see that
the directions of twisting of the M\"{o}bius band  given by $p_1$
and $p_2$ are opposite.

In case of
Arnold surface ${\gA}_{+}$,
consider the Morse-function $f :  {\mathcal C}_{0} \ti ]-1,+1[ \approx
{\mathcal C}_{t,t \in ]-1,1[}
     \to  ]-1,+1[$, $f({\mathcal C}_t)=t$.

Obviously, up to isotopy, one can identify ${\mathcal C}_1$
 with any curve ${\mathcal C}_{1 -\e}$
 $\e >0$ small.
Let us fix such $\e>0$.

Consider the descending one-manifolds $D_t$ of $p_1$ and  $D'_t$ of $p_2$.
Each of them reaches the boundary of the quotient curve
${\mathcal C}_{1-\e} \bk conj$ and
defines a normal bundle $N$ of the real set point  $\RRR {\mathcal C}_{1-\e}$
of ${\mathcal C}_{1 - \e}$.
Choose a smooth tangent vector field $V$ on $\RRR P^2$ such that on
$\RRR {\mathcal C}_{1 - \e}$ it is to
tangent to $\RRR {\mathcal C}_{1 - \e}$ and directed according to the
orientation of $\RRR {\mathcal C}_{1- \e}$.
In such a way,
for $x \in \RRR {\mathcal C}_{1 - \e}$, $N(x)=i{V}(x)$
is directed inside the half of $\CCC {\mathcal C}_{1- \e}$
which induces orientation on $\RRR {\mathcal C}_{1 - \e}$.
Recall that multiplication by $i$ makes a real vector normal
to the real plane
and leaves any vector tangent to $\RRR {\mathcal C}_{1 - \e}$
tangent to $\CCC {\mathcal C}_{1 - \e}$.
Extend the tangent vector field $V$ to a tangent vector field of
${\mathcal C}_{ 1- \e + s}$,$0<s<\e, s+e <1 $, in a neighborhood of
 $\RRR {\mathcal C}_{1 -\e} \subset \RRR P^2_{+}$
 (where
 $\RRR P^2_{+} = \{ x \in \RRR P^2| C_{1-\e}(x) \ge 0 \}$
 with $C_{1-\e}$ is a polynomial giving ${\mathcal C}_{1- \e}$.)

Since complex conjugation acts on gradient-trajectories of $D_t$
(resp $D'_t$) as symmetry of  center $p_1$ (resp, $p_2$),
the perturbation of ${\mathcal C}_0$ to
${\mathcal C}_{1-\e}$
implies on Arnold surface ${\gA}_{+}$
one connected sum with the fibre bundle
of the fibration with base space
$\RRR {\mathcal C}_{1-\e}$ and fiber defined for any point
$x \in \RRR {\mathcal C}_{1-\e}$
as the gradient-trajectory through $x$
quotiented by the action of complex conjugation.

Namely, it adds $ \krest \RRR P^2  \krest \overline{\RRR P^2}$
on  Arnold surface ${\gA}_+$.
Obviously, an argument  analogous to this last one may
be also used to describe  Arnold surfaces ${\gA}_-$.

Q.E.D

\vskip0.1in
From the Lemmas  \mrf{l:lars1} and \mrf{l:lars2},
we shall get topological effect of
perturbations of crossings in the recursive construction of
Harnack curves on Arnold surfaces and thus the description
given in Theorem \mrf{t:theo2}.
We shall use the results of proposition \mrf{p:prop7b}.
and (according to the notations introduced in proposition \mrf{p:prop7b})
 we shall denote by $S"_{2k}$
the complementary subset of $S'_{2k}$ inside $S$.
It consists of the
$(k-2)$ points of $S_{2k}^+ \bk S'_{2k}$.

Let $A_{2k}=\{a_1,...,a_{2k-1},...,a_{4k-2} \}$
be the set of points perturbed in a maximal simple deformation of $\HH_{2k}$
(i.e the ordered set of points $a_i$
where $a_1,...a_{2k-1}$ are
the crossings of $\HH_{2k-1}\cup L$
and $a_{2k},...,a_{4k-2}$ are the crossings
of $B_{2k;1}$).
Let
$A_{2k+1}=\{a_1,...,a_{2k} \}$
be the set of points perturbed in a maximal simple deformation of $\HH_{2k+1}$
(i.e the set of crossings $a_i$ of $\HH_{2k} \cup L$).

For any $k \ge 1$,
denote ${\gA}^{2k}_-$ and ${\gA}^{2k}_+$ the Arnold surfaces.
of the Harnack curve $\HH_{2k}$ of degree $2k$.

For $k=1$,  it follows easily from the Lemma \mrf{l:lars1}
that:
${\gA}^{2}_- \approx \overline{\RRR P^2}$ and
${\gA}^{2}_+ \approx D^2$

For $k>1$, pairs $(S^4,{\gA}^{2k}_-)$, $(S^4,{\gA}^{2k}_+)$
may be deduced by induction on $k$.

Let $\HH_{2k-2}$ be the Harnack curve of even degree $2k-2\ge 1$
and ${\gA}^{2k-2}_-$ and ${\gA}^{2k-2}_+$  be its Arnold surfaces.

Let $\HH_{2k}$ be the Harnack curve of even degree $2k$.
and ${\gA}^{2k}_-$ and ${\gA}^{2k}_+$  be its Arnold surfaces.

Then,  one can deduce the pair $(S^4,{\gA}^{2k}_-)$,
(resp,$(S^4,{\gA}^{2k}_+)$)
from the pair $(S^4,{\gA}^{2k-2}_-)$
(resp,$(S^4,{\gA}^{2k-2}_+)$).

\vskip0.1in
The Harnack curve $\HH_{2k}$ is obtained from $\HH_{2k-2}$ by intermediate
construction of $\HH_{2k-1}$
(see proposition \mrf{p:prop7} and
theorem \mrf{t:theo1} of this chapter).
\vskip0.1in
From Lemma \mrf{l:lars1},
perturbations of singular crossings  of
$\HH_{2k-2} \cup L$
and then
of $\HH_{2k-1} \cup L$
imply $(2k-2)+(2k-1)$
$\krest\overline{\RRR P^2}$
on ${\gA}^{2k}_-$,
and imply
$(2k-2)+(2k-1)$ $\krest D^2$ on ${\gA }^{2k}_+$.

\vskip0.1in

From the Lemma \mrf{l:lars1} and Lemma \mrf{l:lars2},
keeping the notations of the proposition \mrf{p:prop7},
it follows that
perturbations of singular crossings of
$S'_{2k} \subset S_{2k}^+$
imply $(k-1)$ $\krest \overline {\RRR P^2}$ on ${\gA}^{2k}_-$,
and imply $(k-1)$ $\krest D^2$ on ${\gA}^{2k}_+$.
\vskip0.1in
Furthermore, the perturbations of the $(k-2)$ singular crossings of
$S"_{2k} \subset S_{2k}^+$
imply \\(joined  with
perturbations of singular crossings of
$S'_{2k-2} \subset S_{2k-2}^+$)\\
$(k-2)$ $\krest \RRR P^2$ on ${\gA}^{2k}_-$
\vskip0.1in
From the Lemma \mrf{l:lars2}, it follows that perturbations of singular crossings of\\
$S"_{2k}  \subset S_{2k}^+$
joined with perturbations of singular crossings of
$S'_{2k-2} \subset S_{2k-2}^+$
imply $(k-2)$
$ \krest \RRR P^2  \krest \overline{\RRR P^2}$
on ${\gA}^{2k}_{+}$.

Hence ${\gA}^{2k}_- \approx {\gA}^{2k-2}_-
\krest (k-2)  \krest \RRR P^2  \krest
(5k-4) \overline{\RRR P^2}$
and
${\gA}^{2k}_+ \approx {\gA}^{2k-2}_+
 \krest (k-2) \krest  \RRR P^2  \krest
(k-2) \overline{\RRR P^2}$.

Q.E.D

\part{\bf Perestroika Theory on Harnack curves}

As already noticed in the last section, the classification problem
of pairs $(S^4, {\AA}_{\pm})$
amounts to the classification problem of real algebraic curves
up to conjugate equivariant isotopy of $\CCC P^2$.
\vskip0.1in
Thus, there is  an obvious connection  with classification
of Arnold surfaces up to isotopy of $S^4$
and Hilbert's Sixteen Problem on the arrangements of ovals
of real algebraic curves.
\vskip0.1in

First of all we shall detail the
construction of Harnack curves.
Then, we shall  deduce
from this detailed construction, a construction of
curves ${\AA}_m$ of degree $m$
which provides a
description of the pair $(\CCC P^2, \CCC {\AA}_m)$ up to conj-equivariant
isotopy of $\CCC P^2$ (see Theorem \mrf{t:th2typ1} and
Theorem \mrf{t:th1typ2}).
Moreover, this method gives all
possible arrangements of real connected  components
of  curves and therefore is an advancement in the Hilbert's Sixteen Problem.
\vskip0.1in
After that, we shall deal with
Arnold surfaces defined on any curve of even degree with non-empty
real part.
\chapter{Recursive Morse-Petrovskii's theory of
recursive Harnack curves}
\lb{ch:mpr}
In this section, we shall go on the study of critical
points of Harnack's polynomial initiated in chapter  \ref{ch:Mope}.\\
Recall that given a Harnack polynomial
$B_m(x_0,x_1,x_2)=x_0^m.b_m(x_1/x_0,x_2/x_0)$
we call critical point of $B_m(x_0,x_1,x_2)$
any point $(x_0,y_0)$ such that $b_x(x_0,y_0)=0$
and $b_y(x_0,y_0)=0$.
We shall consider only  Harnack curves
given by regular polynomials.
According to the  classification
of Harnack curves $\HH_m$ up to rigid isotopy  (Theorem \ref{t:rigiso}),
one can assume that any Harnack curve $\HH_m$
results of the recursive construction of Harnack curves $\HH_i$,
$1 \le i \le m$, where the curve $\HH_{i+1}$ is deduced from classical
small perturbation of the union $\HH_i \cup L$ of the curve
$\HH_i$ with a line $L$.
This section until its end is devoted to
settle how
critical points of a Harnack polynomial of degree $m$
may be  associated with crossings of
curves $\HH_{m-i} \cup L$, $1 \le i \le (m-1)$,
and how
a critical point associated with a crossing
of $\HH_{m-i} \cup L$, $1 \le i \le (m-1)$,
varies in the recursive construction of Harnack curves $\HH_{m}$.\\

Given a regular Harnack polynomial $B_m(x_0,x_1,x_2)$
with affine associated polynomial
$b_m(x_1/x_0,x_2/x_0)$, we shall study the pencil of curves
given by polynomial
$x_0^m.(b_m(x_1/x_0,x_2/x_0)-c)$, $ c \in \RRR$.
Since Harnack curves realize the maximal number of real components
curves of the pencil may have,
to each critical point of index 1  of a regular Harnack polynomial
corresponds a "gluing" of real components
in the pencil. These gluings are the subject of this section.\\

The main result of this section is gathered in the
Lemma \mrf{l:l1p7}, Lemma \mrf{l:l2p7} and Lemma \mrf{l:l3p7}
where we describe critical points of index $1$ of
Harnack polynomials.

\vskip0.1in
\section{Preliminaries}

Given a regular polynomial $R(x_0,x_1,x_2)=x_0^m.r(x_1/x_0, x_2/x_0)$,
for any real finite critical point
$(x_0,y_0)$ of $r(x,y)$ ,
$r(x_0,y_0) =c_0$,
there exists $\e >0$ sufficiently small such that
$r^{-1} [c_0 - \e,c_0 + \e]$ contains no critical point other than
$(x_0,y_0)$.\\
We shall call {\it topological meaning} of $(x_0,y_0)$ the homeomorphism
type of the following triad of spaces $(W(x_0,y_0); \pr_0 W(x_0,y_0),
 \pr_1 W(x_0,y_0))$
with
$$ W(x_0,y_0) = \cup_{c \in [c_0 - \e,c_0 + \e]} \{ (x_0:x_1:x_2)
\in \CCC P^2 |  x_0^m(r(x_1/x_0,x_2/x_0)-c))=0 \}$$
$$\pr_0 W(x_0,y_0) =\{ (x_0:x_1:x_2)
\in \CCC P^2 | x_0^m(r(x_1/x_0,x_2/x_0)-(c_0- \e))=0 \}$$
$$\pr_1 W(x_0,y_0) =\{ (x_0:x_1:x_2)
\in \CCC P^2 | x_0^m(r(x_1/x_0,x_2/x_0)-(c_0+ \e))= 0 \}$$
\vskip0.1in

Let $B(x_0,y_0) \subset \CCC P^2$  be a small 4-ball around $(x_0,y_0)$
globally invariant by complex conjugation
such that
$\pr_0 W(x_0,y_0) \cap B(x_0,y_0) \not= \emptyset$
and
$\pr_1 W(x_0,y_0) \cap  B(x_0,y_0) \not= \emptyset$

We shall
call {\it  local topological meaning}
of $(x_0,y_0)$ the homeomorphism
type of the following triad of spaces \\
$(W(x_0,y_0)\cap B(x_0,y_0) ;
 \pr_0 W(x_0,y_0) \cap B(x_0,y_0),
 \pr_1 W(x_0,y_0) \cap  B(x_0,y_0))$

In what follows, we shall associate to each critical point $p$ of index 1
of a Harnack polynomial of degree $m$,
two real branches of the Harnack curve of degree $m$ involved in the
local topological meaning of $p$.
\vskip0.1in

{\bf Terminology}
\vskip0.1in
Recall that we distinguish two ways
to perturb a crossing $p$ of a singular
curve ${\AA}$ defined
in the section \mrf{s:comp}
as {\it perturbation of type 1}
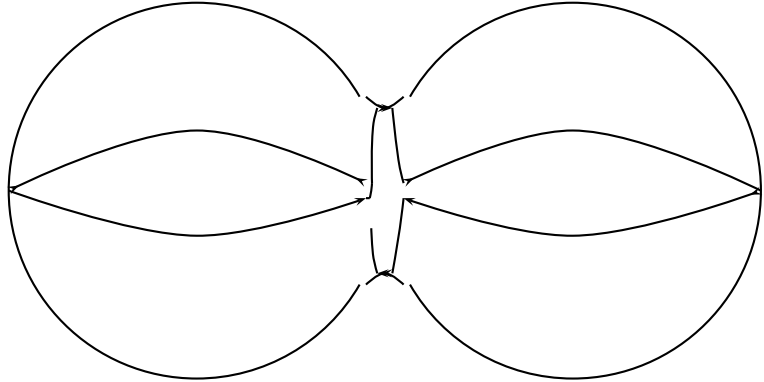
and {\it perturbation of type  2} of $p$.

Thus, local topological meaning of a crossing $(x_0,y_0)$
can be deduced  from  one  topological space
$\pr_1 W(x_0,y_0)\cap B(x_0,y_0)$
or $\pr_0 W(x_0,y_0) \cap B(x_0,y_0)$
and the type $1$ or $2$ of the perturbation  of the point $(x_0,y_0)$.
Let us introduce the notations:
$(\pr_1 W(x_0,y_0) \cap B(x_0,y_0))_i$
and  $(\pr_0 W(x_0,y_0) \cap B(x_0,y_0))_i$
where the subscript  $i$
stands for the type of the perturbation of the crossing.
\vskip0.1in
Although any perturbation of the union of two
projective lines leads to
a conic with orientable real part,
the perturbation of a singular curve of degree $m >2$
of which singular points are crossings may lead to a curve with
non-orientable part.
Nonetheless, in case of a Harnack curve $\HH_m$, any deformation of
a crossing which appears in  a maximal simple deformation of $\HH_m$
(see Proposition \mrf{p:prop7b} and Proposition \mrf{p:prop8})
agrees with a complex orientation of its real part $\RRR \HH_m$.
Besides, it is easy to deduce from
relative orientation and location of real branches of $\RRR \HH_m$
that any such crossing is deformed by a perturbation of type 1.
\vskip0.1in

\section{Critical points and recursive construction of Harnack curves}
\vskip0.1in
\subsection{Introduction}
\vskip0.1in
Let $\HH_{m}$ be the Harnack curve of degree $m$.
In what follows, for any $m>0$, we shall consider only
Harnack polynomials
$B_{m} (x_0,x_1,x_2)=x_0^m.b_m(x_1/x_0,x_2/x_0)$
of degree $m$ and type $\HH^0$.
Since no confusion is possible, we shall call such polynomials
Harnack polynomials.
As already introduced in the preceding section (see proof of
propositions\mrf{p:prop7b} and \mrf{p:prop7}),
one can associate critical points of index $1$ of $B_m$
with crossings of curves $\HH_{m-1} \cup L$.
\vskip0.1in
Formally, we shall say that:

\bede
\mlb{d:uu}
{\it
Let $a$ be a crossing of the curve $\HH_{m-i} \cup L$.
Given a disc $D(a,\e)$ (in the Fubini-Study metric) of radius $\e$
around of $a$ in $\RRR P^2$.

Let $U(a,\e)$ be the neighborhood of $\HH_{m-i+1}$
in $\CCC P^2$
deduced from the action of $U_{\CCC}^3$ on $U(a)$:

$$U(a,\e)=\{ z=<u,p>=(u_0.p_0:u_1.p_1:u_2.p_2) \in \CCC P^2
\mid u=(u_0,u_1,u_2) \in U_{\CCC}^3, p \in U(a) \}$$.

A critical point $(x_0,y_0)$ of index 1 of $B_m$
is associated with a crossing $a$ of the curve
$\HH_{m-i} \cup L$ if there exists $\e_0$ such that:
\been
\item
$a, (x_0,y_0) \in U(a,\e_0)$
\item
the perturbation on the real part of $\HH_m$ involved in the local
topological meaning of $(x_0,y_0)$ is a deformation
on $\CCC \HH_{m} \cap U(a,\e_0)$
\enen}
\end{defi}

\vskip0.1in
As already noticed corollary \mrf{c:gpa}
of the previous section \mrf{s:comp},
without loss of generality, one can consider that Harnack curves
are obtained via the Patchworking method.

Hence, for sake of simplicity, we shall consider
the $T$-Harnack curves introduced in the Chapter \ref{ch:PatchHar}.
For the patchworking construction of Harnack
we can give the following of version definition \ref{d:uu}.\\

Recall that
 $\rho^m : \RRR_+ T_m \ti U_{\CCC}^2 \to \CCC T_{m} \approx \CCC P^2 $
and its restriction \quad \break
$\rho^m|_{ \RRR_+ T_m \ti U_{\RRR}^2}
 : \RRR_+ T_m \ti U_{\RRR}^2 \to \RRR T_{m} \approx \RRR P^2$
denote the natural surjections.

Recall (see definition \mrf{d:def6}) that
given  $\G$  a face of the triangulation
of $T_m$ and $D(p,\e) \in \RRR_+ T_m ^0$
an (euclidian) open 2-disc such that
the moment map $\mu : \CCC T_m  \to T_m$ maps $D(p,\e)$ to a two-disc
$\mu(D(p,\e)$ which contains $\G$   and
intersects only the face $\G$ of the triangulation
of $T_m$;
we call
$U(p)=\rho^m (D(p,\e) \ti U_{\CCC}^2)$
the $\e$-neighborhood of $\CCC \HH_m$ in $\CCC T_m \approx \CCC P^2$
defined from $\G^0$;
we call $\e$-tubular neighborhood
of $\CCC \HH_{m}$ defined from $\G^0$
the $\e$-tubular neighborhood
of $\CCC \HH_{m}$ in
$U(p)$.

\bede
\mlb{d:uu2}
{\it
Let $\G$ be a face of the triangulation of $T_m$ and
$U(p)$ be the $\e$-neighborhood
of $\CCC \HH_{m}$ defined from $\G^0$.
A critical point $(x_0,y_0)$ of index 1 of $B_m$
is associated with a crossing $a$ of the curve
$\HH_{m-i} \cup L$ if :
\been
\item
$a, (x_0,y_0) \in U(p)$
\item
the perturbation on the real part of $\HH_m$ involved in the local
topological meaning of $(x_0,y_0)$ is a deformation
on $\CCC \HH_{m} \cap U(p)$
\enen}
\end{defi}
\vskip0.1in
\bede
\mlb{d:def1p7}
{\it
Let $p_{m}$ (resp $p_{m+j}$, $j >0$) be a
critical point of index 1
of a Harnack polynomial of degree  $m$ (resp $m+j$).
We say that $p_{m+j}$
is equivalent to $p_m$ if
there exist small conj-equivariant open $4$-balls $B(p_{m+j})$
and $B(p_m)$ around  $p_{m+j}$ and $p_m$
with the following properties:
\been
\item
$B(p_{m+j}) \subset B(p_m)$
\item
The triad of spaces
$$( W(p_{m+j} \cap B(p_{m+j}); \pr_0 W(p_{m+j}) \cap B(p_{m+j}),
 \pr_1 W(p_{m+j}) \cap B(p_{m+j}))$$
is the local topological meaning of $p_{m+j}$.
\vskip0.1in
The triad of spaces
$$( W(p_{m} \cap B(p_{m}); \pr_0 W(p_{m}) \cap B(p_{m}),
 \pr_1 W(p_{m}) \cap B(p_{m}))$$
is the local topological meaning of $p_{m+j}$.
\item
Local topological meanings of $p_{m+j}$
and $p_m$ are homeomorphic.
\enen}
\end{defi}
\bere
Let $B_m(x_0,x_1,x_2)$ be a $T$-Harnack polynomial.
It is an easy consequence of the $T$-construction
of  Harnack curves, that one can define $4$-ball of $\CCC P^2$
around any critical point of $B_m(x_0,x_1,x_2)$ as usual $4$-ball of
$(\CCC ^*)^2$.
Indeed, none of the critical points of $B_m(x_0,x_1,x_2)$
belongs to the coordinates axes.
\enre

We shall denote
$p_{m+j}$ is equivalent to $p_{m}$
by
$p_{m+j} \approx p_{m}$.

\bede
\mlb{d:def2p7}
{\it
Given $a$ a crossing of $\HH_{m-1} \cup  L$,
denote by $p_{l}$ the critical point of index 1
of a Harnack polynomial of degree $l \ge m$
associated with $a$.
\been
\item
We say that $p_m$ is
the simple point of $\HH_{m}$ associated with $a$.
\item
If for any $j \ge 1$,
$p_{m+j}$  is equivalent to $p_{m}$,
 we say that
a simple point of $\HH_{m+j}$ with $j \ge 0$
is associated with $a$.
\item
If for any integers $l$,$l'$
,$0 < l \not= l' \le k$,
$p_{m+l}$ is not equivalent to $p_{m+l'}$;
for any $l \le k$,
we say that a $l$-point
$(p_{m},p_{m+1},..., p_{m+l})$ of $\HH_{m+l}$
is associated with $a$.
\item
If $k$ is the smallest integer such that:
for any integers $l,l'$,\\ $0 < l \not = l' \le k$,
$p_{m+l}$ is not equivalent to $p_{m+l'}$;
for any $l \ge k$,
we say that a $k$-point
$(p_{m},p_{m+1},..., p_{m+k})$ of $\HH_{m+l}$
is associated with $a$.
\vskip0.1in
\enen}
\end{defi}

\subsection{Critical points of index $1$ of Harnack polynomials}
\vskip0.1in
We shall now proceed to the study of
critical points of index $1$ of Harnack polynomials.
We shall work
with the notations introduced in the chapter \ref{ch:Mope}.
Given $B_m(x_0,x_1,x_2)= x_0^m.b_m(x_1/x_0,x_2/x_0)$
a Harnack polynomial of type $\HH^0$, we
denote by $S_{m}$ the set of critical points of index 1 of
$b_m(x,y)$, and
by $S_{m}^-$ (resp, $S_{m}^+$) the subset of $S_m$
consisting respectively of critical points of index 1
with negative (resp, positive) critical value.

In case $m=2k$,
we distinguish the two subsets
$S'_{2k}$, $S"_{2k}$ of $S_{2k}^+$ with the properties
 $S'_{2k} \cup S"_{2k}= S_{2k}^+$
 $S'_{2k} \cap S"_{2k}= \emptyset$.
Given $B_{2k}(x_0,x_1,x_2)$ a Harnack polynomial of type $\HH^0$
and $b_{2k}(x_1/x_0,x_2/x_0)$ its affine associated polynomial.
The set $S'_{2k}$  is
constituted by the $c'_1(B_{2k})$
critical points $(x_0,y_0)$ of index 1 $b_{2k}(x_0,y_0)=c_0$, $c_0 >0$
with the property that as $c$ increases from $c_0 - \e$ to $c_0 + \e$
the number of real connected components
of ${\mathcal A}_c= \{ (x_0:x_1:x_2) \in \RRR P^2 | x_0^{2k}.(b_{2k}-c)=0 \}$
intersecting
the line at infinity increases by 1.
The subset $S"_{2k}$  denotes
the complementary set of $S'_{2k}$ inside $S^+_{2k}$.
These sets were already
under consideration in propositions \mrf{p:crit2}, \mrf{p:prop7b} and
theorem \mrf{t:theo2}.
\vskip0.1in
Let us start by the Lemma \mrf{l:l1p7} in which we prove that critical
points of negative critical value of Harnack polynomials of type $\HH^0$
are simple  points in the recursive construction of Harnack curves.
\vskip0.1in

\bele
\mlb{l:l1p7}
\vskip0.1in
{\it
Let $\HH_m$ be a Harnack curve of degree $m$ obtained via the
patchworking method
 and given by a Harnack polynomial of type $\HH^0$.
\been
\item
In case of odd $m$, let $\G$ be a face of the triangulation of
$T_m$ which belongs to $l_{m-1}$.
\item
In case of even $m$, let $\G$ given by vertices
$(c,d+1) (c+1,d)$ with $c=0$ or $d=0$.
\enen
Let $U(p)$ be the $\e$-neighborhood
of $\CCC \HH_m$ in $\CCC P^2$ defined from $\G^0$.
Denote $a$ the unique crossing
of $\HH_{m-1} \cup L$ which belongs to $U(p)$.
Denote $p_m$ the critical point of $\HH_m$ associated with $a$,
and let $B(p_m) \subset U(p)$ a conj-equivariant $4$-ball
around $p_m$.
Then,
\been
\item
\lb{i:l11}
the critical point $p_m$ of $\HH_m$ belongs to $S_m^-$\\
Moreover,
$\HH_{m} \cap B(p_{m}) \approx (\pr_1 W(p_m) \cap B(p_m))_1$
\item
\lb{i:l12}
a simple-point of curves $\HH_{m+j}$, $j \ge 0$,
is associated with $a$
\enen }
\enle
\vskip0.1in
{\bf proof:}
\vskip0.1in
\vskip0.1in
The proof is based on the Petrovskii's theory and propositions
\mrf{p:crit2}, \mrf{p:prop7b} and \mrf{p:prop7}.
Our argumentation is similar to
the one used in the proof of proposition \mrf{p:prop7}.

Assume $\G$ given by vertices
$(c,d+1) (c+1,d)$.
For any $j \ge 0$,
given \\ $B_{m+j}=x_0^{m+j} b_{m+j} (x_1/x_0,x_2/x_0)$
the Harnack polynomial of degree $m+j$, we shall denote
$b_{m+j}^S$ the truncation of $b_{m+j}$ on the monomials
 $x^cy^d$,$x^cy^{d+1}$
 $x^{c+1}y^d$,
 $x^{c+1}y^{d+1}$.

We shall use local description of Harnack
curves, $\HH_{m+j}$, $j \ge 0$,
inside $U(p)$ provided by the patchworking theory.

Recall that for any $j \ge 0$, there exists
an homeomorphism \\ $\tilde{h} :\CCC \HH_{m+j} \cap U(p) \to
 \{(x,y) \in (\CCC^*)^2
| b_{m+j}^S =0 \}\cap U(p)$
such that \\
$\tilde{h}(a)=(x_0,y_0)$ is a critical point of $b_{m+j}^S$.

We shall divide our proof in two parts.
In the first part, we shall consider the curve $\HH_m$ and
study local topological meaning of crossings of $\HH_{m-1} \cup L$.
Then, we shall consider curves $\HH_{m+j}$, $j \ge 1$.
\vskip0.1in
{\bf (\ref{i:l11})}
Let us distinguish the cases $m$ even and $m$ odd.\\
{\bf (\ref{i:l11}).a}
In case of even $m=2k$
\vskip0.1in
Let $\G \in l_{2k-1}$ given by vertices
$(c,d+1) (c+1,d)$ with\\ $c=0$ or $d=0$.
(Namely, $\G$ has vertices $(0,2k-1),
 (1,2k-2)$ or $(2k-1,0), (2k-2,1)$.)
From the patchworking theory,
(see corollary \mrf{c:patch} of the chapter \ref{ch:PatchHar})
$\tilde{h}(a)=(x_0,y_0)$
$x_0>0, y_0<0 $
is a critical point of $b_{2k}^S(x,y)$.

Keeping the notations introduced in
the chapter \ref{ch:Mope}, set $b_{2k}=x_{2k;t}$ with
$t \in ]0,t_{2k}[$.
Then, for fixed $t \in ]0,t_{2k}[$,
let $a_t$ be the crossing of $\HH_{2k-1} \cup L$ contained in $U(p)$.

Therefore,
 $\tilde{h}(a_t)=(x_0,y_0)_t$
$x_0> y_0<0 $
is a critical point of $x_{2k,t}^S(x,y)$.

Set
$x_{2k;t}^S(x,y)=l_t(x,y)+t.k_t(x,y)$
with \\
$l_t(x,y)=
 a_{c,d}  t^{\nu(c,d)}x^cy^d  +
 a_{c+1,d} t^{\nu(c,d+1)}x^{c+1}y^d $;\\
$k_t(x,y)=
 a_{c,d+1} t^{\nu(c+1,d) -1}x^cy^{d+1}  +
-a_{c+1,d+1} t^{\nu(c+1,d+1) -1}x^{c+1}y^{d+1}$.
\vskip0.1in
Then, it is easy to see
that, up to modify the coefficients\\
$a_{c,d},a_{c,d+1},a_{c+1,d},a_{c+1,d}$ if necessary, the point
$\tilde{h}(a_t)=(x_1,y_1)_t$ is
a critical point of the function
$\frac {l_t(x,y)} {k_t(x,y)}$ with negative critical value
$\ge- t>- t_{2k-1}$.

On this assumption, it follows from the equalities\\
$l_t(x,-y)=-l_t(x,y)$, $k_t(x,-y)=k_t(x,y)$,\\
$\frac {\pr l_t} {\pr x}(x,-y) =-\frac {\pr l_t} {\pr x}(x,y)$,
$\frac {\pr l_t} {\pr y} (x,-y) =\frac {\pr l_t} {\pr y} (x,y)$,\\
$\frac {\pr k_t} {\pr x} (x,-y)= \frac {\pr k_t} {\pr x} (x,y)$,
$\frac {\pr k_t} {\pr y} (x,-y) = -\frac {\pr k_t} {\pr y} (x,y)$\\
that $(x_1,- y_1)_t$ is a critical  point of
the function $\frac {l_t(x,y)} {k_t(x,y)}$ with positive critical value
$\le t< t_{2k-1}$.
Therefore, up to modify coefficients of $b_{2k}^S$,
$(x_0,-y_0)$ is a critical point
with negative critical value
of $b_{2k}^S$.
Obviously, (see \cite{Petr}),
such modification does not affect the order and
the topological structure of $b_{2k}^S$.
Therefore, it does not change local topological meaning of
critical points of the Harnack polynomial $B_{2k}$.\\
We shall now proceed to the study of
the Petrovskii's pencil of curves
over $\HH_{2k}$.

Let $d_{2k}$ be the unique  polynomial such that
$b_{2k}= b_{2k}^S + d_{2k}$.
Consider curves of the Petrovskii's pencil over $\HH_{2k}$
outside $d_{2k} =0$,
as level curves of the function
$\frac {b_{2k}^S -c} {d_{2k}}$.

Inside $b_{2k}^S-c=0 \bk d_{2k}=0$,
these curves have critical points the singular
points of $b_{2k}^S-c =0$.

Bringing together  Petrovskii's theory and
the implicit function theorem applied to one-parameter
polynomial \\
$x_0^{2k}.(b_{2k} -c)=
x_0^{2k}.(b_{2k}^S+d_{2k} -c)$
with parameter $c$, it  follows
that a point $p_{2k}$ with
$b_{2k}(p_{2k})=c_0 <0$
is associated with $a$.

When $c=c_0$,
one positive oval touches an other positive oval.
\vskip0.1in

{\bf (\ref{i:l11}).b}
In case of odd $m =2k+1$,
one has to distinguish whether $\G$ intersects or not the
coordinates axes.\\
{\bf (\ref{i:l11}).b.1}

Let $\G \in l_{2k}$ given by vertices
$(c,d+1) (c+1,d)$ with \\ $c=0$ or $d=0$.
 (Namely, $\G$ has vertices $(0,2k),
 (1,2k-1)$ or $(2k,0), (2k-1,1)$.)

According to  the patchworking theory,
(see corollary \mrf{c:patch} of the first section)
$\tilde{h}(a)=(x_0,y_0)$ with
$x_0.y_0>0$
is a critical point of $b_{2k+1}^S(x,y)$.
Set
$\tilde{h}(a)=(x_0,y_0) x_0>0, y_0>0$.

On this assumption, following the previous argumentation,
up to modify coefficients of $b_{2k+1}^S$ if necessary
(without changing the order and
the topological structure of $b_{2k+1}$)
$(x_0,-y_0)$ is a critical point with
negative critical value of $b_{2k+1}^S$.

Therefore, according to Petrovskii's theory,
from an argumentation similar to the one given in even case,
a point $p_{2k+1}$
with $b_{2k+1}(p_{2k+1})=c_0 <0$ is associated with $a$.

When $c=c_0$,
a positive oval touches the one-side component of the curve.\\
{\bf (\ref{i:l11}).b.2}

We shall now consider faces which do not intersect the coordinates axes.
Our argumentation is a slightly modified
version of the previous one.
We shall study the situation inside
$\e$-neighborhood defined from two faces together.
Let $\G \in l_{2k}$  given by vertices
$(c,d+2), (c+1,d+1)$ with $c$ odd,
$c\not=0$ and $d\not=0$ and
$\G' \in l_{2k}$  given by vertices
$ (c+2,d), (c+1,d+1)$.
\vskip0.1in
Consider the convex polygon $K \subset T_{2k+1}$ with vertices\\
$(c+1,d),(c+2,d),(c+2,d+1),(c+1, d+2),(c,d+2),(c,d+1)$.
It is contained in the triangle $T_{2k+1}$ and triangulated by
the triangulation $\tau$ of $T_{2k+1}$.
Denote $U(p)$ the
$\e$-neighborhood of $\CCC \HH_{2k+1}$ defined from $\G^0$.
Denote $U(p')$ the $\e$-neighborhood of $\CCC \HH_{2k+1}$ defined
from ${\G'}^0$.
Let $a$, (resp, $a'$)
be the crossing of
$\HH_{2k} \cup L$
which belongs to $U(p)$,
(resp, $U(p')$).
Denote $b_{2k+1}^S$ (resp, $b_{2k+1}^{S'}$) the truncation
 of $b_{2k+1}$ on the monomials
 $x^cy^{d+1}$,$x^cy^{d+2}$
 $x^{c+1}y^{d+2}$,
 $x^{c+1}y^{d+1}$
(resp,
 $x^{c+1}y^d$,$x^{c+1}y^{d+1}$
 $x^{c+2}y^{d+1}$,
 $x^{c+2}y^{d}$.)

As previously, consider  homeomorphisms \\
$\tilde{h} :\CCC \HH_{2k+1} \cap U(p) \to
 \{(x,y) \in (\CCC^*)^2
| b_{2k+1}^S =0 \}\cap U(p)$
such that $\tilde{h}(a)=(x_0,y_0)$
$x_0>0,y_0>0$
 is a critical point of $b_{2k+1}^S$
\vskip0.1in
and $\tilde{h}' :\CCC \HH_{2k+1} \cap U(p') \to
 \{(x,y) \in (\CCC^*)^2
| b_{2k+1}^{S'}=0 \}\cap U(p')$
such that $\tilde{h}'(a')=(x'_0,y'_0)$
$x_0>0,y_0>0$
is a critical point of $b_{2k+1}^{S'}$.
\vskip0.1in
Keeping the notations introduced in the chapter \ref{ch:Mope},
set $b_{2k+1}=x_{2k;t}$,
$t \in ]0,t_{2k+1}[$.

For fixed $t \in ]0,t_{2k+1}[$, let $a_t$ (resp $a'_t$)
be the crossing of $\HH_{2k} \cup L$ contained in $U(p)$
(resp $U(p')$).
Set
$x_{2k+1;t}^K(x,y)=l_t(x,y)+t.k_t(x,y)$
with \\
$l_t(x,y)=
 a_{c+1,d}  t^{\nu(c,d)}x^{c+1}y^d  +
 a_{c+2,d} t^{\nu(c+2,d)}x^{c+2}y^d +\\
 a_{c+1,d+2}  t^{\nu(c,d)}x^{c+1}y^{d+2}+
 a_{c+2,d+2} t^{\nu(c+2,d)}x^{c+2}y^{d+2}$;\\
$k_t(x,y)=
 a_{c,d+1} t^{\nu(c+1,d) -1}x^cy^{d+1}
 -a_{c+1,d+1} t^{\nu(c+1,d+1) -1}x^{c+1}y^{d+1}+\\
 a_{c+2,d+1} t^{\nu(c+1,d) -1}x^{c+2}y^{d+1}$
\vskip0.1in
Then, it is easy to see
that, up to modify the coefficients
$a_{i,j}$ if necessary,
points
$\tilde{h}(a_t)=(x_0,y_0)_t$
and
$\tilde{h'}(a'_t)=(x'_0,y'_0)_t$
are critical points of the function
$\frac {l_t(x,y)} {k_t(x,y)}$ with negative critical value
$\ge- t>- t_{2k-1}$.
On this assumption, it follows from the equalities \\
$l_t(x,-y)=-l_t(x,y)$, $k_t(x,-y)=k_t(x,y)$,  \\
$\frac {\pr l_t} {\pr x}(x,-y) =-\frac {\pr l_t} {\pr x}(x,y)$,
$\frac {\pr l_t} {\pr y} (x,-y) =\frac {\pr l_t} {\pr y} (x,y)$,\\
$\frac {\pr k_t} {\pr x} (x,-y)= \frac {\pr k_t} {\pr x} (x,y)$,
$\frac {\pr k_t} {\pr y} (x,-y) = -\frac {\pr k_t} {\pr y} (x,y)$
that points $(x_0,- y_0)_t$ and $(x'_0,-y'_0)_t$
are critical  point of
the function $\frac {l_t(x,y)} {k_t(x,y)}$ with positive critical value
$\le t< t_{2k-1}$.
Therefore, up to modify the coefficients of $b_{2k+1}^K$,
$(x_0,-y_0)$ and $(x'_0,-y'_0)$ are critical points
with negative critical value of $b_{2k+1}^K$.
Obviously, (see \cite{Petr}),
such modification does not modify the order and
the topological structure of $b_{2k+1}^K$.
Therefore, it does not change local topological meanings of
critical points of the Harnack polynomial $B_{2k+1}$.\\
We shall now proceed to the study of
the Petrovskii's pencil of curves
over $\HH_{2k+1}$.
Let $U(K^0)$ be the subset $\rho^{2k+1}(\RRR_+ K^0  \ti U_{\CCC}^2)$
of $\CCC P^2$.\\
Obviously, $U(p) \cup U(p') \subset U(K^0)$.
According to patchworking theory,
the truncation $b_{2k+1}^K$ of $b_{2k+1}$ is $\e$-sufficient
for $b_{2k+1}$ in $U(K^0)$.
Therefore,\\
$\tilde{h} :\CCC \HH_{2k+1} \cap U(p) \to
 \{(x,y) \in (\CCC^*)^2
| b_{2k+1}^S =0 \}\cap U(p)$\\ and
$\tilde{h} :\CCC \HH_{2k+1} \cap U(p') \to
 \{(x,y) \in (\CCC^*)^2
| b_{2k+1}^{S'}=0 \}\cap U(p')$
extend to the homeomorphism
$\tilde{h} :\CCC \HH_{2k+1} \cap U(K^0) \to
 \{(x,y) \in (\CCC^*)^2
| b_{2k+1}^K =0 \}\cap U(K^0)$
such that
 $\tilde{h}(a)=(x_0,y_0)$ $x_0.y_0>0$
 $\tilde{h}(a')=(x'_0,y'_0)$ $x'_0.y'_0>0$.
Let $d_{2k+1}$ be the unique polynomial such that
$b_{2k+1}= b_{2k+1}^K + d_{2k+1}$.
Consider curves of the Petrovskii's pencil over $\HH_{2k+1}$
outside $d_{2k+1} =0$, as level curves of the function
$\frac {b_{2k+1}^K -c} {d_{2k+j}}$.

From an argumentation similar to the previous one, according to
the Petrovskii's theory, it follows that:
a point  $p_{2k+1}$
with $b_{2k+1}(p_{2k+1}) =c_0 <0$
is associated with $a$
and a point
$p'_{2k+1}$, $b_{2k+1}(p'_{2k+1}) =c'_0 <0$
is associated with $a'$.
When $c=c_0$, the one-side component of the curve touches itself.
When $c=c'_0$, the one-side component of the curve touches itself.
\vskip0.1in

{\bf (\ref{i:l12})}
Let $a$ be one singular point
of the curve $\HH_{m-1} \cup L$
under consideration in the part $\ref{i:l11}$ of the proof
and $p_m$ be the critical point of the Harnack polynomial of degree $m$
associated with $a$.
We shall now prove that
for any  such crossing $a$ and any $j \ge 1$,
there exists a critical point $p_{m+j}$ of a Harnack polynomial
$B_{m+j}$ associated with $a$.
The point $p_{m+j}$ is equivalent to
the critical point $p_m$ of $B_m$.

Let $j\ge1$ and $B_{m+j}$ be a Harnack polynomial of degree $m+j$

As previously, consider
the homeomorphism \\ $\tilde{h} :\CCC \HH_{m+j} \cap U(p) \to
 \{(x,y) \in (\CCC^*)^2
| b_{m+j}^S =0 \}\cap U(p)$.
\vskip0.1in
{\bf (\ref{i:l12}).a}

Assume $\G$ intersects the coordinates axes
($\G \in l_{2k-1}$ or $\G \in l_{2k}$)
Let $d_{m+j}$ be the unique  polynomial such that
$b_{m+j}= b_{m+j}^S + d_{m+j}$.
Consider curves of the Petrovskii's pencil over $\HH_{m+j}$
outside $d_{m+j} =0$,
as level curves of the function
$\frac {b_{m+j}^S -c} {d_{m+j}}$.
Inside $b_{m+j}^S-c=0 \bk d_{m+j}=0$,
these curves have critical points the singular
points of $b_{m+j}^S-c =0$.

From an argumentation similar to the previous one, according to the
Petrovskii's theory (Petrovskii's pencil over $\HH_{m+j}$),
it follows easily
that for any $j \ge 1$,
a point $p_{m+j}$
is associated with $a$  and is equivalent to the point $p_m$.
\vskip0.1in

{\bf (\ref{i:l12}).b}
Assume $\G$ does not intersect the coordinates axes
($\G \in  l_{2k}$).
Consider as previously the convex polygon $K$ with vertices\\
$(c+1,d),(c+2,d),(c+2,d+1),(c+1, d+2),(c,d+2),(c,d+1)$.
Let $d_{m+j}$ be the unique  polynomial such that
$b_{m+j}= b_{m+j}^K + d_{m+j}$.

Consider curves of the Petrovskii's pencil over $\HH_{m+j}$
outside $d_{m+j} =0$,
as level curves of the function
$\frac {b_{m+j}^K -c} {d_{m+j}}$.
From an argumentation similar to the previous one, according to the
Petrovskii's theory, it follows that for any $j \ge 1$,
a point $p_{m+j}$
is associated with $a$  and is equivalent to the point $p_m$.
Q.E.D.
\vskip0.1in
Then, consider the case when a critical point $p_m$ of  $S^+_{m}$
is associated  to  a crossing of $\HH_{m-1} \cup L$.
Such critical point appears the first time as critical point of
a Harnack polynomial  of degree $4$.
\vskip0.1in
Recall (see the Proposition \mrf{p:crit2} ) that
for a Harnack polynomial
$B_{2k}(x_0,x_1,x_2)=x_0^{2k}.b_{2k}(x_1/x_0,x_2/x_0)$ of type $\HH^0$
we denote  $S'_{2k}$ the set constituted by the $c'_1(B_{2k})$
critical points $(x_0,y_0)$ of index 1 $b_{2k}(x_0,y_0)=c_0$, $c_0 >0$
with the property that as $c$ increases from $c_0 - \e$ to $c_0 + \e$
the number of real connected components intersecting
the line at infinity increases by $1$.
In the Lemma \mrf{l:l2p7}, we shall study how critical points $S'_{2k}$
of a Harnack polynomial of type $\HH^0$ and degree $2k$ vary
in the recursive  construction of Harnack curves.

\bele
\mlb{l:l2p7}
\vskip0.1in
{\it
Let $\HH_{2k}$ be the Harnack curve of degree $2k$ obtained via the
patchworking method and given by a Harnack polynomial of type $\HH^0$.
Let $\G$ be a face of the triangulation of $T_{2k}$
contained into the line $l_{2k-1}$ and
given by vertices $(c+1,d), (c,d+1)$ with $c$ and $d$ odd.
Let $U(p)$ be the $\e$-neighborhood of $\CCC \HH_m$ defined
from $\G^0$.
Denote $a$ the unique crossing
of $\HH_{2k-1} \cup L$
such that $a \in U(p)$.
\vskip0.1in
Then, a $3$-point $(p_{2k},p_{2k+1},p_{2k+2})$
of Harnack curves $\HH_{2k+j}$ ,$j \ge 3$,
is associated with the point $a$.
\been
\item
\lb{i:l21}
The point $p_{2k} \in S'_{2k} \subset S_{2k}^+$.
Besides, there exists a positive oval
${\mathcal O}$ such  that the point
$p_{2k} \in S_{2k}^+$ has the property:
$${\mathcal O} \not \in \pr_1W (p_{2k}),
{\mathcal O} \in \pr_0W(p_{2k})$$
Moreover,
given $B(p_{2k})$ a conj-equivariant $4$-ball around $p_{2k}$,

 $$\HH_{2k} \cap B(p_{2k}) \approx (\pr_0 W(p_{2k}) \cap B(p_{2k}))_1$$
\item
\lb{i:l22}
There exists a negative oval ${\mathcal  O}$ such  that
the point \\ $p_{2k+1} \in S_{2k+1}^+$ has the property:
$${\mathcal  O} \in \pr_0W(p_{2k+1}),
{\mathcal  O} \not \in \pr_1W(p_{2k+1}) $$
Moreover,
given $B(p_{2k+1})$ a conj-equivariant $4$-ball around $p_{2k+1}$,

 $$\HH_{2k+1} \cap B(p_{2k+1}) \approx
(\pr_0 W (p_{2k+1}) \cap B(p_{2k+1}))_1$$
\item
\lb{i:l23}
There exists a positive oval ${\mathcal  O}$ of $\HH_{2k+2}$ such  that
the point \\ $p_{2k+2} \in S_{2k+2}^-$ has the property:
$${\mathcal  O} \in \pr_1W (p_{2k+2}),
{\mathcal  O} \not \in \pr_0W(p_{2k+2})$$
Moreover,
given $B(p_{2k+2})$ a conj-equivariant $4$-ball around $p_{2k+2}$,

 $$\HH_{2k+2} \cap B(p_{2k+2}) \approx (\pr_1 W(p_{2k+2})
 \cap B(p_{2k+2}))_1$$
\enen }
\enle
\bede
We shall call {\it $3$-point of the first kind}
a $3$-point verifying all assumptions of
Lemma \mrf{l:l2p7}.
\end{defi}

\vskip0.1in
{\bf proof:}
\vskip0.1in

As in the proof of Lemma \mrf{l:l1p7},
we shall use local description of
Harnack curves $\HH_{m+j}$, $j \ge 0$,
inside $U(p)$ provided by the patchworking theory.

Recall that for any $j \ge 0$, there exists
an homeomorphism\\ $\tilde{h} :\CCC \HH_{m+j} \cap U(p) \to
 \{(x,y) \in (\CCC^*)^2
| b_{m+j}^S =0 \}\cap U(p)$
such that $\tilde{h}(a)=(x_0,y_0)$ $x_0.y_0<0$
is a critical point of $b_{m+j}^S$.

Our proof uses the results of chapter \ref{ch:Mope} and
in particular proposition \mrf{p:crit2},
and proposition \mrf{p:prop7b}.
We shall work with notations introduced in the chapter \ref{ch:Mope}.
\vskip0.1in
{\bf (\ref{i:l21})}
Let $\HH_{2k}$ the Harnack curve of degree $2k$ $(k>2)$
and $B_{2k}$ be a Harnack polynomial of degree $2k$.
According to proposition \mrf{p:prop7b},
a point  $p_{2k} \in S'_{2k}$ is associated with $a$.

Set $b_{2k}(p_{2k}) =c_0>0$,
then when $c=c_0$,
the non-empty positive oval touches a positive outer oval.
\vskip0.1in
{\bf (\ref{i:l22})}
Let $\HH _{2k+1}$ be the Harnack curve of
degree $2k+1$
and $B_{2k+1}$ be a Harnack polynomial of degree $2k+1$.

Consider the convex polygon $K \subset T_{2k+1}$ with vertices
$(c,d),(c,d+2),(c+1,d+2),(c+2,d),(c+2,d+1)$.
It is contained in the triangle $T_{2k+1}$
and triangulated by
the triangulation $\tau$ of $T_{2k+1}$.

Let $U(K^0)$ be the subset $\rho^{2k+1}(\RRR_+ K^0  \ti U_{\CCC}^2)$
of $\CCC P^2$. Obviously, $U(p) \subset U(K^0)$.

Denote $b_{2k+1}^K$ the truncation
of $b_{2k+1}$ on the monomials
$x^cy^d$, $x^cy^{d+1}$, $x^{c}y^{d+2}$,
 $x^{c+1}y^{d}$, $x^{c+1}y^{d+1}$,
$x^{c+1}y^{d+2}$, $x^{c+2}y^d$,
 $x^{c+2}y^{d+1}$.
(Namely,
$b_{2k+1}^K =b_{2k+1}^S
+ a_{c,d+2}x^cy^{d+2}
+ a_{c+1,d+2}x^{c+1}y^{d+2}
+ a_{c+2,d}x^{c+2}y^d
+ a_{c+2,d+1}x^{c+2}y^{d+1}$ with
$ a_{c,d+2} > 0, a_{c+1,d+2} > 0, a_{c+2,d} > 0, a_{c+2,d+1} > 0$.

According to the patchworking theory,
the truncation $b_{2k+1}^K$ of $b_{2k+1}$ is $\e$-sufficient
for $b_{2k+1}$ in $U(K^0)$.

Therefore,\\
$\tilde{h} :\CCC \HH_{2k+1} \cap U(p) \to
 \{(x,y) \in (\CCC^*)^2
| b_{2k+1}^S =0 \}\cap U(p)$
extends to the homeomorphism
$\tilde{h} :\CCC \HH_{2k+1} \cap U(K^0) \to
 \{(x,y) \in (\CCC^*)^2
| b_{2k+1}^K =0 \}\cap U(K^0)$
 such that $\tilde{h}(a)=(x_0,y_0)$ $x_0.y_0<0$.

According to proposition \mrf{p:prop7}, one can assume that
$(x_0,-y_0)$, $x_0.-y_0 >0$, $x_0<0,-y_0<0$,
is a critical point of
$b_{2k+1}^S=b_{2k}^S$ with positive critical value.

Thus, up to modify coefficients
of the polynomials $b_{2k+1}^K$
(without  changing the order and
the topological structure of $b_{2k+1}$),
the point
$(x_0,-y_0)$ is a critical point of $b_{2k+1}^K$ with
$b_{2k+1}(x_0,-y_0) >0$.

Hence, we shall proceed to the study of
the Petrovskii's pencil of curves
over $\HH_{2k+1}$.
Let $d_{2k+1}$ be the unique  polynomial such that
$b_{2k+1}= b_{2k+1}^K + d_{2k+1}$.
Consider curves of the Petrovskii's pencil over $\HH_{2k+1}$
outside $d_{2k+1} =0$,
as level curves of the function
$\frac {b_{2k+1}^K -c} {d_{2k+1}}$.

Inside $b_{2k+1}^K-c=0 \bk d_{2k+1}=0$,
these curves have critical points the singular
points of $b_{2k+1}-c =0$.

Bringing together Petrovskii's theory and
the implicit function theorem applied to one-parameter
polynomial
$x_0^{2k+1}.(b_{2k+1} -c)$
with parameter $c$, it  follows
that a point $p_{2k+1}$ with
$b_{2k+1}(p_{2k+1})=c_0 >0$
is associated with $a$.

When $c=c_0$, the one-side component
of the curve
$x_0^{2k+1}.(b_{2k+1} -(c_0)$ touches itself.
Besides, according to Petrovskii's Lemma 2 (3.a) and Lemma 3
as $c$ increases from $c_0$ to $c_0 + \e$
one negative oval ${\mathcal  O}^-$
  (of the curve
$x_0^{2k+1}.(b_{2k+1} -(c_0-\e)=0$) disappears.
(In the pencil of curves,
 ${\mathcal  O}^-$ as oval of $\HH_{2k}$
has been created in the construction of $\HH_{2k+1}$ from
$\HH_{2k}  \cup L$)
\vskip0.1in
{\bf (\ref{i:l23})}
Let $\HH_{2k+2}$ be the Harnack curve of degree $2k+2$
and $B_{2k+2}$ be a Harnack polynomial of degree $2k+2$.

Consider $J$ the square with vertices
$(c,d),(c,d+2),(c+2,d),(c+2,d+2)$.
It is contained in the triangle $T_{2k+2}$
and triangulated by
the triangulation $\tau$ of $T_{2k+2}$.
Let $U(J^0)$ be the  subset $\rho^{2k+2}(\RRR_+ J^0  \ti U_{\CCC}^2)$
of $\CCC P^2$. Obviously, $U(p) \subset U(K^0) \subset U(J^0)$.

Denote
$b_{2k+2}^J$ the truncation of $b_{2k+2}$ on the monomials
$x^cy^d$, $x^cy^{d+1}$, $x^cy^{d+2}$,
$x^{c+1}y^d$, $x^{c+1}y^{d+1}$, $x^{c+1}y^{d+2}$,
$x^{c+2}y^d$, $x^{c+2}y^{d+1}$,$x^{c+2}y^{d+2}$.
(Namely,
$b_{2k+1}^J= b_{2k+1}^K+
+ a_{c+2,d+2}x^{c+2}y^{d+2}$ with $a_{c+2,d+2} >0$,
$b_{2k+1}^J= a_{c,d}x^{c,d}+ a_{c+1,d}x^{c+1}y^d+
+ a_{c,d+1}x^{c,d+1}- a_{c+1,d}x^{c+1}y^d+
+ a_{c,d+2}x^cy^{d+2}
+ a_{c+1,d+2}x^{c+1}y^{d+2}
+ a_{c+2,d}x^{c+2}y^d
+ a_{c+2,d+1}x^{c+2}y^{d+1}$ with
$ a_{c,d+2} > 0, a_{c+1,d+2} > 0, a_{c+2,d} > 0, a_{c+2,d+1} > 0$.

According to the patchworking theory,
the truncation $b_{2k+2}^J$ of $b_{2k+2}$ is $\e$-sufficient
for $b_{2k+2}$ in $U(J^0)$.
Therefore,
$\tilde{h} :\CCC \HH_{2k+2} \cap U(p) \to
 \{(x,y) \in (\CCC^*)^2
| b_{2k+2}^S =0 \}\cap U(p)$ extends to the homeomorphism
$\tilde{h} :\CCC \HH_{2k+2} \cap U(J^0) \to
 \{(x,y) \in (\CCC^*)^2
| b_{2k+2}^J =0 \}\cap U(J^0)$
 such that $\tilde{h}(a)=(x_0,y_0)$ $x_0.y_0<0$.

It is easy to see that, up to modify coefficients
of the polynomials $b_{2k+2}^J$,
(without  changing the order and
the topological structure of $b_{2k+2}$),
the point
$(x_0,-y_0)$, $x_0<0, -y_0 <0$
is a critical point of $b_{2k+2}^J$
with $b_{2k+2} (x_0,-y_0) <0$.

Hence, we shall proceed to the study of
the Petrovskii's pencil of curves
over $\HH_{2k+2}$.

Denote
$d_{2k+2}$  the unique  polynomial such that
$b_{2k+2}= b_{2k+2}^J + d_{2k+2}$.

As previously,
consider curves of the Petrovskii's pencil over $\HH_{2k+2}$
outside $d_{2k+2} =0$,
as level curves of the function
$\frac {b_{2k+2}^J -c} {d_{2k+2}}$.

Inside $(b_{2k+2}^J-c=0) \bk(d_{2k+2}=0)$,
these curves have critical points the singular
points of $b_{2k+2}^J-c =0$.
Bringing together Petrovskii's theory,
the implicit function theorem applied to one-parameter
polynomial
$x_0^{2k+2}.(b_{2k+2} -c)=
x_0^{2k}.(b_{2k+2}^J+d_{2k+2} -c)$
with parameter $c$, it  follows
that a point $p_{2k+2}$ with
$b_{2k+2}(p_{2k+2})=c_0 <0$
is associated with $a$.

Let $b_{2k+2}(p_{2k+2})=c_0$.
When $c=c_0$,
one positive  oval touches another positive oval.
Then, as $c$ decreases from $c_0$ to $c_0 -\e$,
a positive oval disappears.
\vskip0.1in
Let $\HH_{2k+j}$ be the Harnack curve of
degree $2k+j$, $j \ge 3$ and let
$B_{2k+j}$ be a Harnack polynomial of degree $m+j$

As previously, consider
$b_{2k+j}^J$ the truncation of $b_{2k+j}$ on the monomials
$x^cy^d$, $x^cy^{d+1}$, $x^{c+1}y^d$, $x^{c+1}y^{d+1}$.

According to patchworking theory,
the truncation $b_{2k+j}^K$ of $b_{2k+j}$ is $\e$-sufficient
for $b_{2k+j}$ in $U(J^0)$.

Therefore,\\
$\tilde{h} :\CCC \HH_{2k+j} \cap U(p) \to
 \{(x,y) \in (\CCC^*)^2
| b_{2k+j}^S =0 \}\cap U(p)$ extends to the homeomorphism \\
$\tilde{h} :\CCC \HH_{2k+j} \cap U(K^0) \to
 \{(x,y) \in (\CCC^*)^2
| b_{2k+1}^K =0 \}\cap U(K^0)$
 such that $\tilde{h}(a)=(x_0,y_0)$ $x_0.y_0>0$.

Let $d_{2k+j}$ be the unique polynomial such that
$b_{2k+j}= b_{2k+j}^J+d_{2k+j}$.
Consider curves of the Petrovskii's pencil over $\HH_{2k+j}$
outside $d_{2k+j} =0$,
as level curves of the function
$\frac {b_{2k+j}^J -c} {d_{2k+j}}$.
From an argumentation similar to the previous one,
it follows that
a point $p_{2k+j}$ of $S^-_{2k+j}$ equivalent to $p_{2k+2}$
is associated with $a$.

Therefore, a {\it 3-point} of curves $\HH_{2k+j}$, $j \ge 3$,
is associated with $a$.

Q.E.D

Recall (see propositiony \mrf{p:prop7b})
that given $B_{2k}= x_0^{2k}.b_{2k}(x_1/x_0,x_2/x_0)$
a Harnack polynomial of degree $2k$ and
type $\HH^0$,
$S"_{2k}$ denotes $S^+_{2k} \bk S'_{2k}$
the complementary set of  $S'_{2k}$ inside $S^+_{2k}$.
In the Lemma \mrf{l:l3p7}, we shall study how critical points $S"_{2k}$
of a Harnack polynomial of type $\HH^0$ and degree $2k$ vary
in the recursive  construction of Harnack curves.

\bele
\mlb{l:l3p7}
\vskip0.1in
{\it
Let $\HH_{2k}$ be the Harnack curve of degree $2k$ obtained via the
patchworking method and given by a Harnack polynomial of type $\HH^0$.
Let $\G$ be a face of the triangulation of $T_{2k}$
contained into the line $l_{2k-1}$ and
given by vertices $(c+1,d), (c,d+1)$ with $c$ and $d$
even and strictly positive.
Let $U(p)$ be the
$\e$-neighborhood of $\CCC \HH_m$ defined from $\G^0$.
Denote $a$ the unique crossing
of $\HH_{2k-1} \cup L$
such that $a \in U(p)$.
\vskip0.1in
Then, a $3$-point $(p_{2k},p_{2k+1},p_{2k+2})$
of Harnack curves $\HH_{2k+j}$, $j \ge 3$, is associated with $a$.
\been
\item
\lb{i:l31}
There exists a negative oval ${\mathcal  O}$ such  that
the point\\ $p_{2k} \in S"_{2k} \subset S_{2k}^+$ has the property:
${\mathcal  O} \in \pr_0W (p_{2k})$
${\mathcal  O} \not \in \pr_1W(p_{2k})$.\\
Moreover,
given $B(p_{2k})$ a conj-equivariant $4$-ball around $p_{2k}$,

$$\HH_{2k} \cap B(p_{2k}) \approx (\pr_0 W(p_{2k}) \cap B(p_{2k}))_1$$
\item
\lb{i:l32}

There exists a negative oval ${\mathcal  O}$ such  that
the point\\ $p_{2k+1} \in S_{2k+1}^+$ has the property:
${\mathcal  O} \in \pr_0W(p_{2k+1})$
${\mathcal  O} \not \in \pr_1W(p_{2k+1}) $.\\
Moreover,
given $B(p_{2k+1})$ a conj-equivariant $4$-ball around $p_{2k+1}$,

$$\HH_{2k+1} \cap B(p_{2k+1}) \approx (\pr_0 W(p_{2k+1}) \cap B(p_{2k+1}))_1$$
\item
\lb{i:l33}

There exists a negative oval ${\mathcal  O}$ such  that
the point\\ $p_{2k+2} \in S_{2k+2}^+$ has the property:
${\mathcal  O} \in \pr_0W(p_{2k+2}) $
${\mathcal  O} \not \in \pr_1W(p_{2k+2}) $.\\
Moreover,
given $B(p_{2k+2})$ a conj-equivariant $4$-ball around $p_{2k+2}$,

 $$\HH_{2k+2} \cap B(p_{2k+2}) \approx
(\pr_0 W (p_{2k+2}) \cap B(p_{2k+2}))_1 $$
\enen}
\enle
\bede
We shall call {\it $3$-point of the second kind}
a $3$-point verifying all assumptions of Lemma \mrf{l:l3p7}.
\end{defi}
\vskip0.1in
{\bf proof}
\vskip0.1in
Our proof is based on arguments
similar to the one given in the proof of Lemma \mrf{l:l2p7}.

We shall use local description of Harnack curves $\HH_{m+j}$ $j \ge 0$
inside $U(p)$ provided by the patchworking theory.

Recall that for any $j \ge 0$, there exists
an homeomorphism \\ $\tilde{h} :\CCC \HH_{m+j} \cap U(p) \to
 \{(x,y) \in (\CCC^*)^2
| b_{m+j}^S =0 \}\cap U(p)$
such that $\tilde{h}(a)=(x_0,y_0)$ is a critical point of $b_{m+j}^S$.
\vskip0.1in
{\bf (\ref{i:l31})}
Let $\HH_{2k}$ be the Harnack curve of degree $2k$, $(k>3)$
and $B_{2k}$ be Harnack polynomial of degree $2k$.

From proposition \mrf{p:prop7b},
a point $p_{2k}$  of $S^+_{2k}$ is associated with $a$.
Set $b_{2k}(p_{2k}) =c_0$.
Then, when $c=c_0$,
the non-empty positive oval touches a negative (inner) oval;
as $c$  increases from $c_0$
to $c_0 + \e$,  one negative oval
${\mathcal  O}^-$ disappears.
(In the pencil of curves,
${\mathcal  O}^-$ as oval of $\HH_{2k-1}$ has
been created in the construction of $\HH_{2k-1}$ from $\HH_{2k-2}
\cup L$).

{\bf (\ref{i:l32})}\\
Let $\HH _{2k+1}$ be the Harnack curve of
degree $2k+1$ and $B_{2k+1}$ be a Harnack polynomial of degree $2k+1$.

Consider the convex polygon $K \subset T_{2k+1}$ with vertices
$(c,d),(c,d+2),(c+1,d+2),(c+2,d),(c+2,d+1)$.
It is contained in triangle $T_{2k+1}$
and triangulated by
the triangulation $\tau$ of $T_{2k+1}$.

Let $U(K^0)$ be the subset $\rho^{2k+1}(\RRR_+ K^0  \ti U_{\CCC}^2)$
of $\CCC P^2$. Obviously, $U(p) \subset U(K^0)$.

Denote $b_{2k+1}^K$ the truncation
of $b_{2k+1}$ on the monomials
$x^cy^d$, $x^cy^{d+1}$, $x^{c}y^{d+2}$,
 $x^{c+1}y^{d}$, $x^{c+1}y^{d+1}$,
$x^{c+1}y^{d+2}$, $x^{c+2}y^d$,
 $x^{c+2}y^{d+1}$.

According to the patchworking theory,
the truncation $b_{2k+1}^K$ of $b_{2k+1}$ is $\e$-sufficient
for $b_{2k+1}$ in $U(K^0)$.

Therefore,
$\tilde{h} :\CCC \HH_{2k+1} \cap U(p) \to
 \{(x,y) \in (\CCC^*)^2
| b_{2k+1}^S =0 \}\cap U(p)$
extends to the homeomorphism
$\tilde{h} :\CCC \HH_{2k+1} \cap U(K^0) \to
 \{(x,y) \in (\CCC^*)^2
| b_{2k+1}^K =0 \}\cap U(K^0)$
 such that $\tilde{h}(a)=(x_0,y_0)$ $x_0.y_0<0$.
According to proposition \mrf{p:prop7},  one can assume that
$(x_0,-y_0)$ is a critical point of $b_{2k+1}^S=b_{2k}^S$.
It is easy to see that, up to modify coefficients
of the polynomials $b_{2k+1}^K$
(without  changing the order and
the topological structure of $b_{2k+1}$),
the point
$(x_0,-y_0)$ $x_0>0,-y_0>0$
is a critical point of $b_{2k+1}^K$ with
$b_{2k+1}(x_0,-y_0) >0$.

Let $d_{2k+1}$ be the unique  polynomial such that
$b_{2k+1}= b_{2k+1}^K + d_{2k+1}$.
Consider curves of the Petrovskii's pencil over $\HH_{2k+1}$
outside $d_{2k+1} =0$,
as level curves of the function
$\frac {b_{2k+1}^K -c} {d_{2k+1}}$.

Inside $(b_{2k+1}^K-c=0) \bk (d_{2k+1}=0)$,
these curves have critical points the singular
points of $b_{2k+1}-c =0$.

Bringing together Petrovskii's theory and
the implicit function theorem applied to one-parameter
polynomial
$x_0^{2k+1}.(b_{2k+1} -c)$
with parameter $c$, it  follows
that a point $p_{2k+1}$ with
$b_{2k+1}(p_{2k+1})=c_0 >0$
is associated with $a$.

When $c=c_0$, the one-side component
of the curve touches itself.
Then,
as $c$ increases from $c_0$ to $c_0 + \e$
 one negative oval ${\mathcal  O}^-$ disappears.
(In the pencil of curves,
  ${\mathcal  O}^-$ as oval of $\HH_{2k-1}$ has
been created in the construction of $\HH_{2k-1}$ from $\HH_{2k-2}
\cup L$).
\vskip0.1in

{\bf (\ref{i:l33})}
Let $\HH_{2k+2}$ be the Harnack curve of degree $2k+2$ and
$B_{2k+2}$ be the Harnack polynomial of degree $2k+2$

Consider $J$ the square with vertices
$(c,d),(c,d+2),(c+2,d),(c+2,d+2)$.
It is contained in the triangle $T_{2k+2}$
and triangulated by
the triangulation $\tau$ of $T_{2k+2}$.

Let $U(J^0)$ be the  subset $\rho^{2k+2}(\RRR_+ J^0  \ti U_{\CCC}^2)$
of $\CCC P^2$. Obviously, $U(p) \subset U(K^0) \subset U(J^0)$.

Denote
$b_{2k+2}^J$ the truncation of $b_{2k+2}$ on the monomials
$x^cy^d$, $x^cy^{d+1}$, $x^cy^{d+2}$,
$x^{c+1}y^d$, $x^{c+1}y^{d+1}$, $x^{c+1}y^{d+2}$,
$x^{c+2}y^d$, $x^{c+2}y^{d+1}$,$x^{c+2}y^{d+2}$.

According to  patchworking theory,
the truncation $b_{2k+2}^J$ of $b_{2k+2}$ is $\e$-sufficient
for $b_{2k+2}$ in $U(J^0)$.

Therefore,\quad \break
$\tilde{h} :\CCC \HH_{2k+2} \cap U(p) \to
 \{(x,y) \in (\CCC^*)^2
| b_{2k+2}^S =0 \}\cap U(p)$
extends to the homeomorphism
$\tilde{h} :\CCC \HH_{2k+2} \cap U(J^0) \to
 \{(x,y) \in (\CCC^*)^2
| b_{2k+2}^J =0 \}\cap U(J^0)$
 such that $\tilde{h}(a)=(x_0,y_0)$ $x_0.y_0<0$.

It is easy to see that up to modify coefficients
of the polynomials $b_{2k+2}^J$
(without  changing the order and
the topological structure of $b_{2k+2}$),
the point
$(x_0,-y_0)$ $x_0>0, y_0 <0$ is a critical point of $b_{2k+2}^J$ with
$b_{2k+2}^J(x_0,-y_0) >0$.

 According to the Petrovskii's theory,
from an argumentation similar to the previous one,
it  follows
that a point $p_{2k+2}$ with
$b_{2k+2}(p_{2k+2})=c_0 >0$
is associated with $a$.

When $c=c_0$,
the negative oval ${\mathcal  O}^-$ touches another negative oval.
(In the pencil of curves, this last negative oval as oval of $\HH_{2k+1}$
is, in the patchworking scheme, situated in front of ${\mathcal  O}^-$ and
has been created in the construction
of $\HH_{2k+1}$ from $\HH_{2k} \cup L$).

Then,
as $c$ increases from $c_0$ to
$c_0 + \e$, one negative oval disappears.

\vskip0.1in
Let $j \ge 3$ be an integer and $\HH_{2k+j}$ be the Harnack curve of
degree $2k+j$.

Let $B_{2k+j}$ be a Harnack polynomial of degree $m+j$.
As previously, consider
$b_{2k+j}^J$ the truncation of $b_{2k+j}$ on the monomials
$x^cy^d$, $x^cy^{d+1}$, $x^{c+1}y^d$, $x^{c+1}y^{d+1}$.
According to the patchworking theory,
the truncation $b_{2k+j}^K$ of $b_{2k+j}$ is $\e$-sufficient
for $b_{2k+j}$ in $U(J^0)$.
Therefore,
$\tilde{h} :\CCC \HH_{2k+j} \cap U(p) \to
 \{(x,y) \in (\CCC^*)^2
| b_{2k+j}^S =0 \}\cap U(p)$.
extends to the homeomorphism
$\tilde{h} :\CCC \HH_{2k+j} \cap U(K^0) \to
 \{(x,y) \in (\CCC^*)^2
| b_{2k+1}^K =0 \}\cap U(K^0)$
 such that $\tilde{h}(a)=(x_0,y_0)$ $x_0.y_0<0$.

From an argumentation similar to the previous one,
it follows that
a point $p_{2k+j}$ of $S^+_{2k+j}$ equivalent to $p_{2k+2}$
is associated with $a$.

Therefore, a {\it 3-point} is associated with $a$.
Q.E.D
\vskip0.1in

\bere
Given $\HH_{2k}$ the Harnack of degree $2k$.
It follows immediately from the Lemmas
\mrf{l:l1p7}, \mrf{l:l2p7}, \mrf{l:l3p7}
above that
$(k-2)^2$ $3$-points are associated
to crossing of curves $\HH_{2k-3-j} \cup L$
$(L \approx \CCC l_{2k-3 -j})$  $(0 \le j \le (2k-7))$.
\vskip0.1in
Likewise,
given $\HH_{2k+1}$ the Harnack curve  of odd degree $2k+1$,
$(k-2)^2$  $3$-points are associated
to crossing of curves $\HH_{2k-3-j} \cup L$
$(L \approx \CCC l_{2k-3 -j})$  $(0 \le j \le (2k-7))$,
and $(2k-3)$ $2$-points
are associated with crossing points of $\HH_{2k-1} \cup L$.
Hence,
for any Harnack curve  $\HH_m$ of degree $m$
with polynomial $B_m$,
one can recover from the Lemmas above
local topological meaning of critical points
of index 1 of $B_m$.
\enre
\vskip0.1in

\chapter{ Perestroika theory on Harnack curves}
 {\bf -Construction of any curve of degree $m$ with non-empty real part-}
\vskip0.1in

In this chapter, we shall give a method of construction of
any curve of a given degree $m$ which provides
description of  pairs
 $(\CCC P^2, \CCC {\AA}_m)$
up to conj-equivariant isotopy.
In a few words, our method of construction of curves of degree $m$
is based on the definition of a chain of modification on the real
connected components of the Harnack curve $\HH_m$.\\

We shall divide this section into three sections.
Our method of construction of algebraic curve
uses invariants of real algebraic curves similar to the Arnold's invariants
of generic immersion of the circle into the plane.
These invariants
and their analogous for algebraic curves
were respectively introduced in \cite{Ar2} and \cite{Vi2}.
In the first section, we state the problem  of construction of algebraic
curves from the pair
$(\CCC P^2, \CCC \HH_m)$ and recall definitions of these invariants.
Classification of pairs $(\CCC P^2, \CCC {\AA}_m)$
leads us to distinguish
curves with orientable set of real points usually called
curves of type $I$,
and curves  with non-orientable set of real points
usually called curves of type $II$.
In the second section, we state a method
of construction of curves of type $I$.
Then, in the third section, we enlarge the method to
curves of type $II$.
The main result of this section is gathered in
Theorem \mrf{t:th2typ1}
and Theorem \mrf{t:th1typ2}.
This statement (Theorem \mrf{t:th2typ1} and \mrf{t:th1typ2})
is the counterpart
of Theorem \mrf{t:theo1} established for Harnack curves.
Namely, we present any curve $\AA_m$ of degree $m$ with non-empty real part
up to conj-equivariant isotopy of $\CCC P^2$ as follows:
Outside a finite number of $4$-balls $B(a_i)$ globally invariant by
complex conjugation, $\AA_m$ is the union of $m$ non-intersecting projective
lines; inside any $4$-ball $B(a_i)$ it is the perturbation
of a crossing.

\vskip0.1in
\section{ Introduction-Classical Topological Facts-}
\mlb{su:ClasTop}
\vskip0.1in

In any cases, (curves  ${\AA}_m$ with an arbitrary
number of real components and arbitrary type),
we shall bring the problem of description of $(\CCC P^2,\CCC {\AA}_m)$
to the problem of description
of all possible moves of the real components of $\HH_m$.
Let us  recall in introduction
general facts (see \cite{Vi} for example)
which will state the space of moves
providing an arbitrary curve of degree $m$
from the Harnack curve $\HH_m$.

\subsection{Introduction}
\mlb{susu:IntPer}

Passing from polynomials to real set of points, it follows easily
that real algebraic curves form a real projective space of dimension
$\frac {m (m+3)} 2$.
We shall denote this space by the symbol $\RRR {\mathcal C}_m$.
It easy to see that one can trace the real part of a curve of degree $m$
through any
$\frac {m (m+3)} 2$ real points, and that it is uniquely defined for
points in general position.
Moreover,
real algebraic curves can be considered as complex curves of special kind.
Similarly, passing from polynomials to complex set of points,
it follows easily
that complex algebraic curves of degree $m$ form a complex projective space
of dimension $\frac {m (m+3)} 2$.
We shall denote this space by $\CCC {\mathcal C}_m$.
Obviously, $\RRR {\mathcal C}_m$ coincides with the real part of
$\CCC {\mathcal C}_m$.
It is easy to see that through any $\frac {m (m+3)} 2$ points
one can always trace a curve of degree $m$.
Moreover, the set of $\frac {m (m+3)} 2$-points for which
such a curve is unique is open and dense in the space of all
 $\frac {m (m+3)} 2$-points subset of $\CCC P^2$.
\vskip0.1in

Let $\RRR {\mathcal D}_m$ denote the subset of
$\RRR {\mathcal C}_m$ corresponding
to real singular curves.
Let $\CCC {\mathcal D}_m$ denote the subset of
$\CCC {\mathcal C}_m$ corresponding
to singular curves.

We call a  path in the complement
$\RRR {\mathcal C}_m \bk \RRR {\mathcal D}_m$
of the discriminant hypersurface in $\RRR {\mathcal C}_m$
a {\it rigid isotopy} of real points set
of nonsingular curves of degree $m$.
We call a smooth path in the complement
$\CCC {\mathcal C}_m \bk \CCC {\mathcal D}_m$
of the discriminant hypersurface in $\CCC {\mathcal C}_m$
a {\it rigid isotopy} of complex point set
of nonsingular curves of degree $m$.
\vskip0.1in
These definitions give  rise naturally
to the classification problem of
non-singular curves of degree $m$ up to rigid isotopy.
\vskip0.1in
From the complex viewpoint, the rigid isotopy classification problem has
trivial solution:
the complex point sets of any two nonsingular
curves of degree $m$ are rigidly isotopic.
Moreover, they are diffeotopic in $\CCC P^2$.
From the real viewpoint, even if we consider only rigid isotopy,
this property can not be extended.
\vskip0.1in
\vskip0.1in

\subsection{State of the problem}
\mlb{su:statepb}
\vskip0.1in
From now on, one can state
some properties of the sets of points $\CCC {\AA}_m$,
$\RRR {\AA}_m$  of a  curve ${\AA}_m$.
Let $\HH_m$ be the Harnack  curve of degree $m$ and ${\AA}_m$ be a curve
of degree $m$ such that the pair $(\RRR P^2, \RRR {\AA}_m)$
is non-homeomorphic to the pair $(\RRR P^2, \RRR \HH_m)$.
From the elementary topological facts recalled previously, there exists
a smooth path $h:[0,1] \to \CCC {\mathcal C}_m \bk \CCC {\mathcal D}_m$ with
$h(0) = \CCC \HH_m, h(1) =\CCC {\AA}_m$ such that the family $\CCC h(t)$
is a diffeotopy of submanifolds of $\CCC P^2$.
Thus, if there exists
$B \subset \CCC P^2$  globally invariant by
complex conjugation such that
$h_t (\CCC \HH_m) \cap B)  \subset B$, $t \in [0,1]$,
topological spaces $(B \cap \HH_m)$ and
$(B \cap {\AA}_m)$ are isotopic in $B$.

According to corollary \mrf{c:gpa}, assume $\HH_m$
obtained via $T$-inductive construction of Harnack curves.
From the Lemma \mrf{l:l1p7}, Lemma \mrf{l:l2p7} and Lemma \mrf{l:l3p7}
of the previous section,
any critical point $(x_0,y_0)$ of index 1 of $B_m$
is associated with a crossing $a$ of one curve
$\HH_{m-i} \cup L$.
Namely,  there exists
$U(p)=\rho^m(D(p,\e) \ti U_{\CCC}^2) \subset(\CCC^*)^2$
$\e$-neighborhood
of $\CCC \HH_{m}$ defined the interior $\G^0$
of a face $\G$ of the triangulation of $T_m$
such that :
\been
\item
$a, (x_0,y_0) \in U(p)$
\item
the perturbation on the real part of $\HH_m$ involved in the local
topological meaning of $(x_0,y_0)$ is a deformation
on $\CCC \HH_{m} \cap U(p)$
\enen
Let us denote by ${\mathcal P}_m$ the set of points $p$ of $T_m$
with the property that
for any $p \in {\mathcal P}_m$,
the $\e$-neighborhood $U(p)$ verifies the two properties above.
Note that
for any two distinct points
$p , p'$ the subsets $U(p)$ and $U(p')$ have empty-intersection.
\vskip0.1in
It easy to see
that the set ${\mathcal P}_m$ contains $\frac {m(m-1)}2$ points
of $T_m$.\\
Moreover, for any $p \in {\mathcal P}_m$,
$U(p) =\rho^m(D(p,\e) \ti U_{\CCC}^2)$
contains conj-equivariant $4$-balls $B(x_0,y_0)$ and $B(a)$
around $(x_0,y_0)$ and $a$ such that\\
$\CCC \HH_m \cap B((x_0,y_0)) \not= \emptyset$
,$\CCC \HH_m \cap B(a) \not= \emptyset$.
\vskip0.1in
Thus, according to
the density of $\frac {m (m+3)} 2$-uple of points which define
uniquely  a curve in the space of $\frac {m (m+3)} 2$-uple
of points of $\CCC P^2$, if there exists a
diffeotopy of submanifolds of $\CCC P^2$
$h :[0,1] \to \CCC {\mathcal C}_m \bk \CCC {\mathcal D}_m$ with
$h(0) = \CCC \HH_m, h(1) =\CCC {\AA}_m$
such that
$h_t (\CCC \HH_m   \cap U(p))  \subset U(p)$, $t \in [0,1]$,
the curve ${\AA}_m$ is entirely defined from its restriction
on  the union
$\cup_{\{p \in{\mathcal P}_m\}} U(p)$
of $\e$-neighborhood $U(p)$
taken over the $\frac {m(m-1)}2$ points of ${\mathcal P}_m$.

In such a way, we bring
the problem of construction of a curve ${\AA}_m$
to the definition of a  path
$h:[0,1] \to \CCC {\mathcal C}_m \bk \CCC {\mathcal D}_m$ with
$h(0) = \CCC \HH_m, h(1) =\CCC {\AA}_m$.
\vskip0.1in
\subsection{Varieties of irreducible curves of
type $I$, degree $m$, and genus $g$}
\mlb{susu:irred,I}

As already introduced for smooth
curves, curves with orientable real set of points
are called curves of type $I$.

In what follows, we shall consider smooth curves of type $I$ in the
broader set of irreducible curves of type $I$.
Any irreducible curve of type $I$
is such that
the real part of its normalization divides the set
of complex points of its normalization into two halves.
We shall call the images of the halves of the normalization
in the set of points of the curve, the {\it halves } of the curve.
Each of the halves is oriented and induces an orientation
on the real part as on its boundary.
We shall consider irreducible curves of
degree $m$, genus $g$
and type $I$
with a distinguished complex orientation.
Curves of this kind constitute a finite-dimensional
stratified real algebraic variety
we shall denote ${\mathcal C}_{I; g }$.

We shall set
${\mathcal C}_{I,m}= \cup_{ 0 \le g \le \frac {(m-1)(m-2)} 2}
{\mathcal C}_{I;g}$
 the subset of ${\mathcal C}_m$ constituted by all irreducible curves of
degree $m$, genus $g$, $0 \le g \le \frac {(m-1)(m-2)} 2$,
and type $I$.

Passing from polynomials to real (resp, complex)
set of points, it follows easily
that varieties ${\mathcal C}_{I;g}$ define subspaces
$\RRR {\mathcal C}_{I;g}$ (resp,$\CCC C_{I;g}$)
of the space $\RRR {\mathcal C}_m$
(resp, $\CCC {\mathcal C}_m$.)

We shall set
$\RRR {\mathcal C}_{I,m}= \cup_{ 0\le g \le \frac {(m-1)(m-2)} 2}
\RRR {\mathcal C}_{I;g}$
the subset of
$\RRR {\mathcal C}_m$ constituted by real set of points of irreducible
curves of degree $m$ and type $I$,
$\CCC {\mathcal C}_{I,m}= \cup_{ 0\le g \le \frac {(m-1)(m-2)} 2}
\CCC {\mathcal C}_{I;g}$
the subset of
$\CCC {\mathcal C}_m$ constituted by  complex set of points of irreducible
curves of degree $m$ and type $I$.
Obviously, $\RRR {\mathcal C}_{I,m}$ coincides with the real part of
$\CCC {\mathcal C}_{I,m}$.
\vskip0.1in
Let us study the set of singular algebraic curves
in the variety ${\mathcal C}_{I;g}$
of curves of degree $m$, genus $g$ and type $I$.
\vskip0.1in
{\it Generic curves}
\vskip0.1in
Following usual terminology, we call
{\it ordinary double point singularity} of a real
irreducible algebraic curve ${\AA}_m$ a non-degenerate singular point of
${\AA}_m$.
All ordinary double points singularity
are equivalent from the complex viewpoint.
From the real viewpoint, one distinguishes several types of such points.
\been
\item
real double point of intersection of two real branches called crossing
\item
real point of intersection of two complex branches
conjugated to each other called {\it solitary double point}
\item
imaginary double point of intersection of the different halves of
$\CCC {\AA}_m \bk \RRR {\AA}_m$
denoted $\a$-point.
\item
imaginary double point of self-intersection of one of the halves
of $\CCC {\AA}_m \bk \RRR {\AA}_m$
denoted $\b$-point.
\enen

Define a {\it generic curve } as a real irreducible
algebraic curve with only ordinary double singularities.
As it is well known, generic curves in the variety ${\mathcal C}_{I;g}$
constitute a
Zarisky open set in the variety ${\mathcal C}_{I;g}$.

For a generic curve ${\AA}_m$  of degree $m$ and  genus $g$
,$0 < g < \frac {(m-1)(m-2)} 2$,
by {\it smoothing} of its real point set
$\RRR {\AA}_m$
we shall understand a smooth oriented $1$-dimensional submanifold
of $\RRR P^2$ obtained from
$\RRR {\AA}_m$ by the modification at each double point determined
by the complex orientation (see figure 5.1).

\vskip0.1in
By smoothing of its complex point set
we shall understand a smooth oriented $1$-dimensional complex submanifold
of $\CCC P^2$ obtained from $\CCC {\AA}_m$
by the modification at
each real double point as above and
by the modification at
each complex double point.
The two complex branches
which merge in  $\a$-point (resp  $\b$-point) become after smoothing
two branches of different halves of $\CCC {\AA}_m \bk \RRR {\AA}_m$
(resp,
two branches of one of the  halves of $\CCC {\AA}_m \bk \RRR {\AA}_m$)
\vskip0.1in
We shall call {\it smoothing}
of a generic curve ${\AA}_m$,
the smooth $1$-dimensional complex submanifold of $\CCC P^2$
deduced from ${\AA}_m$ by smoothing of
$\RRR {\AA}_m$ and $\CCC {\AA}_m$ as described above.
We shall say that the singular points of ${\AA}_m$ are smoothed.
\vskip0.1in
{\it Discriminant Hypersurface}
\vskip0.1in
Consider now the complement of the set of generic curves of all real
algebraic curves in the variety ${\mathcal C}_{I; g}$.
It can  be considered as a
discriminant hypersurface of which strata consist of curves with only
one singular point which is not an ordinary double point.
One can distinguish six main strata
defined by the type of the singularity of the curves it contains:
\been
\item real cusp,
\item real point of direct ordinary tangency,
\item real point of inverse ordinary tangency,
\item real point of ordinary tangency of two imaginary branches,
\item real ordinary triple point of intersection of three real branches,
\item real ordinary triple point of intersection of a real branch and
two conjugate imaginary branches.
\enen
\vskip0.1in
{\it Perestroika}
\vskip0.1in
A {\it generic  path}  in the variety
${\mathcal C}_{I; g}$
intersects the discriminant hypersurface in
a finite number of points, and these points belong to the main strata.
\vskip0.1in
We call {\it perestroika} a change
experienced by a generic path in the space
${\mathcal C}_{I; g}$ when it goes through
the main strata.
The notion of perestroika was initially
introduced  in the context of
generic immersion of the circle into the plane
(see for example  \cite{Ar2}, \cite{Vi2}).

Given a perestroika,
we call {\it smoothed perestroika} the change obtained by
smoothing the fragments involved in the perestroika.
\vskip0.1in
Moreover, to define how a generic path crosses the strata
of the discriminant hypersurface
one has to specify a coorientation of the strata.
Let us define the {\it positive direction} of a perestroika.
The opposite direction is naturally called {\it negative direction}.
\been
\item
In the case of cusp, define the {\it positive direction}
the direction from curve with one more crossing point  to curve
with one more solitary double point.
\item
In  case of real point of tangency, there is a natural
{\it positive direction}
from curve with less  real double points
to curve with more double points.
\item
In case of triple-point, the coorientation of the stratum can
be defined as follows.
The crossing of a triple point by a path gives rise
to the {\it vanishing triangle}: the dying triangle
which existed
just before the crossing and the new born
existed just after it.
A cyclic order on the sides of the triangle is given by the
order the sides are visited.
In the case of triple-point, the coorientation rule assigns {\it signs}
to triangle.
Let $q$ be the number of the sides of the vanishing triangle whose directions
coincide with that given by the cyclical order.
The sign  of a triangle is $(-1)^q$.
The crossing of the triple point is {\it positive} (resp, {\it negative})
if the newborn triangle is positive (resp, {\it negative})
(and hence the dying one is negative (resp, positive)).
\enen
\vskip0.1in
Let us study perestroikas and smoothed perestroikas in the
positive direction.
Properties of perestroikas in the negative direction can be naturally
deduced from properties in the positive one.
\vskip0.1in
We shall first consider perestroikas in which
no imaginary double point is involved;
namely, cusp perestroika and  triple-point perestroika.
It is obvious that the corresponding  smoothed perestroikas
do not change the complex part. Therefore, we shall consider only
real part.
\been
\item
Cusp Perestroika\\
It is easy to see that smoothed cusp perestroika does not
change real part.
\item
Triple-point perestroika.
\vskip0.1in
Relatives orientation of the three real branches of a real triple-point
gives rise to distinguish two kinds of real triple-point.
Consider vectors at the triple point tangent
to the branches and directed according to their orientations.
If one of the vectors can be presented as linear combination of the two
others, then the triple-point is said {\it weak}
(see figure 5.2);
otherwise it is said {\it strong}( see figure 5.3)
Smoothed weak  triple-point perestroika does not change
real part.

Smoothed strong triple-point perestroika changes
real part as shown in figure 5.(3.b).
\enen
Then,
consider perestroika in which imaginary singularities are involved.

Introduce a function $J$  counting the difference between
the number of (resp, smoothed) $\a$-points and $\b$-points
under (resp, smoothed) perestroika.
Obviously,
if $J(\a)$ and $J(\b)$ are the values of $J$ under one
perestroika $\pi$ in the positive direction
the function $J$ takes values $-J(\a)$, $-J(\b)$
under the perestroika $\pi$ in the negative direction.
\been
\item
As already known,
relative orientation of the real branches
of a real self-tangency point gives rise to  distinguish
real point of direct ordinary tangency and real point of inverse
ordinary tangency.
\been
\item
It is easy to see that smoothed real  direct ordinary
tangency perestroika does not
change smoothed real part.
Nonetheless,
the complex part is changed
in such a way that:\\
$J(\a) = 0, J(\b) = -2$ under
 real  direct ordinary
tangency perestroika.
\item
Smoothed real inverse tangency perestroika
changes the real part (see figure 5.4 ).
Besides,
the complex part is changed in such a way that :
$ J(\a) = -2, J(\b) = 0$ under
real inverse tangency perestroika.
\enen
\item
Perestroika of solitary self tangency  is as follows:
two solitary double points come from the world of imaginary
to form a solitary self-tangency point; then arise
two solitary double points with opposite orientation.
Besides,
the complex part is changed in such a way that:
$ J(\a) = -2, J(\b) = 0$ under
a solitary inverse tangency perestroika.
\item
Perestroika of a triple-point with imaginary  branches changes the real part
as shown in figure 5.5.
Besides, it changes the complex part in such a way
that $J(\a) = +2, J(\b) = -2$.
\enen\section{Construction of Curves of type $I$ }
\mlb{su:typ1}
\vskip0.1in
In this section, we shall consider curves
${\AA}_m$ of degree $m$ and type $I$
with pair $(\RRR P^2, \RRR {\AA}_m)$
non-homeomorphic to the pair $(\RRR P^2, \RRR \HH_m)$.
As already noticed, for any curve ${\AA}_m$ of type $I$,
the pair $(\RRR P^2, \RRR {\AA}_m)$  with orientation  of the real point
set $\RRR {\AA}_m$ provides
the pair $(\CCC P^2, \CCC {\AA}_m)$,
up to conj-equivariant isotopy of $\CCC P^2$.
Bringing together the recursive Morse-Petrovskii's theory
for Harnack's curves  $\HH_m$ of chapter \ref{ch:mpr} and
properties of generic paths
in the varieties ${\mathcal  C}_{I;g}$,
we define in
Proposition \mrf{p:prop1typ1}
a path $S : [0,1] \to \RRR {\mathcal  C}_{m}$
with $S(0)=\RRR \HH_m$, $S(1)=\RRR {\AA}_m$.
Such path is obtained from lifting a generic path
in the space of generic immersion of the circle  into
$\RRR^2$ and is described as a set of moves defined
on the connected real components of $\RRR \HH_m$.\\
In this way, the  classification of $M$-curves of a prescribed degree $m$
amounts to the description of all the possible locations
of the real components of the Harnack curve $\HH _m$.

The first two parts of this section are devoted to the definition of the
lifting. Along a generic path in the space of
generic immersion of the circle  into
$\RRR^2$, three type of events "perestroikas"
(namely, the instateneous triple crossings and the inverse and
 direct self-tangency) may happen.
By means of perestroikas and their counterparts  for
algebraic curves, we bring the problem of the
definition of the path $S$ to
the one of lifting perestroikas in the space $\RRR {\mathcal C}_m$.\\
In the third part,
we give in Theorem \mrf{t:th2typ1} a description of
pairs $(\CCC P^2,\CCC {\AA}_m)$ up to conj-equivariant isotopy
which extends the properties of the Harnack curve $\HH_m$
stated in Theorem \mrf{t:theo1}
to any smooth curve of type $I$.

\subsection{Smoothing generic immersion of the circle into $\RRR^2$}

Recall that by a generic immersion of the circle $S^1$
into the plane $\RRR^2$, one means an immersion
with only ordinary double points of transversal intersection,
namely  without triple points and
without points of self-tangency.
The space of all immersions is an infinite-dimensional manifold which
consists of an infinite countable set of irreducible components
described by H.Whitney \cite{Whi} in 1937.

\vskip0.1in
{\bf Whitney's Theorem (1937)} \cite{Whi}
{\it The space of the immersions of a circle into the plane with the
same Whitney index is pathwise connected}
\vskip0.1in
The {\it Whitney index} of an immersion of an  oriented curve
into the plane is the rotation number of the tangent vector
(the degree of the Gauss map).

Arbitrary differentiable immersion of the circle into the plane
does not admit complexification.
Nonetheless,
we shall  generalize Whitney's Theorem to the case of real algebraic
curves of type $I$.
\vskip0.1in
Let us
recall in introduction the following property of the genus.
\vskip0.1in
{\bf Genus Property}
If a collection of disjoint circles embedded into a closed
orientable surface of genus $g$ does not divide the surface,
the number of circles is at most $g$.
In particular $g+1$ disjoint circles always divide the surface.\\
The set of complex points of any smooth
curve of degree $m$ is an orientable
surface of genus $g_m= \frac {(m-1) (m-2)} 2$.
Therefore, one can easily conjecture the following
one-to-one correspondence
between the $g_m+1$ (resp, $l \le g_m+1$)
real connected components of an $M$-curve (resp, an $M-i$-curve)
${\AA}_m$
and the $g_m$ $1$-handles and the sphere $S^2$:
the sphere $S^2$ and each $1$-handle contains
one real connected component of $\RRR {\AA}_m$
which divides it into two halves
(resp,
the sphere $S^2$ and each $1$-handle contains
a part of a real connected component of $\RRR {\AA}_m$;
such part divides it into two halves.)
\vskip0.1in
Call {\it regular curve},
the smooth oriented submanifold of $\RRR P^2$
deduced from the generic immersion of the circle by modification
at each real double point which is
either the Morse modification in $\RRR ^2$
in the direction coherent to a complex orientation
or the {\it Morse modification at infinity in   $\RRR P^2$}
(in the direction coherent to a complex orientation)
of the double point
which associates to the double point of  $\RRR ^2$
two points of the line at infinity
of $\RRR P^2$.
\vskip0.1in

In Lemma \mrf{l:l1typ1},
we give
a generalization of Whitney's theorem to the case of real algebraic
curves of type $I$.

\bele
\mlb{l:l1typ1}
{\it
Let ${\AA}_m$ be a smooth curve of degree $m$ and type $I$,
then its real point set $\RRR {\AA}_m  \subset \RRR P^2$
is a  regular curve
deduced from a generic immersion of the circle $S^1$
into the plane $\RRR ^2$
with
$\frac {(m-1)(m-2)} 2 \le n \le \frac {(m-1)(m-2)} 2 + [\frac {m} 2]$,
double points and Whitney index $\frac {(m-1)(m-2)} 2 + 1$.
The smoothing is such that
$\frac {(m-1)(m-2)} 2$ double points are
smoothed in $\RRR ^2$ by Morse modification
in the direction coherent to a complex orientation.
The $n - \frac {(m-1)(m-2)} 2$ others double points are
smoothed by Morse modification at infinity in $\RRR P^2$.}

\enle

\vskip0.1in
\bere
One can associate to any smooth curve of type $I$
a three-dimensional rooted tree.\\
Recall \cite{Ar2} that a generic immersion of the circle into
the plane is a {\it tree-like curve} if any of its double points
subdivides it into two disjoints loops.
On the assumption of Lemma \mrf{l:l1typ1},
let ${\AA}_m$ be a smooth curve of degree $m$ and type $I$, and $\phi$ the
corresponding generic immersion of the circle $S^1$
into the plane $\RRR ^2$
with
$\frac {(m-1)(m-2)} 2 \le n \le \frac {(m-1)(m-2)} 2 + [\frac {m} 2]$.
Smooth the double-points of $\phi$
which result from points of $\RRR {\AA}_m$ at infinity in such a way that
it results a curve in $\RRR ^3$
with the property that any of its double points
subdivides it into two disjoints loops; it follows
the rooted-tree associated
to ${\AA}_m$.\\
In case of Harnack curves ${\HH}_m$,
an immediate definition of the root of its rooted-tree is the following:
assume ${\HH}_m$ obtained via patchworking process;\\
for odd $m$, the root consists of the part of the odd component
of the curve which intersects $T_m$;\\
for even $m$, the root consists of the part of the  positive oval
(non-empty for $m \ge 4$)
which intersects $T_m$.
\enre

{\bf proof:}
\vskip0.1in

Our proof makes use of properties of the complex point set
$\CCC {\AA}_m$ embedded in $\CCC P^2$.

Consider the usual handlebody decomposition of
$\CCC P^2= B_0 \cup B_1 \cup B$  where $B_0$,$B_1$,$B$
are respectively 0, 2 and 4 handles.

The balls $B_0$ and $B_1$ meet along an unknotted  solid torus
$S^1 \ti B^2$. The gluing diffeomorphism
$S^1 \ti B^2 \to S^1 \ti B^2$ is given by the $+1$ framing map.
In such a way,
the canonical $\RRR P^2$ can be seen as the union of a M\"{o}bius
band $\mathcal  M$ and the disc $D^2 \subset B$ glued along their boundary.
The M\"{o}bius band $\mathcal  M$ lies in $B_0 \cap B_1 \approx S^1 \ti D^2$
with $\pr {\mathcal  M}$ as the $(2,1)$ torus knot, and
$D^2 \subset B$ as the properly imbedded unknotted disc.
The complex conjugation switches $B_0$ and $B_1$ and lets fix
$\mathcal  M$, it rotates $B$ around
$D^2$.

The set of complex points of ${\AA}_m$
is an orientable
surface of genus $g_m= \frac {(m-1) (m-2)} 2$, i.e it is diffeomorphic to a
sphere $S^2$ with $g_m$ $1$-handles $S^2_{g_m}$.
We shall denote $h$ the diffeomorphism $h :\CCC {\AA}_m \to S^2_{g_m}
\subset \RRR^3$.

Recall that surfaces $S^2_{g_m}$ are constructed
as follows.
From the sphere $S^2$,
$(m-1)$ pairwise non-intersecting open discs are removed,
and then the resulting holes are closed up by
$g_m$ orientable cylinders $\approx S^1 \ti D^1$   connecting
the boundary  circles of the discs removed.

Assume $S^2$ provided with a complex conjugation with fixed
point set a
circle $S^1$ which divides the sphere $S^2$ into two halves.
Without loss of generality, one can assume
that any disc removed from
$S^2$ intersects $S^1$ and the two halves of $S^2$.
In such a way, $\RRR {\AA}_m \subset \RRR P^2$ intersects
each $1$-handle.

Fix $D^2$ the two disc of $\RRR P^2= {\mathcal  M} \cup D^2$
in such a way that the boundary circle of $D^2$ is $S^1$ and
therefore each one handle belongs to the solid torus
$S^1 \ti D^2 \supset {\mathcal  M}$.

Let $D^2_{\e} \supset D^2$ be the disc $D^2$ thickened.
Since the interior  $(D^2_{\e})^0$  of $D^2_{\e}$ is homeomorphic to
$\RRR^2$, one can project
$\RRR {\AA}_m$ (up to homeomorphism $\RRR ^2 \approx (D^2_{\e})^0$)
in a direction perpendicular to $\RRR^2$
onto $\RRR^2$.
We may suppose that the direction of the projection is generic
i.e all points of self-intersection of the image on $\RRR^2$  are double
and the angles of intersection are non-zero.\\
Let $r :\RRR^3 \to \RRR^2$
be the  projection which
maps $h(\RRR {\AA}_m) \subset \RRR ^3$ to $\RRR ^2$.

Consider an
oriented tubular fibration
$N \to \CCC {\AA}_m$.
Since $\CCC {\AA}_m$ is diffeomorphic to a
sphere $S^2$ with $g_m$ $1$-handles, one can
consider the restriction of $\CCC {\AA}_m$
diffeomorphic to each torus $T^2$ given by
the sphere $S^2$ with one of the $g_m$ $1$-handles.
The oriented tubular neighborhood of
$h(\CCC {\AA}_m) \cap T^2$, as oriented tubular neighborhood
of the torus $T^2$, intersects the solid
torus $S^1 \ti D^2 \supset {\mathcal  M}$
with $\pr M \subset  S^1 \ti S^1$ as the $(2,1)$ torus knot.
Hence, since in $\CCC P^2$ each real line
is split by its real part into two
halves lines conjugate to each others,
and $2$ disjoint circles always divide the torus;
the real part of the restriction of
$h(\CCC {\AA}_m) \cap T^2$ belongs to the
boundary of the M\"{o}bius band in such a way
that its  projection to $\RRR ^2$ gives one crossing.
Besides these doubles, some double-points of
 $r(h (\RRR \AA_m))$ may
result either from two points of $\RRR {\AA}_m$
which belong to two different handles or from
two points which belong respectively
to a $1$-handle and to $S^2$.
It is easy  to deduce from the relative location
of the sphere $S^2$ and the $1$-handles that
the projection of $\RRR {\AA}_m$
leads to an even
number of these double-points of
$r(h (\RRR \AA_m))$.
Therefore,
the projection $r(h(\RRR {\AA}_m))$ contains  $g_m (mod 2)$ crossings
and
$\RRR {\AA}_m$ can be seen as a regular curve.

By use of
the Whitney expression of the index $ind$ of a plane
curve with $n$ double points - $ind =\sum \e_i \pm 1$
where the summation is over the set of $n$ double points, $\e_i$ is a sign
associated to each double point, and the term $\pm$ depends on the
orientations of the curve and of the plane -
we shall deduce
the Whitney index of $r (h(\RRR {\AA}_m))$.

According to the definition of $\RRR P^2$,
each real part of $h(\RRR {\AA}_m)$ contained in a $1$-handle
has the same projection to $\RRR ^2$
as the $(2,1)$ torus knot on the torus defined from
the sphere $S^2$ with the $1$-handle.

Therefore,
according to the gluing diffeomorphism\\
$S^1 \ti B^2 \to S^1 \ti B^2$,
given an orientation on the
circle $S^1$ induced by one of the halves
of $S^2$,
the projection of each real part of $h(\RRR {\AA}_m)$
contained in a $1$-handle
gives a crossing with sign $+1$ of
the plane $\RRR^2$.
Besides, any $2$ crossings
which result from the projection of points of two different handles
or of points of a $1$-handle and the sphere $S^2$
have opposite Whitney index.
Hence, $r(h(\RRR {\AA}_m))$ has Whitney index $g_m+1$.

Moreover, since ${\AA}_m$ is of
type $I$,
the real part $\RRR {\AA}_m$ is deduced
from $r(h(\RRR {\AA}_m))$
by smoothings which are
either
standard Morse modifications in $\RRR ^2$ at the double point
determined by the complex  orientation, or smoothings at infinity
in $\RRR P^2 \bk \RRR^2$
which associate to a crossing its two pre-image
under $r^{-1}$.

Two real branches of $r(h(\RRR {\AA}_m))$
which intersect in a crossing
become after smoothing two oriented non-intersecting real branches
which belong either to two different connected components
of $\RRR {\AA}_m$ or to one branch which lies as a part of
the $(2,1)$ torus knot
on the boundary M\"{o}bius ${\mathcal  M} \subset \RRR P^2$.
According to the Rokhlin's orientation formula, for any $m > 2$,
$\RRR {\AA}_m$
has at least two connected components.

Therefore, for $m > 2$, $g_m$ crossings are smoothed
in $\RRR ^2$  by Morse Modification
in the direction coherent to the complex orientation
and become two oriented non-intersecting real branches
of distinct connected components of $\RRR {\AA}_m$ in $\RRR P^2$.
The other double-points are smoothed in such a way that
it results two points of
$\RRR {\AA}_m$ in $\RRR P^2 \bk \RRR ^2$.\\
According to the Bezout's theorem,
since each of these
$n-g_m$ double-points
results from two points which belong to  $\RRR P^2 \bk \RRR ^2$
and the number of points in the intersection
of $\RRR {\AA}_m$
with the real line  $\RRR P^2  \bk \RRR ^2$
is less or equal to $m$ and congruent to $m (mod 2)$,
the difference $n-g_m$ is less or equal to
$[\frac {m} 2] $.
It is obvious that one can assume, without loss of generality,
that these
$[\frac {m} 2] $ double points results from
points of $\RRR {\AA}_m$
which belong to two different handles or from
two points which belong respectively
to a $1$-handle and to $S^2$.

Q.E.D.
\vskip0.1in
Given a generic immersion of the circle $S^1$
into the plane $\RRR ^2$ and  ${\mathcal  S}$ the set of its singular
points.
We shall call
{\it partially regular curve}
the smooth oriented submanifold of $\RRR P^2$
deduced from a generic immersion of the circle by modification
of a set ${\mathcal  K} \subset {\mathcal  S}$,
${\mathcal  K} \not = {\mathcal  S}$,
of real double points where modification
at each real double point of ${\mathcal  K}$ is
either the Morse modification in $\RRR ^2$
(in the direction coherent to a complex orientation)
or the Morse modification at infinity of the double point
which associates to the double point of  $\RRR ^2$
two points of the real line $\RRR P^2  \bk \RRR^2$.

\vskip0.1in

The next Lemma  enlarges the preceding statement
to generic curves of degree $m$, type $I$ and genus $g$,
and more generally to
singular curves of degree $m$ and type $I$
with non-degenerate singular points.
\vskip0.1in

\bele
\mlb{l:l2typ1}
{\it
Let ${\AA}_m$ be a singular curve of degree $m$ and type $I$
with non-degenerate singular points,
then its real point set $\RRR {\AA}_m  \subset \RRR P^2$ is a
partially regular curve
deduced from a generic immersion of the circle $S^1$
into the plane $\RRR ^2$
with
$\frac {(m-1)(m-2)} 2 \le n \le \frac {(m-1)(m-2)} 2 + [\frac {m} 2]$,
double points and Whitney index
$\frac {(m-1)(m-2)} 2 + 1$ if and only if its set of singular points
consists of at most
$\frac {(m-1)(m-2)} 2$ crossings.}
\enle
\vskip0.1in

{\bf proof:}
It follows from an argument similar to the one of the proof of
Lemma \mrf{l:l1typ1}.
Q.E.D
\vskip0.1in
\subsection{Lifting of a generic path in the space of generic
 immersion of the circle into the plane  to $\RRR {\mathcal C}_m$ }
\mlb{susu:lift}
\vskip0.1in

According to the Lemma \mrf{l:l1typ1},
one can consider $\RRR \HH_m$ and $\RRR {\AA}_m$
as smoothed immersions of the circle into
the plane with the same Whitney index.
Let us denote
$\RRR \tilde{\HH}_m$, $\RRR \tilde{\AA}_m$
the corresponding immersions.
From the Whitney's Theorem,
there exists a path $\tilde{h}$  which connects
$\RRR \tilde{\HH}_m$ and $\RRR \tilde{\AA}_m$.
In this section,  we shall in Proposition \mrf{p:prop1typ1} and
Theorem \mrf{t:th1typ1} define
a path $S$ in $\RRR {\mathcal  C}_m$
$S(0)=\RRR \HH_m$,\\ $S(1)=\RRR {\AA}_m$
from lifting a path $\tilde{h}$  which connects
$\RRR \tilde{\HH}_m$ and $\RRR \tilde{\AA}_m$ in the space of immersions
of the circle into $\RRR^2$.

By means of the Lemmas \mrf{l:l1typ1} and \mrf{l:l2typ1},
we shall develop properties provided by the genus property and
interpret the path
$\tilde{h}:
\RRR \tilde{\HH}_m \to \RRR \tilde{\AA}_m$
as a change of relative position
of the $1$-handles of
$(\CCC \HH_m) \approx S^2_{g_m}$,
$(\CCC \HH_m) \supset \RRR \HH_m$
which gives
$(\CCC {\AA}_m) \approx  S^2_{g_m}$,
$(\CCC {\AA}_m) \supset \RRR {\AA}_m$.
The path $S$  appears as the track on $\RRR P^2$
of a diffeotopy $h_t$ of $\CCC P^2$, $t \in [0,1]$,
$h(0)= \CCC \HH_m$, $h(1)=\CCC {\AA}_m$.
\vskip0.1in
{\it (\mrf{susu:lift}).A
 Lifting Perestroikas in the space $\RRR {\mathcal C}_{m}$}\\

Recall that given a perestroika along a path in the space
${\mathcal C}_{I,g}$
we call {\it smoothed perestroika} the change obtained by
smoothing the fragments involved in the perestroika.
We shall call {\it diffeotopic perestroika}
of perestroika defined along a path $h$
in the space  ${\mathcal  C}_{I,g})$ of irreducible
curves of type $I$ and genus $g$ and degree $ 1 \le n \le m$,
the change experienced by a curve
in ${\mathcal C}_{I,g}$ with the following property.\\
For any  point $a$ which participated in a perestroika (along the path $h$)
let $U(a)$ be a small $\e$-neighborhood
 around $a$ of the complex point
set of the initial curve embedded in $\CCC P^2$,
(one can choose also $U(a)$ as follows:
$U(a) = \{ z=<u,p>=(u_0.p_0:u_1.p_1:u_2.p_2) \in \CCC P^2 \mid
u=(u_0,u_1,u_2) \in U_{\CCC}^3, p=(p_0:p_1:p_2) \in D(a,\e) \}$
where $D(a,\e)$ denotes a small disc
(in the Fubini-Study metric)  of radius $\e$ around $a$
in $\RRR P^2$.
Restrictions to $U(a)$ of  the  complex point set
of the initial curve
and of
the  complex point set
which results from  the diffeotopic perestroika
are diffeotopic.
(As will become clear later, in most cases, given a curve  ${\mathcal A}_m$
the change experienced by the curve
${\mathcal A}_m$
after just one diffeotopic perestroika
does not lead to a curve.)

Given a diffeotopic perestroika, we call
{\it smoothed diffeotopic perestroika}
the change obtained by smoothing the fragments involved in
the diffeotopic perestroika (i.e fragments inside $U(a)$
for any $a$ which participated in the perestroika).
\vskip0.1in

{\it
 (\mrf{susu:lift}).B
   Chain of Diffeotopic Perestroika}

Let us recall that up to regular deformation of its polynomial
(see Theorem \mrf{t:bli}),
curve $\HH_m$
results from the recursive construction of Harnack curves
$\HH_i$, $1 \le i \le m$, where  $\HH_{i+1}$ is deduced from classical
small perturbation of the union $\HH_i \cup L$ of the curve $\HH_i$
with a line $L$.

Hence, when searching to
interpret the path
$\tilde{h}:
\RRR \tilde{\HH}_m \to \RRR \tilde{\AA}_m$
as a change of relative position
of the $1$-handles of
$(\CCC \HH_m) \approx S^2_{g_m}$,
$(\CCC \HH_m) \supset \RRR \HH_m$,
one can proceed by induction
and search for $1 \le j \le (m-1)$ to move the  union
of the $g_{j}$ 1-handles of
$\CCC \HH_{j} \approx S^2_{g_j}$  with
the $(j-1)$ more 1-handles of
$\CCC \HH_{j+1} \approx S^2_{g_{j}+(j-1)}$
in such  a way that at the end of the process
it results
$(\CCC {\AA}_m) \approx  S^2_{g_m}$,
$(\CCC {\AA}_m) \supset \RRR {\AA}_m$.

Using the recursive Morse-Petroskii's theory of chapter \ref{ch:mpr}
as a tool to describe the lifting of $\tilde{h}$ as a track of a diffeotopy
$h_t$ of $\CCC P^2$, $t \in [0,1]$,
$h(0)= \CCC \HH_m$, $h(1)=\CCC {\AA}_m$,
we shall need {\it a global order} on the set of crossings of curves
$\HH_i \cup L$ and the associated $k$-uple.\\

{\it a)  global order}\\

According to Theorem \mrf{t:theo1}, for any
Harnack curve $\HH_m$ of degree $m$,
there exists a finite number $I$
($I= \frac {m(m+1)} 2  +  \Sigma_{k=2}^{k=[m/2]} (2k-3)$)
of disjoint $4$-balls $B(a_i)$ of radius $\e$ (in the Fubini-Study metric)
invariant by complex conjugation centered in points $a_i$
of $\RRR P^2$ such that
up to conj-equivariant isotopy of $\CCC P^2$,
$\HH_m \bk \cup_{i \in I} B(a_i)= \cup_{i=1}^m L_i \bk
  \cup_{i=1}^I B(a_i) $
 where $L_1,..., L_m$ are $m$ distinct projective lines
 with
 $$L_i \bk \cup_{i =1}^I B(a_i) \cap
  L_j \bk \cup_{i=1}^I  B(a_i) = \emptyset$$
  for any $i \not=j$, $1\le i, j \le m$.\\

The proof of Theorem \mrf{t:theo1} is based on an induction on the degree $m$
of the curve $\HH_m$ which provides a natural order on the lines $L_j$,
$1 \le j \le m$.
Since any $4$-ball $B(a_i)$)
intersects exactly $2$ lines $L_{j-1}$,
$L_{j}$
the natural order on the set of lines $L_j$
extends to an order on the set of points $a_i$.
We shall say that
$a_i$ has order $j$, $2 \le j \le m$, if and only if
$B(a_i)$ intersects the lines $L_{j-1}$ and $L_{j}$.
(In other words, in the inductive construction,
any crossing of $\HH_j \cup L$ has order $j$;
and according to Definition \mrf{d:uu},
the first point of a
$k$-uple has the order of the crossing to which
it is associated.) \\

Such order extends on $3$-uple (resp the $2$-uple)
as follows.
Given $a_j$ of global order $j$ the first point of a
$3$-point $(a_j,a_{j+1},a_{j+2})$,
($j$ even)
associated to a
$a_{j+1}$  has global order $j+1$,
$a_{j+2}$  has global order $j+2$.

Consider the set $a_j$ $1 \le j \le J=\frac {m(m+1)} 2$
of crossings of  curves $\HH_{i} \cup L$ , $1 \le i \le (m-1)$
in the inductive construction of $\HH_m$.
It may be easily easily extracted from the proof of Theorem \mrf{t:theo1}
(see also Proposition \mrf{p:prop8} and Definition \mrf{d:uu})
that
$\cup_{i=1}^I B(a_i) \supset \cup_{j=1}^{\frac {m(m+1)} 2} U(a_j)$
where $B(a_i)$ is the $4$-ball around $a_i$ of radius $\e$ and
$U(a_j)$ is the $\e$-neighborhood around $a_j$.
(According to terminology of section \mrf{susu:cji}
(see proposition \mrf{p:prop8})
the set of points $\cup_{i \in I} a_i = \cup_{n=1}^m A_n$
where $A_n$ denotes the set of points perturbed in a maximal
simple of deformation of $\HH_n$.)\\

For sake of simplicity,
according to corollary \mrf{c:gpa},
we shall consider the $T$-inductive construction of Harnack curves
and $e$-tubular neighborhoods $U(p)$ (see definition \mrf{d:uu2})
defined from faces $\G \subset l_{j}$
, $l_{j}=\{ (x,y) \in (\RRR^+)^2 | x +y =j \}$.
As already introduced in subsection \mrf{su:statepb},
${\mathcal P}_m$ denotes the set of points $p$
with the property:
$\cup_{p \in {\mathcal P}_m}  U(p) = \cup_{a_i \in I} U(a_i)$.
(In such a way, we have the following equivalent definition of
the global order for the patchworking construction.
Any point $a$ of the set $a_i \in I$
has order $j+1$, $2 \le j +1 \le m$, if and only if $a \in  U(p)$
where
$U(p)$ is the $\e$-tubular neighborhood
defined from the interior $\G^0$ of a face $\G \subset l_{j}$
, $l_{j}=\{ (x,y) \in (\RRR^+)^2 | x +y =j \}$.
This global order extends on $k$-uple
as previously explained.)\\
Let us extend the global order on the set points $a_i, i \in I$
to the set $\cup_{a_i \in I} U(a_i)$.
We shall say that $U(a_i)$ has order $j$ if $a_i$
has global order $j$.
Equivalently, in the pathchworking scheme, $U(p)=U(a_i)$
has order $j$ if it
is defined from the interior $\G^0$ of a face $\G \subset l_{j-1}$.

\vskip0.1in
{\it b) Covering of $\CCC P^2$}\\

The union $\cup_{i \in I} U(a_i)$  does not cover $\CCC P^2$.
We shall consider consider the family
$\tilde{U}(a_j)$, $\cup_{j=1}^J
\tilde{U}(a_j) \supset \cup_{i=1}^I B(a_i)$,
of neighborhoods of points $a_j$, $j \in J$,
verifying the following properties:

- any $\tilde{U}(a_l)$  of order $j$
intersects two neighborhoods $\tilde{U}(a_k)$  of order $j-1$,

\begin{eqnarray}
\lb{e:p1}
&{\cup_{j=2}^m \cup_{a_i~of~order~j} \tilde{U}(a_i)}& =\CCC P^2\nonumber\\
\lb{e:p2}
-for~j=2: &\nonumber\\
&\cup_{a_i~of~order~2} \tilde{U}(a_i)&
\supset \CCC L_2 \bk  \cup_{a_{i} of order 3} \tilde{U}(a_i)\nonumber\\
&\cup_{a_i~of~order~2} \tilde{U}(a_i)&
\supset \CCC L_{1}
\nonumber\\
\lb{e:p3}
-for~3 \le j <m: &\nonumber\\
&\cup_{a_i~of~order~j} \tilde{U}(a_i)&
\supset \CCC L_j \bk  \cup_{a_{i}~of~order~j+1} \tilde{U}(a_i)
\nonumber\\
&\cup_{a_i~of~order~j} \tilde{U}(a_i)&
\supset \CCC L_{j-1} \bk  \cup_{a_{i}~of~order~j-1} \tilde{U}(a_i)
\nonumber\\
\lb{e:p4}
-for~j =m: &\nonumber\\
&\cup_{a_i~of~order~m} \tilde{U}(a_i)&
\supset \CCC L_m\nonumber\\
&\cup_{a_i~of~order~m} \tilde{U}(a_i)&
\supset \CCC L_{m-1} \bk  \cup_{a_{m}~of~order~m-1} \tilde{U}(a_i)\nonumber
\end{eqnarray}

Given $\HH_m$ the Harnack curve of degree $m$
and ${\AA}_m$ a smooth curve of type $I$,
we shall
in the Proposition \mrf{p:prop1typ1},
define a path
$$S :[0,1] \to \RRR {\mathcal  C}_{m}$$
$S(0)=\RRR \HH_m$,
$S(1)=\RRR {\AA}_m$
described locally
in the opens $\tilde{U}(a_i)$ and sequentially
by induction on the global order on the set points
$a_i$,$i \in I$ and the associated $k$-uple
by diffeotopic perestroikas of perestroikas
in the spaces ${\mathcal C}_{I,g}$ of
irreducible curves of degree $1 \le n \le m$, type $I$ and genus $g$.
Besides,
any singular point which participated in diffeotopic
perestroika belongs to $\cup_{i=1}^I B(a_i)$.

\bepr
\mlb{p:prop1typ1}
{\it
Let $\HH_m$ be the Harnack curve of degree $m$
and ${\AA}_m$ be a smooth curve of type $I$.
There exists a path
$$S :[0,1] \to \RRR {\mathcal  C}_{m}, S(0)=\RRR \HH_m, S(1)=\RRR {\AA}_m$$
which crosses the discriminant hypersurface $\RRR {\mathcal D}_m$.
Up to conj-equivariant isotopy of $\CCC P^2$,
the path  may be sequentially defined in the open $\tilde{U}(a_i)$
by induction on the global order by smoothed
diffeotopic perestroikas of
diffeotopic perestroikas  on the spaces
${\mathcal C}_{I,g}$ of curves of type $I$, degree $1 \le n \le m$,
and genus $g$ as follows:
\vskip0.1in

The sequence of diffeotopic perestroika is defined as follows:\\
Assume $\HH_m$ obtained via the patchworking method.
\been
\item
one can assume that only
points $\cup_{i \in I} a_i$
are double-points involved in a diffeotopic perestroika.

\item
any diffeotopic perestroika
of which  real ordinary double-point are points of order at most $j$
is a diffeotopic perestroika in the space ${\mathcal C}_{I,g}$ of curves
of degree $j$, type $I$ and genus $g$.
\item
a diffeotopic perestroika defines double-points involved in next
as follows.
As  branches involved in
the topological meaning of the critical point associated to
a real ordinary double-point of order $j$ involved in the
diffeotopic perestroika
move under the diffeotopic perestroika they define
real ordinary double-point of order $j+1$ involved in a next
diffeotopic
perestroika.
in the space ${\mathcal C}_{I,g}$ of curves
of degree $j+1$, type $I$ and genus $g$.\\
Any imaginary point which participated to
such sequence of diffeotopic perestroika belongs to
the intersection of a neighborhood
 $\tilde{U}(a_j)$
of a point order $j$ and a neighborhood of a point
$\tilde{U}(a_{j+1})$ of order $j+1$.
Besides, imaginary points participated in such a way that if one
$\a$-point (or $\b$-point) appears (resp, disappears)
after a diffeotopic perestroika on real points of order $j$, then
it disappears (resp, appears) after a
diffeotopic perestroika
on real points of order $j+1$.
\enen}
\enpr
\vskip0.1in

{\bf proof:}
\vskip0.1in
Let us  explain the method of our proof.
\vskip0.1in
We shall prove that
up to slightly modify the coefficients of the polynomial
giving the curve ${\AA}_m$,
one can always assume that there exists a diffeotopy $h_t$
of $\CCC P^2$
$h(0)= \CCC \HH_m$, $h(1)=\CCC {\AA}_m$
with the property
$h_t(\CCC {\AA}_m   \cap U(p))  \subset U(p)$
for any $p \in {\mathcal  P}_m$.

Our argumentation is based on
the Lemma \mrf{l:l1typ1} and its proof.
According to the Lemma \mrf{l:l1typ1},
one can consider $\RRR \HH_m$ and $\RRR {\AA}_m$
as smoothed immersions of the circle into
the plane with the same Whitney index.
We shall denote
$\RRR \tilde{\HH}_m$, $\RRR \tilde{\AA}_m$
the corresponding immersions.
From the Whitney's theorem,
there exists a path $\tilde{h}$  which connects
$\RRR \tilde{\HH}_m$ and $\RRR \tilde{\AA}_m$.

Therefore, the definition of the path $S$ in $\RRR {\mathcal  C}_m$
$S(0)=\RRR \HH_m$,\\ $S(1)=\RRR {\AA}_m$
may be  deduced from a lifting of the path $\tilde{h}$  which connects
$\RRR \tilde{\HH}_m$ and $\RRR \tilde{\AA}_m$ in the space of immersions
of the circle into $\RRR^2$.

Since any smooth curve is irreducible, and any reducible polynomial
is the product of a finite number of reducible polynomials,
we shall, according to Lemma \mrf{l:l1typ1} and Lemma \mrf{l:l2typ1},
lift the path
$\tilde{h}: \RRR \tilde{\HH_m} \to \RRR \tilde{{\AA}}_m$
in  the space of immersion of the circle into the plane
from smoothed diffeotopic perestroikas of
diffeotopic perestroikas in the  spaces
${\mathcal  C}_{I;g}$ of curves of degree $n \le m$, type $I$ and genus
$0 \le g \le \frac {(n-1).(n-2)} 2$.

For any diffeotopic perestroika $\pi$ experienced along a generic path in
${\mathcal  C}_{I;g}$,
locally a description of the real components of the resulting curve
is deduced  from double points involved in the diffeotopic perestroika.
Thus, any smoothed diffeotopic perestroika $\pi$ experienced along $S$,
may be defined from its smoothed double points.\\
From the detailed construction of $\HH_m$ given  in
chapter \ref{ch:mpr}, we shall get smoothed points of smoothed
diffeotopic perestroika experienced along the path $S$.\\
As
already introduced in the proof of Lemma \mrf{l:l1typ1},
we shall consider $h$ the diffeomorphism
$h:\CCC \HH_m \to S^2_{g_m} \subset \RRR^3$
and $r: \RRR^3 \to \RRR^2$ the projection which maps $h(\RRR \HH_m)$
to $r(h(\RRR \HH_m)=\RRR \tilde{\HH}_m \subset \RRR^2$.
We shall prove that only points of
$\RRR \tilde{\HH}_m$
 smoothed in $\RRR^2$ to give $\RRR \HH_m$
may lift to smoothed double points
of an irreducible curve of degree $n$, $1 \le n \le m$
type $I$ and genus $g$ involved in smoothed diffeotopic
perestroika.
Thus, double-points of $\RRR \tilde{\HH}_m$ smoothed at infinity
are involved along the path
$\tilde{h}:
\RRR \tilde{\HH}_m \to \RRR \tilde{\AA}_m$ as the relative position
of the $1$-handles of
$(\CCC \HH_m) \approx S^2_{g_m}$,
$(\CCC \HH_m) \supset \RRR \HH_m$
changes to give
$(\CCC {\AA}_m) \approx  S^2_{g_m}$,
$(\CCC {\AA}_m) \supset \RRR {\AA}_m$.

In such a way,
using the path provided
by  Lemma \mrf{l:l1typ1} and
Whitney's theorem, according to the Lemma \mrf{l:l2typ1} and the
fact that $\CCC \HH_{m}$ and $\CCC {\AA}_m$
are diffeotopic, we shall characterize the track $S$
on $\RRR P^2$
of a diffeotopy $h_t$ of $\CCC P^2$, $t \in [0,1]$,
$h(0)= \CCC \HH_m$, $h(1)=\CCC {\AA}_m$.
\vskip0.1in

We shall divide our proof into three parts:

\been
\item In the first part, we
define double-points of curves along $S$.
\item
The second part deals with
smoothed diffeotopic perestroikas along $S$.
\item
In the third part, we give a method which provides
any curve ${\AA}_m$ of degree $m$ and type $I$
from $\HH_m$, that is $\RRR {\AA}_m$ from
$\RRR {\HH}_m$.
\enen

\vskip0.1in
The first two parts make use of the idea that locally
a real branch does not differ from another real branch.
In these parts, we  consider
the situation more locally than globally and fix definitions
we shall need in the third part.
In the third part, we give a method of organizing these local
moves in such a way that globally it results a real algebraic curve
in $\CCC P^2$.
\vskip0.1in
The exposition of the main ideas of our proof is now finished, and we
shall proceed to precise arguments.
\vskip0.1in

{\bf 1. Double points of curves along $S$.}

Let us start our definition of double points of curves along
the path $S$ we have planed to construct by  the following  remark.
Assume $\HH_m$ constructed via the patchworking method.
Given $\HH_m$,
and its projection $r (h (\RRR \HH_m))= \RRR \tilde{\HH_m}$,
it is easy to deduce from
the Lemma \mrf{l:l1typ1}  and its proof and Lemma \mrf{l:l2typ1},
that for $m > 2$ only the $g_m$ points of $\RRR \tilde{\HH}_m$
smoothed in $\RRR^2$ to give $\RRR {\HH}_m$  may lift to real
double points of curves along $S$.

\vskip0.1in
Denote \quad \break
$\rho^{m} :(T_{m} \ti U_{\CCC}^2) \to \CCC T_{m} \approx
\CCC P^2$  the natural surjection.
According to the  chapter \ref{ch:Mope} (see Proposition \mrf{p:prop8} and
Theorem \mrf{t:theo1}),
the Harnack curve $\HH_m$ of degree $m$ is
obtained from recursive perturbation of curves $\HH_{m-i} \cup L$, $0<i<m$,
( with $L \approx \CCC l_{m-i} = \rho^m(l_{m-i} \ti U_{\CCC}^2)$)
in such way that outside a finite number $I$
of $4$-ball $B(a_i) \subset \CCC P^2$
there exists a conj-equivariant isotopy $j_t$ of subset of $\CCC P^2$
which maps:
\been
\item
$\HH_{m-i} \bk \cup_{i=1}^{I} B(a_i)$
onto a part of the curve $\HH_{m-i+1}$.\\
(in the patchworking  scheme, such part is contained in the restriction
$\rho ^{m}(T_{m-i} \ti U_{\CCC}^2)$ of
$\CCC P^2 \approx \CCC T_{m}$.)
\item
$(\HH_{m-i} \bk j_1(\HH_{m-i+1} \bk \cup_{i=1}^{I} B(a_i))$ onto
a part of the projective line $L$.\\
(in the patchworking scheme, such part  is contained in the restriction
$\rho ^{m}(D_{m-i+1,m-i} \ti U_{\CCC}^2 )$ of
$\CCC P^2 \approx \CCC T_{m}$.)
\enen
Consequently, (see  Theorem \mrf{t:theo1}),
outside a finite number of $4$-balls $B(a_i)$, the curve
$\HH_m$, is up to conj-equivariant isotopy of $\CCC P^2$,
the union of $m$ projective lines
minus their intersections with $4$-balls $B(a_i)$.
\vskip0.1in
Since projective lines are diffeotopic to circles,
outside a finite number of $4$-balls $B(a_i)$,
components of $\RRR \HH_m$  can be considered
as components of $m$ concentric circles.
Thus, smoothed real points involved in diffeotopic perestroikas
along $S$ belong to the union
of real components of $\HH_m$ contained in
$\rho^m (D_{m-i +1, m-i}\cup D_{m-i, m-i -1} \ti U_{\CCC}^2)$
where
$D_{m-i+1, m-i}$ and $D_{m-i, m-i-1}$
 are bands of $T_m$
(for $j \in \{ m-i-1, m-i \}$,
$D_{j+1, j } =\{ x \ge 0, y \ge 0 | j \le x+ y \le j +1 \}$)
with intersection
$D_{m-i+1, m-i}\cap D_{m-i, m-i-1}=l_{m-i}$.\\

Therefore, locally any smoothed real point
can be considered as
the smoothed intersection of two real branches of
$ \RRR \HH_m \cap
\rho^m (D_{m-i +1, m-i}\cup D_{m-i, m-i -1} \ti U_{\RRR}^2)$
of which intersection point belongs to a neighborhood of
$\rho^m (l_{m-i} \ti U_{\RRR}^2)$.

Consequently, it is easy to deduce from relative orientation  and
relative location of real branches of
$\RRR \HH_m
\cap \rho^m (D_{m-i +1, m-i}\cup D_{m-i, m-i -1} \ti U_{\CCC}^2)$
that
smoothed real points involved in diffeotopic perestroika along $S$
belong to $\rho^m(D(p) \ti U_{\RRR}^2)= U(p) \cap \RRR ^2$
with $p \in {\mathcal  P}_m$ and that any $U(p)$ contains one
such smoothed double-point.
\vskip0.1in

In such a way, according to
Lemma \mrf{l:l1typ1} and Lemma \mrf{l:l2typ1},
smoothed real points along $S$
are smoothed crossings which belong to the set \\
$\cup_{p \in {\mathcal  P}_m }\rho^m( D(p)\ti U_{\RRR}^2)
= \cup_{p \in {\mathcal  P}_m } U(p) \cap \RRR^2$.
\vskip0.1in
\vskip0.1in
\vskip0.1in

{\bf 2.Smoothed diffeotopic perestroika on
smoothed double points of curves along $S$.}

Recall that
from the previous Morse-Petrovskii's study of Harnack curves
(see chapter \ref{ch:mpr})
for any $p \in {\mathcal  P}_m$,
there exists a crossing $a \in \HH_{m-i} \cup L$
with $a \in U(p)= \rho^m (D(p) \ti U_{\CCC}^2)$
such that the crossing
$a \in {\mathcal  P}_m$ is associated to:
\been
\item
either
a simple-point of $\HH_{m}$
\item
either
a $3$-point $(p_{m-i +1}, p_{m-i +2}, p_{m-i +3})$
(where $p_{m-i +3} \approx p_m$
of $\HH_{m}$
\item
or
in case $m$ odd a $2$-point
$(p_{m-1}, p_{m})$ of $\HH_{m}$
\enen
\vskip0.1 in

As recalled previously, cusp perestroika
and weak triple point perestroika
do not change real and complex smoothed parts.
Therefore, there is not need to consider them.
Besides, since smoothed real points of
generic curves along $S$ are smoothed crossings,
no smoothed imaginary self-tangency perestroika
is undergone along $S$. (This restriction is not
in contradiction with  the fact that
${\AA}_m$ is an arbitrary curve
of degree $m$ and type $I$
since $\HH_m$ is an $M$-curve.)

Obviously, in $U(p)$ any change implied by a smoothed
diffeotopic perestroika
is a perturbation in which,
besides real branches of the crossing $a \in U(p)$ of $\HH_j \cup L$,
 real branches
involved in the local
topological meaning of simple point
and of any point of a $3$-point (or $2$-point in case $m$ odd)
participated. \\
We shall study the situations
inside the four $4$-balls with non-empty real part
and globally invariant by complex conjugation
of $U(p) \subset (\CCC^*)^2$  simultaneously.

We shall distinguish
perestroikas in which only simple points
are involved and perestroikas in which $k$-points are involved.

{\bf 2.a Smoothed diffeotopic perestroikas defined on simple points.}

Let us consider
smoothed diffeotopic perestroikas defined on simple points.

Recall that
given $a$ a crossing of $\HH_{j} \cup L$,
we say
that $p_l$ the critical point of $\HH_l$ $j+1 \le l \le m$,
associated to $a$ is
simple point if:
\been
\item
$l \ge j+1$,
\item
any critical point $p_{k}$,
$j+1 \le k \le l$
is equivalent to $p_{j+1}$.
\enen
\vskip0.1in

Given $a$ a crossing of $\HH_{j}\cup L$ ($L \approx \CCC l_{j}$)
$1 \le j \le m$
denote $(x_0,y_0)$
the simple point of $\HH_m$,
associated the crossing $a$.
Let $p \in {\mathcal  P}_m$, such that
$ a, (x_0,y_0) \in \rho^m (D(p) \ti U_{\RRR}^2)= U(p) \cap (\RRR^*)^2$
Real branches in $U(p)$ are the real branches involved in
the  topological meaning of $(x_0,y_0)$ and
the real branches which result from the perturbation of
the crossing $a$.
Let us remark that if the relative position of real
branches which result from the perturbation of
the crossing $a$ is preserved, then relative position of the real
real branches involved in
the  topological meaning of $(x_0,y_0)$ is also preserved.
Otherwise, according to Lemmas \mrf{l:l1p7},  \mrf{l:l2p7},  \mrf{l:l3p7}
and Lemma  \mrf{l:l1typ1},
it would lead with contradiction  with the fact that we consider
only curves of type $I$ that is only Morse modification in
the direction coherent  with complex orientation.

As already noticed, in most perestroikas
imaginary points are involved.
Considering together imaginary and real points
involved in these perestroikas,
we shall define smoothed diffeotopic perestroikas
of perestroika  along a path $h$
in the space  ${\mathcal  C}_{I,g})$ (see section \mrf{susu:irred,I})
and get relative position of real branches in $U(p)$ after
such diffeotopic perestroikas.
\vskip0.1in
\bele

\mlb{l: Lemma Aty1}
Let $\pi$ be either a
direct or an inverse self tangency
smoothed diffeotopic perestroika or
a triple point with imaginary branches smoothed
diffeotopic perestroika along the path
$$S :[0,1] \to \RRR {\mathcal  C}_{m}$$
$S(0)=\RRR \HH_m$,
$S(1)=\RRR {\AA}_m$.
Let
$c \in\rho^m(D(p)\ti U_{\RRR}^2)=U(p) \cap (\RRR^2)$,
with $p \in {\mathcal  P}_m$
(which may be $a$ or $(x_0,y_0)$) a real point involved in $\pi$.

Then,
under $\pi$ relative position of real branches in $U(p)$ changes
as follows.\\
One can always assume that
the crossing $a \in \HH_{j} \cup L$
is the real (smoothed) point involved in the perestroika,
the relative position of real branches involved in
the local topological meaning of $(x_0,y_0)$ is changed as follows.
Draw a dotted (unknotted)
line between the two real branches involved in the
local topological meaning of the critical point $(x_0,y_0)$
The diffeotopic perestroika moves the real component
around $a$ according to the perestroika and
leaves the two real branches linked by
the unknotted dotted line.
\enle
\vskip0.1in
{\bf proof:}

(a)\\
Let $\pi$ be a
direct self tangency perestroika
or inverse self tangency perestroika
in the negative direction.

Denote $(x_0,y_0)$
the critical point of $B_{j+1}$ associated to the crossing $a$
such that\\
$ a, (x_0,y_0) \in \rho^m (D(p) \ti U_{\RRR}^2)= U(p) \cap (\RRR^*)^2$.

Denote  $c, c'$ the two real points involved
in $\pi$.\\
Let $p, p' \in {\mathcal  P}_m$
such that
$c \in \rho^m(D(p) \ti U_{\RRR}^2)=U(p) \cap \RRR^2$,\\
$c' \in \rho^m(D(p') \ti U_{\RRR}^2)=U(p') \cap \RRR^2$
Then, imaginary  points $i,i'$ involved in $\pi$
are such that:
$i \in\rho^m(D(p) \ti U_{\CCC}^2)=U(p)$,\\
$i'\in \rho^m(D(p') \ti U_{\CCC}^2)=U(p')$.

In an $\e$-neighborhood $U(p)$, $p \in {\mathcal P}_m$,
the complex point set of the Harnack curve
$\HH_{j+1}$ can be considered  as the image of a smooth  section of a
tubular fibration
$N \to (\CCC V\bk \{a\}))$ where
$V = (\HH_{j+1} \cup L) \cap U(p)$
(i.e, from Morse Lemma, in a neighborhood of $a$, $V$ looks like
the intersection of two real lines in the point $a$)
and $N$ is an $\e$-tubular neighborhood of $V$.

Assume $c=a$, (resp $c=(x_0,y_0)$),
relative location of the real branches in $U(p)$
is changed after diffeotopic perestroika
in such a way that after diffeotopic perestroika
the complex point set
is  the image of a smooth  section of a
tubular fibration
$N \to (\CCC V\bk \{a\}))$.
Therefore, since $(x_0,y_0)$
(resp, $a$) is not a real point involved in the
perestroika, the
complex point $i$ lies in complex 4-ball around $(x_0,y_0)$
which contains real branches involved in
the topological meaning of $(x_0,y_0)$ (resp,  which results from
the perturbation of the singular crossing $a$).
Hence, it follows the relative position
of the real branches in $U(p)$.
Obviously, the relative position
of the real branches inside $U(p')$ may be deduced from the
same argument.
\vskip0.1in
(b)\\
Consider now $\pi$ a triple point perestroika with
imaginary branches.
Obviously, triple point with imaginary branches
perestroika in the positive (resp, negative) direction
may appear only after direct (resp, inverse) self tangency perestroika.
Such perestroika provides two imaginary smoothed points,
each of
which belongs to an $\e$-neighborhood $U(p)$
$p \in {\mathcal P}_m$.

Assume an oval ${\mathcal O}$
is involved in a smoothed
triple point perestroika with imaginary branches.
Then,
${\mathcal O}$ intersects
$\cup_{ p \in {\mathcal P}_m } U(p)$
in four of its neighborhood $U(p)$.
Each of these four neighborhoods $U(p)$ contains one imaginary point
involved in the perestroika:
two of them contain one of the two
imaginary points
required before the perestroika,
and the two others contain one of the two imaginary points
which appear after the perestroika.

Q.E.D

\vskip0.1in

{\bf 2.b Smoothed diffeotopic perestroikas defined on 2-points
and 3-points.} \\

We shall now consider crossings
associated to  $3$-point of $\HH_m$ and $2$-point
in case of odd $m$ in the construction of $\HH_m$.
Let us prove the following Lemma.
\vskip0.1in
\bele
\mlb{l: Lemma Bty1}
\vskip0.1in
A smoothed  strong real triple point
diffeotopic perestroika $\pi$
may be experienced  along the path
$$S :[0,1] \to \RRR {\mathcal  C}_{m}$$
$S(0)=\RRR \HH_m$,
$S(1)=\RRR {\AA}_m$  only
in a neighborhood which contains a crossing $a \in \HH_{j} \cup L$
associated to a triple-point or
a double-point of the curve $\HH_m$.
\enle
\vskip0.1in
{\bf proof:}

Recall that given
$a$ a crossing of one curve $\HH_{j-1} \cup L$
(with  $j$ even) associated to a $3$-uple $(p_j,p_{j+1},p_{j+2})$
of $\HH_{j+2}$, $j+2 \le m$,
we consider the following natural order on the $3$-uple
$(p_j,p_{j+1},p_{j+2})$:
$p_j$ has order $j$,
$p_{j+1}$ has  order $j+1$,
$p_{j+2}$  has order $j+2$.
\vskip0.1in

{\bf 1}Let us first consider the $2$-uple $(p_j,p_{j+1})$.\\
{\bf 1.a}
Let $b_{j}(p_{j})=c_0$
with $(p_j,p_{j+1})$
a $2$-uple of the first or the second kind.
It follows from Lemma \mrf{l:l2p7} of chapter \ref{ch:mpr}
that
in case of $2$-uple of the first kind:
when $c=c_0$,
the non-empty positive oval touches a positive
outer oval;
and from  Lemma \mrf{l:l3p7} of chapter \ref{ch:mpr}
that in case of $2$-uple of the second kind:
when $c=c_0$,
the non-empty positive oval touches a negative
outer oval.\\
{\bf 1.b}
Let $b_{j+1}(p_{j+1}) =c_0$
with $(p_j,p_{j+1})$
a $2$-uple of the first or the second kind.
When $c=c_0$, the one-side component
of the curve $\HH_{j+1}$ touches itself.
Then,
as $c$ increases from $c_0$ to $c_0 + \e$,
one negative oval ${\mathcal  O}^-$ disappears.

Besides, when consider together local topological meaning
of $p_j$ and $p_{j+1}$ three real branches are involved
which belong  respectively to a part
of $\HH_{j-1}$ isotopic to the real line
 $L_{j-1} \approx \RRR l_{j-2}$,
a part of $\HH_{j}$
isotopic to the real line $L_j \approx \RRR l_{j-1}$,
a part of $\HH_{j+1}$
isotopic to the real line $L_{j+1} \approx \RRR l_{j}$.
Besides,
the branches which are respectively
isotopic to a part of the line $L_j \approx \RRR l_{j-1}$ and
isotopic to a part of the line $L_{j+1} \approx \RRR l_{j}$ belong
to an oval ${\mathcal  O}$ (negative in case of $2$-uple of the first kind
and positive in case of $2$-uple of the second kind)
which  results from the perturbation of
$\HH_{j} \cup L$, $L_{j+1} \approx \CCC l_j$, around two crossings.

Consider $\CCC \HH_{j+1}$ inside
 $\rho^m( D_{j, j-1} \cup  D_{j+1, j} \ti U_{\RRR}^2) \subset \RRR P^2$
where \\ $D_{j+1, j }
 =\{ x \ge 0, y \ge 0 | j \le x+ y \le j +1 \}$,\\
$D_{j, j -1} =\{ x \ge 0, y \ge 0 | j-1 \le x+ y \le j +1 \}$.

It  is an easy consequence of the Lemma \mrf{l:l2p7} and
\mrf{l:l3p7} of chapter \ref{ch:mpr}
that one can consider the change from $p_j$ to $p_{j+1}$
as a "jump" of the branch
of $\HH_{j-1}$ isotopic to the real line $L \approx \RRR l_{j-2}$
from
 $\rho^m( D_{j, j-1} \ti U_{\RRR}^2)$
to
 $\rho^m( D_{j+1, j} \ti U_{\RRR}^2)$, namely
as a smoothed  triple point perestroika.

\vskip0.1in
{\bf 2}
Consider now $3$-uple  $(p_j,p_{j+1}, p_{j+2})$ of the first
and  second kind.\\
Let $b_{j+2}(p_{j+2}) =c_0$,\\
{\bf 2.a}\\
in case of $3$-uple of the first kind, it follows
from Lemma \mrf{l:l2p7} of chapter \ref{ch:mpr}, that when $c=c_0$
one positive outer oval of
the curve $\HH_{j+2}$
touches another outer positive oval;\\
{\bf 2.b}\\
in case of $3$-uple of the first kind, it follows
from Lemma \mrf{l:l3p7}  of chapter \ref{ch:mpr}, that  when $c=c_0$
one  negative inner oval  of
the curve $\HH_{j+2}$
touches another negative oval.
\vskip0.1in

It is obvious that the  previous description of the strong triple
point perestroika defined on the $2$-uple $(p_j,p_{j+1})$
extends to the definition of a strong triple point
perestroika defined on the $3$-uple $(p_j,p_{j+1},p_{j+2})$ in
such a way that
relative position and orientation
of the two  real branches involved in the topological meaning of $p_{j+2}$
does not change under the perestroika.

\vskip0.1in
Therefore, the Lemma is straightforward consequence of the
previous study of the $2$-uple and $3$-uple.
The  varieties $\CCC {\HH}_m$ and $\CCC {\AA}_m$ are diffeotopic;
under the perestroika the diffeotopy
is preserved locally in any neighborhood $U(p)$, $p \in {\mathcal  P}_m$.
\vskip0.1in
Q.E.D

\vskip0.1in
As explained in section \mrf{susu:irred,I},
any smoothed perestroika $\pi$
required a given location and orientation of real branches.
In particular,
smoothed triple-point perestroika is possible only
in a neighborhood of an oval and of three real branches close to this oval
with relative location and orientation as described in
section \mrf{susu:irred,I}.
Namely, the real branches involved  in the local
topological meaning of $p_{j}$ may be
involved in a triple point
diffeotopic perestroika  only if they have been
reproached by previous perestroikas.
It is easy to deduce
from  the relative location  of branches involved in the
topological meaning of critical point
that in case of $2$-uple of the first kind,
the singular situation  $b_{j+1}(p_{j+1}) =c_0$
is possible if and only if two positive
outer ovals coalesce with the one-side component of the curve.
In case of $2$-uple of the second kind,
the singular situation  $b_{j+1}(p_{j+1}) =c_0$
is possible if and only if two negative inner ovals
coalesce with the one-side component of the curve.
We shall precise this remark in
the third part.

{\bf 3. Method which provides $(\CCC P^2,\CCC {\AA}_m)$
from $(\CCC P^2,\CCC \HH_m)$}
\vskip0.1in

Denote ${\mathcal  S}_{I,m}$
the set
${\mathcal  C}_{I,m} \bk {\mathcal  D}_m $ of smooth
curves of degree $m$ and type $I$.

We  shall define a family of moves
$\phi :\HH_{m} \to  {\mathcal  S}_{I,m} $
with the property
 that in these moves
the set $\cup_{i=1}^I a_i=\cup_{n=1}^m A_n$
of crossing-points
perturbed in the construction of $\HH_{m}$ participate as follows.
Moves of real components of $\HH_{m}$ correspond
to an other choice of perturbation on these crossings.
\vskip0.1in
Briefly the method of our proof is the following:

We shall apply
successively Lemma
\mrf{l: Lemma Aty1} and \mrf{l: Lemma Bty1}
inside the sets:
$$\rho^{m}((D_{1,0} \cup D_{2,1}) \ti U_{\CCC}^2)$$

$$\rho^{m}((D_{2,,1} \cup D_{3,2}) \ti U_{\CCC}^2)$$
$$\rho^{m}((D_{3,2} \cup D_{2,1}) \ti U_{\CCC}^2)$$
$$...$$
$$\rho^{m}((D_{m-i+1,m-i} \cup D_{m-i+2,m-i+1}) \ti U_{\CCC}^2)$$
$$\rho^{m}((D_{m-i+3,m-i+2 } \cup D_{mm-i+2,m-i+2}) \ti U_{\CCC}^2)$$
$$...$$
$$\rho^{m}((D_{m-1,m-2 } \cup D_{m-2,m-3}) \ti U_{\CCC}^2)$$
$$\rho^{m}((D_{m,m-1} \cup D_{m-1,m-2}) \ti U_{\CCC}^2)$$
to describe the track $S$ on $\RRR P^2$
of a diffeotopy $h_t$ of $\CCC P^2$, $t \in [0,1]$,
$h(0)= \CCC \HH_m$, $h(1)=\CCC {\AA}_m$ with
${\AA}_m$ a curve of type $I$ and degree $m$.

\vskip0.1in
We shall now proceed to precise argument.
\vskip0.1in

{\it 3.a Preliminaries}
\vskip0.1in
Recall that, according to chapter \ref{ch:mpr},
for any point $a$ perturbed in the construction of $\HH_m$,
there exists $p \in {\mathcal  P}_m$
(where $\sharp {\mathcal  P }_m = \frac {m (m-1)} 2 $)
such that $a \in U(p)$ where
$U(p)$ is the $\e$-tubular neighborhood
defined from the interior $\G^0$ of a face $\G$
of the triangulation of $T_m$.
According to Lemma \mrf{l:l1typ1}  $\RRR \HH_m$
(resp, $\RRR {\AA}_m$) is  a smoothed
generic immersion of the  circle $S^1$
into the plane $\RRR ^2$,
$\RRR \tilde{\HH}_m$ (resp, $\RRR \tilde{{\AA}}_m$)
with  $n \ge \frac {(m-1)(m-2)} 2$ double points and Whitney index
$\frac {(m-1)(m-2)} 2 + 1$.
The smoothing is such that
$\frac {(m-1)(m-2)} 2$ double points are
smoothed in $\RRR ^2$ by Morse modification
in the direction coherent to the complex orientation.
We shall
denote by ${\mathcal  G}_m$ the set of the
$g_m =\frac {(m-1)(m-2)} 2$ double points of $\RRR \tilde{\HH}_m$
smoothed in $\RRR ^2$ by Morse modification
in the direction coherent to the complex orientation.

According to the part 1, we shall assume
that the set ${\mathcal  G}_m$
of double-points of
$\tilde{\RRR \HH_m}$
has the following properties:
\been
\item
\lb{i:id}
any two distinct points of  ${\mathcal  G}_m$
belong to two distinct
neighborhoods $U(p)$ with $p \in {\mathcal  P}_m$
\item
\lb{i:ide}
there exists a
bijective correspondence  between
points of ${\mathcal  G}_m$
and  the $g_m$ ovals
${\mathcal  O}$ of $\HH_m$ which do not
intersect the line at infinity.
\enen
\vskip0.1in
Properties (\ref{i:id}) and (\ref{i:ide})
are straightforward
consequences of the patchworking construction
(see for example Corollary \mrf{c:patch} and its proof)
that any oval ${\mathcal  O}$  of $\HH_m$ which does not
intersect the line at infinity,
intersects two  neighborhoods  $U(p)$
defined from  a face $\G \subset l_{j}$ ,$1\le j \le m-1$,
$l_{j}=\{ (x,y) \in (\RRR^+)^2 | x +y =j \}$, namely
two neighborhoods $U(p)$
of global order $j+1$.

Besides, an oval ${\mathcal  O}$ may intersect:
\been
\item
either
$4$ neighborhoods $U(p)$ of the union
$\cup_{p \in {\mathcal  P}_m} U(p)$ such that
two neighborhoods $U(p)$ are
of global order $j$, $3 \le j \le m-1$,
the two others are of global order $j-1$, $j+1$
\item
or $3$ neighborhoods $U(p)$ of the union
$\cup_{p \in {\mathcal  P}_m} U(p)$ such that
two neighborhoods $U(p)$ are
of global order $m$, and the other is of global order $m-1$.

\enen

Therefore, one can choose
the set ${\mathcal  G}_m$ such that
any point of  ${\mathcal  G}_m$
belongs to
a neighborhood $U(p)$, $p \in {\mathcal  P}_m$
of global order $j$, $2 \le j \le m-1$.
\vskip0.1in
Consider
the family $\tilde{U}(p) \supset U(p)$, $p \in {\mathcal  P}_m$
$\cup_{p \in {\mathcal  P}_m} \tilde{U}(p)=\CCC P^2$
verifying the properties
(\ref{e:p1}), (\ref{e:p2}), (\ref{e:p3}), (\ref{e:p4}).
Using the patchworking construction, according  to
the moment map
$\mu: \CCC T_m \to T_m$ and the natural surjection
$\rho^m: \RRR^+ T_m \ti U_{\CCC}^2 \to \CCC T_m$
(see the preliminary section ), one can  choose for example:
\been
\item
 for any $U(p)$  of global order $2 \le j \le m-1$,\\
$\tilde{U}(p)=\rho^m( S \ti U_{\CCC}^2) $
where $S$ is the  unique square $S$ with vertices
 $(c,d), (c+1,d), (c,d+1), (c+1,d+1)$, $c+d+1=j$ such that
 $U(p) \subset \tilde{U(p)}$.

\item
for any $U(p)$  of global order $m$,
$\tilde{U}(p)=\rho^m( K \ti U_{\CCC}^2)$
where  $K$ is the  polygon  with vertices
$(c,d), (c+2,d), (c,d+2)$, $c+d+1=m-1$ such that
 $U(p) \subset  \tilde{U}(p)$.
\enen
\vskip0.1in
{\it 3.b Statement of the Method}
\vskip0.1in

We shall now deduce a lifting of the path  $\tilde{h}$
in the space of immersions of the circle into $\RRR ^2$ with index
$g_m= \frac {(m-1) (m-2)} 2$ in the space $\RRR {\mathcal  C}_{m}$
and describe the topological pair $(\CCC P^2, \CCC {\AA}_m)$.
Such lifting is defined by smoothed diffeotopic perestroikas
on $\RRR \HH_m$ inside $\cup_{p \in {\mathcal P}_m} \tilde{U}(p)$.
\vskip0.1in

Recall that,
according to parts 1 and 2 of our proof,
we shall consider only
direct or inverse self tangency
perestroika, triple point with imaginary branches perestroika,
strong real triple point
perestroika.
Obviously,
in any smoothed  diffeotopic perestroika $\pi$ undergone along the lifting
$S$ of $\tilde{h}:\RRR \tilde{\HH}_m \to \RRR  \tilde {\AA}_m$,
at least one oval ${\mathcal  O}$
of $\RRR \HH_m$ which does not intersect the line
at infinity is involved.
As already noticed, in most perestroikas
imaginary points are involved.
We shall consider the function $J$ counting increment and decrement
of $\a$-point (double point of different halves)
and $\b$-point (double point of one of the halves)
along the path $S: [0,1] \to \RRR  {\mathcal  C}_{m}$.
\vskip0.1in
Let ${\mathcal  O}$ be an oval of $\HH_m$.
It intersects $\cup_{p \in {\mathcal P}_m} \tilde{U}(p)$
in four neighborhood $\tilde{U}(p)$, $p \in {\mathcal P}_m$:
 $U(p_1)$ of global
order $i$ with $2 \le i \le m-1$,
, $U(p_2)$ and $U(p_3)$
of global order $i+1$, $U(p_4)$ of global order $i+2$.

We shall, from the  Lemma \mrf{l: Lemma Aty1}
and \mrf{l: Lemma Bty1},
get moves of the real branches contained in $U(p_i)$, $1 \le i \le 4$,
under any smoothed diffeotopic perestroika $\pi$.

\vskip0.1in
We shall prove that smoothed diffeotopic perestroikas
on the real components of $\HH_m$
of which smoothed points are
points of global order $2 \le j \le m-1 $
( in  the union of neighborhoods
$\tilde{U}(p)$, $p \in {\mathcal  P}_m$
of global order $j$, $2 \le j \le m-1$), induce smoothed
diffeotopic perestroikas
on the real component of $\HH_m$ in  neighborhoods
$\tilde{U}(p)$ $p \in {\mathcal  P}_m$
of global order $j+1$.\\
We shall call this process the step $j$ of the  lifting,
and define
the lifting of $\tilde{h}: \RRR \tilde{\HH}_m \to \RRR \tilde{{\AA}}_m$
by induction on the global order on the set of points $a_i, i \in I$
and the associated $k$-uple.
\vskip0.1in

\vskip0.1in
1) Let us first consider
direct and inverse self-tangency
and triple point with imaginary branches diffeotopic perestroika.
According to Lemma \mrf{l: Lemma Aty1}
one can assume,
without loss of generality
that the smoothed real points involved
in triple point with imaginary branches diffeotopic perestroika
are the two crossings which belong to the neighborhoods $U(p_2)$, $U(p_3)$
of global order $i \le m-1$.
In such a way,
according to  the
Lemma \mrf{l: Lemma Aty1}, we get
description of relative position
of real branches in  $U(p_2)$ and $U(p_3)$
and the smoothed points of global order $j+1$ which may be involved
in a next diffeotopic perestroika is deduced

Besides, along the path $S$, as the curve $\HH_m$ experiences
various diffeotopic perestroikas
$J(\a)$ and $J(\b)$ increases or decreases.
It follows from the
Lemma \mrf{l: Lemma Aty1} that imaginary points involved
in a perestroika $\pi$
belong to neighborhoods $U(p_i)$, $1 \le i \le 4$,
in such a way that diffeotopic  perestroikas
in neighborhoods
of global order $j$ provide or smooth the
imaginary points involved in  diffeotopic perestroika
in neighborhoods $\tilde{U}(p)$ of order $j+1$.
\vskip0.1in
2)Consider now $\pi$  a strong triple diffeotopic perestroika.

From the
Lemma \mrf{l: Lemma Bty1},
 strong triple point smoothed diffeotopic perestroika $\pi$
requires a neighborhood $\tilde{U}(p)$
which contains a crossing $a$ associated to
a triple-point or
a double-point of the curve $\HH_m$.
It is an easy consequence of the
Lemma \mrf{l: Lemma Bty1},
that, according to our notations,
the $k$-uple ($k \in \{ 2,3 \}$)
is of the form $(p_{i},p_{i+1},p_{i+2})$, $i$ even
(resp, $(p_{i},p_{i+1})$, $i$ even)
with $p_{i} \in U(p_1)$
belongs to
the neighborhood $\tilde{U}(p_1)$ of global
order $i$ with $i+1 \le j \le m-1$.
As already noticed in the proof of Lemma \mrf{l: Lemma Bty1},
the restriction of $3$-uple
$(p_{i},p_{i+1},p_{i+2})$ to
$(p_{i},p_{i+1})$ is sufficient to define $\pi$;
 relative position  of real branches involved
in the  topological meaning of $p_{i+2}$ under $\pi$
follows by induction.\\
Besides,the diffeotopic perestroika $\pi$  is possible
only if orientation
and relative location of the real branches
in the neighborhoods $U(p_1)$, $U(p_2)$, $U(p_3)$, $U(p_4)$
is one of the two required (see section \mrf{susu:irred,I}).
It follows from the
Lemma \mrf{l: Lemma Bty1}
except around real branches involved in $\pi$,
the  situation remains the same
  inside
$\cup_{ 1 \le  i  \le 4} \tilde{U}(p)$
 under diffeotopic perestroika.

In such a way, at the end of the step $j$
of the lifting, we get an irreducible curve of degree $j$
and type $I$ inside
$\rho^{j}(T_{j} \ti U_{\CCC}^2)$.
This irreducible curve of degree $j$
has no real double-point singularity but
may have an even number of imaginary  double-point
singularity.
The step $j+1$ of the lifting can be interpreted as
the perturbation of
the union of this curve of degree $j$ with a line
-namely, according to the Lemma \mrf{l:l1typ1} and its proof
as the way to move the  union
of the $g_{j}$ 1-handles of
$\CCC \HH_{j} \approx S^2_{g_j}$  with
the $(j-1)$ more 1-handles of
$\CCC \HH_{j+1} \approx S^2_{g_{j}+(j-1)}$-.
\vskip0.1in

We have considered only
neighborhoods $\tilde{U}(p)$ of global order $j$, $2 \le j \le m-1$, and
$\sharp{\mathcal  P}_m -g_m=m-1$.
Nonetheless,
the lifting of $\tilde{h} : \RRR \tilde{\HH}_m \to \RRR \tilde{{\AA}}_m$
is entirely defined  by induction on the global order, and
we get the topological pair
$(\CCC P^2, \CCC {\AA}_m)$
from the method proposed.
Indeed,
consider the $m-2$ ovals ${\mathcal  O}$  which intersect
two neighborhoods $\tilde{U}(p)$
of global order $m$, and one neighborhood $\tilde{U}(p)$
of global order $m-1$.
(The union of
these $m-2$ ovals intersects any neighborhood $\tilde{U}(p)$,
$p \in {\mathcal  P}_m$
of global order $m$.)
Any of these  m-2 ovals ${\mathcal  O}$ is involved in a
diffeotopic perestroika
implied (in the chain of diffeotopic perestroika) by a
diffeotopic perestroika on
neighborhood of order $m-1$.
Hence,
the lifting of $\tilde{h}: \RRR \tilde{\HH}_m \to \RRR \tilde{\AA}_m$
follows by induction.

Consequently, we get description of $\CCC {\AA}_m$ in any
neighborhood $\tilde{U}(p)$, $p \in {\mathcal  P}_m$
and thus in ${\cup_{p \in {\mathcal  P}_m} \tilde{U}(p)}= \CCC P^2$.

According to Lemma \mrf{l:l1typ1} and Lemma \mrf{l:l2typ1},
we have defined a lifting of the path
$\tilde{h} : \RRR \tilde{\HH}_m \to \RRR \tilde{{\AA}}_m$
in the space of generic immersion of
the space of immersions of a circle into the plane to
$S: \RRR \HH_m \to \RRR {\AA}_m$ which is
the track on $\RRR P^2$
of a diffeotopy $h_t$ of $\CCC P^2$, $t \in [0,1]$,
$h(0)= \CCC \HH_m$, $h(1)=\CCC {\AA}_m$.
Q.E.D
\vskip0.1in

\subsection{ Description of pairs $(\CCC P^2,\CCC {\AA}_m)$
up to conj-equivariant isotopy }

In this section, we shall first give a combinatorial method
to describe pairs $(\CCC P^2,\CCC {\AA}_m)$, up to conj-equivariant
isotopy, where ${\AA}_m$ is a smooth
curve of degree $m$ and type $I$.
Then,
we extend in Theorem \mrf{t:th2typ1}
the description of the pair $(\CCC  P^2, \CCC \HH_m)$
stated in Theorem \mrf{t:theo1} for Harnack curve
to any smooth curve of type $I$.

{\bf Combinatorial Description of pairs $(\CCC P^2,\CCC {\AA}_m)$}
\vskip0.1in
In this part, we shall detail the method given in the third part of the
proof of the Proposition \mrf{p:prop1typ1}
and describe
the family of moves
$\phi:\HH_{m} \to  {\mathcal  S}_{I,m} $
where ${\mathcal  S}_{I,m}={\mathcal  C}_{I,m} \bk {\mathcal  D}_m$.
\vskip0.1in
For any curve ${\AA}_m$ of degree $m$ and type $I$,
we have defined a lifting of the path
$\tilde{h}:\RRR \tilde{\HH}_m \to \RRR \tilde{\AA}_m$
which is the track $S$
on $\RRR P^2$
of a diffeotopy $h_t$ of $\CCC P^2$, $t \in [0,1]$,
$h(0)= \CCC \HH_m$, $h(1)=\CCC {\AA}_m$.

(Obviously,
the path is not necessarily lifted continuously in the space
of real algebraic curves in the sense that we can not affirm
that any immersion of the circle into $\RRR ^2$ with index $g_m$
along an arbitrary  path from $\RRR \tilde{\HH}_m$ to
$\RRR \tilde{\AA}_m$
lift to a curve of degree  $m$ and type $I$.)
\vskip0.1in

In this part, we shall describe
the family of moves
$\phi:\HH_{m} \to  {\mathcal  S}_{I,m} $

We shall work with the terminology introduced in section \mrf{susu:cji},
(see Proposition \mrf{p:prop8}, and also Theorem \mrf{t:theo1}
and consider the set  of points $\cup_{i \in I} a_i = \cup_{n=1}^m A_n$
where $A_n$ denotes the set of points perturbed in a maximal
simple of deformation of $\HH_n$.
Recall that
$A_{2k}=\{a_1,...,a_{2k-1},...,a_{4k-2} \}$ where
$a_1,...a_{2k-1}$ are
the crossings  of $\HH_{2k-1}\cup L$
and $a_{2k},...,a_{4k-2}$ are the crossings
of $\HH_{2k; 1}$ where $\HH_{2k;t}$, $t \in [0,1]$
is a maximal simple deformation of $\HH_{2k}$;
$A_{2k+1}=\{a_1,...,a_{2k} \}$
where $a_1,...a_{2k}$ are
the  crossings of $\HH_{2k} \cup L$.
In such a way, we assign to the set
$A{2k}=\{a_1,...,a_{2k-1},...,a_{4k-2} \}$ the global order $2k$;
and to the set
$A_{2k+1}=\{a_1,...,a_{2k} \}$ the global order $2k+1$.

\vskip0.1in

Fold the real point set of $\HH_m$ in such a way
for any $ 1 < j <m$, the two real points of $\HH_{j+1}$
resulting from the perturbation of one crossing of  $\HH_{j} \cup L$
are glued each other in the initial crossing.
Then, fold the curve $\HH_m$ in any
other points perturbed in its construction;
namely, glue the two branches involved in the local topological meaning
of any critical point $a_i$, $2k \le i \le  4k-2$ of Harnack polynomials
$B_{2k}$, $2k \le m$.\\
The curve ${\AA}_m$ results from desingularization
of all these points with multiplicity  $2$ compatible with the way $S$
goes through the discriminant hypersurface.
It is equivalent to unfold in a way different from the one
given $\HH_m$.
Restriction on the unfolding is given by
the set of smoothed
diffeotopic perestroika experienced along the path $S$ from
$\RRR \HH_m$ to  $\RRR {\AA}_m$.
Along the path only
triple-point with imaginary branches,
strong triple point,
inverse and direct self-tangency diffeotopic perestroika
are encountered.
\vskip0.1in

Moreover, according to the third part of the
proof of the Proposition \mrf{p:prop1typ1},

\been
\item
desingularization (unfolding) are done in
increasing global order on the set of points perturbed
in the construction of $\HH_m$.
Locally, any choice of a desingularization  is compatible
with one diffeotopic perestroika or lets fix  the real components of $\HH_m$.
\item
Imaginary points are involved in
one diffeotopic perestroika $\pi$ along the path $S$
in such a way that if  $J_{\pi}(\gamma)= 2$,
  $\gamma \in \{ \a, \b \}$,
then it provides imaginary points
involved in
-at least one and at most two-
next diffeotopic perestroikas $\pi'$ along $S$
such that $J_{\pi'} (\gamma) =-2$. \\
From the
Lemma \mrf{l: Lemma Aty1},  any  real or imaginary point involved in
a diffeotopic perestroika  belongs to
a neighborhood
$\tilde{U}(p) \supset U(p)=\rho^m ( D(p) \ti U_{\CCC}^2)$,
$p \in {\mathcal  P}_m$,
and the next perestroikas are provided by move of real branches
in $U(p)$.
Consequently, in the chain of moves,
imaginary singular points are smoothed in such a way a the end of
the procedure one gets a curve without singular points.
\enen

\vskip0.1in

Let us summarize  properties of
the path $S:[0,1] \to \RRR {\mathcal  C}_{m}$
in the following Theorem \mrf{t:th1typ1}.

Recall that given the triangle $T_m$,
the real projective space $\RRR P^2$ may be
deduced  from the square $T_m^*$
made of $T_m$ and it symmetric copies $T_{m,x}=s_x(T_m)$, $T_{m,y}=s_y(T_m)$
$T_{m,xy}=s(T_m)$ where $s_x, s_y, s=s_x \circ s_y$ are reflections
with respect to the coordinates axes.
Given a triangle $T$ of the set $T_m$,$T_{m,x}$, $T_{m,y}$, $T_{m,xy}$
$s(T)$ denotes the symmetric copy of $T$.
\vskip0.1in
\beth
\mlb{t:th1typ1}
\vskip0.1in
{\it Let ${\AA}_m$ be a curve of degree $m$
 and type $I$.
Let $\HH_m$ be the Harnack curve of degree $m$.
Then, up to conj-equivariant isotopy of $\CCC P^2$,
there exists a path  $S :[0,1] \to \RRR {\mathcal  C}_{m}$
with $S(0)=\RRR \HH_m$, $S(1)=\RRR {\AA}_m$
which crosses
the discriminant hypersurface $\RRR {\mathcal  D}_m$
in curves with crossings, real points of intersection
of a real branch and two conjugate imaginary branches,
strong triple-points,
direct and inverse self-tangency points.

Curves along $S$
may be deduced one from another from
smoothed diffeoopic perestroikas with the property
that any smoothed diffeotopic perestroika
along $S$ defines double-points involved in
next smoothed diffeotopic perestroikas along $S$.
\vskip0.1in
Furthermore,
\been
\item
any real ordinary double point of curves along $S$ is
a crossing which
belongs to  the set of points
perturbed in the construction of $\HH_m$
\item
any  double-point involved in a diffeotopic perestroika
is a crossing which belongs to the set of points
perturbed in the whole construction of $\HH_m$.
\enen
Besides, smoothed diffeotopic perestroikas may be
described in the patchworking scheme
as follows:
Let $T \in \{T_m,T_{m,x},T_{m,y},T_{m,xy} \}$;
\been
\item
inverse self tangency  diffeotopic perestroika
in the positive direction
acts as follows:
an oval ${\mathcal O}$ of $\HH_m$ contained in a  triangle  $T$
becomes an oval $s({\mathcal O})$
of the symmetric triangle $s(T)$,
\vskip0.1in
in the negative direction,
according to the definition of the
self tangency diffeotopic perestroika
$s({\mathcal O})$ disappears in
two real branches with opposite orientation.

Moreover, in the positive (resp negative)
direction, it requires (resp, smooths)
imaginary points
provided by a triple point with imaginary diffeotopic perestroika in the
negative direction.

\item

direct  self tangency diffeotopic perestroika
acts as follows:
it does not change the real part
but in the negative direction it provides
the imaginary points of a triple-point diffeotopic perestroika
with imaginary  conjugate branches in the positive direction;
\vskip0.1in
in the positive direction it smooths
the imaginary  point
of a triple-point diffeotopic perestroika
with imaginary conjugate branches in the negative direction.
\item
triple-point with imaginary conjugate branches
diffeotopic perestroika is as follows:
\vskip0.1in
an oval ${\mathcal O}$ of $\HH_m$ contained in a triangle  $T$
becomes an oval $s({\mathcal O})$
of the symmetric triangle $s(T)$.
\item
strong triple point perestroika
may be experienced only in an open which
contains
a $k$-point (in case $m$ even : $k=3$;
in case $m$ odd: $k=3$ or $k=2$) of $\HH_m$.
\enen
\vskip0.1in}
\enth

\bere
\vskip0.1in
\been
\item
One can easily deduce the maximal nest curve of degree $m=2k$
(which are $L$-curves and thus  have standard Arnold surfaces)
from the Harnack curve $\HH_{2k}$.\\
Indeed, choose an  orientation on  the  recursive line
$\RRR l_{2k -j}$ opposite to the one  chosen  in
the construction of $\HH_m$.
A slight perturbation of the resulting curves given by
standard modification on real double point singularity
compatible with the orientation of the resulting
curve leads to the maximal nest curve.\\
Besides,
according to the Theorem \mrf{t:th1typ1}, there
exists a path from
$\HH_{2k}$ to the maximal nest curves $<1> ^{k}$.
Such path may be described by smoothed diffeotopic perestroika,
some of which (for $k > 2$) are smoothed strong triple point diffeotopic
perestroikas.
\item
Obviously, from the method summarized in theorem \mrf{t:th1typ1},
one can deduce whether there exist  curves
of a degree $m$ and type $I$ with a given real scheme.

For example, the non-existence of the curve of degree $7$ and
real scheme $J \cup <1<14>>$ is easy to deduce.
Indeed, it can be easily proved that the existence of such curve
is in contradiction with the possible situations inside
the three neighborhoods $U(p)$, $p \in {\mathcal  P}_7$
of global order $6$ which contain crossings
associated to $2$-uple.
(More generally,
one can notice that the $2k-3$ $2$-uples of $\HH_{2k+1}$
define possible intersection of the one-side component of $\HH_{2k+1}$
with the boundary of the M\"{o}bius band  of $\RRR P^2$ embedded in
$\CCC P^2$. )
\enen
\enre
\vskip0.1in

\vskip0.1in
{\bf  Description of pairs $(\CCC P^2,\CCC {\AA}_m)$
up to conj-equivariant isotopy }
\vskip0.1in
We shall now in Theorem \mrf{t:th2typ1}
extend the results of Theorem \mrf{t:theo1}
to any smooth curve of type $I$.\\
According to Theorem \mrf{t:theo1}
there exists a finite number $I$ of $4$-balls
$B(a_i)$ globally invariant by complex conjugation;
such that,
up to conj-equivariant isotopy of $\CCC P^2$
$\HH_m \bk \cup_{i \in I} B(a_i)= \cup_{i=1}^m L_i \bk
  \cup_{i=1}^I B(a_i)$; inside any $B(a_i)$ $\HH_m$
 is the perturbation of crossing
 centered the in point $a_i$.
The following statement is a consequence
of the definition of $S: [0,1] \to \RRR {\mathcal C}_m$
given in Proposition \mrf{p:prop1typ1}
and the fact that
any real singular point which participated in diffeotopic
perestroika along $S$ belongs to $\cup_{i=1}^I a_i$.

\beth
\mlb{t:th2typ1}
{\it
Let ${\AA}_m$ be a curve of degree $m$ and type $I$.
There exists a finite number $I$
($I=1+2 ...m +  \Sigma_{k=2}^{k=[m/2]} 2k-3$)
 of disjoints $4$-balls $B(a_i)$
invariant by complex conjugation and centered in points
$a_i$ of $\RRR P^2$ such that,
up to conj-equivariant isotopy of $\CCC P^2$ :
\been
\item
 ${\AA}_m \bk \cup_{i \in I}B(a_i) =
 \cup_{i=1}^m L_i\bk \cup_{i \in I}B(a_i) $
where $L_1$,...,$L_m$ are $m$ projective lines
with $L_i  \bk \cup_{i \in I} B(a_i) \cap
L_j  \bk \cup _{i \in I} B(a_i)=\emptyset$
for any $i \not = j$,$1 \le i,j \le m$.
\item
situations inside $4$-balls $B(a_i)$ are covered by
perturbation of type $1$ or type $2$ of the crossing $a_i$.
\enen}
\enth
\vskip0.1in
{\bf proof:}
\vskip0.1in

Our proof is based on the proof of Proposition \mrf{p:prop1typ1},
Theorem \mrf{t:th1typ1}
and Theorem \mrf{t:theo1}.

According to Theorem \mrf{t:theo1},
there exists a finite number $I$ of disjoints $4$-balls $B(a_i)$
invariant by complex conjugation and centered in points
$a_i$ of $\RRR P^2$ such that
, up to conj-equivariant isotopy of $\CCC P^2$ :
\been
\item
 $ \HH_m \bk \cup_{i \in I}B(a_i) =
 \cup_{i=1}^m L_i\bk \cup_{i \in I}B(a_i) $
where $L_1$,...,$L_m$ are $m$ projective lines
with $L_i  \bk \cup_{i \in I} B(a_i) \cap
L_j  \bk \cup _{i \in I} B(a_i)=\emptyset$
for any $i \not = j$,$1 \le i,j \le m$.
\item
situations inside $4$-balls $B(a_i)$ are perturbations of type 1
of crossing.
\enen

Let $\CCC \HH_m^+$ be the half of $\RRR \HH_m$
which induces
orientation on the real part of $\RRR \HH_m$.
The conj-equivariant isotopy brings  $\CCC \HH_m^+ \bk
   \cup_{i=1}^I B(a_i)$
 to halves $\CCC L_i^+  \bk \cup_{i=1}^I B(a_i)$
of lines $L_i \bk  \cup_{i=1}^I B(a_i)$, $1 \le i \le m$,
 which induce an orientation on the real part $\RRR L_i$.

The path $S:[0,1] \to \RRR {\mathcal  C}_{m,I}$
$S(0)= \CCC \HH_m$, $S(1)=\CCC {\AA}_m$
can be seen as the track
on $\RRR P^2$
of a diffeotopy $h_t$ of $\CCC P^2$, $t \in [0,1]$,
$h(0)= \CCC \HH_m$, $h(1)=\CCC {\AA}_m$ with
$h_t( \CCC \HH_m) \cap U(p) \subset U(p)$ and may be described
as a family moves defined by smoothed diffeotopic perestroikas.

According to the proof of proposition \mrf{p:prop1typ1},
any real double point involved in a diffeotopic perestroika belongs to
$\cup_{i \in I}B(a_i)$.

Besides, nonetheless relative location of any part of
$\HH_m \bk \cup_{i \in I}B(a_i)$ conj-equivariant isotopic
to a part $L_i\bk \cup_{i \in I}B(a_i)$
is changed under a diffeotopic perestroika,
its orientation remains the same
before and after diffeotopic perestroika.

Therefore, for any line $L_i$,
 the half $\CCC L_i^+ \bk \cup_{i \in I}B(a_i)$
of  $L_i \bk \cup_{i \in I}B(a_i)$
which induces orientation on $\RRR L_i
 \bk \cup_{i \in I}B(a_i)$
is up to conj-equivariant isotopy a part of the half
$\CCC {\AA}_m^+$
of $\CCC {\AA}_m$ which induces orientation on $\RRR {\AA}_m$.

Hence, outside $\cup_{i=1}^I B(a_i)$, the curve ${\AA}_m$
is union of $m$ lines minus their intersection with
$\cup_{i=1}^I B(a_i)$.

Moreover,
situations inside $4$-balls $B(a_i)$ are covered by two cases
which are locally perturbations of type $1$ or type $2$ of
a crossing.

\vskip0.1in
Q.E.D
\section{Construction of Curves of type $II$}
\mlb{su:typ2}
\vskip0.1in

As already known, if an $M$-curve is constructed
curves with fewer components are easily constructed.
Besides, $M$-curves, as  curves of type $I$,
have been studied in the
section  \mrf{su:typ1}.
We shall now introduce a  method inspired from the preceding one,
which will give all $M-i$-curves of degree $m$
and therefore any curves of type $II$ from the Harnack curve $\HH_m$.
\vskip0.1in
In the previous section,
it was essential
that all curves have orientable real point set.
In order to  consider
curves with  non-orientable real set of points
we shall modify the argument
introduced in the preceding part.
We shall  call {\it perturbation
in the non-coherent direction of a crossing} of a curve
the desingularization
which induces non-orientability of the real part of the resulting curve.
\vskip0.1in
For a generic curve ${\AA}_m$  of degree $m$ and  genre $g$,
$0 < g < \frac {(m-1)(m-2)} 2$,
by {\it smoothing in the coherent direction}
of one double point of its
real point set $\RRR {\AA}_m$, we shall understand
the Morse modification in $\RRR ^2$ in the  direction
coherent to a complex orientation of ${\AA}_m$,
by {\it smoothing in the non-coherent direction}
of one double  point of its real point set $\RRR {\AA}_m$,
we shall understand the Morse modification in $\RRR ^2$
at the double point
in the direction non-coherent to
a complex orientation of ${\AA}_m$,
i.e not compatible with an orientation
of the real point set of the resulting curve.
\vskip0.1in
Let us start our study of curves of type $II$
by a Lemma analogous to the Lemma \mrf{l:l1typ1}
stated for curves of type $I$.
\vskip0.1in

Call {\it smoothed generic immersion of the circle $S^1$
into the plane $\RRR ^2$}, the smooth submanifold of $\RRR P^2$
deduced from the generic immersion of the circle by modification
at each real double point which is
either a Morse modification in $\RRR ^2$
coherent or non-coherent to the complex orientation,
or the Morse modification in $\RRR P^2$
which associates to the double point of  $\RRR ^2$
two points of the line at infinity of $\RRR P^2$.

A generalization of Whitney's theorem to the case of real algebraic
curves  will be provided by the following Lemma:
\vskip0.1in
\bele
\mlb{l:l1typ2}
\vskip0.1in
{\it
Let ${\AA}_m$ be a smooth curve of degree $m$ and type $II$
with non-empty real points set $\RRR {\AA}_m$
then  $\RRR {\AA}_m  \subset \RRR P^2$
is a smoothed generic immersion of the circle $S^1$
into the plane $\RRR ^2$
with \\
$\frac {(m-1)(m-2)} 2 \le n \le \frac {(m-1)(m-2)} 2 + [\frac {m} 2]$,
double points and Whitney index $\frac {(m-1)(m-2)} 2 + 1$.
The smoothing is such that
$\frac {(m-1)(m-2)} 2$ double points are
smoothed in $\RRR ^2$ by Morse modification
and at least one Morse modification is in the non
coherent direction.}
\enle
\vskip0.1in
\bere
One can associate to any smooth curve of type $II$
a three-dimensional  rooted tree.
\enre

{\bf proof:}
\vskip0.1in
Our proof makes use of properties of the complex point set
$\CCC {\AA}_m$ embedded in $\CCC P^2$.
It is deduced from an argument similar to the one given for
curves of type $I$.

Consider the usual handlebody decomposition of
$\CCC P^2= B_0 \cup B_1 \cup B$  where $B_0$,$B_1$,$B$
are respectively 0,2 and 4 handles.

The ball $B_0$ and $B_1$ meet along an unknotted  solid torus
$S^1 \ti B^2$.
The gluing diffeomorphism
$S^1 \ti B^2 \to S^1 \ti B^2$ is given by the $+1$ framing map.
In such a way,
the canonical $\RRR P^2$ can be seen as the union of a M\"{o}bius
band $\mathcal  M$ and the disc $D^2 \subset B$ glued along their boundary.
The M\"{o}bius band $\mathcal  M$ lies in $B_0 \cap B_1 \approx S^1 \ti D^2$
with $\pr {\mathcal  M}$ as the $(2,1)$ torus knot, and
$D^2 \subset B$ as the properly imbedded unknotted disc.
The complex conjugation switches $B_0$ and $B_1$ and lets fix
$\mathcal  M$, it rotates $B$ around
$D^2$.

The set of complex points
of ${\AA}_m$
is an orientable
surface of genus
$g_m= \frac {(m-1) (m-2)} 2$, i.e it is diffeomorphic to a
sphere $S^2$ with
$g_m$ $1$-handles $S^2_{g_m}$.
We shall denote $h$ the diffeomorphism
$h :\CCC {\AA}_m \to S^2_{g_m}
\subset \RRR^3$.

Assume $S^2$ provided with a complex conjugation with fixed
point set a circle $S^1$ which divides the sphere $S^2$ into two halves.

Without loss of generality, one can assume
that any disc removed from
$S^2$ and then closed up to give $S^2_{g_m}$
intersects $S^1$ and the two halves of $S^2$.
In such a way, $\RRR {\AA}_m \subset \RRR P^2$ intersects
each $1$-handle.

Fix $D^2$ the two disc of $\RRR P^2= {\mathcal  M} \cup D^2$
in such a way that the boundary circle of $D^2$ is $S^1$ and
therefore each one handle belongs to the solid torus
$S^1 \ti D^2 \supset {\mathcal  M}$.

Since the interior  $(D^2)^0$  of $D^2$ is homeomorphic to
$\RRR^2$, one can project
$\RRR {\AA}_m$
(up to homeomorphism $\RRR ^2 \approx (D^2)^0$)
in a direction perpendicular to $\RRR^2$
onto $\RRR^2$.
We may suppose that the direction of the projection is generic
i.e all points of self-intersection of the image on $\RRR^2$  are double
and the angles of intersection are non-zero.

Let $r :\RRR^3 \to \RRR^2$
be the  projection which
maps $h(\RRR {\AA}_m) \subset \RRR ^3$ to $\RRR ^2$.
Consider an
oriented tubular fibration
$N \to \CCC {\AA}_m$.
Since $\CCC {\AA}_m$ is diffeomorphic to a
sphere $S^2$ with $g_m$ $1$-handles, one can
consider the restriction of $\CCC {\AA}_m$
diffeomorphic to the torus $T^2$ given by
the sphere $S^2$ with one of $g_m$ $1$-handle.

The oriented tubular neighborhood of
$h(\CCC {\AA}_m) \cap T^2$, as oriented
tubular neighborhood of the torus $T^2$, intersects the solid
torus\\ $S^1 \ti D^2 \supset {\mathcal  M}$
with $\pr M \subset  S^1 \ti S^1$ as the $(2,1)$ torus knot.

Hence, since in $\CCC P^2$, each real line
is split by its real part into two
halves lines conjugate to each others,
and $2$ disjoint circles always divide the torus,
one can assume, without loss of generality, that the
real part of the restriction of $h(\CCC {\AA}_m) \cap T^2$ belongs to the
boundary of the M\"{o}bius band in such a way
that its  projection to $\RRR ^2$ gives one crossing.

Besides, some double-points of $r(h(\RRR \HH_m))$  may
result either from two points of $\RRR {\AA}_m$
which belong to two different handles or
one point which belongs to a $1$-handle and the other point
belongs  to $S^2$.
It easy  to see that, in both cases, the projection leads to an even
number of such double-points.

Hence, from an  argumentation similar to the one given
in case of curves of type $I$,
it follows that $\RRR {\AA}_m$
is the smoothed immersion of a generic immersion of the circle $S^1$
into the plane $\RRR ^2$
with $n$,
$\frac {(m-1)(m-2)} 2 \le n \le \frac {(m-1)(m-2)} 2 + [\frac {m} 2]$,
double points and Whitney index
$\frac {(m-1)(m-2)} 2 + 1$.

Besides, it is obvious, since ${\AA}_m$ is of type $II$,
that at least one-double point of $r(h(\RRR \tilde {\AA}_m))$
is smoothed in the non-coherent direction (in other words,
at least two circles of the $g_m+1$ which
divide the surface $S^2_{g_m}$ glue and disappear in one circle.)
Q.E.D.

\vskip0.1in
Given a generic immersion of the circle $S^1$
into the plane $\RRR ^2$ and  ${\mathcal  S}$ the set of its singular
points.
We shall call
{\it partially smoothed} generic immersion of the circle $S^1$
into the plane $\RRR ^2$
the singular real part
deduced from the generic immersion of the circle by modification
of a set ${\mathcal  K} \subset {\mathcal  S}$, ${\mathcal  K} \not = {\mathcal  S}$
of real double points.
\vskip0.1in
The following Lemma, analogue to the Lemma \mrf{l:l2typ1}
of the previous part,  enlarges the preceding statement
to singular curves.
\vskip0.1in
\bele
\mlb{l:l2typ2}
\vskip0.1in
{\it
Let ${\AA}_m$ be a singular curve of degree $m$ and type $II$
with non-empty real points set $\RRR {\AA}_m$
and non-degenerate singular points,
then its real point set $\RRR {\AA}_m  \subset \RRR P^2$ is a
partially smoothed generic immersion of the circle $S^1$
into the plane $\RRR ^2$
with  $n \ge \frac {(m-1)(m-2)} 2$ double points and Whitney index
$\frac {(m-1)(m-2)} 2 + 1$ if and only if its set of singular points
consists of at most
$\frac {(m-1)(m-2)} 2$ crossings and at least one crossing is smoothed
in $\RRR^2$ in the non-coherent direction.}
\enle
\vskip0.1in
{\bf proof:}
\vskip0.1in
It follows from an argument similar to the one used in
the proof of Lemma \mrf{l:l1typ2}.
Q.E.D
\vskip0.1in
Although
the stratification of the set singular curves in the variety
${\mathcal  C}_{I;g}$
may not be extended to curves of type $II$, we shall lift the path
provided by Whitney's theorem and Lemma \mrf{l:l1typ2}
to algebraic curves of type $II$.

As previously,
we shall consider the space $\RRR {\mathcal  C}_m$ of real algebraic curves
of degree $m$ and its
subset $\RRR {\mathcal  D}_m$ constituted
by real singular algebraic curves of degree $m$.
The set $\RRR {\mathcal  D}_m$ has
an open every dense subset which consist of curves with only one singular
point which is a non-degenerate double point (i.e
a solitary real double-point, or a crossing).
This subset is called {\it principal part} of the set $\RRR {\mathcal  D}_m$.
A  generic path in $\RRR {\mathcal  C}_m$
intersects $\RRR {\mathcal  D}_m$  only in its principal part and only
transversally.
When a singular curve occurs in a generic one-parameter
family curves, the moving curve is passing through a Morse
modification.
\vskip0.1in
Here are  simple properties of Morse modifications which
motivate our study:
\been
\item
Under a Morse modification a curve of type $I$
can turn only to a curve of type $II$ and the
number of its real components decreases.
\item
A complex orientation of a non-singular curve turns
into an orientation of the singular curve which appears at the moment
of the modification:
the orientations defined
by complex orientation of the curve
 on the two arcs
which approach each other and  merge
should be coherent with an orientation of the singular curve.
\enen

We shall generalize the previous method giving any smooth curve
${\AA}_m$ of type  $I$  to smooth curves of type $II$ and deduce
curves of type $II$ from smoothing (including smoothing in a direction
non-coherent to a complex orientation)
generic curves of type $I$, genus  $g$ and type $I$.

\vskip0.1in

Assume $\HH_m$ the Harnack of degree $m$ obtained via the
patchworking method.
For any point $p \in {\mathcal  P}_m$,
the subset
$U(p)=\rho^m(D(p,\e) \ti U_{\CCC}^2) \subset(\CCC^*)^2$
intersects $(\RRR^*)^2$
in four discs.
We shall
denote by ${\mathcal  B}_m$ the set constituted by the union
of the four $2$-discs of $U(p) \cap (\RRR^*)^2$
taken over the points $p \in {\mathcal  P}_m$.

\vskip0.1in
\bepr
\mlb{p:prop1typ2}
{\it
Let $\HH_m$ be the Harnack curve of degree $m$
and ${\AA}_m$ be a smooth curve of type $II$.
Then, up to conj-equivariant isotopy of $\CCC P^2$,
there exists a path
$$S :[0,1] \to \RRR {\mathcal  C}_{m}$$
$S(0)=\RRR \HH_m$,
$S(1)=\RRR {\AA}_m$
which intersects the discriminant hypersurface $\RRR {\mathcal  D}_m$
in such a way that:
\been
\item
any real ordinary double point of curves along $S$ is
a crossing which belongs to a real $2$-disc of ${\mathcal  B}_m$.
\item
the curve ${\AA}_m$ is deduced from smoothing
the real part of a generic curve of degree $m$
of which double points are crossings.
At least one point is smoothed in a direction non-coherent
with a complex orientation of the real part of the curve.
\enen}
\enpr

{\bf proof:}
\vskip0.1in
As already done in case of curves of type $I$,
using the path provided
by  Lemma \mrf{l:l1typ2} and
Whitney's theorem, according to the Lemma \mrf{l:l2typ2},
and the fact that $\CCC \HH_{m}$ and $\CCC {\AA}_m$
are diffeotopic, we shall characterize the track $S$
on $\RRR P^2$
of a diffeotopy $h_t$ of $\CCC P^2$, $t \in [0,1]$,
$h(0)= \CCC \HH_m$, $h(1)=\CCC {\AA}_m$.
\vskip0.1in

{\it  Method which provides  $(\CCC P^2,\CCC{\AA}_m)$ from
$(\CCC P^2,\CCC \HH_m)$.}
\vskip0.1in

Our method is a slightly modified version of the
preceding one where any curve ${\AA}_m$ of type $I$
is obtained from $\HH_m$.
\vskip0.1in
We shall prove that
up to modify the coefficients of the polynomial
giving the curve ${\AA}_m$,
one can always assume that there exists a diffeotopy $h_t$
of $\CCC P^2$
$h(0)= \CCC \HH_m$, $h(1)=\CCC {\AA}_m$
with the property
$h_t(\CCC {\AA}_m   \cap U(p))  \subset U(p)$ for any $p \in {\mathcal  P}_m$.
\vskip0.1in
We shall refer to the method of section \mrf{su:typ1}
and stress only the modification.
\vskip0.1in
According to the Lemma \mrf{l:l1typ2},
consider $\RRR \HH_m$ and $\RRR {\AA}_m$
as smoothed immersions of the circle into
the plane with the same Whitney index.
Denote
$\RRR \tilde{\HH}_m$, $\RRR \tilde{\AA}_m$
the corresponding immersions.
From the Whitney's theorem (\cite{Whi}),
there exists a path  $\tilde{h}$ which connects
$\RRR \tilde{\HH}_m$ and $\RRR \tilde{\AA}_m$.
Thus,
the definition of the path $S$ is reduced to the definition of
a lifting of the path  which connects
$\RRR \tilde{\HH}_m$ and $\RRR \tilde{\AA}_m$ in
the space of immersions
of the circle into $\RRR^2$ to a path in the space of curve of degree $m$
and type $II$.

Define now the {\it smoothed perestroika} of a perestroika,
the change obtained by smoothing the fragments involved in the perestroika
where smoothing before perestroika
are taken  only in a direction  coherent
and smoothing after perestroika
may be taken in a direction non-coherent.

Then, as in the proof of  proposition  \mrf{p:prop1typ1}
of section \mrf{su:typ1}
where curves of type $I$ were under consideration,
since any smooth curve is irreducible, and any reducible polynomial
is the product of a finite reducible one,
according to Lemma \mrf{l:l1typ2} and Lemma \mrf{l:l2typ2},
we may lift the path $\tilde{h}$
in  the space of immersion of the circle into the plane
from smoothed perestroikas in the  space
${\mathcal  C}_{I;g}$ of curve of degree $n \le m$, type $I$ and genus
$0 \le g \le \frac {(n-1).(n-2)} 2$.

In such a way, using the path provided
by  Lemma \mrf{l:l1typ2} and
Whitney's theorem, according to the Lemma \mrf{l:l2typ2} and the
fact that $\CCC \HH_{m}$ and $\CCC {\AA}_m$
are diffeotopic, we characterize the track $S$
on $\RRR P^2$
of a diffeotopy $h_t$ of $\CCC P^2$, $t \in [0,1]$,
$h(0)= \CCC \HH_m$, $h(1)=\CCC {\AA}_m$.
\vskip0.1in
Q.E.D
\vskip0.1in

\bere
 As already known, if two  $M-2$-curves
are obtained one from the other by a deformation through a double
nondegenerate point, one of them is of type $I$
and the other is of type $II$.
Since, any $M-1$-curve is of type $II$,
according to our method
it remains to say that two  deformations in the non-coherent
direction may turn a curve of type $I$ to a curve of type $I$.
\enre

As a curve of degree $m$ moves as a point of $\RRR {\mathcal  C}_m$ along
an arc which has a transversal intersection with the principal part,
then the set of real points of this curve undergoes either
a Morse modification of index $0$ or $2$
(the curves acquires a solitary singular point, which then becomes
a new oval, or else one of the ovals contracts to point)
or Morse modification of index $0$, an oval contracts to a point
or a Morse modification of index $1$
( two real arcs of the curve approach one other and merge,
then diverge in their modified form.)

As already introduced, we
denote  $A_{m+1}$ the set of crossings
perturbed in the construction
of $\HH_{m+1}$ from $\HH _m$ introduced in proposition \mrf{p:prop7}.

In case $m+1=2k$,
$A_{2k}=\{a_1,...,a_{2k-1},...,a_{4k-2} \}$
where  $a_1,...a_{2k-1}$ are
crossings  of $\HH_{2k-1}\cup L$
and $a_{2k},...,a_{4k-2}$ are crossings
of $\tilde{X}_{2k; \tau}$.

In case $m+1=2k +1$,
$A_{2k+1}=\{a_1,...,a_{2k} \}$
where $a_1,...a_{2k}$ are
crossings of $\HH_{2k}\cup L$.

We may characterize Morse modifications on a curves of
degree $m$ as follows:

\vskip0.1in

\begin{coro}
{\it Any Morse modification on a curve of degree $m$
may be described, up to conj-equivariant isotopy,
by perturbations of crossings
of a subset of  $\cup_{n=1}^m A_m$}
\end{coro}
\vskip0.1in
{\bf proof:}
It is a straightforward  consequence of proposition \mrf{p:prop1typ2}.
Q.E.D.
\newpage
{\bf  Description of pairs $(\CCC P^2,\CCC {\AA}_m)$
up to conj-equivariant isotopy }
\vskip0.1in
\beth
\mlb{t:th1typ2}
\vskip0.1in
{\it Let ${\AA}_m$ be a curve of degree $m$ and type $II$
with non-empty real part.
There exists a finite number $I$ of disjoints $4$-balls $B(a_i)$
invariant by complex conjugation and centered in points
$a_i$ of $\RRR P^2$ such that,
up to conj-equivariant isotopy of $\CCC P^2$ :
\been
\item
 ${\AA}_m \bk \cup_{i \in I}B(a_i) =
 \cup_{i=1}^m L_i\bk \cup_{i \in I}B(a_i) $
where $L_1$,...,$L_m$ are $m$ projective lines
with $L_i  \bk \cup_{i \in I} B(a_i) \cap
L_j  \bk \cup _{i \in I} B(a_i)=\emptyset$
for any $i \not = j$,$1 \le i,j \le m$.
\item
situations inside $4$-balls $B(a_i)$
are covered by  perturbation of type $1$ or type $2$
of the crossing $a_i$.
\enen}
\enth

\vskip0.1in
{\bf proof:}
\vskip0.1in

Our proof is based on the previous method and the
theorem \mrf{t:theo1}.
The argument is analogous to the one of the proof of
theorem \mrf{t:th2typ1} where curves of type $I$ were under consideration.
Q.E.D

\chapter{Arnold surfaces of curves of even degree\\
-Proof of the Rokhlin's conjecture-}
\lb{ch:arn}
\vskip0.1in
In this section, we shall prove the Rokhlin's conjecture:
\vskip0.1in
{\bf Theorem} (Rokhlin's conjecture)
\vskip0.1in

{\it Arnold surfaces $\gA$ are standard for all  curves of
even degree with non-empty real part.}
\vskip0.1in

{\bf  proof:}
\vskip0.1in
The proof is based on the Livingston's theorem
and the following statement
(deduced from  theorem \mrf{t:th2typ1}
and theorem \mrf{t:th1typ2} )

Given ${\AA}_m$ a curve of degree $m$ and type $I$ or type $II$
with non-empty real part.
There exists a finite number $I$ of disjoints $4$-balls $B(a_i)$
invariant by complex conjugation and centered in points
$a_i$ of $\RRR P^2$ such that
up to conj-equivariant isotopy of $\CCC P^2$:
\been
\item
 ${\AA}_m \bk \cup_{i \in I}B(a_i) =
 \cup_{i=1}^m L_i\bk \cup_{i \in I}B(a_i) $
where $L_1$,...,$L_m$ are $m$ projective lines
with $L_i  \bk \cup_{i \in I} B(a_i) \cap
L_j  \bk \cup _{i \in I} B(a_i)=\emptyset$
for any $i \not = j$,$1 \le i,j \le m$.
\item
situations inside any $4$-balls $B(a_i)$ are covered by
perturbation of type $1$ or type $2$ of the crossing $a_i$.
\enen

It follows from an argumentation similar to the one
given in Chapter \ref{ch:ArnHar} (see theorem \mrf{t:theo2}
where it is stated that
Arnold surfaces of Harnack curves $\HH_{2k}$ are standard surfaces in
$S^4$)
that any Arnold surface of a curve of
even degree $2k$ with non-empty real part is standard.

Q.E.D
\newpage

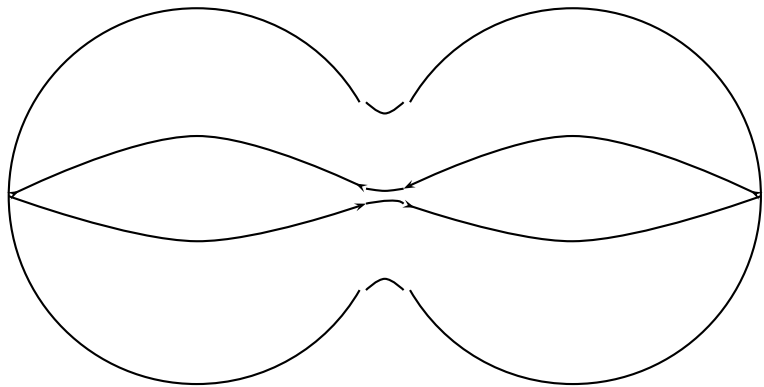
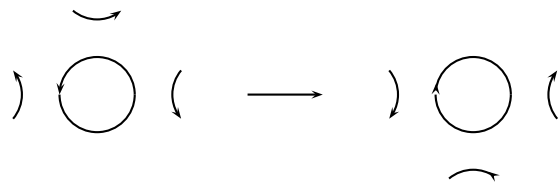
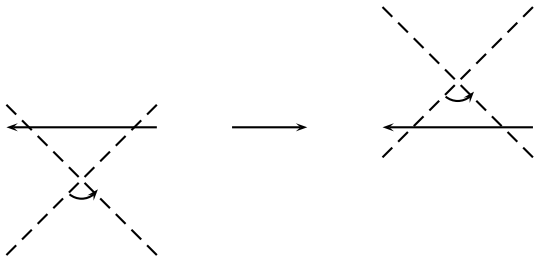
\begin{pspicture}(0,1)(0,6)
\rput(1,6){{\it figure 1.1:perturbation of type 1}}

\psarc(2.5,2.5){2.5}{30}{330}
\psarc(7.5,2.5){2.5}{210}{150}

\pscurve{>-<}(0,2.5)(2.5,3.3)(4.75,2.6)

\pscurve{->}(0,2.5)(2.5,1.9)(4.75,2.4)
\pscurve{>->}(10,2.5)(7.5,1.9)(5.25,2.4)

\pscurve{-<}(10,2.5)(7.5,3.3)(5.25,2.6)

\pscurve(4.75,3.75)(5,3.6)(5.25,3.75)
\pscurve(4.75,1.25)(5,1.4)(5.25,1.25)
\pscurve (4.82,2)(4.85,1.6)(4.9,1.4)

\psline{<-<}(4.9,1.4)(5.1,1.4)

\pscurve(5.1,1.4)(5.2,2)(5.25,2.4)
\psline(4.75,2.4)(4.8,2.4)
\pscurve(4.8,2.4)(4.82,2.5)(4.83,2.6)
\pscurve(4.83,2.6)(4.85,3.4)(4.9,3.6)
\psline{>->}(4.9,3.6)(5.1,3.6)
\pscurve(5.1,3.6)(5.2,2.8)(5.25,2.6)
\end{pspicture}

\rput(1,-6){{\it figure 1.2: perturbation of type 2}}

\begin{pspicture}(0,9)(0,14)
\psarc(2.5,2.5){2.5}{30}{330}
\psarc(7.5,2.5){2.5}{210}{150}
\pscurve{>-<}(0,2.5)(2.5,3.3)(4.75,2.6)

\pscurve{->}(0,2.5)(2.5,1.9)(4.75,2.4)

\pscurve{-<}(10,2.5)(7.5,1.9)(5.25,2.4)

\pscurve{>->}(10,2.5)(7.5,3.3)(5.25,2.6)

\pscurve (4.75,3.75)(5,3.6)(5.25,3.75)
\pscurve (4.75,1.25)(5,1.4)(5.25,1.25)
\pscurve(4.75,2.6)(5,2.57)(5.25,2.6)
\pscurve(4.75,2.4)(5.2,2.43)(5.25,2.4)
\end{pspicture}

\newpage

\begin{pspicture}(-5,-5)(5,5)
\rput(0,3){figure 3.1}
\psarc{-<}(-1,1){1}{270}{450}
\rput(0.3,1){$p_1$}
\psarc{->}(3,1){1}{90}{270}
\rput(1.8,1){$p_2$}

\psarc{-<}(1,1){1}{90}{270}

\psarc{-<}(1,1){1}{270}{450}
\end{pspicture}

\newpage

\begin{pspicture}(0,1)(0,5)
\rput(0,3){Figure 5.1-Smoothings}
\psline{>->}(-1,-1)(1,1)
\psline{>->}(-1,1)(1,-1)
\psarc{>->}(0,1){0.5}{180}{360}
\psarc{<-<}(0,-1){0.5}{0}{180}
\end{pspicture}
\begin{pspicture}(-5,0)(0,1)
\psline[linestyle=dashed](-1,-1)(1,1)
\psline[linestyle=dashed](-1,1)(1,-1)
\psarc{->}(0,0){0.5}{30}{150}
\psline{->}(3,0)(4,0)
\psarc{->}(5,0){0.25}{0}{180}
\psarc(5,0){0.25}{180}{360}
\end{pspicture}
\begin{pspicture}(0,5)(0,10)
\rput(0,3){Figure 5.2-Weak triple point Perestroika}
\psline{->}(-1.5,0)(1.5,0)
\psline{->}(-1.5,-0.5)(0.5,1.5)
\psline{->}(-0.5,1.5)(1.5,-0.5)
\psline{->}(2,0)(3,0)
\psline{->}(4,0)(7,0)
\psline{->}(4,0.5)(5.5,-1.5)
\psline{->}(5,-1.5)(7,0.5)
\end{pspicture}
\begin{pspicture}(0,9)(0,14)
\rput(0,3){Figure 5.(3.a)-Strong triple point Perestroika}
\psline{<-}(-1.5,0)(1.5,0)
\psline{->}(-1.5,-0.5)(0.5,1.5)
\psline{->}(-0.5,1.5)(1.5,-0.5)
\psline{->}(2,0)(3,0)
\psline{<-}(4,0)(7,0)
\psline{->}(4,0.5)(5.5,-1.5)
\psline{->}(5,-1.5)(7,0.5)
\end{pspicture}
\begin{pspicture}(0,13)(0,18)

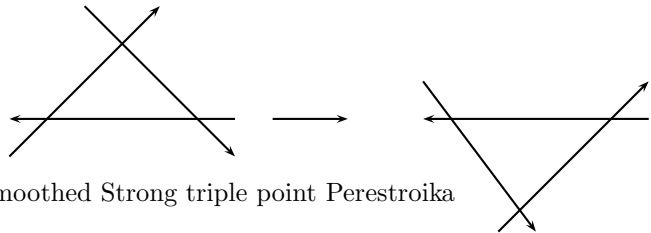
\rput(0,3){Figure 5.(3.b)-Smoothed Strong triple point Perestroika}
\psarc{->}(0,0){0.5}{0}{180}
\psarc(0,0){0.5}{180}{360}

\psarc{->}(-1.5,0){0.5}{320}{400}
\psarc{->}(1.5,0){0.5}{140}{220}
\psarc{->}(0,1.5){0.5}{230}{310}

\psline{->}(2,0)(3,0)
\psarc{-<}(5,0){0.5}{0}{180}
\psarc(5,0){0.5}{180}{360}

\psarc{<-}(3.5,0){0.5}{320}{400}
\psarc{<-}(6.5,0){0.5}{140}{220}
\psarc{>-}(5,-1.5){0.5}{50}{130}
\end{pspicture}
\newpage
\begin{pspicture}(0,1)(0,6)
\rput(0,3){Figure 5.4-Smoothed real inverse tangency Perestroika}
\psarc{->}(-1.5,0){1}{320}{400}
\psarc{->}(1.5,0){1}{140}{220}
\psline{->}(2,0)(3,0)
\psarc{<-}(4,2){1}{230}{310}
\psarc{->}(4,0){0.5}{0}{180}
\psarc(4,0){0.5}{180}{360}
\psarc{<-}(4,-2){1}{50}{130}
\end{pspicture}
\begin{pspicture}(0,9)(0,14)
\rput(0,3){Figure 5.5-Triple-point with imaginary branches Perestroika}
\psline{<-}(-1,0.7)(1,0.7)
\psline[linestyle=dashed](-1,-1)(1,1)
\psline[linestyle=dashed](-1,1)(1,-1)
\psarc{->}(0,0){0.25}{230}{330}
\psline{->}(2,0.7)(3,0.7)
\psline{<-}(4,0.7)(6,0.7)
\psline[linestyle=dashed](4,0.3)(6,2.3)
\psline[linestyle=dashed](4,2.3)(6,0.3)
\psarc{->}(5,1.3){0.25}{230}{330}
\end{pspicture}
\newpage

\begin{pspicture}(0,9)(0,14)
\rput(1,5){Morse modification in $\RRR ^2$ }
\rput(1,4){in the direction coherent and non-coherent
to a complex orientation}
\psline{>->}(-2,-2)(2,2)
\psline{>->}(-2,2)(2,-2)
\psarc{>->}(0,2){1}{180}{360}
\psarc{<-<}(0,-2){1}{0}{180}
\psarc(2,0){1}{90}{270}
\psarc(-2,0){1}{270}{450}
\end{pspicture}
\newpage

\begin{pspicture}(-6,6)(18,-18)
\rput(0,-16){Harnack T-curve of degree 10}

\psline(1,0)(0,1)
\psline(2,0)(0,2)
\psline(3,0)(0,3)
\psline(4,0)(0,4)
\psline(5,0)(0,5)
\psline(6,0)(0,6)
\psline(7,0)(0,7)
\psline(8,0)(0,8)
\psline(9,0)(0,9)
\psline(10,0)(0,10)

\psline(-1,0)(0,1)
\psline(-2,0)(0,2)
\psline(-3,0)(0,3)
\psline(-4,0)(0,4)
\psline(-5,0)(0,5)
\psline(-6,0)(0,6)
\psline(-7,0)(0,7)
\psline(-8,0)(0,8)
\psline(-9,0)(0,9)
\psline(-10,0)(0,10)

\psline(-1,0)(0,-1)
\psline(-2,0)(0,-2)
\psline(-3,0)(0,-3)
\psline(-4,0)(0,-4)
\psline(-5,0)(0,-5)
\psline(-6,0)(0,-6)
\psline(-7,0)(0,-7)
\psline(-8,0)(0,-8)
\psline(-9,0)(0,-9)
\psline(-10,0)(0,-10)

\psline(1,0)(0,-1)
\psline(2,0)(0,-2)
\psline(3,0)(0,-3)
\psline(4,0)(0,-4)
\psline(5,0)(0,-5)
\psline(6,0)(0,-6)
\psline(7,0)(0,-7)
\psline(8,0)(0,-8)
\psline(9,0)(0,-9)
\psline(10,0)(0,-10)

\psline(-10,0)(10,0)
\psline(-9,1)(9,1)
\psline(-8,2)(8,2)
\psline(-7,3)(7,3)
\psline(-6,4)(6,4)
\psline(-5,5)(5,5)
\psline(-4,6)(4,6)
\psline(-3,7)(3,7)
\psline(-2,8)(2,8)
\psline(-1,9)(1,9)

\psline(-9,-1)(9,-1)
\psline(-8,-2)(8,-2)
\psline(-7,-3)(7,-3)
\psline(-6,-4)(6,-4)
\psline(-5,-5)(5,-5)
\psline(-4,-6)(4,-6)
\psline(-3,-7)(3,-7)
\psline(-2,-8)(2,-8)
\psline(-1,-9)(1,-9)

\psline(0,-10)(0,10)
\psline(1,-9)(1,9)
\psline(2,-8)(2,8)
\psline(3,-7)(3,7)
\psline(4,-6)(4,6)
\psline(5,-5)(5,5)
\psline(6,-4)(6,4)
\psline(7,-3)(7,3)
\psline(8,-2)(8,2)
\psline(9,-1)(9,1)

\psline(-1,-9)(-1,9)
\psline(-2,-8)(-2,8)
\psline(-3,-7)(-3,7)
\psline(-4,-6)(-4,6)
\psline(-5,-5)(-5,5)
\psline(-6,-4)(-6,4)
\psline(-7,-3)(-7,3)
\psline(-8,-2)(-8,2)
\psline(-9,-1)(-9,1)

\psline(0.5,-0.5)(0.5,0)
\psline(0.5,0)(0,0.5)
\psline(0,0.5)(-0.5,0.5)

\psline(-0.5,0.5)(-0.5,1)
\psline(-0.5,1)(0,1.5)

\psline(0,1.5)(0,1.5)
\psline(0,1.5)(0.5,1.5)
\psline(0.5,1.5)(0.5,2)
\psline(0.5,2)(0,2.5)
\psline(0,2.5)(-0.5,2.5)
\psline(-0.5,3)(-0.5,2.5)
\psline(-0.5,3)(0,3.5)
\psline(0,3.5)(0.5,3.5)
\psline(0.5,3.5)(0.5,4)

\psline(0,4.5)(0.5,4)
\psline(0,4.5)(-0.5,4.5)
\psline(-0.5,5)(-0.5,4.5)
\psline(-0.5,5)(0,5.5)
\psline(0,5.5)(0.5,5.5)
\psline(0.5,5.5)(0.5,6)

\psline(0,6.5)(0.5,6)
\psline(0,6.5)(-0.5,6.5)
\psline(-0.5,7)(-0.5,6.5)
\psline(-0.5,7)(0,7.5)
\psline(0,7.5)(0.5,7.5)
\psline(0.5,7.5)(0.5,8)

\psline(0,8.5)(0.5,8)
\psline(0,8.5)(-0.5,8.5)
\psline(-0.5,9)(-0.5,8.5)
\psline(-0.5,9)(0,9.5)
\psline(0,9.5)(0.5,9.5)

\psline(0,1.5)(-0.5,1)

\psline(0,1.5)(-0.5,1)

\psline(-0.5,-1)(-1,-0.5)
\psline(1,-0.5)(0.5,-0.5)
\psline(-0.5,0.5)(-0.5,1)

\psline(1.5,0)(1,-0.5)
\psline(1.5,0)(1.5,0.5)
\psline(1.5,0.5)(2,0.5)
\psline(2,0.5)(2.5,0)
\psline(2.5,0)(2.5,-0.5)
\psline(2.5,-0.5)(3,-0.5)
\psline(3,-0.5)(3.5,0)

\psline(3.5,0)(3,-0.5)
\psline(3.5,0)(3.5,0.5)
\psline(3.5,0.5)(4,0.5)
\psline(4,0.5)(4.5,0)
\psline(4.5,0)(4.5,-0.5)
\psline(4.5,-0.5)(5,-0.5)
\psline(5,-0.5)(5.5,0)

\psline(5.5,0)(5,-0.5)
\psline(5.5,0)(5.5,0.5)
\psline(5.5,0.5)(6,0.5)
\psline(6,0.5)(6.5,0)
\psline(6.5,0)(6.5,-0.5)
\psline(6.5,-0.5)(7,-0.5)
\psline(7,-0.5)(7.5,0)

\psline(7.5,0)(7,-0.5)
\psline(7.5,0)(7.5,0.5)
\psline(7.5,0.5)(8,0.5)
\psline(8,0.5)(8.5,0)
\psline(8.5,0)(8.5,-0.5)
\psline(8.5,-0.5)(9,-0.5)
\psline(9,-0.5)(9.5,0)
\psline(9.5,0)(9.5,0.5)

\psline(2,1.5)(1.5,2)
\psline(1.5,2)(1.5,2.5)
\psline(2,1.5)(2.5,1.5)
\psline(2.5,1.5)(2.5,2)
\psline(2,2.5)(2.5,2)
\psline(2,2.5)(1.5,2.5)

\psline(2,3.5)(1.5,4)
\psline(1.5,4)(1.5,4.5)
\psline(2,3.5)(2.5,3.5)
\psline(2.5,3.5)(2.5,4)
\psline(2,4.5)(2.5,4)
\psline(2,4.5)(1.5,4.5)

\psline(2,5.5)(1.5,6)
\psline(1.5,6)(1.5,6.5)
\psline(2,5.5)(2.5,5.5)
\psline(2.5,5.5)(2.5,6)
\psline(2,6.5)(2.5,6)
\psline(2,6.5)(1.5,6.5)

\psline(4,1.5)(3.5,2)
\psline(3.5,2)(3.5,2.5)
\psline(4,1.5)(4.5,1.5)
\psline(4.5,1.5)(4.5,2)
\psline(4,2.5)(4.5,2)
\psline(4,2.5)(3.5,2.5)

\psline(6,1.5)(5.5,2)
\psline(5.5,2)(5.5,2.5)
\psline(6,1.5)(6.5,1.5)
\psline(6.5,1.5)(6.5,2)
\psline(6,2.5)(6.5,2)
\psline(6,2.5)(5.5,2.5)

\psline(4,3.5)(3.5,4)
\psline(3.5,4)(3.5,4.5)
\psline(4,3.5)(4.5,3.5)
\psline(4.5,3.5)(4.5,4)
\psline(4,4.5)(4.5,4)
\psline(4,4.5)(3.5,4.5)

\psline(1.5,8)(1.5,8.5)
\psline(1.5,8)(2,7.5)
\psline(2,7.5)(2.5,7.5)

\psline(3.5,6)(3.5,6.5)
\psline(3.5,6)(4,5.5)
\psline(4,5.5)(4.5,5.5)

\psline(5.5,4)(5.5,4.5)
\psline(5.5,4)(6,3.5)
\psline(6,3.5)(6.5,3.5)

\psline(7.5,2)(7.5,2.5)
\psline(7.5,2)(8,1.5)
\psline(8,1.5)(8.5,1.5)

\psline(-1,-0.5)(-0.5,-1)
\psline(-0.5,-1)(-0.5,-1.5)
\psline(-1,-0.5)(-1.5,-0.5)
\psline(-1.5,-0.5)(-1.5,-1)
\psline(-1,-1.5)(-1.5,-1)
\psline(-1,-1.5)(-0.5,-1.5)

\psline(-1,-6.5)(-0.5,-7)
\psline(-0.5,-7)(-0.5,-7.5)
\psline(-1,-6.5)(-1.5,-6.5)
\psline(-1.5,-6.5)(-1.5,-7)
\psline(-1,-7.5)(-1.5,-7)
\psline(-1,-7.5)(-0.5,-7.5)

\psline(1,-7.5)(0.5,-8)
\psline(0.5,-8)(0.5,-8.5)
\psline(1,-7.5)(1.5,-7.5)
\psline(1.5,-7.5)(1.5,-8)
\psline(1,-8.5)(1.5,-8)
\psline(1,-8.5)(0.5,-8.5)

\psline(-1,-6.5)(-0.5,-7)
\psline(-0.5,-7)(-0.5,-7.5)
\psline(-1,-6.5)(-1.5,-6.5)
\psline(-1.5,-6.5)(-1.5,-7)
\psline(-1,-7.5)(-1.5,-7)
\psline(-1,-7.5)(-0.5,-7.5)

\psline(1,-1.5)(0.5,-2)
\psline(0.5,-2)(0.5,-2.5)
\psline(1,-1.5)(1.5,-1.5)
\psline(1.5,-1.5)(1.5,-2)
\psline(1,-2.5)(1.5,-2)
\psline(1,-2.5)(0.5,-2.5)

\psline(-1,-2.5)(-0.5,-3)
\psline(-0.5,-3)(-0.5,-3.5)
\psline(-1,-2.5)(-1.5,-2.5)
\psline(-1.5,-2.5)(-1.5,-3)
\psline(-1,-3.5)(-1.5,-3)
\psline(-1,-3.5)(-0.5,-3.5)

\psline(1,-3.5)(0.5,-4)
\psline(0.5,-4)(0.5,-4.5)
\psline(1,-3.5)(1.5,-3.5)
\psline(1.5,-3.5)(1.5,-4)
\psline(1,-4.5)(1.5,-4)
\psline(1,-4.5)(0.5,-4.5)

\psline(-1,-4.5)(-0.5,-5)
\psline(-0.5,-5)(-0.5,-5.5)
\psline(-1,-4.5)(-1.5,-4.5)
\psline(-1.5,-4.5)(-1.5,-5)
\psline(-1,-5.5)(-1.5,-5)
\psline(-1,-5.5)(-0.5,-5.5)

\psline(1,-5.5)(0.5,-6)
\psline(0.5,-6)(0.5,-6.5)
\psline(1,-5.5)(1.5,-5.5)
\psline(1.5,-5.5)(1.5,-6)
\psline(1,-6.5)(1.5,-6)
\psline(1,-6.5)(0.5,-6.5)

\psline(-3,-0.5)(-2.5,-1)
\psline(-2.5,-1)(-2.5,-1.5)
\psline(-3,-0.5)(-3.5,-0.5)
\psline(-3.5,-0.5)(-3.5,-1)
\psline(-3,-1.5)(-3.5,-1)
\psline(-3,-1.5)(-2.5,-1.5)

\psline(-5,-0.5)(-4.5,-1)
\psline(-4.5,-1)(-4.5,-1.5)
\psline(-5,-0.5)(-5.5,-0.5)
\psline(-5.5,-0.5)(-5.5,-1)
\psline(-5,-1.5)(-5.5,-1)
\psline(-5,-1.5)(-4.5,-1.5)

\psline(-3,-2.5)(-2.5,-3)
\psline(-2.5,-3)(-2.5,-3.5)
\psline(-3,-2.5)(-3.5,-2.5)
\psline(-3.5,-2.5)(-3.5,-3)
\psline(-3,-3.5)(-3.5,-3)
\psline(-3,-3.5)(-2.5,-3.5)

\psline(-3,-4.5)(-2.5,-5)
\psline(-2.5,-5)(-2.5,-5.5)
\psline(-3,-4.5)(-3.5,-4.5)
\psline(-3.5,-4.5)(-3.5,-5)
\psline(-3,-5.5)(-3.5,-5)
\psline(-3,-5.5)(-2.5,-5.5)

\psline(-7,-0.5)(-6.5,-1)
\psline(-6.5,-1)(-6.5,-1.5)
\psline(-7,-0.5)(-7.5,-0.5)
\psline(-7.5,-0.5)(-7.5,-1)
\psline(-7,-1.5)(-7.5,-1)
\psline(-7,-1.5)(-6.5,-1.5)

\psline(-5,-2.5)(-4.5,-3)
\psline(-4.5,-3)(-4.5,-3.5)
\psline(-5,-2.5)(-5.5,-2.5)
\psline(-5.5,-2.5)(-5.5,-3)
\psline(-5,-3.5)(-5.5,-3)
\psline(-5,-3.5)(-4.5,-3.5)

\psline(-5,-0.5)(-4.5,-1)
\psline(-4.5,-1)(-4.5,-1.5)
\psline(-5,-0.5)(-5.5,-0.5)
\psline(-5.5,-0.5)(-5.5,-1)
\psline(-5,-1.5)(-5.5,-1)
\psline(-5,-1.5)(-4.5,-1.5)

\psline(-3,-2.5)(-2.5,-3)
\psline(-2.5,-3)(-2.5,-3.5)
\psline(-3,-2.5)(-3.5,-2.5)
\psline(-3.5,-2.5)(-3.5,-3)
\psline(-3,-3.5)(-3.5,-3)
\psline(-3,-3.5)(-2.5,-3.5)

\psline(-3,-4.5)(-2.5,-5)
\psline(-2.5,-5)(-2.5,-5.5)
\psline(-3,-4.5)(-3.5,-4.5)
\psline(-3.5,-4.5)(-3.5,-5)
\psline(-3,-5.5)(-3.5,-5)
\psline(-3,-5.5)(-2.5,-5.5)

\psline(-3,-4.5)(-2.5,-5)
\psline(-2.5,-5)(-2.5,-5.5)
\psline(-3,-4.5)(-3.5,-4.5)
\psline(-3.5,-4.5)(-3.5,-5)
\psline(-3,-5.5)(-3.5,-5)
\psline(-3,-5.5)(-2.5,-5.5)

\psline(-7,-0.5)(-6.5,-1)
\psline(-6.5,-1)(-6.5,-1.5)
\psline(-7,-0.5)(-7.5,-0.5)
\psline(-7.5,-0.5)(-7.5,-1)
\psline(-7,-1.5)(-7.5,-1)
\psline(-7,-1.5)(-6.5,-1.5)

\psline(-5,-2.5)(-4.5,-3)
\psline(-4.5,-3)(-4.5,-3.5)
\psline(-5,-2.5)(-5.5,-2.5)
\psline(-5.5,-2.5)(-5.5,-3)
\psline(-5,-3.5)(-5.5,-3)
\psline(-5,-3.5)(-4.5,-3.5)

\psline(-2.5,1)(-2.5,0.5)
\psline(-2.5,1)(-2,1.5)
\psline(-2,1.5)(-1.5,1.5)
\psline(-1.5,1.5)(-1.5,1)
\psline(-1.5,1)(-2,0.5)
\psline(-2,0.5)(-2.5,0.5)

\psline(-4.5,1)(-4.5,0.5)
\psline(-4.5,1)(-4,1.5)
\psline(-4,1.5)(-3.5,1.5)
\psline(-3.5,1.5)(-3.5,1)
\psline(-3.5,1)(-4,0.5)
\psline(-4,0.5)(-4.5,0.5)

\psline(-6.5,1)(-6.5,0.5)
\psline(-6.5,1)(-6,1.5)
\psline(-6,1.5)(-5.5,1.5)
\psline(-5.5,1.5)(-5.5,1)
\psline(-5.5,1)(-6,0.5)
\psline(-6,0.5)(-6.5,0.5)

\psline(-8.5,1)(-8.5,0.5)
\psline(-8.5,1)(-8,1.5)
\psline(-8,1.5)(-7.5,1.5)
\psline(-7.5,1.5)(-7.5,1)
\psline(-7.5,1)(-8,0.5)
\psline(-8,0.5)(-8.5,0.5)

\psline(-2.5,3)(-2.5,2.5)
\psline(-2.5,3)(-2,3.5)
\psline(-2,3.5)(-1.5,3.5)
\psline(-1.5,3.5)(-1.5,3)
\psline(-1.5,3)(-2,2.5)
\psline(-2,2.5)(-2.5,2.5)

\psline(-2.5,5)(-2.5,4.5)
\psline(-2.5,5)(-2,5.5)
\psline(-2,5.5)(-1.5,5.5)
\psline(-1.5,5.5)(-1.5,5)
\psline(-1.5,5)(-2,4.5)
\psline(-2,4.5)(-2.5,4.5)

\psline(-2.5,7)(-2.5,6.5)
\psline(-2.5,7)(-2,7.5)
\psline(-2,7.5)(-1.5,7.5)
\psline(-1.5,7.5)(-1.5,7)
\psline(-1.5,7)(-2,6.5)
\psline(-2,6.5)(-2.5,6.5)

\psline(3,-1.5)(2.5,-2)
\psline(2.5,-2)(2.5,-2.5)
\psline(3,-1.5)(3.5,-1.5)
\psline(3.5,-1.5)(3.5,-2)
\psline(3,-2.5)(3.5,-2)
\psline(3,-2.5)(2.5,-2.5)

\psline(5,-1.5)(4.5,-2)
\psline(4.5,-2)(4.5,-2.5)
\psline(5,-1.5)(5.5,-1.5)
\psline(5.5,-1.5)(5.5,-2)
\psline(5,-2.5)(5.5,-2)
\psline(5,-2.5)(4.5,-2.5)

\psline(7,-1.5)(6.5,-2)
\psline(6.5,-2)(6.5,-2.5)
\psline(7,-1.5)(7.5,-1.5)
\psline(7.5,-1.5)(7.5,-2)
\psline(7,-2.5)(7.5,-2)
\psline(7,-2.5)(6.5,-2.5)

\psline(3,-3.5)(2.5,-4)
\psline(2.5,-4)(2.5,-4.5)
\psline(3,-3.5)(3.5,-3.5)
\psline(3.5,-3.5)(3.5,-4)
\psline(3,-4.5)(3.5,-4)
\psline(3,-4.5)(2.5,-4.5)

\psline(3,-5.5)(2.5,-6)
\psline(2.5,-6)(2.5,-6.5)
\psline(3,-5.5)(3.5,-5.5)
\psline(3.5,-5.5)(3.5,-6)
\psline(3,-6.5)(3.5,-6)
\psline(3,-6.5)(2.5,-6.5)

\psline(5,-3.5)(4.5,-4)
\psline(4.5,-4)(4.5,-4.5)
\psline(5,-3.5)(5.5,-3.5)
\psline(5.5,-3.5)(5.5,-4)
\psline(5,-4.5)(5.5,-4)
\psline(5,-4.5)(4.5,-4.5)

\psline(-4.5,3)(-4.5,2.5)
\psline(-4.5,3)(-4,3.5)
\psline(-4,3.5)(-3.5,3.5)
\psline(-3.5,3.5)(-3.5,3)
\psline(-3.5,3)(-4,2.5)
\psline(-4,2.5)(-4.5,2.5)

\psline(-4.5,5)(-4.5,4.5)
\psline(-4.5,5)(-4,5.5)
\psline(-4,5.5)(-3.5,5.5)
\psline(-3.5,5.5)(-3.5,5)
\psline(-3.5,5)(-4,4.5)
\psline(-4,4.5)(-4.5,4.5)

\psline(-6.5,3)(-6.5,2.5)
\psline(-6.5,3)(-6,3.5)
\psline(-6,3.5)(-5.5,3.5)
\psline(-5.5,3.5)(-5.5,3)
\psline(-5.5,3)(-6,2.5)
\psline(-6,2.5)(-6.5,2.5)

\psline(-9.5,-0.5)(-9,-0.5)
\psline(-9,-0.5)(-8.5,-1)
\psline(-8.5,-1.5)(-8.5,-1)

\psline(-7.5,-2.5)(-7,-2.5)
\psline(-7,-2.5)(-6.5,-3)
\psline(-6.5,-3.5)(-6.5,-3)

\psline(-5.5,-4.5)(-5,-4.5)
\psline(-5,-4.5)(-4.5,-5)
\psline(-4.5,-5.5)(-4.5,-5)

\psline(-3.5,-6.5)(-3,-6.5)
\psline(-3,-6.5)(-2.5,-7)
\psline(-2.5,-7.5)(-2.5,-7)

\psline(-1.5,-8.5)(-1,-8.5)
\psline(-1,-8.5)(-0.5,-9)
\psline(-0.5,-9.5)(-0.5,-9)

\end{pspicture}

\end{document}